\documentclass[a4paper,11pt,oneside]{article}
\usepackage[utf8x]{inputenc}
\usepackage{graphicx}
\usepackage{enumerate}
\usepackage{todonotes}
\usepackage{colonequals}
\usepackage[nottoc,notlot,notlof]{tocbibind}   
\usepackage[square, numbers, sort]{natbib}
\usepackage{subcaption}
\usepackage{layouts}

\usepackage{amsmath}
\usepackage{amsthm}
\usepackage{amssymb}
\usepackage{mathtools}
\usepackage{authblk}
\usepackage[shortlabels]{enumitem}
\usepackage{fullpage}
\usepackage{esint}

\usepackage{algorithm}
\usepackage{algpseudocode}

\usepackage{tikz}
\newcommand*\circled[1]{\tikz[baseline=(char.base)]{
            \node[shape=circle,draw,inner sep=1pt] (char) {#1};}}
\usetikzlibrary{shapes.geometric, arrows, arrows.meta}

\usepackage{empheq}

\newcommand{\CC}{C\nolinebreak\hspace{-.05em}\raisebox{.4ex}{\tiny\bf +}\nolinebreak\hspace{-.10em}\raisebox{.4ex}{\tiny\bf +}}
\def\CC{{C\nolinebreak[4]\hspace{-.05em}\raisebox{.4ex}{\tiny\bf ++}}}

\usepackage{hyperref}
\hypersetup{colorlinks, linkcolor=blue, urlcolor=red, filecolor=green, citecolor=blue}

\theoremstyle{plain}
\newtheorem{theorem}{Theorem}[section]
\newtheorem{corollary}[theorem]{Corollary}
\newtheorem{lemma}[theorem]{Lemma}
\newtheorem{proposition}[theorem]{Proposition}

\newtheorem{definition}[theorem]{Definition}
\newtheorem{assumption}[theorem]{Assumption}

\theoremstyle{remark}
\makeatletter
    \def\@endtheorem{\hfill$\diamond$\endtrivlist\@endpefalse }
\makeatother
\newtheorem{remark}[theorem]{Remark}

\newcommand{\N}{\mathbb{N}}
\newcommand{\R}{\mathbb{R}}

\newcommand{\Z}{\mathbb{Z}}
\newcommand{\dd}{\mathop{}\!\mathrm{d}}
\DeclareMathOperator*{\argmin}{arg\,min}
\newcommand\e{\varepsilon}
\newcommand\esssup{\mathop{\operatorname{ess\,sup}}}

\newcommand\dist{\operatorname{dist}}
\newcommand\sym{\operatorname{sym}}
\newcommand{\SO}[1]{\operatorname{SO}(#1)}

\newcommand\loc{{\operatorname{loc}}}

\newcommand\II{{\operatorname{I\!I}}}
\newcommand{\step}[2]{\medskip\noindent\textbf{Step #1~}(#2)\newline\noindent}
\newcommand{\reftoproof}[1]{\noindent(See Section~\ref{#1} for the proof.)\smallskip}
\newcommand{\substep}[1]{\medskip\noindent\textit{Substep #1. }}
\newcommand\deform{v}
\newcommand\minvec{\xi}

\newcommand{\wto}{\rightharpoonup}
\newcommand{\wtto}{\stackrel{2}{\rightharpoonup}}
\newcommand{\stto}{\stackrel{2}{\longrightarrow}}
\newcommand{\meanh}[1]{\langle #1\rangle_{\rm h}}
\DeclarePairedDelimiter{\abs}{\lvert}{\rvert}
\DeclarePairedDelimiter{\norm}{\lVert}{\rVert}
\def\uloc{{\rm uloc}}
\def\HH{{\mathbf H}}
\def\eff{{\rm eff}}
\def\iso{{\rm iso}}

\begin{document}
\begin{center}
  {\LARGE A homogenized bending theory for prestrained plates}
  \medskip

  \normalsize
  Klaus B\"ohnlein\footnote{klaus.boehnlein@tu-dresden.de}\textsuperscript{*}, Stefan Neukamm\footnote{stefan.neukamm@tu-dresden.de}\textsuperscript{*}, David Padilla-Garza\footnote{david.padilla-garza@tu-dresden.de}\textsuperscript{*} and
  Oliver Sander\footnote{oliver.sander@tu-dresden.de}\textsuperscript{*} \par \bigskip

  \textsuperscript{*}Faculty of Mathematics, Technische Universit\"at Dresden\par \bigskip

  \today
\end{center}

\begin{abstract}

  The presence of prestrain can have a tremendous effect on the mechanical behavior of slender structures.
  Prestrained elastic plates show spontaneous bending in equilibrium---a property that makes such objects relevant for the fabrication of active and functional materials.
  In this paper we study microheterogeneous, prestrained plates
  that feature non-flat equilibrium shapes.
  Our goal is to understand the relation between the properties of the prestrained microstructure
  and the global shape of the plate in mechanical equilibrium.
  To this end, we consider a three-dimensional, nonlinear elasticity model that describes a periodic material that occupies a domain with small thickness.
  We consider a spatially periodic prestrain described in the form of a multiplicative decomposition of the deformation gradient.
  By simultaneous homogenization and dimension reduction, we rigorously derive an effective plate model
  as a $\Gamma$-limit for vanishing thickness and period.
  That limit has the form of a nonlinear bending energy with an emergent spontaneous curvature term.

  The homogenized properties of the bending model (bending stiffness and spontaneous curvature) are characterized by corrector problems.
  For a model composite---a prestrained laminate composed of isotropic materials---we investigate the dependence of the homogenized properties on the parameters of the model composite.
  Secondly, we investigate the relation between the parameters of the model composite and the set of shapes with minimal bending  energy.
  Our study reveals a rather complex dependence of these shapes on the composite parameters.
  For instance, the curvature and principal directions of these shapes depend  on the parameters
  in a nonlinear and discontinuous way; for certain parameter regions we observe uniqueness and non-uniqueness of the shapes.
  We also observe size effects: The geometries of the shapes depend on the aspect ratio between the plate thickness and the composite period.
  As a second  application of our theory we study a problem of shape programming: We prove that any target shape (parametrized by a bending deformation) can be obtained (up to a small tolerance) as an energy minimizer of a composite plate,
  which is simple in the sense that the plate consists of only finitely many grains that are filled with a parametrized composite with a single degree of freedom.
\end{abstract}
  \smallskip


  \textbf{MSC-2020:} 74B20 35B27 74Q05

\tableofcontents

\section{Introduction}
\begin{figure}[h]
  \begin{center}
    \includegraphics{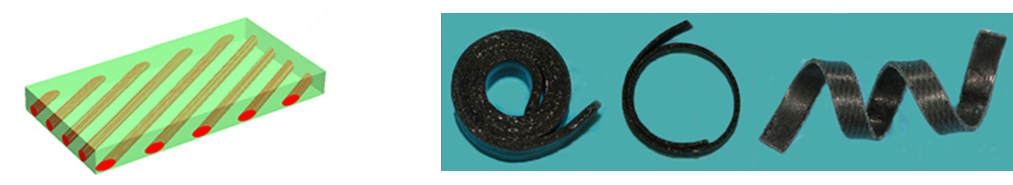}
  \end{center}
  \caption{Design and experimentally observed cylindrical, equilibrium shapes of 3d-printed. composite bilayer plates. Reproduced from \cite{ge2013active}, with the permission of AIP Publishing.}
  \label{fig:expstrip}
\end{figure}

\paragraph{General motivation.}
Natural and synthetic elastic materials often are prestrained. For slender structures, the presence of prestrain may have a huge impact on the mechanical behavior: Plates and films with prestrain often exhibit a complex equilibrium shape due to spontaneous bending, wrinkling and symmetry breaking \cite{klein2007shaping, ware2015voxelated, van2018programming}. 
Prestrain can be the result of different physical mechanisms (e.g.,~swelling or de-swelling of gels \cite{ionov2013biomimetic,tanaka1979kinetics}, thermal expansion~\cite{sigmund1997design,gibiansky1997thermal}, or the nematic-elastic coupling in liquid crystal elastomers \cite{warner2007liquid}), and
can be triggered by different stimuli---a property that makes such materials interesting for the fabrication of functional materials with a controlled shape change; see~\cite{van2018programming} for a review.
New manufacturing techniques such as additive manufacturing even enable the design of microstructured, prestrained materials whose functionality results from a complex interplay between the geometry, the material properties, and the prestrain distribution on a small length scale. An example is a self-assembling cube shown in~\cite{ge2013active}, whose functionality is due to a sandwich-type prestrained composite plate designed with a fibred microstructure, see Figure~\ref{fig:expstrip}.

The development of reliable models and simulation methods that are able to predict the macroscopic behavior
based on the specification of the material on the small scale is an important part of understanding such materials and subject of ongoing research.

\paragraph{Scope and main results of the paper.} In this paper we are interested in the effective elastic behavior of prestrained composite plates. In particular, we seek to understand the \emph{microstructure--shape relation}, i.e., the relation between the mechanical properties and prestrain distribution of the plate on the small length scale, and the emergent macroscopic equilibrium shape.
The starting point of the analysis is the following energy functional of three-dimensional, nonlinear elasticity:
\begin{equation}\label{eq:I}
  \mathcal E^{\varepsilon,h}(\deform)\colonequals\frac1h\int_{\Omega_h} W_\e\big(x,\nabla \deform(x) \,A^{-1}_{\varepsilon,h}(x)\big)\,dx,\qquad \Omega_h\colonequals S\times\Big(-\tfrac h2,\tfrac h2\Big).
\end{equation}
Here, $\Omega_h$ denotes the reference configuration of the three-dimensional plate, $S\subset\R^2$ is the midsurface,
$0<h\ll 1$ denotes the plate thickness, and $\deform:\Omega_h\to\R^3$ is the deformation of the plate.
The stored energy density function $W_{\e}(x,F)$ is assumed to be frame indifferent with a single, non-degenerate energy well at $\SO 3$ for almost all $x \in \Omega_h$. It describes a heterogeneous composite with a microstructure that
oscillates locally periodically in in-plane directions on a length scale $0<\e\ll 1$ (see Section~\ref{sec:2.2} for details).
In addition, \eqref{eq:I} models a microheterogeneous prestrain based on a multiplicative decomposition of the deformation gradient $\nabla \deform$ with the matrix field $A_{\e,h}:\Omega_h\to \R^{3\times 3}$ as in, e.g., \cite[Section~2.1]{bauer2019derivation}.
Like the stored energy function $W_\varepsilon$ itself, we assume that $A_{\e,h}$ oscillates in in-plane directions.

As explained in \cite{van2018programming}, there are two different mechanisms for prestrain-induced shape changes of plates, namely, the ``buckling strategy'' and the ``bending strategy''. In this paper, we are interested in the latter.
As is well-known, bending of plates can be driven by gradients of the prestrain along the thickness of the plate with a magnitude comparable to the thickness. Therefore, we assume that $A_{\varepsilon,h}=I_{3\times 3}+h B_{\e,h}$ with $\|B_{\e,h}\|_{L^\infty(\Omega_h)}\lesssim 1$ uniformly in $\e$ and $h$.

The first mathematical problem that we address is the rigorous derivation of a homogenized,
nonlinear plate model with an effective prestrain as a $\Gamma$-limit of $\mathcal E^{\varepsilon,h}$
when both parameters $\e$ and $h$ tend to $0$ simultaneously (Theorem~\ref{T1}).  In the special case
of a globally periodic microstructure, the derived model is given by a bending energy with an effective prestrain:
\begin{equation}\label{eq:II}
  \mathcal I^\gamma_{\hom}:H^2_{\iso}(S;\R^3)\to\R,\qquad \mathcal I^\gamma_{\hom}(\deform)=\int_S Q^\gamma_{\hom}\big(\II_\deform(x')-B_{\rm eff}^\gamma\big)\,dx',\qquad x'\colonequals(x_1,x_2).
\end{equation}
Here, $H^2_{\iso}(S;\R^3)$ is the nonlinear space of bending deformations, i.e., the set of all $\deform\in H^2(S;\R^3)$ satisfying the isometry constraint $(\nabla'\deform)^\top\nabla'\deform=I_{2\times 2}$ where $\nabla' \colonequals (\partial_1,\partial_2)$.
We denote by $\II_\deform\colonequals\nabla'\deform^\top\nabla'(\partial_1\deform\wedge\partial_2\deform)$ the second fundamental form associated with $\deform$ and note that it captures the curvature of the deformed plate. The effective bending moduli are described by means of a positive-definite quadratic form $Q^\gamma_{\hom}:\R^{2\times 2}_{\sym}\to[0,\infty)$, and the effective prestrain of the two-dimensional  plate is described by a matrix $B^\gamma_{\eff}\in\R^{2\times 2}_{\sym}$. Both $Q^\gamma_{\hom}$ and $B^\gamma_{\eff}$ can be derived from $W_{\e}$ and $A_{\e,h}$ by homogenization formulas that require to solve certain corrector problems.
These corrector problems are the equilibrium equations of linear elasticity posed on a representative volume with periodic boundary conditions, see Proposition~\ref{P:1}.

It turns out that the precise form of the limiting energy $\mathcal I^\gamma_{\hom}$ is sensitive
to the relative scaling of the plate thickness $h$ and the size $\varepsilon$ of the microstructure:
In order to capture this size effect we introduce the parameter $\gamma\in(0,\infty)$ and we shall assume
that $\frac{h}{\e}\to \gamma$ as $(\e,h)\to 0$. As indicated by the notation, $\gamma$ enters
the formulas for the effective quantities $Q^\gamma_{\hom}$ and $B^\gamma_{\eff}$. We remark that our result,
Theorem~\ref{T1}, includes the more general case of a \emph{locally periodic composite}
(defined in Assumption~\ref{ass:W} below), which leads to a $x'$-dependence of $Q^\gamma_{\hom}$ and $B^\gamma_{\eff}$. Furthermore, we discuss displacement boundary conditions in Theorem~\ref{T2}.
\medskip

\begin{figure}
  \includegraphics[width=10cm]{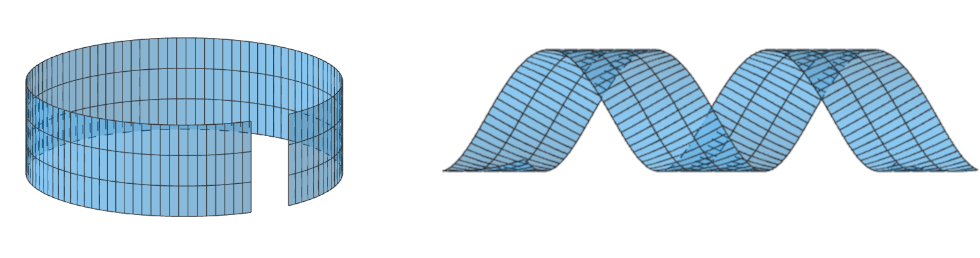}
  \centering
  \caption{Shapes with minimial energy predicted by our theory. A composite plate with prestrained fibres in one of the layers is considered. The design parameters are the volume fraction of the fibres and the angle that they form with the longitudinal axis of the strip. See Section~\ref{S:textures1} for details. }
  \label{fig:cylindrical_strip_teaser}
\end{figure}

The standard theory of $\Gamma$-convergence implies that (almost) minimizers of the scaled global energy $\frac1{h^2}\mathcal E^{\varepsilon, h}$ converge (up to subsequences) to minimizers of the plate energy $\mathcal I^\gamma_{\hom}$. The minimizers of the latter thus capture the effective equilibrium shapes of the three-dimensional plate for $0<\e,h\ll 1$.
In Section~\ref{S:micro-shape} we investigate the minimizers of $\mathcal I^\gamma_{\hom}$ and their dependence on the three-dimensional composite, i.e., on $W_\e$ and $B_{\e,h}$. Understanding this relation can be seen as two steps:
\begin{enumerate}[(a)]
\item (Microstructure--properties relation). In Section~\ref{S:properties} we first investigate the map $(W_\e, B_{\e,h})\mapsto (Q_{\hom}^\gamma, B^\gamma_{\eff})$, which for each composite requires to solve a set of three corrector problems of the type~\eqref{eq:corrector_equation}.
While this can be done numerically (as we plan to do in a forthcoming paper), here we focus mostly
on analytic results: In Lemma~\ref{L:orth1} we prove that isotropic composites with isotropic prestrain and a reflection symmetric geometry lead to orthotropic $Q_{\hom}^\gamma$ and diagonal $B^\gamma_{\eff}$. Moreover, in Lemma~\ref{S:ex:C3} we obtain explicit formulas for $Q_{\hom}^\gamma$ and $B^\gamma_{\eff}$ in the case of a parametrized, isotropic laminate, see Figure~\ref{sketchofmicrostructure} for a schematic visualization.

\item (Properties--shape relation). Once the effective coefficients of $Q_{\hom}^\gamma$ and $B^\gamma_{\eff}$ are known,
minimizing the energy functional \eqref{eq:II} determines the equilibrium shape of the plate.
In the general case, there is no hope for explicit formulas for the minimizers. However, when $Q_{\hom}^\gamma$ and $B^\gamma_{\eff}$ are independent of $x'\in S$, then free minimizers correspond to cylindrical shapes with constant fundamental form. In this case, the minimization of the functional $\mathcal I^\gamma_{\hom}$ simplifies to an algebraic minimization
problem~\eqref{L:ex:1:eq1}, which can be solved without the need for solving a nonlinear partial differential equation.
Lemma~\ref{S:ex:ex1} establishes a trichotomy result for the set of  minimizers in the case when $Q^\gamma_{\hom}$ is orthotropic and yields a way to evaluate minimizers algorithmically.
\end{enumerate}
By combining both steps we recover the desired \emph{microstructure--shape relation}.
In Sections~\ref{S:properties} and \ref{S:micro-shape} we illustrate this procedure on the level of a parametrized, isotropic, two-component laminate.
In Section~\ref{S:properties} we focus on the microstructure--properties relation.
We obtain explicit formulas that relate the parameters of the composite (in particular, the volume fractions of the components, the stiffness contrast, and the strength of the prestrain) to the effective quantities $(Q_{\hom}^\gamma,B^\gamma_{\eff})$.
Based on this, in Section~\ref{S:micro-shape}, we explore the parameter dependence of the geometry of shapes with minimal bending energy; see Figure~\ref{fig:cylindrical_strip_teaser} for an example.

Our detailed parameter study shows that the geometries of shapes with minimal bending energy
depend on the parameters in a rather complex, and partially counterintuitive way. In particular, we observe a discontinuous dependence of the geometry on the parameters: E.g., the bending direction and the sign of the curvature may jump when perturbing the volume fraction of the components of the composite. We also observe a size effect: Qualitative and quantitative properties of the set of shapes with minimal energy depend on the scale ratio $\gamma$. Furthermore, we observe a break of symmetry: By changing the volume fraction we may transition from a situation where the set of shapes with minimal energy are a rotationally symmetric one-parameter family to a situation where the set consists of a unique shape. 
\medskip

In Section~\ref{S:textures} we turn to the problem of shape programming: In Theorem~\ref{T:shape} we prove that, in a nutshell, any shape that can be parametrized by a bending deformation $\deform\in H^2_{\iso}(S;\R^3)$ can be approximated by low energy deformations of a prestrained, three-dimensional composite plate with a simple design. Here, simple means that the plate is partitioned into finitely many grains and each grain is filled by realizations of a prescribed, one-parameter family of a periodic, prestrained composite.
\bigskip

\paragraph{A brief survey of the methods and previous results.} Our $\Gamma$-convergence result, Theorem~\ref{T1}, is concerned with
simultaneous dimension reduction and homogenization in nonlinear elasticity. The asymptotics $h\to 0$ correspond to the derivation of a bending plate energy in the spirit of the seminal papers \cite{friesecke2002theorem, friesecke2006hierarchy}.
Our analysis heavily relies on the geometric rigidity estimate of \cite{friesecke2002theorem} and the
corresponding method to prove compactness for sequences of 3d-deformations with finite bending energy.
On the other hand, the limit $\e\to0$ describes the homogenization of microstructure. For the analysis of the simultaneous $\Gamma$-limit $(\e,h)\to0$ we follow the general strategy introduced by the second author in the case of rods \cite{neukamm2010homogenization,neukamm2012rigorous}. It relies on the fact that to leading order, the nonlinear energy can be written as a convex, quadratic functional of the scaled nonlinear strain. This allows using methods from convex homogenization, in particular, the notion of two-scale convergence \cite{nguetseng1989general, allaire1992homogenization, visintin2007two}, and a representation of effective quantities with the help of corrector problems. Results on simultaneous homogenization and dimension reduction for plates in the von-Kármán and bending regimes have been obtained by the second author together with Vel{\v{c}}i{\'c} and Hornung in \cite{neukamm2013derivation, hornung2014derivation}, see also \cite{neukamm2015homogenization, cherdantsev2015bending, bukal2017simultaneous, hornung2018stochastic} for related works.
\medskip

Extending these ideas, in the present paper we consider materials with a microheterogeneous prestrain, whose magnitude scales with the thickness $h$ of the plate. The derivation of bending theories for prestrained plates (without homogenization) has first been studied by Schmidt \cite{schmidt2007plate, schmidt2007minimal}; for related results in the context of nematic plates see \cite{agostiniani2017shape, agostiniani2017dimension, agostiniani2020rigorous}. In these works (as in our paper) the prestrain is modeled by a multiplicative decomposition $F_\text{el} A^{-1}$ of the deformation gradient $F_\text{el}$ with a factor $A$ that is close to identity, i.e., $A=I+h B$. As a consequence, admissible plate deformations of the limiting model are isometries for the Euclidean metric, and the multiplicative prestrain turns into a linearized, additive one.
While the prestrain in~\cite{schmidt2007plate} is a macroscopic quantity, in the present paper
we consider a microheterogeneous prestrain leading to an effective, homogenized prestrain
in the limit model. A similar extension has been considered by the second author in \cite{bauer2019derivation}
for the case of nonlinear rods. In contrast to previous results on dimension reduction for plates,
our derivation, which invokes homogenization, requires a precise characterization of two-scale limits of the scaled nonlinear strain along sequences with finite bending energy. This is achieved in Proposition~\ref{P:chartwoscale}, which, in particular, affirmatively identifies the two-scale limits of the nonlinear strain in flat regions---a problem that remained open in \cite{hornung2014derivation}.
The proof of Proposition~\ref{P:chartwoscale} is based on a wrinkling ansatz introduced by the third author in \cite{padilla2020dimension}.
\medskip

Minimizers of bending energies for plates with prestrain have been studied first in the spatially homogeneous case with 2d-prestrains that are multiples of the identity matrix~\cite{schmidt2007minimal}.
In particular, \cite[Lemma 3.1]{schmidt2007minimal} contains the convenient observation that in the homogenous case, free minimizers are cylindrical. As usual in homogenization, non-isotropy of effective quantities typically emerges even for composites consisting of isotropic constituents materials.
For the model \eqref{eq:II} this means that the quadratic form $Q^\gamma_{\hom}$ is typically
not isotropic, and that $B^\gamma_{\rm eff}$ is not a multiple of the identity. Therefore, compared to \cite{schmidt2007minimal}, in our case the structure of minimizers is considerably richer.
In particular, with Lemma~\ref{L:char} we classify the sets of minimizers into three families.
In the spatially heterogeneous case or in the case of prescribed displacement boundary conditions, bending deformations with minimal energy are not necessarily cylindrical and explicit formulas are not available.

The numerical minimization of the energy \eqref{eq:II} is highly nontrivial.
Most works in the literature consider only the case without prestrain, i.e., $B_\eff^\gamma = 0$.
The first difficulty is the discretization of the space $H^2_\text{iso}(S;\R^3)$
of bending isometries, i.e., of deformations $\deform$
in $H^2$ satisfying $(\nabla' \deform)^\top \nabla' \deform = I_{2 \times 2}$.  No fully conforming discretizations
seem to exist in the literature.  Nonconforming discretizations based on the MINI- and Crouzeix--Raviart elements
have been proposed in~\cite{bartels2013finite}.  Alternative discretizations using discrete Kirchhoff triangles
or Discontinuous Galerkin (DG) finite elements appear in~\cite{bartels2013DKT,rumpf2021finite}
and~\cite{bonito2021dg}, respectively.
An approach using spline functions (which are in $H^2(S;\R^3)$, but are not pointwise isometries)
appears in~\cite{mohan2021minimal}.
Convergence results are given, e.g., in~\cite{bartels2017bilayer}.

Prestrain is included in a few models, but the attention has been restricted so far
only to the isotropic case, where $Q_{\hom}^\gamma(G)=|G|^2$ and $B_{\eff}^\gamma=\rho I$
for some $\rho\in\R$, see  \cite{bonito2021numerical,bonito2022ldg,bartels2017bilayer}.
In~\cite{GNPG1}, a model with a more general effective prestrain has been considered
in the context of liquid crystal elastomer plates.

The second difficulty is the minimization of the non-convex energy functional~\eqref{eq:II},
which is a challenge even without the prestrain. The works of \citeauthor{bartels2013DKT} use
different numerical gradient flows together with linearizations of the isometry constraint
\cite{bartels2013DKT,bartels2013finite, bartels2017bilayer,bonito2021dg}.
This leads to a (controllable) algebraic violation of the constraints beyond the one
introduced by the discretization.
\citet{rumpf2021finite} use a Lagrange multiplier formulation and a Newton method instead,
and preserve the exact isometry constraints at the grid vertices.
Note that we do not numerically minimize~\eqref{eq:II} in the present manuscript.
Situations requiring such approaches will be the subject of a later paper.
\medskip

Finally, we remark that prestrained plates have also been intensively studied from the perspective
of non-Euclidean elasticity, see \cite{bhattacharya2016plates, lewicka2011scaling, lewicka2020dimension, lewicka2010foppl, lewicka2010shell, lewicka2011matching, lewicka2017plates, lewicka2020quantitative}. In this context,
the reference configuration is assumed to be a Riemannian manifold and the factor $A$ in
the decomposition $F_{\rm el}A^{-1}$ is viewed as the square root of the metric. As observed
in \cite{bhattacharya2016plates} there is an interesting interplay between the critical scaling
of the energy with respect to~the thickness of the plate $h$ and the curvature of the metric $\sqrt{A}$. More specifically,
the minimum energy (per volume) is of order at most $h^{2}$ (as in our case) if and only if there exists an isometric immersion
with finite bending energy of the metric $\sqrt{A}$. Recent works have also considered scaling regimes different from the bending one. For instance, scaling by $h^4$ leads to von-K\'arm\'an plate models, see \cite{de2021hierarchy, de2020energy}. We note that the condition that
the minimum energy scales at most like $h^{4}$ 
is intimately linked to the structure of the Riemann curvature tensor of the metric $\sqrt{A}$. 
Moreover, the membrane scaling has been considered in~\cite{plucinsky2018patterning} and \cite{plucinsky2018actuation} with applications to nematic liquid crystal elastomer plates and nonisometric origami.
\medskip

\paragraph{Organization of the paper.} In Section~\ref{sec:2.2} we introduce the three-dimensional model. In Section~\ref{sec:2:3} we present the limit plate model and state the $\Gamma$-convergence result. Section~\ref{S:correctors} contains the definition of the effective quantities via homogenization formulas and their characterization with the help of correctors. Section~\ref{S:nonlinstrain} establishes a characterization of the two-scale structure of limits of the nonlinear strain and establishes strong two-scale convergence for the nonlinear strain for minimizing sequences. Section~\ref{S:properties} is devoted to the study of the microstructure--properties relation, and Section~\ref{S:micro-shape} to the properties--shape and microstructure-shape relations.
In Section~\ref{S:textures} we present our result on shape programming. All proofs are presented in Section~\ref{S:proofs}. In the appendix, in Appendix~\ref{A:lambda} and \ref{A:twoscale} we discuss various properties of two-scale convergence for grained microstructures---a variant of two-scale convergence that we introduce in this paper. 


\section{Setup of the model and derivation of the prestrained plate model}

We derive the plate model by simultaneous homogenization and dimension reduction of
a three-dimensional model. The proofs for all results stated in this section
are collected in Section~\ref{S:proofs}.

\subsection{The three-dimensional model and assumptions on the material law and microstructure}
\label{sec:2.2}

We denote by $\Omega_h\colonequals S\times (-\frac h2,\frac h2)$ the reference configuration
of a three-dimensional plate with thickness $h>0$, where $S\subset \R^2$ is an open, bounded
and connected Lipschitz domain. We call $\Omega\colonequals S\times(-\frac12,\frac12)$ the corresponding domain of unit thickness. We use the shorthand notation $x=(x',x_3)$ with $x'\colonequals(x_1,x_2)$ and denote the \textit{scaled deformation gradient} by $\nabla_h\colonequals\big(\nabla',\tfrac{1}{h}\partial_3\big)$ where $\nabla'\colonequals(\partial_1,\partial_2)$. Moreover, we write $I_{d\times d}$ to denote the unit matrix in $\R^{d\times d}$.

By rescaling \eqref{eq:I} and specializing to the case
\begin{equation*}
A_{\e,h}^{-1}(x',hx_3) = I_{3\times 3}-hB_{\e,h}(x)
\end{equation*}
(i.e., $A_{\e,h}\approx I+h B_{\e,h}$)
we obtain the energy functional $\mathcal I^{\e,h}:L^2(\Omega;\R^3)\to[0,+\infty]$,
\begin{equation}\label{def:ene}
  \mathcal I^{\e,h}(\deform)\colonequals
  \begin{cases}
    \displaystyle \frac1{h^2}\int_{\Omega} W_\e\Big(x,\nabla_h \deform(x)(I_{3\times 3}-hB_{\e,h}(x))\Big)\dd x&\text{if }\deform\in H^1(\Omega;\R^3),\\
    +\infty&\text{otherwise.}
  \end{cases}
\end{equation}
We shall study the $\Gamma$-limit of $\mathcal I^{\e,h}$ as both small parameters $\e$ and $h$
converge to $0$ simultaneously, and as already mentioned in the introduction, it turns out that the obtained $\Gamma$-limit will depend on the limit of the ratio $h/\e$.
To capture this size effect, we make the following assumption:
\begin{assumption}[Relative scaling of $h$ and $\e$]\label{A:gamma}
  There exists a number $\gamma\in (0,\infty)$ and a monotone function $\e:(0,\infty)\to (0,\infty)$
  such that $\lim_{h\downarrow 0}\e(h)=0$ and $\lim_{h\downarrow 0}\frac{h}{\e(h)}=\gamma$.
\end{assumption}
The case $\gamma\ll 1$ corresponds to the situation of a plate that is thin compared to
the typical size of the microstructure, while $\gamma\gg1$ corresponds to the case of a microstructure
that is very fine not only compared to the macroscopic dimensions of the problem, but also to
the small thickness of the plate. Note that Assumption~\ref{A:gamma} excludes the extreme cases
$\gamma=0$ (i.e., $h \ll\e \ll 1$) and $\gamma=\infty$ (i.e., $\e\ll h\ll 1$).
We comment on these cases in  Remark~\ref{rem:casesgamma}.

Next, we state our assumptions on the material law, the microstructure of the composite, and the microstructure of the prestrain. Following \cite{bauer2019derivation} we describe prestrained elastic materials by combining
\begin{itemize}
\item a geometrically nonlinear, stored energy function  that describes a non-prestrained, elastic material with a \emph{stress-free, nondegenerate} reference state at $\SO 3$,
\end{itemize}
with
\begin{itemize}
\item a \emph{multiplicative decomposition} of the deformation into an elastic part and a prestrain that is
of the order of the plate's thickness $h$ and locally periodic (in in-plane directions) on the scale $\varepsilon$.
\end{itemize}
In \cite{eckart1940thermodynamics,L69,kroner1959allgemeine}, such a multiplicative decomposition
of the deformation has been introduced in the context of finite strain elastoplasticity.

The stored energy functions we consider have to have certain standard properties.
We collect appropriate functions and their linearizations in so-called \emph{material classes}:
\begin{definition}[Nonlinear material class]\label{D1}
  Let $0<\alpha\leq\beta$, $\rho>0$, and let $r:[0,\infty)\to[0,\infty]$ denote a monotone function satisfying $\lim_{\delta\to0}r(\delta)=0$.
\begin{itemize}
\item The class $\mathcal W(\alpha,\beta,\rho,r)$ consists of all measurable functions $W:\R^{3\times 3}\to[0,+\infty]$ that
\begin{enumerate}[label=({W}\arabic*)]
 \item \label{item:nonlinear_material_w1}
  are frame indifferent: $W(RF)=W(F)$ for all $F\in\R^{3\times 3}$, $R\in \SO 3$;
 \item are non-degenerate:
 \begin{alignat*}{2}
  W(F) &\geq \alpha \dist^2(F,\SO 3) & \qquad & \text{for all $F\in\R^{3\times 3}$,}\\
  W(F) &\leq \beta \dist^2(F,\SO 3)  &        & \text{for all $F\in\R^{3\times 3}$ with $\dist^2(F,\SO 3)\leq \rho$;}
 \end{alignat*}
 \item \label{item:nonlinear_material_w3}
  are minimal at $I_{3\times 3}$: $W(I_{3\times 3})=0$;
 \item \label{item:nonlinear_material_w4}
  admit a quadratic expansion at $I_{3\times 3}$:
  For each $W$ there exists a quadratic form $Q:\R^{3\times 3}\to\R$ such that
  \begin{equation*}
   |W(I_{3\times 3}+G)-Q(G)|\leq |G|^2r(|G|)\qquad\text{for all $G\in\R^{3\times 3}$.}
  \end{equation*}
\end{enumerate}
\item The class $\mathcal Q(\alpha,\beta)$ consists of all quadratic forms $Q$ on $\R^{3\times 3}$ such that
  \begin{equation*}
    \alpha|\sym G|^2\leq Q(G)\leq\beta|\sym G|^2
    \qquad
    \text{for all }G\in\R^{3\times 3},
  \end{equation*}
  where $\sym{G} \colonequals \frac{1}{2}(G + G^\top)$ is the symmetric part of $G$.
  We associate with each $Q \in \mathcal{Q}(\alpha,\beta)$ the fourth-order tensor $\mathbb L\in\operatorname{Lin}(\R^{3\times 3},\R^{3\times 3})$ defined by the polarization identity
  \begin{equation}\label{eq:polarization}
    \mathbb L F:G\colonequals\frac12\big(Q(F+G)-Q(F)-Q(G)\big),
  \end{equation}
  where $:$ denotes the standard scalar product in $\R^{3\times 3}$.
\end{itemize}
\end{definition}
Properties~\ref{item:nonlinear_material_w1}--\ref{item:nonlinear_material_w4} are standard assumptions in the context of dimension reduction. In particular,
stored energy functions of class $\mathcal W(\alpha,\beta,\rho,r)$ can be linearized at the identity (see, e.g., \cite{padilla2020dimension, MN11, GN11, neukamm2010homogenization})
and the elastic moduli of the linearized model are given by the quadratic form $Q$
of~\ref{item:nonlinear_material_w4}. Furthermore, for any stored energy function $W$ we have
by~\cite[Lemma~2.7]{neukamm2012rigorous}
\begin{equation*}
  W\in\mathcal W(\alpha,\beta,\rho,r)\qquad\implies\qquad Q\in\mathcal Q(\alpha,\beta),
\end{equation*}
(which motivates the definition of the class $\mathcal{Q}(\alpha,\beta)$),
and thus $Q$ is a bounded and positive definite quadratic form on symmetric matrices in $\R^{3\times 3}$.

\bigskip

For the plate model we consider particular stored energy functions that are elements of $\mathcal{W}$.

\begin{assumption}[Nonlinear material]\label{ass:nonlinear}
  For all $\e>0$ and $\alpha$, $\beta$, $\rho$, $r$ as in Definition~\ref{D1},
  the elastic energy density $W_\e:\Omega\to[0,+\infty]$ of~\eqref{def:ene}
  is a Borel function such that $W_\e(x,\cdot)\in\mathcal W(\alpha,\beta,\rho,r)$ for almost every $x\in\Omega$.
\end{assumption}
Additionally, we shall assume that the microstructure of the composite is locally periodic, that is,
we consider countably many open subsets $S_j$, called ``grains'' that partition $S$ up to a null set,
and assume that on each $S_j\times(-\tfrac12,\tfrac12)\subset\Omega$ the composite features a laterally periodic microstructure, possibly with a different reference lattice in each grain (Figure~\ref{grains}). This leads to the following definition:
\begin{figure}
    \centering
    \includegraphics[scale=0.5]{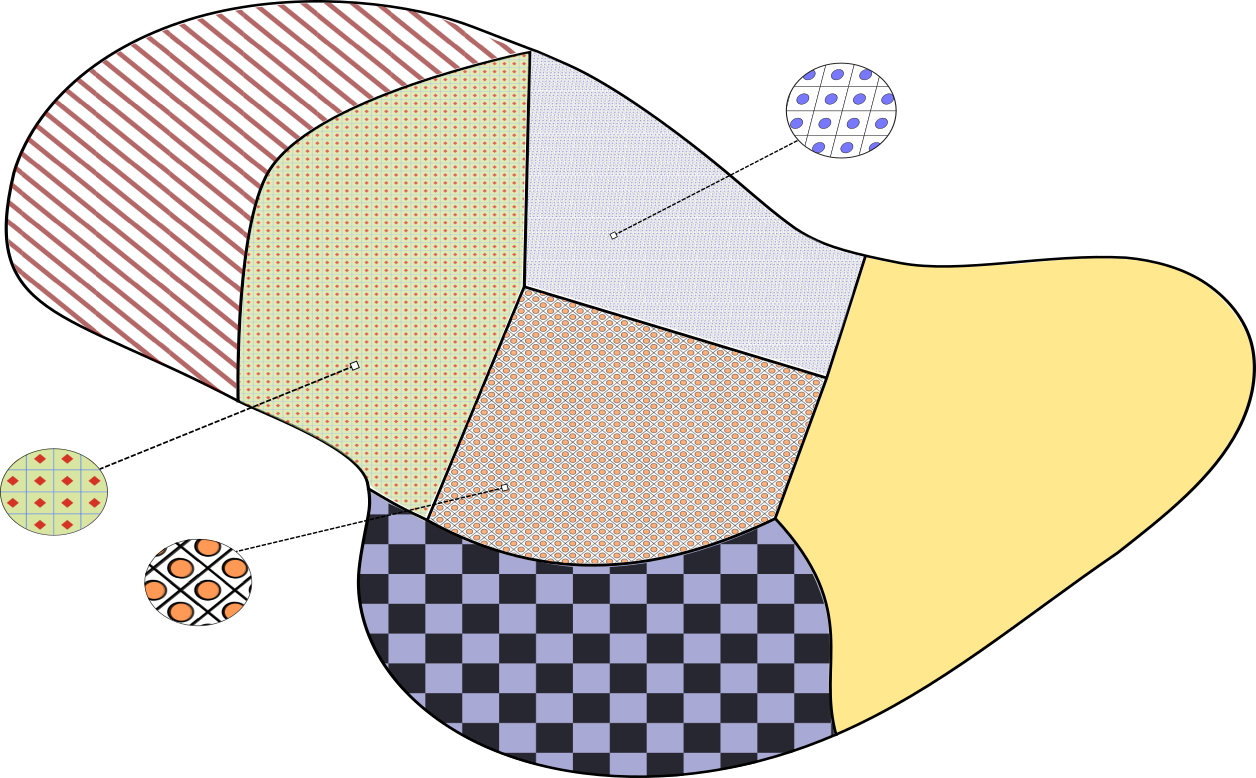}
    \caption{Schematic picture of a microstructure with a different reference lattice and composite in each grain}
    \label{grains}
\end{figure}
\begin{definition}[$\Lambda$-periodicity, grain structure, local periodicity]\label{ass:pcperiodicity}\mbox{}
  \begin{enumerate}[(i)]
  \item \label{ass:pcperiodicity:item1} Let $\Lambda\in\R^{2\times 2}$ be invertible. A measurable function $\varphi:\R^2\to\R$ is called
  \emph{$\Lambda$-periodic} if $\varphi(\cdot+\tau)=\varphi(\cdot)$ for all $\tau\in\Lambda\Z^2$.

  \item A \emph{grain structure} is a finite or countable family  $\{S_j,\Lambda_j\}_{j\in J}$ consisting of open, disjoint subsets $S_j$ of $S$ and matrices $\Lambda_j\in \R^{2\times 2}$ such that $S\setminus(\bigcup_{j\in J}S_j)$ is a null set and
    \begin{equation}\label{ass:pcperidicity:eq}
      \frac1{C_J}\leq \Lambda_j^\top\Lambda_j\leq C_J
      \qquad \text{in the sense of quadratic forms},
    \end{equation}
    for all $j\in J$ with a constant $C_J>0$ independent of $j\in J$.
  \item A measurable function $\varphi:\Omega\times\R^2\to\R$ is called \emph{locally periodic (subordinate to the grain structure} $\{S_j,\Lambda_j\}_{j\in J}$), if
    \begin{equation}\label{eq:locper}
      \forall j\in J\,:\quad \varphi(x,\cdot)\text{ is }\Lambda_j\text{-periodic for a.e.~$x\in S_j\times(-\tfrac12,\tfrac12)$.}
    \end{equation}
  \end{enumerate}
\end{definition}
Note that by \eqref{ass:pcperidicity:eq} the geometry of the local lattices of periodicity is uniformly controlled.
Our core assumption on the composite's microstructure is now
that $W_\e$ describes a composite that is locally periodic on scale $\e>0$. %
Likewise, we shall assume that the prestrain is locally periodic and satisfies a smallness condition.
\newcounter{saveenum}
\begin{assumption}[Local periodicity of the composite and prestrain] \label{ass:W}
  Let $\{S_j,\Lambda_j\}_{j\in J}$ be a grain structure as in Definition~\ref{ass:pcperiodicity}. Let Assumption~\ref{ass:nonlinear} be satisfied and  suppose that there exists $Q:\Omega\times\R^2\times\R^{3\times 3}\to\R$
  with $Q(x,y,\cdot)\in\mathcal Q(\alpha,\beta)$ for a.e.~$x\in\Omega$, $y\in\R^2$, such that the following properties hold:
  \begin{enumerate}[(i)]
  \item \label{ass:W:1}(Local periodicity of $Q$). \label{item:local_periodicity_ii}
    The map $\mathbb L$ associated with $Q$ via \eqref{eq:polarization} is locally periodic subordinate to $\{S_j,\Lambda_j\}_{j\in J}$, and for each grain $S_j$, $j\in J$, the map $S_j\ni x'\mapsto \mathbb L(x',\cdot)\in L^\infty\big((-\tfrac12,\tfrac12)\times\R^2;\R^{3\times 3\times 3\times 3}\big)$ is continuous.
  \item \label{ass:W:2} (Strong two-scale approximation of $Q_\e$). The quadratic term $Q_\e(x,\cdot)$ in the expansion of $W_\e(x,\cdot)$ (cf.~\ref{item:nonlinear_material_w4} of Definition~\ref{D1}) satisfies
    \begin{equation*}
      \limsup_{\e\to0}\esssup_{x\in S_j\times(-\tfrac12,\tfrac12)}\max_{G\in\R^{3\times 3}\atop |G|=1} \abs[\Big]{Q_\e(x,G)-Q(x,\tfrac{x'}\e,G)} = 0
      \qquad \text{for all $j\in J$}.
    \end{equation*}
   \setcounter{saveenum}{\value{enumi}}
  \end{enumerate}
  Furthermore, we suppose that for all $\e,h>0$ the prestrain $B_{\e,h}:\Omega\to\R^{3\times 3}_{\sym}$ of \eqref{def:ene} is measurable and we suppose that there exists a Borel function $B:\Omega\times\R^2\to\R^{3\times 3}_{\sym}$ such that:
  \begin{enumerate}[(i)]
  \setcounter{enumi}{\value{saveenum}}
  \item \label{ass:W:3} (Local periodicity of $B$). The function $B$ is locally periodic subordinate to $\{S_j,\Lambda_j\}_{j\in J}$, and
    \begin{equation*}
      \sum_{j\in J}\int_{S_j}\fint_{\Box_{\Lambda_j}}|B(x',x_3,y)|^2\dd (x_3,y)\dd x'<\infty,
    \end{equation*}
    where $\Box_{\Lambda_j}:=(-\frac12,\frac12)\times Y_{\Lambda_j}$ and $Y_{\Lambda_j}:=\Lambda_j[-\frac12,\frac12)^2$. Moreover, we write $\fint_{\Box_{\Lambda_j}} \colonequals \frac{1}{\abs{\Box_{\Lambda_j}}} \int_{\Box_{\Lambda_j}}$ for the integral mean.

  \item  \label{ass:W:4}  (Strong two-scale approximation and boundedness of $B_{\e(h),h}$).
   We have
    \begin{equation*}
      \limsup_{h\to 0}\int_\Omega\Big|B_{\e(h),h}(x)-B\big(x,\tfrac{x'}{\e(h)}\big)\Big|^2\dd x=0
      \qquad \text{and} \qquad
      \limsup_{h\to 0}\sqrt h\|B_{\e(h),h}\|_{L^\infty(\Omega)}<\infty,
    \end{equation*}
    where $\e(\cdot)$ is as in Assumption~\ref{A:gamma}.
  \end{enumerate}
\end{assumption}
The main reason for considering not only the $\Z^2$-periodic case but the more general
and flexible structure of Assumption~\ref{ass:W} is our application to shape programming presented in Section~\ref{S:textures}.
Note that by Assumption~\ref{ass:W} we do not assume that $W_\e$ is itself locally periodic.
We only suppose that the quadratic term $Q_\e$ in the expansion of $W_\e$ is close to a
locally periodic quadratic form (scaled by $\e$).
\begin{remark}[Examples for locally periodic composites]\label{R:Ex1}\mbox{}
  \begin{enumerate}[(a)]
  \item ($\Z^2$-periodic case). A special case of a locally periodic composite in the sense of Assumption~\ref{ass:W} is the \emph{$\Z^2$-periodic case}, where the partition consists of the single grain $S$ and the lattice of periodicity is everywhere the same and given by $\Lambda=I_{2\times 2}$. More explicitly, in this case, we may consider a Borel function $W:(-\frac12,\frac12)\times\R^2\times\R^{3\times 3}\to[0,\infty]$ such that $W(x_3,y,\cdot)\in\mathcal W(\alpha,\beta,\rho,r)$ for a.e.\ $y\in\R^2$,
  $x_3\in(-\frac12,\frac12)$, and $\R^2\ni y\to W(x_3,y,F)$ is $\Z^2$-periodic for a.e. $x_3\in(-\frac12,\frac12)$ and all $F\in\R^{3\times 3}$. Then the family  of scaled stored energy functions $W_\e(x',x_3,F)\colonequals W(x_3,\tfrac{x'}{\e},F)$ satisfies Assumption~\ref{ass:W}. Note that in this example $W_\e(x',x_3,F)$ is itself $\e\Z^2$-periodic w.r.t.~$x'\in\R^2$.
  \item A simple example of a composite that is locally periodic subordinate to a grain structure $\{S_j,\Lambda_j\}_{j\in J}$ is given by $W_\e(x',x_3,F)\colonequals\sum_{j\in J}1_{S_j}(x')W(x_3,\Lambda_j^{-1}\tfrac{x'}{\e},F)$ where $W$ is as above. Then the family $W_\e$ satisfies Assumption~\ref{ass:W}.
  \end{enumerate}
\end{remark}

\begin{remark}[Examples of prestrains]\mbox{}
\begin{enumerate}[(a)]
\item Under reasonable assumptions, polymer hydrogels, i.e., networks of hydrophilic rubber molecules,
 can be modeled by considering a prestrain of the form $B = \frac{\lambda - 1}{h} I_{3\times 3}$,
 where $\lambda$ is a material parameter depending on the free swelling factor \cite{flory1943statistical}.

\item  Nematic liquid crystal elastomers are
  materials consisting of a polymer network inscribed with liquid crystals.
  These materials feature a coupling between the liquid crystal orientation and elastic properties that can be modeled by considering a prestrain of the form $B = \overline{r}(\frac{1}{3} I_{3\times 3} -  n \otimes n)$, where $n : S \to \R^3$ is a unit vector field describing the local orientation of the liquid crystals, and $\overline{r}$ is a material parameter (see \cite{warner2007liquid}).

\item In finite-strain thermoelasticity (see \cite{vujovsevic2002finite}), prestrains of the form $ B = ( \overline{\beta} - \overline{\alpha})n_{0} \otimes n_{0} - \overline{\beta} I_{3\times 3}$ are considered to describe materials that (after a change of temperature) stretch in a direction  $n_{0}: S \to \mathcal{S}^{2}$ with factor $\overline\alpha$, and contract in directions orthogonal to $n_0$ with factor $\overline\beta$.
\end{enumerate}
\end{remark}

\subsection{The limit model and convergence results}\label{sec:2:3}

We now discuss the model that results from letting $(\e, h) \to 0$ simultaneously.
The limit energy can be written as the sum of two contributions. The first one
is a homogenized bending energy $\mathcal I_{\hom}^\gamma:L^2(S;\R^3)\to[0,\infty]$ that includes
an effective prestrain $B_{\eff}^\gamma\in L^2(S;\R^{2\times 2}_\text{sym})$ that captures the impact of the prestrain $B$ on the macroscale:
\begin{equation}\label{limit:eq1}
  \mathcal I_{\hom}^{\gamma}(\deform)
   \colonequals
  \begin{dcases}
    \int_{S}Q_{\hom}^\gamma(x',\II_{\deform}-B_{\eff}^\gamma(x'))\dd x' & \text{if $\deform\in H^2_{\iso}(S;\R^3)$},\\
    +\infty&\text{else},
  \end{dcases}
\end{equation}
where
\begin{equation*}
  \II_\deform\colonequals\nabla'\deform^\top\nabla' b_{\deform}\ \text{ and }\ b_{\deform}\colonequals(\partial_1\deform\wedge\partial_2\deform)
\end{equation*}
denote the second fundamental form (expressed in local coordinates) of the surface parametrized by $\deform$ and the surface normal, respectively.
Above,  the quadratic form $Q_{\hom}^\gamma$ describes the homogenized elastic moduli of the composite.
The second contribution $\mathcal I_\text{res}^\gamma$
(defined in Definition~\ref{D:eff} below) is a residual energy that is independent of the deformation $u$,
and is quadratic in the prestrain $B$. Both effective quantities---the quadratic form $Q_{\hom}^\gamma$ describing
the homogenized elastic moduli of the composite, and the effective prestrain $B_{\eff}^\gamma$---only depend on the linearized material law $Q$, the prestrain $B$, and the scale ratio $\gamma$.
More specifically,
\begin{itemize}
\item $Q_{\hom}^\gamma(x',\cdot):\R^{2\times 2}_{\sym}\to\R$ is a positive definite quadratic form given by the homogenization formula of Definition~\ref{def:Qgamma} below,
\item $B_{\eff}^\gamma(x')\in\R^{2\times 2}_{\sym}$ is given by the averaging formula of Definition~\ref{D:eff} below,
\item $\mathcal I^\gamma_{\rm res}(B)\geq 0$ is defined in Definition~\ref{D:eff}.
\end{itemize}
Both $Q_{\hom}^\gamma$ and $B_{\eff}^\gamma$ can be evaluated for any $x'\in S$
by solving linear corrector problems. In Section~\ref{S:correctors} we present an algorithm
to evaluate these quantities, and we show numerical experiments in Section~\ref{S:properties}.
\medskip

Our first result establishes $\Gamma$-convergence for $(h,\e)\to 0$:

\begin{theorem}[$\Gamma$-convergence]\label{T1}
  Let Assumptions~\ref{A:gamma} and~\ref{ass:W} be satisfied. Then the following statements hold:
  \begin{enumerate}[(\alph*)]
  \item \label{item:T1:compactness}
   (Compactness). Let $(\deform_h)\subset L^2(\Omega;\R^3)$ be a sequence with equibounded energy, i.e.,
    \begin{equation}\label{eq:equibounded}
      \limsup_{h\to0}\mathcal I^{\e(h),h}(\deform_h)<\infty.
    \end{equation}
    Then
    \begin{equation}\label{eq:FBE}
      \limsup_{h\to0}\frac{1}{h^2}\int_\Omega \dist^2\big(\nabla_h\deform_h(x),\SO 3\big)\dd x<\infty,
    \end{equation}
    and there exists $\deform\in H^2_{\iso}(S;\R^3)$ and a subsequence (not relabeled) such that
    \begin{subequations}\label{T1:conv}
      \begin{alignat}{2}
      \deform_h-\fint_{\Omega}\deform_h\dd x & \to \deform & \qquad & \text{in $L^2(\Omega)$},\\
      \shortintertext{and}
      \nabla_h \deform_h & \to (\nabla'\deform,b_\deform) && \text{in $L^2(\Omega)$}.
    \end{alignat}
    \end{subequations}
    Here and below, $b_\deform\colonequals\partial_1\deform\wedge\partial_2\deform$ is the surface normal vector,
    and $\wedge$ denotes the vector product.
  \item \label{item:T1:lower_bound} (Lower bound). If $(\deform_h)\subset L^2(\Omega;\R^3)$ is a sequence with $\deform_h-\fint_{\Omega}\deform_h\dd x\to \deform$ in $L^2(\Omega)$, then
    \begin{equation*}
      \liminf_{h\to0} \mathcal I^{\e(h),h}(\deform_h)\geq \mathcal I_{\hom}^\gamma(\deform)+\mathcal I_{\rm res}^\gamma(B).
    \end{equation*}
  \item \label{item:T1:recovery_sequence} (Recovery sequence). For any $\deform\in H^2_{\iso}(S;\R^3)$ there exists a sequence $(\deform_h)\subset W^{1,\infty}(\Omega;\R^3)$ with $\deform_h\to \deform$ strongly in $H^1(\Omega;\R^3)$ such that
    \begin{equation}\label{T1:c}
      \lim_{h\to0}\mathcal I^{\e(h),h}(\deform_h)=\mathcal I_{\hom}^\gamma(\deform)+\mathcal I_{\rm res}^\gamma(B),
    \end{equation}
    and, in addition,
    \begin{equation}\label{T1:equnif}
      \limsup_{h\to 0}h^\frac12 \norm[\big]{ \nabla_h \deform_h-(\nabla'\deform,b_\deform) }_{L^\infty(\Omega)}=0.
    \end{equation}
  \end{enumerate}
\end{theorem}%
\noindent
(See Sections~\ref{SS:compactness},  \ref{SS:lowerbound},
and \ref{SS:upperbound} for the proof of (a),(b), and (c), respectively.)
\smallskip

\begin{remark}[Identification of functions defined on $S$ with their canonical extension to $\Omega$]\label{R:extension}
  In the paper we tacitly identify functions defined on $S$, say $\deform:S\to\R^3$, with their  canonical extension to $\Omega$, namely, $\Omega\ni(x',x_3)\mapsto \deform(x')$. This clarifies the meaning of statements such as \eqref{T1:conv}.
\end{remark}

The main points of Theorem~\ref{T1} are the parts~\ref{item:T1:lower_bound} and~\ref{item:T1:recovery_sequence}.
The implication \eqref{eq:FBE} $\implies$ \eqref{T1:conv}, which is the main point of
Part~\ref{item:T1:compactness}, has already been proven by~\citet{friesecke2002theorem}
using their celebrated geometric rigidity estimate.
\begin{remark}[The cases $\gamma = \infty$ and $\gamma = 0$]\label{rem:casesgamma}
  In Theorem~\ref{T1} we assume that $\gamma\in(0,\infty)$, which means that $h$ and $\e(h)$ converge to $0$ with the same order. The result changes for $\gamma=\infty$ and $\gamma=0$, respectively. We note that Theorem~\ref{T1} can be extended  to the case $\gamma = \infty$ (i.e.~when $\e\ll h$) based on the methods developed in this paper and \cite{hornung2014derivation}. In contrast, the case $\gamma = 0$, which corresponds to $h\ll \e$,  is more subtle. Even in the case without prestrain (i.e., $B = 0$) it is not fully understood, and the resulting $\Gamma$-limit can be qualitatively different depending on the relative scaling between $h^2$ and $\varepsilon$:
  The cases $h^{2} \ll \varepsilon \ll h$ and $\varepsilon \ll h^{2}$ are treated in \cite{cherdantsev2015bending} and \cite{velvcic2015derivation}, respectively.
  The sequential limit $\e\to 0$ after $h\to 0$ is discussed in \cite{neukamm2015homogenization}.
\end{remark}
%

%

We also establish a variant of the theorem for plates with displacement boundary conditions on straight parts of the boundary.
For the precise formulation of the resulting two-dimensional model we introduce a set $\mathcal A_{BC}$
of bending deformations with the appropriate boundary conditions. Note that the following definition introduces additional assumptions on the geometry of~$S$:
\begin{definition}[Two-dimensional displacement boundary conditions]\label{ass:BC}
  Let $S$ be a convex and bounded Lipschitz domain with piece-wise $C^1$-boundary. We consider $k_{BC}\in\N$ relatively open, non-empty lines segments $\mathcal L_i\subset\partial S$ and a reference isometry $\deform_{BC}\in H^2_{\iso}(S;\R^3)$ such that for $i=1,\ldots,k_{BC}$  (the trace of) $\nabla' \deform_{BC}$ is constant on $\mathcal L_{i}\cap\partial S$. We introduce a space of bending deformations that satisfy the following displacement boundary conditions on the line-segments:
  \begin{equation*}
    \mathcal A_{BC}\colonequals\Big\{\deform\in H^2_{\iso}(S;\R^3)\,:\,\deform=\deform_{BC}\text{ and }\nabla' \deform=\nabla'\deform_{BC}\text{ on }\mathcal L_1\cup\ldots\cup\mathcal L_{k_{BC}}\,\Big\}.
  \end{equation*}
\end{definition}
We shall see that the boundary conditions of Definition~\ref{ass:BC} emerge from sequences of 3d-deformations $(\deform_h)$ with finite bending energy that satisfy boundary conditions of the form
\begin{equation}\label{eq:BC3d}
  \deform_h=(1-h\delta)\deform_{BC}+hx_3b_{\deform_{BC}}
  \qquad\text{on }(\mathcal L_1\cup\ldots\cup\mathcal L_{k_{BC}})\times \Big(\!\!-\frac12,\frac12 \Big),
\end{equation}
where $\delta\in\R$ denotes a fixed parameter.
More precisely, we recall the following
extension of the compactness part~\ref{item:T1:compactness} of Theorem~\ref{T1},
proved in~\cite[Theorem~2.9 (a)]{GNPG1}:
\begin{lemma}[Emergence of two-dimensional boundary conditions]\label{Lbc}
  Let $S$ and $\mathcal A_{BC}$ be as in Definition~\ref{ass:BC}. Consider a sequence $(\deform_h)\subset H^1(\Omega;\R^3)$ satisfying \eqref{eq:FBE} and \eqref{eq:BC3d} for some fixed $\delta\in\R$. Then there exists a subsequence (not relabeled) and $\deform\in \mathcal A_{BC}$ such that $\deform_h\to \deform$ in $L^2(\Omega)$ and $\nabla_h\deform_h\to(\nabla'\deform,b_\deform)$ in $L^2(\Omega)$.
\end{lemma}
We finally show that we can construct recovery sequences that feature the boundary conditions~\eqref{eq:BC3d}:
\begin{theorem}[Recovery sequences subject to displacement boundary conditions]\label{T2}
  Consider the setting of Theorem~\ref{T1} and assume boundary conditions as in Definition~\ref{ass:BC}.
  If additionally
  \begin{equation}\label{ass:unifB}
    C_B
    \colonequals
    \sup_{j\in J} \; \esssup_{x'\in S_j}\fint_{\Box_{\Lambda_j}}|B(x',x_3,y)|^2\dd (x_3,y)
    <
    \infty,
  \end{equation}
  then there exists $\delta>0$ (only depending on $C_B$) such that for any $\deform\in H^2_{\iso}(S;\R^3)$
  there exists a sequence $(\deform_h)\subset H^1(\Omega;\R^3)$ with $\deform_h\to \deform$ strongly in $H^1(\Omega;\R^3)$ satisfying \eqref{T1:c}, \eqref{T1:equnif}, and the displacement boundary conditions \eqref{eq:BC3d}.
\end{theorem}
\reftoproof{SS:upperbound}

The argument for Theorem~\ref{T2} is based on recent results obtained by us together with \citet{GNPG1}. In particular, there we prove for $S$ and $\mathcal A_{BC}$ as in Definition~\ref{ass:BC} the approximation property
  \begin{equation}\label{ass:BC:dense}
    \mathcal A_{BC}\cap C^\infty(\bar S;\R^3)\text{ is dense in }(\mathcal A_{BC},\|\cdot\|_{H^2(S;\R^3)}),
  \end{equation}
  see \cite[Proposition~2.11]{GNPG1}.
  It extends previous results by Pakzad~\cite{pakzad2004sobolev} and Hornung~\cite{hornung2011approximation} on the approximation of isometries by smooth isometries to the case of affine boundary conditions on line segments.
\begin{remark}[Convergence of almost minimizers]\label{R:convmin1}
  Theorem~\ref{T1} together with Theorem~\ref{T2} and Lemma~\ref{Lbc} implies that the sequence of functionals
  \begin{equation*}
    \mathcal I_h:L^2(\Omega;\R^3)\to[0,\infty],\qquad \mathcal I_h(\deform)\colonequals
    \begin{cases}
      \mathcal I^{\e(h),h}(\deform)&\text{if $\deform$ satisfies \eqref{eq:BC3d}},\\
      +\infty      &\text{else},
    \end{cases}
  \end{equation*}
  $\Gamma$-converges in $L^2(\Omega)$ for $h\to 0$ to the functional
  \begin{equation*}
    \mathcal I:L^2(\Omega;\R^3)\to[0,\infty],\qquad \mathcal I(\deform)\colonequals
    \begin{cases}
      \mathcal I^\gamma_{\hom}(\deform)+\mathcal I^\gamma_{\rm res}(B)&\text{if }\deform\in\mathcal A_{BC},\\
      +\infty      &\text{else.}
    \end{cases}
  \end{equation*}

  Thanks to the boundary condition and strict convexity of \eqref{limit:eq1} viewed
  as a function of $\II_\deform$,
  the limit functional admits a minimizer in the set $\mathcal A_{BC}$. Furthermore, by the imposed boundary conditions, the sequence $(\mathcal I_h)$ is equicoercive in $L^2(\Omega;\R^3)$, and thus standard arguments from the theory of $\Gamma$-convergence imply that any sequence of almost minimizers $(\deform_h)$, i.e.,
  functions for which
  \begin{equation*}
    \mathcal I_h(\deform_h)\leq \inf\mathcal I_h+h,
  \end{equation*}
  converges (modulo a sub-sequence, not relabeled) to $\deform_{*},$ a minimizer of $\mathcal I^{\gamma}_{\hom}$ in the set $\mathcal{A}_{BC}$. In view of the compactness part~\ref{item:T1:compactness} of Theorem~\ref{T1} we even get
  \begin{equation*}
    \nabla_h\deform_h\to (\nabla'\deform_*,b_{\deform_*})\qquad\text{strongly in }L^2(\Omega).
  \end{equation*}
  Furthermore, in Section~\ref{S:strongtwoscale} below we shall also see that the associated sequence of nonlinear strains strongly two-scale converges to a limit that is completely determined by $\deform_*$, cf.~Remark~\ref{R:strongconvergenceforalmostminimizers} in connection with Proposition~\ref{strongconvergenceforalmostminimizers}.
\end{remark}

\subsection{The homogenization formula and corrector problems}\label{S:correctors}
This section presents the definitions of the homogenized quadratic form $Q_{\hom}^\gamma$, the effective prestrain $B_{\eff}^\gamma$, and the residual energy $\mathcal I^\gamma_{\rm res}$.
The definition of $Q_{\hom}^\gamma(x',\cdot)$
and $B^\gamma_{\eff}(x')$ for a fixed material point $x'\in S$,
only depends on the quadratic energy density
\begin{equation*}
  \Big(-\frac12,\frac12\Big)\times\R^2\times\R^{3\times 3}\ni(x_3,y,G)\mapsto Q(x',x_3,y,G)
\end{equation*}
and the prestrain tensor
\begin{equation*}
  \Big(-\frac12,\frac12\Big)\times\R^2\ni(x_3,y)\mapsto B(x',x_3,y)
\end{equation*}
of Assumption~\ref{ass:W}---both only at the position $x'$.
In view of this, we first present a \emph{local} definition of the effective quantities for $x'\in S$ fixed, and discuss their continuity properties.
Secondly, we present the global definition of the effective quantities and of the residual energy $\mathcal I^\gamma_{\rm res}$, see Definition~\ref{def:effquant} below.

\paragraph{Local definition of the effective quantities.}
For the local definition it is useful to introduce the following terminology:
\begin{definition}[$\Lambda$-periodic, admissible quadratic form and prestrain]\label{D:admissible}
  Let $\Lambda\in\R^{2\times 2}$ be invertible and $0<\alpha\leq \beta<\infty$.
  \begin{enumerate}[(i)]
  \item
    A Borel function $Q:(-\frac12,\frac12)\times\R^2\times\R^{3\times 3}\to\R$ is called a
    \emph{$\Lambda$-periodic, admissible quadratic form}, if $Q(x_3,y,\cdot)$ is a quadratic form
    of class $\mathcal Q(\alpha,\beta)$ for a.e.\ $(x_3,y) \in (-\frac{1}{2},\frac{1}{2}) \times \R^2$,
    and $Q(x_3,\cdot,G)$ is $\Lambda$-periodic for all $G\in\R^{3\times 3}$ and a.e.\ $x_3$.

   \item A Borel function $B:(-\frac12,\frac12)\times\R^2\to\R^{3\times 3}_{\sym}$ is called a \emph{$\Lambda$-periodic, admissible prestrain}, if $B$ is square integrable and $B(x_3,\cdot)$ is $\Lambda$-periodic for a.e.\ $x_3$.
  \end{enumerate}
\end{definition}
Note that in this definition, the variable $y$ lives in the space $\R^2$ where the periodic
microstructure is defined, and not in the macroscopic space $S \subset \R^2$.
In the following we frequently use the notation $\Box_{\Lambda}=(-\tfrac12,\tfrac12)\times Y_\Lambda$ and $Y_\Lambda=\Lambda[-\frac12,\frac12)^2$ that we introduced in Assumption~\ref{ass:W}~\ref{ass:W:3}.
By $\Lambda$-periodicity, $Q$ and $B$ can be restricted to $\Box_\Lambda$ without loss of information.
\begin{remark}[Reference cell of periodicity and representative volume element]\label{rem:rve}
  If we consider $Q:\Omega\times\R^2\times\R^{3\times 3}\to\R$ and $B:\Omega\times\R^2\to\R^{3\times 3}_{\sym}$ from Assumption~\ref{ass:W}, then the functions $Q(x',\cdot)$ and $B(x',\cdot)$ are $\Lambda$-periodic and admissible in the sense of Definition~\ref{D:admissible}. (More precisely, we have to choose $\Lambda=\Lambda_j$ if $x_j\in S_j$ for some $j\in J$).
  The set $Y_\Lambda$ is the reference cell of periodicity, see Appendix~\ref{A:lambda}.
  From the perspective of homogenization, the set $\Box_{\Lambda}$
  can be viewed as a three-dimensional \emph{representative volume element} that is attached to a (macroscopic) position $x'$ in the midsurface $S$---the restriction of $Q(x',\cdot,G)$ and $B(x',\cdot)$ to $\Box_\Lambda$ represents the microstructure of the composite at position $x'$.
\end{remark}

We start by presenting a variational definition for the effective quantities $Q^\gamma_{\hom}$
and $B^\gamma_{\eff}$ of Theorem~\ref{T1} associated to a general $\Lambda$-periodic and admissible
quadratic form $Q$ and prestrain $B$.
After that, we shall establish a representation of these quantities based on correctors.
This yields a convenient computational scheme for evaluating $Q_{\hom}^\gamma$ and $B_{\eff}^\gamma$,
see Proposition~\ref{P:1} below.
We then show that the effective quantities $Q^\gamma_{\hom}$ and $B^\gamma_{\eff}$ continuously depend on $Q$ and $B$.
\smallskip

The definition of the homogenized quadratic form involves the function space
$H^1_{\gamma}(\Box_\Lambda;\R^3)$, which consists of all local $H^1$-functions $\varphi:(-\tfrac12,\tfrac12)\times\R^2\to\R^3$
that are $\Lambda$-periodic in the second argument (the $y$-variable); see~\eqref{def:Hgammaper} for the precise definition.
%
\begin{definition}[Homogenized quadratic form]\label{def:Qgamma}
  Given a $\Lambda$-periodic, admissible quadratic form~$Q$, we define the associated homogenized quadratic form $Q^\gamma_{\hom}:\R^{2\times 2}_{\sym}\to[0,\infty)$ as
  \begin{equation}\label{eq:def:min}
    Q^\gamma_{\hom}(G)\colonequals\inf_{M,\varphi}\fint_{\Box_\Lambda}Q\big(x_3,y,\iota(x_3G+M)+\sym\nabla_\gamma\varphi\big)\dd (x_3,y),
  \end{equation}
  where the infimum is taken over all $M\in\R^{2\times 2}_{\sym}$ and $\varphi\in H^1_\gamma(\Box_\Lambda;\R^3)$, and where $\nabla_\gamma\colonequals(\nabla_y,\frac1\gamma\partial_{3})$. Above, we denote by $\iota(G)$ the unique $3\times 3$-matrix whose upper-left $2\times 2$-block is equal to $G\in\R^{2\times 2}$ and whose third column and row are zero.

\end{definition}
\begin{remark}[Special cases]
 For $\Lambda=I_{2\times 2}$ we recover the known case of a $\Z^2$-periodic composite.
 In that situation, $Q^\gamma_{\hom}(G)$ coincides with the formula derived in \cite{hornung2014derivation}.
 In the spatially homogeneous case where $Q$ is independent of $x_3$ and $y$, we recover the formula
 of~\cite{friesecke2002theorem},
  \begin{equation*}
    Q^\gamma_{\hom}(G)=\frac{1}{12}\min_{d\in\R^3}Q\big(\iota(G)+d\otimes e_3\big).
  \end{equation*}
\end{remark}
Next, we turn to the definition of the effective prestrain $B^\gamma_{\eff}$.
We adapt the scheme in \cite{bauer2019derivation}, which is based on orthogonal projections in the Hilbert space of functions $G:(-\frac12,\frac12)\times\R^2\to\R^{3\times 3}_{\sym}$, $(x_3,y)\mapsto G(x_3,y)$ that are $\Lambda$-periodic in $y$ and square integrable on $\Box_\Lambda$. Specifically, let $L^2(\Lambda;\R^{3\times 3}_{\rm sym})$ be the space of $\Lambda$-periodic functions
in $L^2_{\loc}(\R^2;\R^{3\times 3}_{\rm sym})$, see Appendix~\ref{A:twoscale}, and consider the Hilbert space
\begin{equation*}
  \HH_{\Lambda}
  \colonequals
  L^2\big((-\tfrac12,\tfrac12);L^2(\Lambda;\R^{3\times 3}_{\rm sym})\big),
  \qquad
  \Big(G,G'\Big)_{\Lambda}\colonequals\fint_{\Box_\Lambda}\mathbb L G:G'\dd (x_3,y).
\end{equation*}
Note that the induced norm satisfies
\begin{equation}
\label{def:innprod}
  \|G\|^2_{\Lambda}=\fint_{\Box_\Lambda}Q(x_3,y,G)\dd (x_3,y).
\end{equation}
We further consider the subspaces
\begin{align*}
 \HH^\gamma_\text{rel,$\Lambda$}
 & \colonequals
 \Big\{\iota(M)+\sym\nabla_\gamma\varphi\,:\,M\in\R^{2\times 2}_{\sym},\,\varphi\in H^1_{\gamma}(\Box_\Lambda;\R^3)\Big\}\\
 \shortintertext{and}
 \HH^\gamma_{\Lambda}
 & \colonequals
 \Big\{\iota(x_3 G)+\chi\,:\,G\in\R^{2\times 2}_{\rm sym},\,\chi\in\HH^\gamma_{{\rm rel},\Lambda}\Big\}
\end{align*}
and note that they
are closed subspaces of $\HH_{\Lambda}$. Here and below, we understand $\iota(x_3G)$ as the map
in $\HH_\Lambda$ defined by $(x_3,y)\mapsto \iota(x_3G)$.
The closedness of the subspaces can be seen by using Korn's inequality in form of Lemma~\ref{L:korn-lambda} in combination
with Poincar\'e's inequality. We  denote by $\HH^{\gamma,\perp}_{{\rm rel},\Lambda}$ the orthogonal
complement of $\HH^\gamma_{{\rm rel},\Lambda}$ in $\HH^\gamma_{\Lambda}$, and write $P^{\gamma,\perp}_{{\rm rel},\Lambda}$
for the orthogonal projection from $\HH_{\Lambda}$ onto $\HH^{\gamma,\perp}_{{\rm rel},\Lambda}$.
Similarly, we write $P^{\gamma}_{\Lambda}$ for the orthogonal projection from $\HH_\Lambda$
onto $\HH^\gamma_{\Lambda}$.

\begin{remark}[Interpretation of $\HH^\gamma_\text{rel, $\Lambda$}$ and $\HH^\gamma_{\Lambda}$]\label{R:HH}
  Note that we have the orthogonal decomposition
  \begin{equation*}
    \HH_\Lambda = \HH^\gamma_\Lambda\oplus \big(\HH^\gamma_\Lambda\big)^\perp=\HH^\gamma_{{\rm rel},\Lambda}\oplus \HH^{\gamma,\perp}_{{\rm rel},\Lambda}\oplus \big(\HH^\gamma_\Lambda\big)^\perp
  \end{equation*}
  The spaces $\HH^\gamma_\text{rel, $\Lambda$}$ and $\HH^\gamma_{\Lambda}$ naturally show up in the two-scale analysis of the non-linear strain: Indeed, as we shall prove below in Proposition~\ref{P:chartwoscale}, when considering a sequence of deformations $(\deform_h)$ in $H^1(S;\R^3)$ with finite bending energy and limit $\deform\in H^2_{\iso}(S;\R^3)$, the associated sequence of nonlinear strains $E_h(\deform_h) \colonequals \frac{\sqrt{\nabla_h\deform_h^\top\nabla_h\deform_h}-I_{3\times 3}}{h}$ weakly two-scale converges (up to a subsequence) to a limit $E$, and for a.e.~$x'\in S$, the limiting strain takes the form $E(x',\cdot)=\iota(x_3\II_\deform(x'))+\chi$ for some $\chi\in\HH^\gamma_{\text{rel,$\Lambda$}}$. The field $\chi$ can be interpreted as a corrector that captures the oscillations on the scale $\e(h)$ that emerge along the selected subsequence of  $\big(E_h(y_h)\big)$. Note that in the definition of the effective quadratic form $Q^\gamma_{\hom}$ we relax the local energy by infinimizing over all $\iota(M)+\sym\nabla_\gamma\varphi=\chi\in\HH^{\gamma}_{\text{rel,$\Lambda$}}$, see~\eqref{eq:def:min}.
  With help of the scalar product $(\cdot,\cdot)_\Lambda$, the infimization can be rephrased as an orthogonal projection: $Q^\gamma_{\hom}(G)=\|P^{\gamma,\perp}_{{\rm rel},\Lambda}\iota(x_3 G)\|^2_\Lambda$.
\end{remark}
Next, we turn to the definition of the effective prestrain $B^\gamma_{\eff}$. We first note that $B$ can be decomposed into two parts:
$P^\gamma_\Lambda B$ (the orthogonal projection of $B$ onto $\HH^\gamma_\Lambda$),
and the orthogonal complement, $(I-P^\gamma_\Lambda)B$.
The energy contribution of the latter is captured by the residual energy introduced below.
As we shall see, only $P^\gamma_\Lambda B$  interacts with the deformation.
We define $B^\gamma_{\eff}$ in such a way that we can express the energy contribution associated with $P^\gamma_\Lambda B$ in the form $Q_{\hom}^\gamma(B^\gamma_{\eff})$.
For this purpose, in the following lemma, we introduce an operator $\mathbf E_{\Lambda}^\gamma$.
%
\begin{lemma}\label{L:E}
  Let $\gamma\in(0,\infty)$, and let $Q$ be $\Lambda$-periodic and admissible in the sense of Definition~\ref{D:admissible}. Then the map
  \begin{equation*}
    \mathbf E_{\Lambda}^\gamma:\R^{2\times 2}_{\sym}\to \HH^{\gamma,\perp}_{{\rm rel},\Lambda},\qquad \mathbf E_{\Lambda}^\gamma(G)\colonequals P^{\gamma,\perp}_{{\rm rel},\Lambda}\big(\iota(x_3G)\big)
  \end{equation*}
  is a linear isomorphism, and
  \begin{equation*}
    \sqrt{\tfrac\alpha{12}}|G|\leq \bigg(\fint_{\Box_\Lambda}|\mathbf E^\gamma_{\Lambda}(G)|^2\bigg)^\frac12\leq \sqrt{\tfrac\beta{12}}|G|,
  \end{equation*}
  for all $G\in \R^{2\times 2}_{\sym}$.
\end{lemma}
\reftoproof{SS:homogenization}

\begin{remark}
  In view of the definition of $\mathbf E_{\Lambda}^\gamma$ and of $Q^\gamma_{\hom}$ in~Definition~\ref{def:Qgamma}, we have $Q^\gamma_{\hom}(G)=\fint_{\Box_\Lambda}Q(x_3,y,\mathbf E_{\Lambda}^\gamma(G))\dd (x_3,y)$ for all $G\in\R^{2\times 2}_{\sym}$.
  Furthermore, the proof of Lemma~\ref{L:E} reveals that
  \begin{equation}\label{eq:bounds}
    \frac{\alpha}{12}|G|^2\leq Q^\gamma_{\hom}(G)\leq \frac{\beta}{12}|G|^2\qquad\text{for all }G\in\R^{2\times 2}_{\sym}.
  \end{equation}
\end{remark}
\begin{definition}[Effective prestrain]\label{D:eff}
    Let $\gamma\in(0,\infty)$, and let $Q,B$ be $\Lambda$-periodic and admissible in the sense of Definition~\ref{D:admissible}. We define the effective prestrain $B_{\eff}^\gamma\in\R^{2\times 2}_{\sym}$ associated with $Q$ and $B$ by
  \begin{equation}\label{eq:def:beff}
    B_{\eff}^\gamma\colonequals(\mathbf E_{\Lambda}^\gamma)^{-1}\Big(P^{\gamma,\perp}_{{\rm rel},\Lambda}(\sym B)\Big).
  \end{equation}
\end{definition}
Next, we represent $Q^\gamma_{\hom}$ and $B^\gamma_{\eff}$ by means of three corrector problems.
These are linear Korn-elliptic partial differential equations with domain $\Box_\Lambda$ subject to
periodic boundary conditions in the $y$-variable.
The weak formulation~\eqref{eq:corrector_equation} (below) of these corrector problems can be
phrased as a variational problem in the Hilbert space $H^1_{\gamma}(\Box_\Lambda;\R^3)$
(defined in \eqref{def:Hgammaper}),
and the coefficients are given by $\mathbb L$, the symmetric fourth-order tensor
obtained in~\eqref{eq:polarization} from $Q$ via polarization.
\begin{lemma}[Existence of a corrector]\label{L:corrector}
  Let $Q$ and $B$ be $\Lambda$-periodic and admissible in the sense of Definition~\ref{D:admissible},
  and assume that $\frac{1}{C_\Lambda}\leq \Lambda^\top\Lambda\leq C_\Lambda$. Let $G\in \R^{2\times 2}_{\sym}$. Then there exists a unique pair
  \begin{equation*}
    M_G\in\R^{2\times 2}_{\sym},\qquad \varphi_G\in H^1_{\gamma}(\Box_\Lambda;\R^3)\text{ with }\int_{\Box_\Lambda}\varphi_G\dd (x_3,y) =0,
  \end{equation*}
  solving the corrector problem
  \begin{equation}\label{eq:corrector_equation}
    \int_{\Box_\Lambda}\mathbb L\big(\iota(x_3 G+M_G)+\sym(\nabla_\gamma\varphi_G)\big):\big(\iota(M')+\sym(\nabla_\gamma\varphi')\big)\dd (x_3,y)=0
  \end{equation}
  for all $\varphi'\in H^1_{\gamma}(\Box_\Lambda;\R^3)$ and $M'\in\R^{2\times 2}_{\sym}$.  Moreover, there exists a constant $C=C(\alpha,\beta,\gamma,C_\Lambda)$ such that
  \begin{equation}\label{eq:corrector_apriori}
    |M_G|^2+\fint_{\Box_\Lambda}|\nabla_\gamma\varphi_G|^2\dd (x_3,y) \leq C|G|^2.
  \end{equation}
  We call $(M_G,\varphi_G)$ the corrector associated with $G$.
\end{lemma}
\reftoproof{SS:homogenization}
\begin{remark}\label{R:corrector}
  We note that \eqref{eq:corrector_equation} is the Euler-Lagrange equation of the minimization problem in the definition of $Q^\gamma_{\hom}(G)$ in \eqref{eq:def:min}. In particular, $(M_G,\varphi_G)$ is the unique minimizer of \eqref{eq:def:min} (modulo an additive constant for $\varphi_G$).
\end{remark}
\begin{proposition}[Representation via correctors]\label{P:1}
 Let $Q$ and $B$ be $\Lambda$-periodic and admissible in the sense of Definition~\ref{D:admissible},
 and assume that $\frac{1}{C_\Lambda}\leq \Lambda^\top\Lambda\leq C_\Lambda$. Let
 $G_1$, $G_2$, $G_3$ be an orthonormal basis of $\R^{2\times 2}_{\sym}$ and denote by $(M_{G_i},\varphi_{G_i})$ the corrector associated with $G_i$ in the sense of Lemma~\ref{L:corrector}.
  \begin{enumerate}[(a)]
  \item \label{P:1:b} (Representation of $Q^\gamma_\textnormal{hom}$).
    The matrix  $\widehat Q\in\R^{3\times 3}$ defined by
    \begin{equation*}
      \widehat Q_{ik}\colonequals\fint_{\Box_\Lambda}\mathbb L\big(\iota(x_3 G_i+M_{G_i})+\sym(\nabla_\gamma\varphi_{G_i})\big):\iota(x_3 G_k)\dd (x_3,y)
    \end{equation*}
    is symmetric and positive definite, and we have
    \begin{equation}\label{P:1:coercivityQhat}
      \frac{\alpha}{12} I_{3\times 3}\leq \widehat Q\leq \frac{\beta}{12} I_{3\times 3},
    \end{equation}
    in the sense of quadratic forms. Moreover, for all $G\in\R^{2\times 2}_{\rm sym}$ we have the representation
    \begin{equation*}
      Q_{\rm hom}^\gamma(G)=\sum_{i,j=1}^3\widehat Q_{ij}\widehat G_i\widehat G_j,
    \end{equation*}
    where $\widehat{G}_1$, $\widehat{G}_2$, $\widehat{G}_3$ are the coefficients of $G$
    with respect to the basis $G_1$, $G_2$, $G_3$.

  \item  \label{P:1:c} (Representation of $B^\gamma_\textnormal{eff}$).
    Define $\widehat B\in\R^3$ by
    \begin{equation*}
      \widehat B_i
      \colonequals
      \fint_{\Box_\Lambda}\mathbb L\big(\iota(x_3 G_i+M_{G_i})+\sym(\nabla_\gamma\varphi_{G_i})\big):B\dd (x_3,y),\qquad i=1,2,3.
    \end{equation*}
    Then we have $B_{\eff}^\gamma=\sum_{i=1}^3\big(\widehat Q^{-1}\widehat B\big)_iG_i$.
  \end{enumerate}
\end{proposition} 
\reftoproof{SS:homogenization}
The following lemma shows that the correctors and the effective quantities associated with an admissible  pair $(Q,B)$ continuously depends on $(Q,B)$:
\begin{lemma}[Continuity]\label{L:cont}
  Consider a sequence of $\Lambda$-periodic and admissible pairs $(Q_n,B_n)$, $n\in\N\cup\{\infty\}$, and assume that for $n\to\infty$,
  \begin{alignat*}{2}
    Q_n(x_3,y,G) & \; \to \; Q_\infty(x_3,y,G) & \qquad & \text{for all $G\in\R^{3\times 3}$ and a.e.\ $(x_3,y)$},\\
    B_n          & \; \to \; B_\infty          &        & \text{strongly in $L^2(\Box_\Lambda)$}.
  \end{alignat*}
  Denote by $Q^{\gamma}_{\hom,n}$, $B^{\gamma}_{\hom,n}$, and $(M_{n,i},\varphi_{n,i})$ the effective quantities and correctors associated with $(Q_n,B_n)$ in the sense of Proposition~\ref{P:1}. Then for $n\to\infty$,
  \begin{alignat}{2}
    Q^{\gamma}_{\hom,n}(G) & \; \to \; Q^{\gamma}_{\hom,\infty}(G) & \qquad & \text{for all $G\in\R^{2\times 2}_{\sym}$},\\
    B^{\gamma}_{\eff,n}    & \; \to \; B^{\gamma}_{\eff,\infty} & & \text{in $\R^{2\times 2}_{\sym}$},\\\label{pf:L:cont:eq1}
    M_{n,i}                & \; \to \; M_{\infty,i} & & \text{in $\R^{2\times 2}_{\sym}$},\\\label{pf:L:cont:eq2}
    \varphi_{n,i}          & \; \to \; \varphi_{\infty,i}          &        & \text{strongly in $H^1_{\gamma}(\Box_\Lambda;\R^3)$}.
  \end{alignat}
\end{lemma}
\reftoproof{SS:homogenization}

\paragraph{Global definition of the effective quantities.}
We present the \emph{global} definition of the effective quantities and introduce the residual energy associated with a \emph{locally periodic} composite:
\begin{definition}[Homogenized coefficients, effective prestrain and residual energy]\label{def:effquant}
  For $(Q,B)$ as in Assumption~\ref{ass:W} define the \emph{homogenized quadratic form} $Q^\gamma_{\hom}:S\times\R^{2\times 2}_{\sym}\to[0,\infty)$ and \emph{the effective prestrain} $B^\gamma_{\eff}:S\to\R^{2\times 2}_{\sym}$ as follows: For  all $x'\in S_j$, $j\in J$, we define $Q^\gamma_{\hom}(x',G)$ and $B^\gamma_{\eff}(x')$ by \eqref{eq:def:min} and \eqref{eq:def:beff} applied with $\Lambda=\Lambda_j$, $Q=Q(x',\cdot)$, and $B=B(x',\cdot)$. Furthermore, define the \emph{residual energy} as
  \begin{equation}\label{D:Ires}
    \mathcal I_{\rm res}^\gamma(B)\colonequals\sum_{j\in J}\int_{S_j}\fint_{\Box_{\Lambda_j}}Q\Big(x',x_3,y,(I-P^\gamma_{\Lambda_j})(\sym B(x',\cdot))\Big)\dd (x_3,y)\dd x'.
  \end{equation}
\end{definition}
\begin{remark}[Regularity in $x'$]
  Assumption~\ref{ass:W}~\ref{item:local_periodicity_ii} in connection with Lemma~\ref{L:cont}
  implies that $Q^{\gamma}_{\hom}(\cdot,G)$ is continuous in each of the grains $S_j$, $j\in J$, and that $B^\gamma_{\eff}$ is measurable.
  In particular, we deduce that the integral in $\mathcal I_{\hom}^\gamma$ (which requires measurability of the integrand) is well-defined.
\end{remark}

Let us anticipate that in the $\Gamma$-convergence proof, we shall first obtain as a $\Gamma$-limit the ``abstract'' functional $\widetilde{\mathcal I}^\gamma:H^2_{\iso}(S;\R^3)\to\R$,
  \begin{equation}
  \label{def:Igamma}
    \widetilde{\mathcal I}^{\gamma}(\deform)
    \colonequals
    \min_{M,\varphi}\sum_{j\in J}\int_{S_j}\fint_{\Box_{\Lambda_j}}Q\Big(x',x_3,y,\iota(x_3\II_{v}+M)+\sym\nabla_\gamma\varphi-B\Big)\dd (x_3,y)\dd x',
  \end{equation}
  where the minimization is over all corrector pairs $(M,\varphi)$ with $M\in L^2(S;\R^{2\times 2}_{\sym})$
  and $\varphi\in L^2\big(S;H^1_{\gamma,\uloc}\big)$ that are locally periodic in the sense of \eqref{eq:locper}.
  The reason why the relaxation w.r.t.~$M$ and $\varphi$ occurs can be explained by means of the two-scale structure of the limiting strain, which we analyze in Proposition~\ref{P:chartwoscale} in the next section.
  The following lemma, whose proof is based on the projection scheme introduced above, shows that the abstract $\Gamma$-limit decomposes as claimed in Theorem~\ref{T1}:
%
  \begin{lemma}[Representation of the abstract $\Gamma$-limit]\label{L:3.3}
    Let Assumption~\ref{ass:W} be satisfied and let $\widetilde{\mathcal I}^\gamma:H^2_{\iso}(S;\R^3)\to\R$ be defined by \eqref{def:Igamma}.
Then
 \begin{equation*}
   \widetilde{\mathcal I}^{\gamma}(\deform) = \mathcal I^\gamma_{\hom}(\deform)+\mathcal I^\gamma_{\rm res}(B),
 \end{equation*}
 with $\mathcal I^\gamma_{\hom}$ and $\mathcal I^\gamma_{\rm res}$ defined
 in~\eqref{limit:eq1} and~\eqref{D:Ires}, respectively.
\end{lemma}
\reftoproof{SS:homogenization}

\section{Two-scale limits of nonlinear strain}\label{S:nonlinstrain}
As in previous works on simultaneous homogenization and dimension reduction
\cite{neukamm2010homogenization,neukamm2012rigorous,hornung2014derivation},
it is crucial to have a precise understanding of the oscillatory behavior of the nonlinear strain
\begin{equation}\label{D:nonlinestrain}
  E_h(\deform_h)\colonequals\frac{\sqrt{(\nabla_h\deform_h)^\top\nabla_h \deform_h}-I_{3\times 3}}{h}
\end{equation}
for sequences of deformations $(\deform_h)\subseteq H^1(\Omega;\R^3)$ with finite bending energy in the sense of~\eqref{eq:FBE}.
In this section, we give a precise characterization of weak two-scale limits of the
nonlinear strain $E_h$ along sequences of deformations with finite bending energy (Proposition~\ref{P:chartwoscale}).
Furthermore, we prove that $E_h$ strongly two-scale converges along sequences of deformations with converging energy (Proposition~\ref{strongconvergenceforalmostminimizers}).

\subsection{Characterization of sequences with finite bending energy}
In~\cite{friesecke2002theorem} it is shown that any sequence $(\deform_h)$ satisfying \eqref{eq:FBE}
admits a subsequence (also called $(\deform_h)$) such that there is a bending deformation
$\deform\in H^2_{\iso}(S;\R^3)$ for which
\begin{alignat}{2}\label{eq:uhconv}
  \deform_h-\fint_{\Omega}\deform_h & \to \deform  & \qquad & \text{strongly in $L^2(\Omega)$},\\
\shortintertext{and vector fields $M\in L^2(S;\R^{2\times 2}_{\sym})$ and $d\in L^2(\Omega;\R^3)$ such that}
  E_h(\deform_h) & \wto \iota(x_3\II_\deform+M)+\sym(d\otimes e_3) && \text{weakly in $L^2(\Omega)$}.\label{eq:weaknonlinear}
\end{alignat}

%
When adding homogenization to the game, the identification of the weak limit of the
nonlinear strain
is not sufficient; we need to resolve the two-scale structure of \eqref{eq:weaknonlinear}.
For this we appeal to the following variant of two-scale convergence. It is a rather straightforward extension of the notion introduced in \cite{neukamm2012rigorous, neukamm2010homogenization} to the locally periodic setting that we consider here.
In the following definition, $L^2_{\uloc}(\R^2)$ denotes the space of uniform locally $p$-integrable functions, see~\eqref{D:uloc}.
Likewise, $L^2(\Lambda_j)$ denotes the space $\Lambda_j$-periodic functions in $L^2_{\loc}(\R^2)$, see Appendix~\ref{A:lambda}. Also recall the notation $Y_{\Lambda_j}=\Lambda_j[-\frac12,\frac12)^2$ and $\Box_{\Lambda_j}=(-\frac12,\frac12)\times Y_{\Lambda_j}$.
\begin{definition}[Two-scale convergence for locally periodic functions]\label{D:twoscale}
  Let $\{S_j,\Lambda_j\}_{j\in J}$ be a grain structure in the sense of Definition~\ref{ass:pcperiodicity}, and suppose that $h\mapsto\e(h)$ satisfies Assumption~\ref{A:gamma}. We say that a sequence $(\varphi_h)\subset L^2(\Omega)$ weakly two-scale converges in $L^2$  as $h\to0$ to a function $\varphi\in L^2(\Omega;L^2_{\uloc}(\R^2))$ if $(\varphi_h)$ is bounded in $L^2(\Omega)$, and
  \begin{enumerate}[(i)]
  \item  $\varphi$ is locally periodic in the sense of \eqref{eq:locper}, and
  \item for all $j\in J$ and $\psi\in C_c^\infty\big(S_j\times(-\tfrac 12,\tfrac 12);C(\Lambda_j)\big)$,
    \begin{equation*}
      \lim_{h\to0}\int_\Omega \varphi_h(x)\psi\big(x,\tfrac{x'}{\e(h)}\big)\dd x
      =
      \int_{S_j}\fint_{\Box_{\Lambda_j}}\varphi(x',x_3,y)\psi(x',x_3,y)\dd(x_3,y)\dd x'.
    \end{equation*}
    Here the space of two-scale test-functions  $C_c^\infty\big(S_j\times(-\tfrac 12,\tfrac 12);C(\Lambda_j)\big)$  consists of all smooth $C(\Lambda_j)$--valued functions with support compactly contained in $S_j\times(-\frac12,\frac12)$ 
    where $C(\Lambda_j)$ denotes the space of continuous, $\Lambda_j$-periodic functions, see Appendix~\ref{A:lambda}.
  \end{enumerate}
  We say that $(\varphi_h)$ strongly two-scale converges to $\varphi$ if additionally
  \begin{equation*}
    \int_\Omega|\varphi_h|^2\dd x
    \;\to\;
    \sum_{j\in J}\int_{S_j}\fint_{\Box_{\Lambda_j}}|\varphi(x',x_3,y)|^2\dd (x_3,y)\dd x'.
  \end{equation*}
  We write $\varphi_h\stackrel{2}{\wto} \varphi$ and $\varphi_h\stackrel{2}{\longrightarrow} \varphi$
  in $L^2$ for weak and strong two-scale convergence in $L^2$, respectively.
\end{definition}
Appendix~\ref{A:twoscale} lists some properties of this notion of two-scale convergence. In particular, Proposition~\ref{P:twoscale:gradrecov} in Appendix~\ref{A:twoscale} shows that two-scale limits of bounded sequences
of scaled gradients $\nabla_h \deform_h$ can be written in the form $(\nabla'\deform(x'),0)+\nabla_\gamma\varphi(x',x_3,y)$ where $\nabla_\gamma=(\partial_{y_1},\partial_{y_2},\tfrac1\gamma\partial_3)$ with a
\emph{macroscopic} function $\deform\in H^1(S;\R^3)$ and a correction $\varphi$ that is a locally periodic function
in $L^2(S; H^1_{\gamma,\uloc})$; see \eqref{def:Hgamma} for the definition of the space $H^1_{\gamma,\uloc}$ (it contains all
 functions defined on $(-\frac12,\frac12)\times\R^2$ with values in $\R^3$ whose $H^1$-norm on cubes $(0,z)+(-\frac{1}{2},\frac{1}{2})^3$ can be bounded uniformly in $z\in\R^2$).
This and further two-scale convergence methods are used to establish the next result, which identifies the structure of two-scale limits of the nonlinear strain.
\begin{proposition}[Characterization of the two-scale limiting strain]\label{P:chartwoscale}
  Let $\{S_j,\Lambda_j\}_{j\in J}$ be a grain structure in the sense of Definition~\ref{ass:pcperiodicity}, and suppose that $\e(h)$ is as in Assumption~\ref{A:gamma} for some $\gamma\in(0,\infty)$.
  \begin{enumerate}[(a)]
  \item \label{item:char_twoscale_a} Let $(\deform_h)\subseteq H^1(\Omega;\R^3)$ be a sequence of finite bending energy, satisfying~\eqref{eq:FBE},
  with limit $\deform\in H^2_{\iso}(S;\R^3)$ in the sense of \eqref{eq:uhconv}. Then, up to a subsequence,
    \begin{equation}\label{eq:twoscaleE}
      E_h(\deform_h)\,\wtto\, \iota(x_3\II_\deform+M)+\sym\nabla_\gamma\varphi,
    \end{equation}
    for a matrix field $M\in L^2(S;\R^{2\times 2}_{\rm sym})$ and a corrector $\varphi\in L^2(S;H^1_{\gamma,\uloc})$ that is locally periodic in the sense of \eqref{eq:locper}.

  \item \label{item:char_twoscale_b}
  For all $\deform\in H^2_{\iso}(S;\R^3)$, $M\in L^2(S;\R^{2\times 2}_{\sym})$, and any corrector $\varphi\in L^2(S;H^1_{\gamma,\uloc})$ satisfying~\eqref{eq:locper}, there exists a sequence $(\deform_h)$ in $H^1(\Omega;\R^3)$ such that
  \begin{alignat}{2}
   \nonumber
   \deform_h
   & \to
   \deform &\qquad& \text{strongly in $L^2(\Omega;\R^3)$}\\
   \label{eq:recov:twoscale}
   E_h(\deform_h)
   & \stto
   \iota(x_3 \II_\deform+M)+\sym\nabla_\gamma\varphi && \text{strongly two-scale in $L^2$}.
  \end{alignat}
  Furthermore,
  \begin{equation*}
   \limsup\limits_{h\to 0}h\|E_h(\deform_h)\|_{L^{\infty}} = 0
   \qquad\text{and}\qquad
   \lim_{h\to 0}\|\det(\nabla_h\deform_h)-1\|_{L^{\infty}}=0.
  \end{equation*}

    Finally, if $\deform$ satisfies boundary conditions in the sense that $v\in\mathcal A_{BC}$ (see Definition~\ref{ass:BC}), and if
    \begin{equation*}
      M(x')+\delta I_{2\times 2}\geq 0
      \qquad
      \text{in the sense of quadratic forms},
    \end{equation*}
    for some $\delta>0$  and a.e.\ $x'\in S$,
    then there exists a sequence $(\deform_h)$ satisfying \eqref{eq:recov:twoscale} and the boundary condition \eqref{eq:BC3d}.
  \end{enumerate}
\end{proposition}
\reftoproof{SS:chracterization}

Proposition~\ref{P:chartwoscale} is the key ingredient for determining
the $\Gamma$-limit of $\mathcal I^{\e(h),h}$, cf.~\eqref{def:ene}:
In the proof of Theorem~\ref{T1}, with help of Proposition~\ref{P:chartwoscale}
we shall first establish $\Gamma$-convergence of $\mathcal I^{\e(h),h}$ to the functional $\widetilde{\mathcal I}^\gamma$, whose definition (see  \eqref{def:Igamma}) is closely related to \eqref{eq:twoscaleE}.
Indeed, the argument of the quadratic form under the integral in \eqref{def:Igamma} is precisely the difference of the right-hand side of \eqref{eq:twoscaleE} and the prestrain $B$. Furthermore, $\widetilde{\mathcal I}^\gamma(\deform)$ is then obtained by minimizing out  $M$ and $\varphi$, i.e., the only terms of the right-hand side of \eqref{eq:twoscaleE} that are not uniquely determined by the deformation $\deform$.
\smallskip

Proposition~\ref{P:chartwoscale} extends a previous result
in \cite{hornung2014derivation} in various directions: First, the construction of Part~\ref{item:char_twoscale_b}
takes boundary conditions of the form of Definition~\ref{ass:BC} into account.
Secondly, ``grained'' composites, i.e., composites with a periodic microstructure
whose reference lattice changes from grain to grain, are considered. Thirdly and most importantly,
Proposition~\ref{P:chartwoscale} closes a gap in the characterization of the limiting strain on regions where the second fundamental form vanishes.
More precisely,
 \cite{hornung2014derivation} proves Part~\ref{item:char_twoscale_a} in the single grain case $J=1$,
 $\Lambda_1=I_{2\times 2}$, and establishes Part~\ref{item:char_twoscale_b} only under the additional assumption that
the matrix $M$ in \eqref{eq:twoscaleE} vanishes on all flat parts of the bending deformation $u$, i.e.,
\begin{equation}\label{eq:flatcondition}
  M=0\text{ a.e. in }\big\{x'\in S\,:\,\II_\deform(x')=0\big\}.
\end{equation}
This assumption is critical for the construction in \cite{hornung2014derivation}.
Nevertheless, in the case without prestrain, the construction of \cite{hornung2014derivation}
is sufficient for deriving the $\Gamma$-limit.
In contrast, when a prestrain is present, it is necessary to treat the case where \eqref{eq:flatcondition} is not satisfied
to retrieve the $\Gamma$-limit.
In the case without homogenization, this has recently been achieved by the third author in \cite{padilla2020dimension}
using a flexibility result for isometric immersions that was inspired by \cite{lewicka2017convex}.
In our proof of Proposition~\ref{P:chartwoscale}\ref{item:char_twoscale_b}, we combine
this method with the two-scale ansatz of \cite{hornung2014derivation} and thus give a complete characterization
of the two-scale limits of the nonlinear strain.

\subsection{Strong two-scale convergence of the nonlinear-strain for almost-minimizing sequences}\label{S:strongtwoscale}

Theorem~\ref{T1} implies that a sequence of almost-minimizers $(\deform_h)$ converges (after possibly extracting a subsequence) to a minimizer of the limiting energy. The next result shows that the nonlinear strain $E_h(\deform_h)$ defined in~\eqref{D:nonlinestrain}
strongly two-scale converges along such a sequence to a two-scale strain that is uniquely determined
by the second fundamental form of the limiting deformation. More precisely, if the limiting deformation
is given by $\deform_*\in H^2_{\iso}(S;\R^3)$, then the limiting strain takes the form
\begin{equation}\label{eq:twoscalestrain}
  E_*(x,y) = \iota\big(x_3\II_{\deform_*}(x')+M_*(x')\big)+\sym\nabla_\gamma\varphi_*(x,y),
\end{equation}
where
\begin{align*}
  M_* \in L^2(S;\R^{2\times 2}_{\sym})\text{ and }
  \varphi_* \in L^2\big(S;H^1_{\gamma, \uloc}\big)
\end{align*}
satisfy, for all $j\in J$ and a.e.~$x'\in S_j$, the weak corrector problem
\begin{multline}\label{eq:weakcor}
  \int_{\Box_{\Lambda_j}}\mathbb L(x',\cdot)\big(\iota(x_3 \II_{\deform_*}(x')+M_*(x'))+\sym(\nabla_\gamma\varphi_*(x',\cdot))\big):\big(\iota(M')+\sym(\nabla_\gamma\varphi')\big)\dd (x_3,y)=0, \\
  \text{for all $M'\in\R^{2\times 2}_{\sym}$ and $\varphi'\in H^1_{\gamma}(\Box_{\Lambda_j};\R^3)$}
\end{multline}
subject to the mixed boundary and zero-mean conditions
\begin{equation*}
  \varphi_*(x',\cdot)\in H^1_\gamma(\Box_{\Lambda_j};\R^3)\quad\text{ and }\quad\fint_{\Box_{\Lambda_j}}\varphi_*(x',x_3,y)\dd (x_3,y)=0.
\end{equation*}
\begin{remark}[Corrector representation of $E_*$]
  By linearity of equation~\eqref{eq:weakcor}, for each $x'\in S$ we have the corrector representation
  \begin{equation*}
    \big(M_*(x'),\varphi_*(x',\cdot)\big)
    =
     \sum_{i=1}^{3} \widehat{\II_{\deform_*}(x')}_{i} \left(M_{i},\varphi_{i} \right),
  \end{equation*}
  where the $(M_i,\varphi_{i})$, $i=1,2,3$ are the correctors of Proposition~\ref{P:1} (applied with $Q=Q(x',\cdot)$
  and $B=B(x',\cdot)$), and where $\widehat{\II_{\deform_*}(x')}\in\R^3$ denotes the coefficient vector of $\II_{\deform_*}(x')\in\R^{2\times 2}_{\sym}$
  with respect to the basis $G_1,G_2,G_3$ of Proposition~\ref{P:1}.
\end{remark}
Next, we establish strong two-scale convergence of the nonlinear strain for sequences of deformations whose energy is converging. In particular, it applies to sequences of (almost) minimizers.
\begin{proposition}\label{strongconvergenceforalmostminimizers}
  Let Assumptions~\ref{A:gamma} and~\ref{ass:W} be satisfied. Consider a sequence $(\deform_h)$ in $H^1(\Omega; \mathbb{R}^{3})$ that strongly converges in $L^2(\Omega)$ to some $\deform_*\in H^{2}_{\iso}(S;\R^3)$. Assume that the energy converges in the sense that
  \begin{equation}\label{eq:strongconvergenceforalmostminimizers:1}
    \lim_{h \to 0} \mathcal I^{\e(h),h}(\deform_{h}) = \mathcal I^\gamma_{\hom}(\deform_*)+\mathcal I^\gamma_{\rm res}(B).
  \end{equation}
  Then
  \begin{equation*}
    E_h(\deform_h)\stto E_*\qquad\text{strongly two-scale in }L^2,
  \end{equation*}
  where $E_*$ is defined in~\eqref{eq:twoscalestrain}.
\end{proposition}

Proposition~\ref{strongconvergenceforalmostminimizers}, proved in Section~\ref{SS:strongconvergence},
extends \cite[Theorem 7.1]{friesecke2002theorem}, where strong ``single-scale'' convergence of $E_h(\deform_h)$
was established in the case without homogenization. Also related is Theorem~7.5.1 of~\cite{neukamm2010homogenization},
which shows the result for two-scale homogenization and rods.
\begin{remark}[Application to almost-minimizers]\label{R:strongconvergenceforalmostminimizers}
 Proposition~\ref{strongconvergenceforalmostminimizers} applies in particular to almost-minimizers:
 When $(\deform_h)$ is a sequence of almost-minimizers, that is, a sequence such that
\begin{equation*}
\lim_{h \to 0} \left( \mathcal I^{\e(h),h}(\deform_{h}) - \inf_{\deform\in H^1(\Omega;\R^3)} \mathcal I^{\e(h),h}(\deform) \right) =0,
\end{equation*}
then Theorem~\ref{T1} implies that a subsequence of $(\deform_h)$ converges strongly in $H^{1}$
to a minimizer of $\mathcal I^{\gamma}(\cdot)+\mathcal I^\gamma_{\rm res}(B)$. Hence,
Proposition~\ref{strongconvergenceforalmostminimizers} implies that (along the same subsequence)
$E_h(\deform_h)\stto E_*$ strongly two-scale in $L^2$.
\end{remark}

\section{The microstructure--properties relations}
\label{S:properties}

In this and the following section we investigate the relation between the microstructure of the composite material
and the effective shape of the plate in equilibrium, i.e., we want to predict the geometry of a minimizer $\deform_*\in H^2_{\rm iso}(S;\R^3)$ of $\mathcal I^\gamma_{\hom}$
based on knowledge of the microstructure of the composite and of the prestrain. We focus here on ``free'' minimizers, i.e., we do not take boundary conditions or external loads into account. The relationship can then be split into two parts:
\begin{enumerate}[(Q1)]
\item \label{item:relations_q1}
  How do the effective stiffness $Q_{\rm hom}^\gamma$ and prestrain $B_{\rm eff}^\gamma$
  depend on the microstructure of the composite and its prestrain?
\item \label{item:relations_q2}
  How do the minimizers of $\mathcal I^\gamma_{\hom}$ depend on $Q_{\rm hom}^\gamma$ and $B_{\rm eff}^\gamma$?
\end{enumerate}
Note that in the general case, both problems can only be studied by numerical simulations:
Question~\ref{item:relations_q1} requires solving a system of linear PDEs that admits
closed-form solutions only in special cases. Question~\ref{item:relations_q2} is even more delicate
since it generally requires solving a nonlinear, singular PDE. Solving PDEs, however, can be mostly avoided when studying
the homogeneous case, i.e., when $Q_{\hom}^\gamma$ and $B_{\rm eff}^\gamma$
are independent of $x'\in S$.  In this case it is known that free minimizers have a constant second fundamental form;
see~\cite{schmidt2007minimal} and Lemma~\ref{L:char:cylindrical} below.
The problem then reduces to an algebraic minimization problem. This is the case that the present paper
focuses on. The numerical
investigation of Question~\ref{item:relations_q2} in the case of spatial heterogeneity will be the subject of a further paper.
\smallskip

In this section we investigate Question~\ref{item:relations_q1}, and thus study
the connection between the microstructure (the linearized elastic energy $Q$, and the prestrain $B$) and the effective properties of the thin plate ($Q^{\gamma}_{\rm hom}$ and $B^{\gamma}_{\rm eff}$).
In Section~\ref{S:ex:0}, we first  introduce a class of examples in which some symmetries of $Q^{\gamma}_{\rm hom}$ and $B^{\gamma}_{\rm eff}$ can be deduced from the symmetries of $Q$ and $B$. Later, we  further specialize that class
to one for which explicit formulas for $Q^\gamma_\text{hom}$ and $B^\gamma_\text{eff}$ are available.
Section~\ref{S:ex:2} then numerically explores these formulas.
The discussion of Question~\ref{item:relations_q2} shall follow in Section~\ref{S:micro-shape}.

\subsection{The case of orthotropic effective stiffness}
\label{S:ex:0}

  In this section we introduce a special class of composites that feature simplified formulas for the effective quantities. We shall discuss a series of examples that become progressively more specific. Our final example consists of a parametrized laminate that is composed of two isotropic materials with zero-Poisson ratio. The upshot of the example is that it features a closed-form expression for $Q^\gamma_{\hom}$ and $B^\gamma_{\eff}$ that can be evaluate at low computational cost --- an advantage that we shall exploit in our numerical studies.

  In the following, we use the representation of a matrix $G\in \R^{2\times 2}_{\sym}$ as
  \begin{equation}
    G=\sum_{i=1}^3\widehat G_iG_i,\quad\text{where }    G_1\colonequals e_1\otimes e_1,
    \,
    G_2\colonequals e_2\otimes e_2,
    \,
    G_3\colonequals \frac1{\sqrt 2}(e_1\otimes e_2+e_2\otimes e_1),
    \label{canonicalBasis}
  \end{equation}
  denotes the canonical orthonormal basis of $\R^{2\times 2}_{\sym}$. We call $\widehat G_1,\widehat G_2,\widehat G_3\in\R$ the coefficients of $G$. We start by introducing the following notion of orthotropocity for $Q^\gamma_{\rm hom}$, and note that all examples in this section shall feature this material symmetry.

\begin{definition}[Orthotropicity]\label{def:orthotropicity}
 We call a quadratic form $Q^{\gamma}_{\rm hom}: \R^{2\times 2}_{\sym} \to \R$ \emph{orthotropic}, if there exist $q_1,q_2,q_{12},q_3 \in \mathbb{R}$ such that
  \begin{equation*}
    Q^{\gamma}_{\hom}(G)=\big(\widehat G_1^2q_1+\widehat G_1\widehat G_2q_{12}+\widehat G_2^2q_2\big)+\widehat G_3^2q_3.
  \end{equation*}
  We call $q_1,q_2,q_{12},q_3$ the coefficients of $Q^\gamma_{\hom}$.
\end{definition}
In the rest of this section we consider composites that are $\Lambda$-periodic with $\Lambda=I_{2\times 2}$. For convenience we introduce for the representative volume element the shorthand $\Box\colonequals \Box_{\Lambda}=(-\frac12,\frac12)\times[-\frac12,\frac12)^2$.
The following lemma yields a sufficient condition on the quadratic form $Q$ for orthotropicity.
\begin{lemma}[A sufficient condition for orthotropicity]\label{L:orth1}
  Let the quadratic form $Q$ be admissible in the sense of Definition~\ref{D:admissible},
  and assume that $\Lambda = I_{2\times 2}$. Suppose that $Q$ is independent of $y_2$ and takes the form
  \begin{equation}\label{formulaforQ}
    Q(x_3,y,G)=\lambda(x_3,y_1){\rm tr}(G)^2+2\mu(x_3,y_1)|\sym G|^2
  \end{equation}
  with Lamé coefficients $\lambda$, $\mu$ satisfying the symmetry conditions
  \begin{equation*}
    \lambda(x_3,y_1)=\lambda(-x_3,y_1)=\lambda(x_3,-y_1)
    \qquad \text{and} \qquad
    \mu(x_3,y_1)=\mu(-x_3,y_1)=\mu(x_3,-y_1)
  \end{equation*}
  for all $x_3\in \left(-\frac12,\frac12 \right)$ and a.e.\ $y_1\in(-\frac12,\frac12)$.
  Then $Q^{\gamma}_{\rm hom}$ associated with $Q$ via Definition~\ref{def:Qgamma} is orthotropic. 
  Moreover, the corrector pair $(M_i,\varphi_i)$ associated with $G_i$ via \eqref{eq:corrector_equation} satisfies the following properties:
  \begin{align}
    \label{L:orth1:a}&M_1=M_2=M_3=0,\\
    \intertext{and a.e.~in~$\Box_\Lambda$ we have,}
    \label{L:orth1:c}&\partial_{y_2}\varphi_1=\partial_{y_2}\varphi_2=\partial_{y_2}\varphi_3=\varphi_1\cdot e_2=\varphi_2\cdot e_2=\varphi_3\cdot e_1=\varphi_3\cdot e_3=0,\\
    \label{eq:trace}&\mathbb L\big(\iota(x_3G_3)+\sym\nabla_\gamma\varphi_3\big):(e_j\otimes e_j)=0\qquad\text{for $j=1,2,3$}.
  \end{align}
  Above, we denote by $e_1,e_2,e_3$ the canonical basis of $\R^3$.
\end{lemma}
\reftoproof{SS:Formulas}

Next, we specialize the situation further to a class of examples for which we have explicit solutions to the corrector problem \eqref{eq:corrector_equation}:
In addition to the assumptions of Lemma~\ref{L:orth1}, in the following result we consider a composite in which the elastic law described by $Q$ is additionally independent of $x_3$ and thus describes a laminate. Furthermore, we assume that the components of the composite have a Poisson ratio that vanishes (i.e., $\lambda=0$). We shall also derive simplified formulas for $B^\gamma_{\eff}$ in the case of a prestrain $B$ that is rotationally invariant, i.e., of the form  $B=\rho I_{2\times 2}$ for a scalar function $\rho$.
\begin{lemma}[Laminate of isotropic materials with vanishing Poisson ratio]\label{S:ex:ex1}
  Let $(Q,B)$ be admissible in the sense of Definition~\ref{D:admissible},
  and assume that $\Lambda = I_{2\times 2}$. Suppose that $Q$ takes the form
  \begin{equation*}
    Q(x_3,y,G)=2\mu(y_1)|\sym G|^2,
  \end{equation*}
and assume that $\mu(-y_1)=\mu(y_1)$ for a.e. $y_1\in (-\frac12,\frac12)$.
  We denote the harmonic and arithmetic means of $\mu$ by
  \begin{subequations}
  \begin{align}
   \label{muGamma1}
    \meanh{\mu} & \colonequals \Big( \int_{-\frac12}^{\frac12} \frac{1}{\mu} \,dy_1 \Big)^{-1}
    \qquad \text{and} \qquad
    \overline{\mu}\colonequals \int_{-\frac12}^{\frac12}\mu \,dy_1,
  \end{align}
  respectively, and introduce for $\gamma\in(0,\infty)$ the following weighted average
  \begin{align}
   \label{muGamma2}
    \mu_{\gamma}
    & \colonequals
      \min_{w \in\mathcal H}\fint_{\Box}\mu\Big(\Big(\sqrt{12} x_3+\partial_{y_1}w\Big)^2+\Big(\tfrac{1}{\gamma}\partial_3w\Big)^2\Big) \,d(x_3,y),
  \end{align}
\end{subequations}
where
\begin{equation*}
  \mathcal H\colonequals\Big\{w\in H^1_{\gamma}(\Box;\R)\,:\,\fint_\Box w=0,\,\partial_{y_2}w=0\text{ a.e.~in~$\Box$}\,\Big\}.
\end{equation*}
 Let $Q^\gamma_{\hom}$ and $B^\gamma_{\eff}$ denote the effective quantities associated with $(Q,B)$ defined via Definitions~\ref{def:Qgamma} and \ref{D:eff}.
  Then the following properties hold:
  \begin{enumerate}[(a)]
  \item We have
    \begin{equation}\label{S:ex:muGammaProp1}
      \meanh{\mu}\leq \mu_\gamma \leq  \overline{\mu} \quad \text{for all $\gamma\in (0,\infty)$}.
    \end{equation}
    Furthermore,  the map $(0,\infty)\ni\gamma\mapsto \mu_\gamma$ is continuous and monotonically decreasing, and satisfies
    \begin{equation}\label{S:ex:muGammaProp2}
      \lim\limits_{\gamma\to 0}\mu_\gamma= \overline{\mu}
      \qquad \text{and} \qquad
      \lim\limits_{\gamma\to \infty}\mu_\gamma= \meanh \mu.
    \end{equation}
    Furthermore, if $\mu$ is non-constant, then $\gamma\mapsto\mu_\gamma$ is strictly monotone.
  \item $Q^\gamma_{\hom}$ is orthotropic with coefficients
    \begin{equation}
      \label{eq:q_hom_laminate_coefficients}
      q_{1} = \frac{1}{6}\meanh\mu, \quad q_{2}= \frac{1}{6}\overline\mu, \quad q_{3}=\frac{1}{6}\mu_\gamma, \quad q_{12}=0.
    \end{equation}
  \item Assume that $B(x_3,y)=\rho(x_3,y_1)I_{3\times 3}$ for a scalar function $\rho$. Then the coefficients of $B^\gamma_{\eff}$ are given by
    \begin{equation}\label{eq:Beff}
      \widehat B^\gamma_{{\rm eff},1} = 12 \fint_{\Box}\rho\, x_3\,d(x_3,y),\quad
      \widehat B^\gamma_{{\rm eff},2} = \frac{12}{\overline{\mu}}\fint_{\Box}\mu\rho x_3\,d(x_3,y),\quad
      \widehat B_{{\rm eff},3}^\gamma = 0.
    \end{equation}
  \end{enumerate}
\end{lemma}
\reftoproof{SS:Formulas}

  \begin{remark}[Physical interpretation]
    The harmonic mean $\meanh{\mu}$ and the arithmetic mean $\overline\mu$ are typical averages in homogenization of laminates: By assumption, the laminate under consideration oscillates in the $y_1$-coordinate, while it is constant with respect to $y_2$. As a consequence, the corrector associated with $G_2$ vanishes and leads to the arithmetic mean for the effective coefficient $q_2$. On the other hand, the corrector associated with $G_1$ oscillates and leads to the harmonic mean for $q_1$. Furthermore, since $\fint_\Box|\sqrt{12}x_3|^2\,d(x_3,y)=1$, the quantity $\mu_\gamma$ is a weighted average that interpolates between the arithmetic and harmonic mean.
\end{remark}

\begin{figure}
  \centering
  \begin{subfigure}[c]{0.45\textwidth}
    \def\svgscale{0.5}
    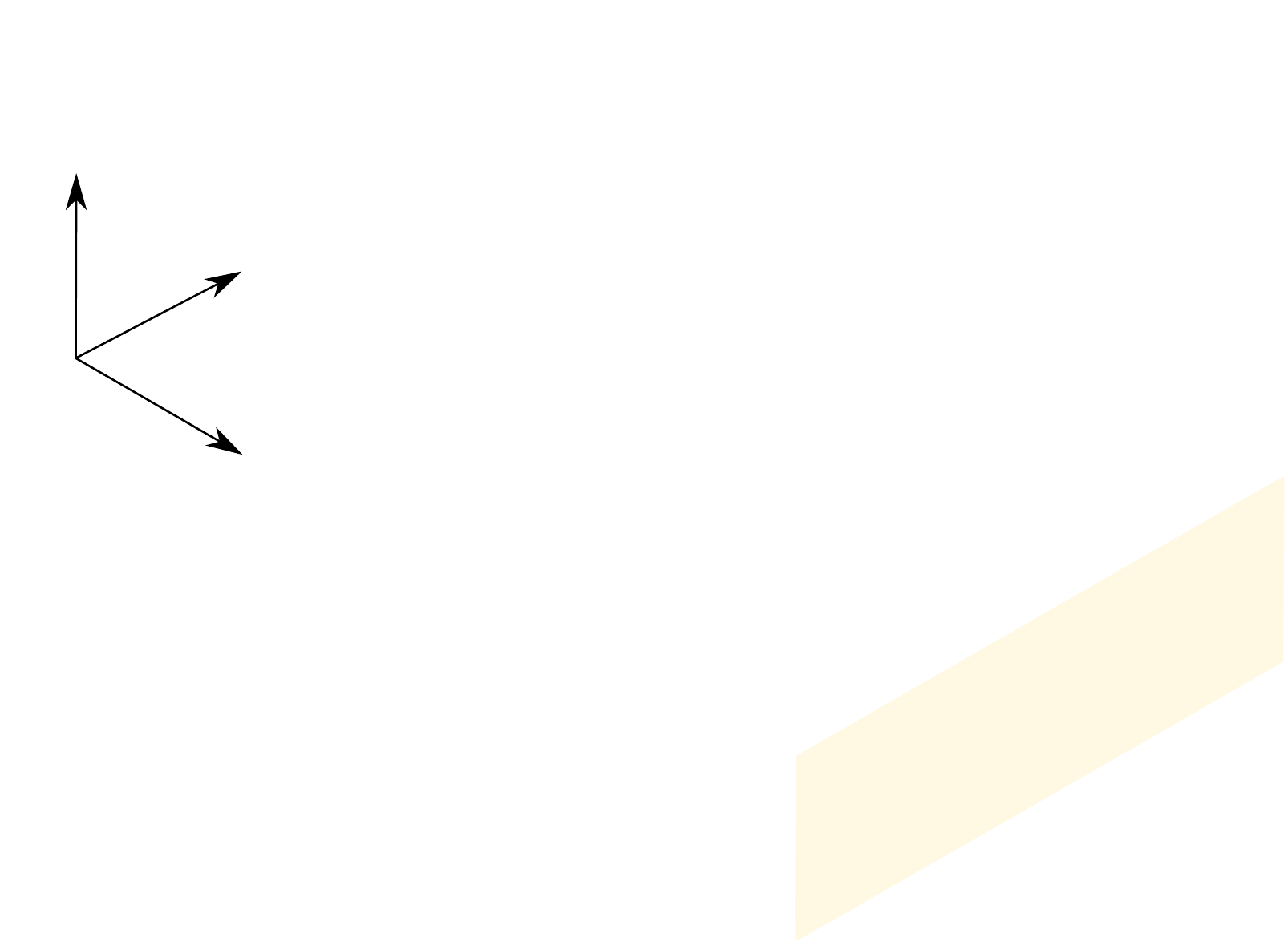 
    \subcaption{Prestrain strength $\rho(x_3,y_1)$}
    \label{sketchofmicrostructure_a}
  \end{subfigure} \hfill%
  \begin{subfigure}[c]{0.45\textwidth}
    \def\svgscale{0.5}
\begingroup%
  \makeatletter%
  \providecommand\color[2][]{%
    \errmessage{(Inkscape) Color is used for the text in Inkscape, but the package 'color.sty' is not loaded}%
    \renewcommand\color[2][]{}%
  }%
  \providecommand\transparent[1]{%
    \errmessage{(Inkscape) Transparency is used (non-zero) for the text in Inkscape, but the package 'transparent.sty' is not loaded}%
    \renewcommand\transparent[1]{}%
  }%
  \providecommand\rotatebox[2]{#2}%
  \newcommand*\fsize{\dimexpr\f@size pt\relax}%
  \newcommand*\lineheight[1]{\fontsize{\fsize}{#1\fsize}\selectfont}%
  \ifx\svgwidth\undefined%
    \setlength{\unitlength}{354.33202874bp}%
    \ifx\svgscale\undefined%
      \relax%
    \else%
      \setlength{\unitlength}{\unitlength * \real{\svgscale}}%
    \fi%
  \else%
    \setlength{\unitlength}{\svgwidth}%
  \fi%
  \global\let\svgwidth\undefined%
  \global\let\svgscale\undefined%
  \makeatother%
  \begin{picture}(1,0.95866539)%
    \lineheight{1}%
    \setlength\tabcolsep{0pt}%
    \put(0,0){\includegraphics[width=\unitlength,page=1]{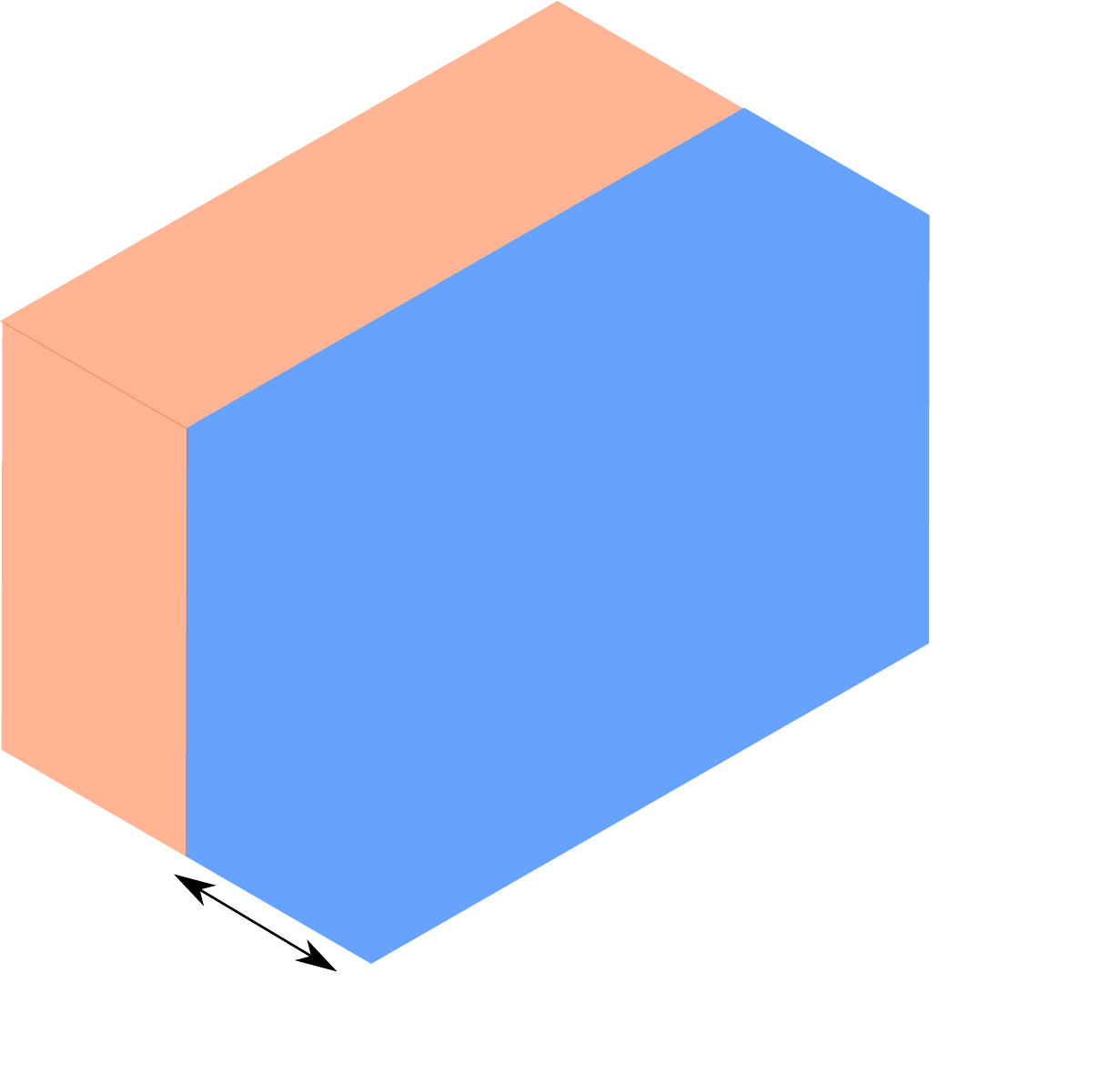}}%
    \put(0.1789901,0.08710228){\color[rgb]{0,0,0}\rotatebox{-32.37269844}{\makebox(0,0)[lt]{\lineheight{1.25}\smash{\begin{tabular}[t]{l}$\theta$\\\end{tabular}}}}}%
    \put(0,0){\includegraphics[width=\unitlength,page=2]{parLam-part2.pdf}}%
    \put(0.58029231,0.52280628){\color[rgb]{0,0,0}\rotatebox{30.2779555}{\makebox(0,0)[lt]{\lineheight{1.25}\smash{\begin{tabular}[t]{l}$ \mu_1$\\\end{tabular}}}}}%
    \put(0.25053368,0.71244845){\color[rgb]{0,0,0}\rotatebox{30.2779555}{\makebox(0,0)[lt]{\lineheight{1.25}\smash{\begin{tabular}[t]{l}$\mu_1$\\\end{tabular}}}}}%
    \put(0.42709436,0.62234149){\color[rgb]{0,0,0}\rotatebox{30.2779555}{\makebox(0,0)[lt]{\lineheight{1.25}\smash{\begin{tabular}[t]{l}$\mu_2$\\\end{tabular}}}}}%
    \put(0,0){\includegraphics[width=\unitlength,page=3]{parLam-part2.pdf}}%
  \end{picture}%
\endgroup%
 
    \subcaption{Lamé parameter $\mu(x_3,y_1)$}
  \end{subfigure}
  \caption{Schematic view of the microstructure of Lemma~\ref{S:ex:C3}}%
  \label{sketchofmicrostructure}
\end{figure}


We will now specialize the particular case of Lemma~\ref{S:ex:ex1} even further:
We divide the representative volume element of the composite $\Box=(-\frac12,\frac12)\times[-\frac12,\frac12)^2$ into two regions in the $x_{3}$-axis (top and bottom) and two regions
in the $y_{1}$-axis (middle and its complement).  The middle region has width $\theta$.
The remaining Lamé parameter $\mu$ takes one value in the middle and another value in the complement region. The prestrain strength $\rho$ takes one value
in the top, complement region, another value in the bottom, middle region, and is $0$ elsewhere. The setting is illustrated in Figure~\ref{sketchofmicrostructure}, and formalized
in the following lemma.

\begin{lemma}[A parametrized laminate]\label{S:ex:C3}
  For parameters $\mu_1>0$, $\rho_1\in\R$, $\theta\in[0,1],\, \theta_{\rho}\in\R,$ and $\theta_{\mu}>0$, we consider the situation of Lemma \ref{S:ex:ex1} where $\rho$ and $\mu$ are defined by
  \begin{equation*}
    \rho(x_3,y_1)
    \colonequals
    \begin{cases}
      \rho_1&\text{if }|y_1|>\frac\theta 2\text{ and }x_3>0,\\
      \rho_2&\text{if }|y_1|<\frac\theta 2\text{ and }x_3<0,\\
      0&\text{else,}
    \end{cases}
    \qquad \text{and} \qquad
    \mu(x_3,y_1)
    \colonequals
    \begin{cases}
      \mu_1&\text{if  } |y_1| > \frac\theta 2,\\
      \mu_2&\text{else},\\
    \end{cases}
  \end{equation*}
 where $\rho_2 \colonequals \theta_{\rho}\rho_1$ and $\mu_2\colonequals \theta_{\mu}\mu_1$.
Then we have
  \begin{equation}\label{S:ex:coeff}
   q_{1}= {\mu_1}\frac{\theta_{\mu}}{6(\theta+(1-\theta)\theta_{\mu})},\qquad q_{2}=\frac{\mu_1}{6}\big((1-\theta)+\theta\theta_{\mu}\big),
  \end{equation}
  and
  \begin{align}
   \label{S:ex:Prestrain}
    \widehat B^\gamma_{{\rm eff},1} = \frac{3\rho_1}{2}\big(1-\theta(1+\theta_{\rho})\big),\quad
    \widehat B^\gamma_{{\rm eff},2} = \frac{3\rho_1}{2} \left( \frac{1-\theta(1+\theta_{\mu}\theta_{\rho})}{1-\theta+\theta\theta_{\mu}} \right), \quad
    \widehat B^\gamma_{{\rm eff},3} = 0.
  \end{align}
\end{lemma}
\reftoproof{SS:Formulas}
%
\begin{table}[h!]
\begin{empheq}[box=\fbox]{align*}
\theta      &: \text{ volume fraction of the components} \\ 
\theta_\mu  &: \text{ stiffness ratio} \\ 
\theta_\rho &: \text{ prestrain contrast} \\
\mu_1       &: \text{ material stiffness }\\ 
\rho_1      &: \text{ strength of the isotropic prestrain}
\end{empheq}
\caption{Parameters for the parametrized laminate introduced in Lemma~\ref{S:ex:C3} }
\label{parList}
\end{table}
\subsection{Numerical computation of the effective quantities}
\label{S:ex:2}
In this section we answer Question~\ref{item:relations_q1}
for the parametrized laminate introduced in Lemma~\ref{S:ex:C3} by numerically exploring the parameter dependence of the effective quantities $Q^\gamma_{\hom}$ and $B^\gamma_{\eff}$. Recall that by Lemma~\ref{S:ex:C3}, $Q^\gamma_{\hom}$ is orthotropic with coefficients $q_1,q_2$ given by \eqref{S:ex:coeff}, $q_{12}=0$, and $q_3(\gamma) = \frac{1}{6} \mu_\gamma$ with $\mu_\gamma$
given by \eqref{muGamma2}.
The coefficients of $B^\gamma_{\eff}$ are given by \eqref{S:ex:Prestrain}.
These coefficients (and thus $Q^\gamma_{\hom}$ and $B^\gamma_{\eff}$) depend on the parameters of the model
shown in Table~\ref{parList}. As can be seen by a close look at the formulas for the coefficients, the following scaling properties hold:
\begin{equation}\label{eq:QBinvariance}
  \begin{aligned}
  Q^\gamma_{\hom}(G;\theta,\theta_\mu,\theta_\rho,\mu_1,\rho_1)=\,&\,\mu_1 Q^\gamma_{\hom}(G;\theta,\theta_\mu,1,1,1),\\
  B^\gamma_{\eff}(\theta,\theta_\mu,\theta_\rho,\mu_1,\rho_1)=\,&\,\rho_1B^{\gamma=1}_{\eff}(\theta,\theta_\mu,\theta_\rho,1,1).
\end{aligned}
\end{equation}
In particular, $Q^\gamma_{\hom}$
does not depend on $\theta_\rho$ and $\rho_1$, and $B^\gamma_{\eff}$
does not depend on $\mu_1$ and the scaling parameter $\gamma$.
In view of this, we set $\mu_1=\rho_1=1$ in the following.
With the idea in mind that the volume fraction $\theta$ can be controlled during the fabrication of the composite, we shall mainly focus on the functional dependence of $Q^\gamma_{\hom}$ and $B^\gamma_{\eff}$ on $\theta$ for various values of $\theta_\mu$ and $\theta_\rho$.
\smallskip
\begin{figure}[t]
 \captionsetup{type=figure}\addtocounter{figure}{-1}
\begin{minipage}[t]{.5\textwidth} 
    \begin{subfigure}[c]{0.9\textwidth}
	 \includegraphics{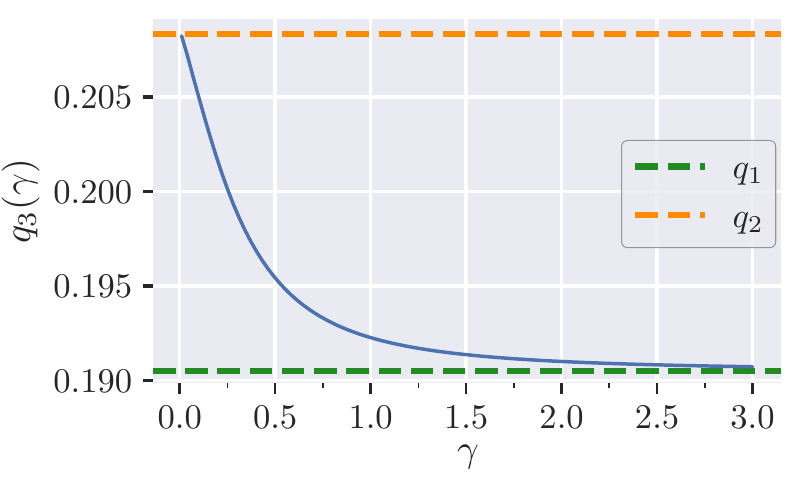}
    \end{subfigure} \hfill 
\end{minipage} \hfill
\begin{minipage}[t]{.4\textwidth}
\vspace{-\baselineskip}
\centering
\begin{tabular}[t]{c|c|c|c}
    Parameters &    $\theta_{\mu}$  &  $\theta_{\rho}$  & $\theta $      \\ \hline
  			   &   		2			&     2				&  $\frac{1}{2}$ \\
\end{tabular}
    \captionof{figure}{Numerical approximation of $q_3(\gamma) = \frac{1}{6}\mu_\gamma$ }
\label{q3-plot}
\end{minipage} 
\end{figure}%
In view of the definition of the coefficients (cf.~\eqref{S:ex:coeff}, \eqref{eq:q_hom_laminate_coefficients}, and \eqref{S:ex:Prestrain}) all except the coefficient $q_3(\gamma)$ are given by explicit formulas, which can be evaluated without computational effort. 
In contrast, for $q_3(\gamma)$ we need to solve the Euler-Lagrange equation for the minimization problem~\eqref{muGamma2}. It is a two-dimensional, elliptic system that we solve using a finite element method implemented
in \CC \ using the \textsc{Dune} libraries \cite{DUNE, sander2020dune}. We also note that $q_3(\gamma)$ is the only coefficient that depends on the scaling parameter $\gamma$.
Figure~\ref{q3-plot} shows the value of $q_3(\gamma)$ obtained numerically for
different values of the relative scaling parameter $\gamma$. 
We clearly observe property~\eqref{S:ex:muGammaProp1}, i.e., $q_{1} < q_{3}(\gamma) <  q_{2}$
for $\gamma \in (0,\infty)$, continuity and monotonicity of $\gamma\mapsto q_3(\gamma)$, and the asymptotic behavior
\begin{equation}\label{ex:proxy}
 \lim_{\gamma \to 0} q_{3}(\gamma) = q_{2}
 \qquad \text{and} \qquad
 \lim_{\gamma \to \infty} q_{3}(\gamma) = q_{1},
\end{equation}
as predicted by Lemma~\ref{S:ex:ex1}. Furthermore, $q_3(\gamma)$ appears to be a convex function with a rapid decay for values of $\gamma$ close to zero.
Let us anticipate that in the next section we shall appeal to the asymptotic behavior \eqref{ex:proxy} to avoid the computational effort required to solve \eqref{muGamma2} numerically: in view of  \eqref{ex:proxy}, in the extreme cases $0< \gamma \ll 1$ and $ \gamma \gg 1 $, we may use  $q_{2}$ and $q_{1}$ as approximations for $q_3(\gamma)$, respectively. 
\begin{figure}[H]
\includegraphics[trim={0 0.25cm 0 0.25cm}, clip]{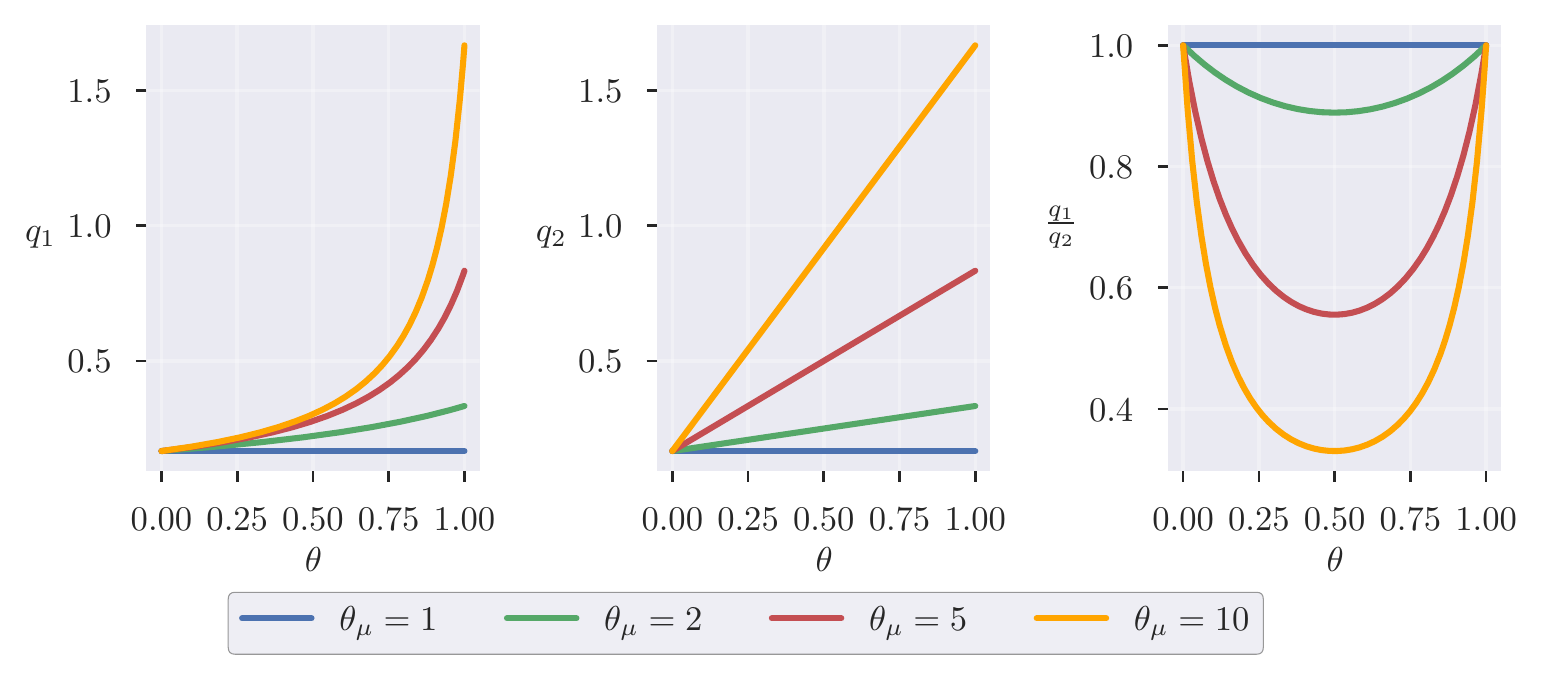} 
\caption{Coefficients $q_1, q_2$ and the ratio $q_1/q_2$  as functions of the
 volume fraction $\theta$ of the constituents, for different values of $\theta_{\mu}$}
\label{plot-q1q2}
\end{figure}
Figure~\ref{plot-q1q2} shows the elastic moduli $q_1, q_2$, and their ratio $\frac{q_1}{q_2}$
as functions of the volume fraction~$\theta$ for different values of the stiffness ratio $\theta_\mu$.
The ratio $\frac{q_1}{q_2}$ serves as a measure for the in-plane anisotropy of the effective material.
As expected from \eqref{S:ex:coeff} we observe a nonlinear dependence of $q_1$ on
the volume fraction $\theta$, and a linear dependence of $q_2$ on $\theta$.
Note that for $\theta_\mu=1$ the laminate reduces to a homogeneous material.
By increasing the stiffness ratio $\theta_\mu$, the slopes of $q_1$ and $q_2$ as functions
of $\theta$ increase while the ratio $\frac{q_1}{q_2}$ decreases. We note that as $\theta_{\mu}$ increases, so does the possible anisotropy. Additionally, for all $\theta_\mu > 1$ the ratio $\frac{q_1}{q_2}$ assumes its minimum at $\theta=\frac{1}{2}$.
Next, we study the effective prestrain $B_{\eff}^\gamma$ and its dependence on the laminate parameters. Recall that  $\widehat B^\gamma_{\eff,3}=0$ by Lemma~\ref{S:ex:ex1}. (Here and below, we use the notation introduced in \eqref{canonicalBasis} for the coefficients of $B^\gamma_{\eff}$.)  
\begin{figure}[H]
\begin{center}  
\includegraphics[trim={0 0.25cm 0 0.25cm}, clip]{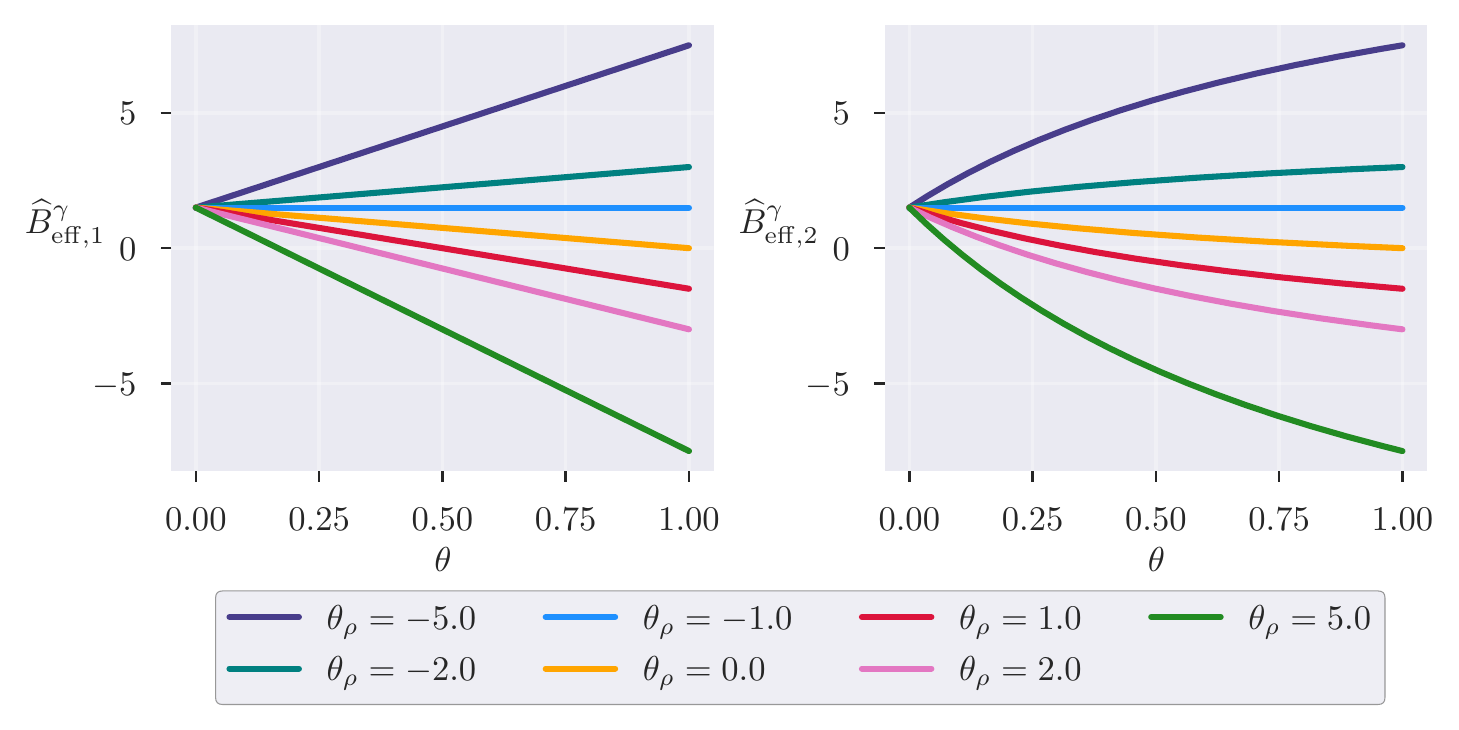}
\end{center}
\vspace{-\baselineskip}
\caption{Effective prestrain components $\widehat B_{\text{eff},1}^{\gamma}$ (left)
   and $\widehat B_{\text{eff},2}^{\gamma}$ (right) for a fixed value of $\theta_{\mu} = 2$
   with different values of $\theta_{\rho}$}
\label{plot-prestrain-thetaRho}
\end{figure}

Figure~\ref{plot-prestrain-thetaRho} displays the effective prestrain coefficients $\widehat B_{\text{eff},1}^{\gamma}$
and $\widehat B_{\text{eff},2}^{\gamma}$ as functions of the volume fraction $\theta$
for different values of the prestrain ratio $\theta_\rho$. Note that the dependence of $\widehat B_{\text{eff},1}^{\gamma}$ on $\theta$ is linear,
while the dependence of $\widehat B_{\text{eff},2}^{\gamma}$ on $\theta$ is nonlinear.
This is somewhat surprising, since the corrector associated to $G_{1}$ is non-zero,
while the corrector associated to $G_{2}$ is zero,
cf.\ the proof of Lemma~\ref{S:ex:ex1}.
Furthermore, note that the effective prestrain is not zero if $\theta=0$, and the sign, as well as the slope of the effective prestrain depends on the prestrain ratio $\theta_\rho$. %
\vspace{\baselineskip}
\begin{figure}[H]
 \captionsetup{type=figure}\addtocounter{figure}{-1}
\begin{minipage}[t]{.5\textwidth} 
\vspace{-\baselineskip}
    \begin{subfigure}[c]{1\textwidth}
\includegraphics{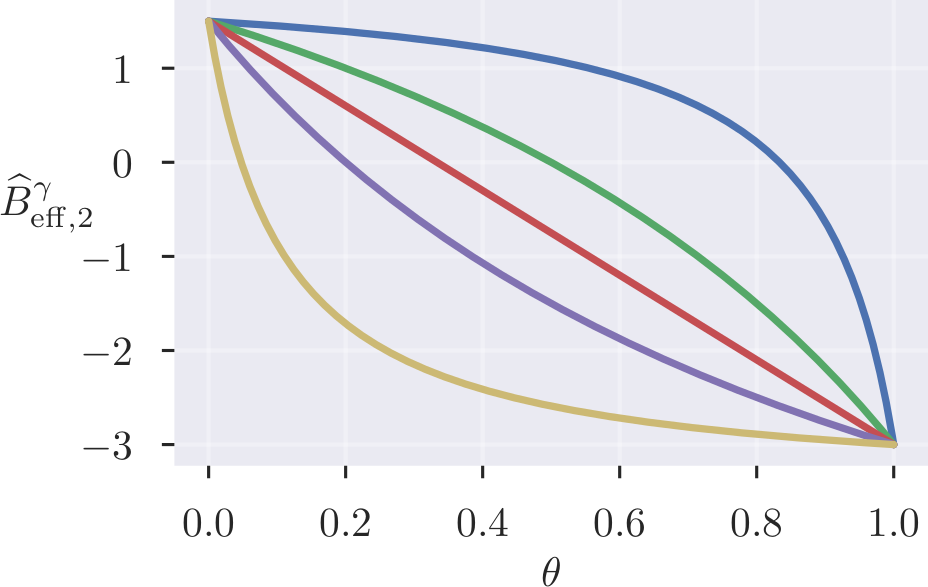} 
    \end{subfigure} \hfill 
\end{minipage} \hfill
\begin{minipage}[t]{.4\textwidth}
\vspace{-\baselineskip}
\centering
\includegraphics{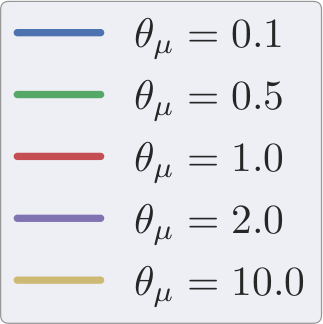}
    \captionof{figure}{Effective prestrain $\widehat B_{\text{eff},2}^{\gamma}$ for a fixed value of $\theta_{\rho}=2$ with different values of $\theta_\mu$}
\label{plot-prestrain-thetaMu}
\end{minipage} 
\end{figure}

In a similar way, Figure~\ref{plot-prestrain-thetaMu} shows the dependence of $\widehat B_{\eff,2}^\gamma$ on the volume fraction $\theta$ but for different stiffness ratios $\theta_\mu$. Note that $\widehat B_{\eff,1}^\gamma$ is independent of $\theta_\mu$, see Lemma~\ref{S:ex:C3}.
Again, the dependence of $\widehat B_{\eff,2}^\gamma$ on $\theta$ is nonlinear except for the homogeneous case $\theta_\mu = 1.0$,
where the dependence becomes affine.



\section{The microstructure--shape relation}\label{S:micro-shape}

In this section, we investigate Question~\ref{item:relations_q2} and combine it with
the results of Section~\ref{S:properties} in order to explore the parameter-dependence of
the shapes of free minimizers of the homogenized energy functional $\mathcal I^\gamma_{\hom}$ for the parametrized laminate material of Lemma~\ref{S:ex:C3}.
Note that unless stated otherwise, with ``minimizer'' we always refer to a \emph{global} minimizer.
As mentioned before, in the present paper, we restrict our analysis to the spatially homogeneous case and thus assume that
\begin{equation}\label{eq:st:homcase}
  \text{$Q_{\hom}^\gamma$ and $B_{\eff}^\gamma$ are independent of $x'\in S$.}
\end{equation}
A sufficient condition for \eqref{eq:st:homcase} is the global periodicity of the composite's microstructure.
Condition~\eqref{eq:st:homcase} allows to simplify the minimization of $\mathcal I^\gamma_{\hom}$ drastically:
Every bending deformation that minimizes $\mathcal I^\gamma_{\hom}$ has a constant fundamental form and thus parametrizes a cylindrical surface with constant curvature. Note that for  $\deform\in H^2_{\iso}(S;\R^3)$, thanks to the isometry constraint $(\nabla'\deform)^\top\nabla'\deform=I_{2\times 2}$, we have
\begin{equation*}
  \II_\deform(x')\in\mathcal G\colonequals \big\{G\in\R^{2\times 2}_{\sym}\,:\,\det G=0\big\}\qquad\text{for a.e. }x'\in S.
\end{equation*}
Thus, under condition \eqref{eq:st:homcase}, minimization of the non-convex integral functional $\mathcal I^\gamma_{\hom}$ reduces to the following \emph{algebraic minimization problem}:
\begin{equation}\label{L:ex:1:eq1}
  \mathcal S
  \colonequals
  \argmin_{G\in\mathcal G}Q_{\hom}^\gamma(G-B_{\rm eff}^\gamma).
\end{equation}
More precisely, we recall the following result from  \cite{schmidt2007minimal}:
\begin{lemma}[\mbox{\cite[Theorem~3.2]{schmidt2007minimal}}]\label{L:char:cylindrical}
  In the situation of Theorem~\ref{T1} suppose that $Q_{\hom}^\gamma$ and $B_{\eff}^\gamma$ are
  independent of $x'\in S$. 
  Then $\deform\in H^2_{\rm iso}(S;\R^3)$ is a minimizer of $\mathcal I^\gamma_{\hom}$,
  if and only if there exists $G$ solving~\eqref{L:ex:1:eq1} such that $\II_\deform=G$
  almost everywhere in $S$, where $\II_\deform$ denotes the second fundamental form of $\deform$.
\end{lemma}
In view of this, we study the dependence of the solution set $\mathcal S$ on the effective coefficients
of $Q^\gamma_{\hom}$ and $B^\gamma_{\eff}$.
Having the parametrized laminate of Lemma~\ref{S:ex:C3} in mind, we focus on the orthotropic case
of Definition~\ref{def:orthotropicity}.
In Section~\ref{S:classification_orthotropic_case} we present a classification result that
describes $\mathcal S$ for a general orthotropic quadratic form and a general effective, diagonal prestrain.
In Section~\ref{S:param2} we then explore the microstructure--shape relation by analyzing the dependence of $\mathcal S$ on the parameters of the parametrized laminate material of Section~\ref{S:ex:2}.

\begin{remark}[Cylindrical surfaces and their parametrization by angle and curvature]\label{R:kappaalpha}
  The geometry of a cylindrical surface can be conveniently parametrized by an angle and a scalar curvature.  We shall use this parametrization in the characterization and visualization of the set~$\mathcal S$ in the upcoming section. Let us first fix our terminology: Recall that $\deform\in H^2_{\iso}(S;\R^3)$ is called  \emph{cylindrical} if $\II_\deform$ is constant, i.e., if there exists
  $G\in\mathcal G$ such that $\II_\deform=G$ a.e.~in $S$.
  We note that for any $G\in\mathcal G$ there exists a deformation $\deform_G\in H^2_{\iso}(S;\R^3)$  with $\II_{\deform_G}=G$.
  To see this, first note that any $G\in\mathcal G$ can represented as
  \begin{equation}\label{D:kappa_e}
    G=\kappa \minvec\otimes \minvec
    \qquad\text{for some unique $\kappa\in\R$ and a vector $\minvec\in\R^2$ with $\abs{\minvec}=1$}.
  \end{equation}
  In the case $G\neq 0$, the line spanned by $\minvec$ is uniquely determined by $G$.
  A direct calculation shows that the map $\deform_G:S\to\R^3$,
  \begin{equation*}
    \deform_G(x')\colonequals\Big(\int_0^{x'\cdot \minvec}\cos(-\kappa s)\,\dd s,\,x'\cdot \minvec^\perp,\,\int_0^{x'\cdot \minvec}\sin(-\kappa s)\,\dd s\Big)^\top,\qquad \minvec^\perp\colonequals
    \begin{pmatrix}
      0&-1\\
      1&0
    \end{pmatrix}\minvec,
  \end{equation*}
  defines a bending deformation with its second fundamental form satisfying  $\II_{\deform_G}=G$.
  In fact, by the rigidity theorem for surfaces (see, e.g.,~\cite[Theorem~3]{ciarlet2002recovery}),
  any parametrized surface $\deform\in H^2_{\iso}(S;\R^3)$ with $\II_\deform=G$ a.e.~in $S$ equals $\deform_G$ modulo a superposition with a Euclidean transformation of $\R^3$.
  Geometrically, $\deform_G$ parametrizes a cylindrical surface, whose (nonzero) principal curvature is given by $\kappa$
  and with associated principal direction (expressed in local coordinates) $\pm \minvec$, see Figure~\ref{fig:angle}.
  %
  For $\kappa<0$ the surface is bent in the direction of the surface normal $\partial_1\deform\wedge\partial_2\deform$.

  For visualizations it is convenient to parametrize the set $\mathcal G\setminus \{0\}$ by
  associating to each $G = \kappa \minvec \otimes \minvec \in \mathcal{G}$ the curvature $\kappa$ and the angle
  \begin{equation}\label{D:angle}
    \alpha(G)\colonequals
    \begin{cases}
      \displaystyle\arctan\frac{\minvec \cdot e_2}{\minvec\cdot e_1}&\text{if }\minvec\cdot e_1\neq 0,\\
      \frac{\pi}{2}&\text{else}.
    \end{cases}
  \end{equation}
  Note that the expression on the right-hand side is the same for $\minvec$ and $-\minvec$ and thus $\alpha(G)$ is well-defined.
  Geometrically, $\alpha(G)$ is the angle required to rotate the line spanned by $e_1$ to the line spanned by $\minvec$
  (in counterclockwise direction).
  The map $\mathcal G\setminus\{0\}\mapsto (\alpha,\kappa)\in (-\frac\pi2,\frac\pi2]\times\R\setminus\{0\}$ is a bijection.
  It is even a homeomorphism if we identify the end points of $(-\frac\pi2,\frac\pi2]$.
\end{remark}

\begin{figure}[H]
 \centering
 \captionsetup{type=figure}\addtocounter{figure}{-1}
\begin{subfigure}[c]{0.3\textwidth}
\includegraphics[scale=0.5,trim={0 2cm 0 2cm},clip]{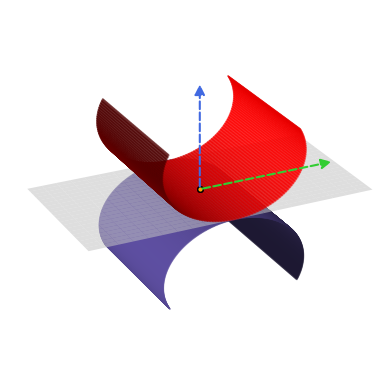}
\subcaption{$\alpha = 0, \ \minvec = \begin{pmatrix}
1 \\ 0   \end{pmatrix}$}
\end{subfigure} \hfill%
\begin{subfigure}[c]{0.3\textwidth}
\includegraphics[scale=0.5,trim={0 2cm 0 2cm},clip]{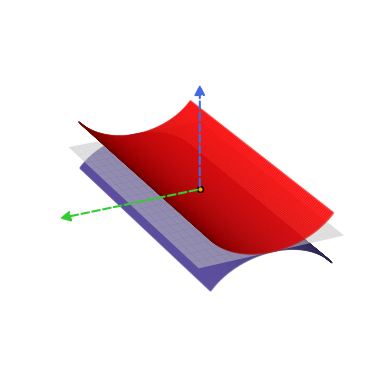}
\subcaption{$\alpha = \frac{\pi}{2}, \ \minvec = \begin{pmatrix}
0 \\ 1   \end{pmatrix}$}
\end{subfigure} \hfill
\begin{subfigure}[c]{0.3\textwidth}
\includegraphics[scale=0.5,trim={0 2cm 0 2cm},clip]{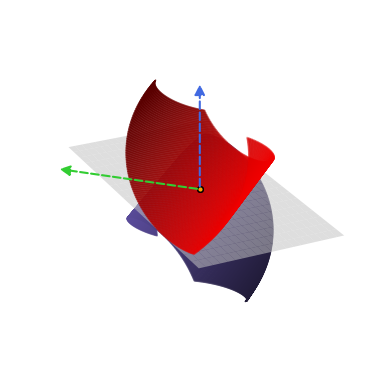}
\subcaption{$\alpha = \frac{\pi}{4}, \ \minvec = \frac{1}{\sqrt{2}}\begin{pmatrix}
1 \\ 1   \end{pmatrix}$}
\end{subfigure}
\par\bigskip \par \bigskip  
\begin{minipage}[t]{.5\textwidth} 
\vspace{-\baselineskip}
    \begin{subfigure}[c]{0.9\textwidth}
    \def\svgwidth{1\textwidth}
    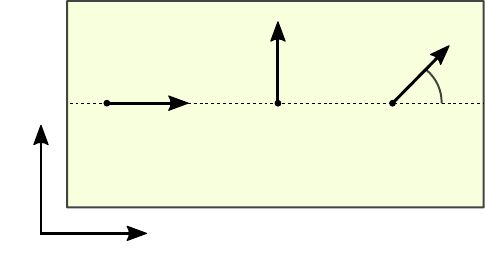
	\subcaption{Principal directions in local coordinates}
    \end{subfigure} \hfill 
\end{minipage}
\begin{minipage}[t]{.45\textwidth}
\vspace{-\baselineskip}
\captionof{figure}{Deformations with constant second fundamental form $\II = \kappa \minvec \otimes \minvec$ and $\abs{\kappa} = 2$.
  Blue color corresponds to positive values of the curvature $\kappa$ and red to negative values of $\kappa$.
  The green vector shows $(\nabla u)\minvec$ while the blue vector represents the surface normal $\partial_1\deform\wedge\partial_2\deform$}
\label{fig:angle}
\end{minipage}    
  
\end{figure}

\subsection{Classification of \texorpdfstring{$\mathcal S$}{S} in the orthotropic case}
\label{S:classification_orthotropic_case}

We analyze the algebraic minimization problem \eqref{L:ex:1:eq1} in the case of
an orthotropic quadratic form $Q^\gamma_\text{hom}$ and a diagonal prestrain $B^\gamma_\text{eff}$. For the upcoming discussion it is not important that $Q^\gamma_\text{hom}$ and $B^\gamma_\text{eff}$ are defined via the homogenization formulas of Section~\ref{S:correctors}. We rather consider a generic quadratic form and prestrain, which we denote by $Q$ and $B$ to simplify the notation.

More precisely, with $G_1,G_2,G_3$ the orthonormal basis of $\R^{2\times 2}_{\sym}$ introduced
in~\eqref{canonicalBasis}, let $Q:\R^{2\times 2}_{\sym}\to\R$ be a positive definite quadratic form, $B\in\R^{2\times 2}_{\sym}$, and assume that
\begin{subequations}\label{eq:ass:orth}
\begin{align}
  &\text{$Q$ is orthotropic in the sense of Definition~\ref{def:orthotropicity} with coefficients $q_1$, $q_2$, $q_3$, $q_{12}$}, \\
  \shortintertext{and}
  &B\text{ is diagonal and invertible, i.e., $\det B\neq 0$}.
\end{align}
\end{subequations}
We remark that in view of \eqref{eq:ass:orth} the positive definiteness of $Q$ is equivalent to
\begin{equation}\label{eq:st:1001}
  q_1,q_2,q_3>0
  \qquad \text{and} \qquad
  -2\sqrt{q_1q_2}<q_{12}<2\sqrt{q_1q_2}.
\end{equation}

Our goal is to determine the set of minimizers
\begin{equation*}
  \mathcal S_{Q,B}\colonequals \argmin_{G\in\mathcal G} Q(G-B).
\end{equation*}
This is a quadratic minimization problem on the non-convex set $\mathcal{G}$, and thus a rich behavior can be expected. We first note that by the positive definiteness of $Q$ and in view of the assumption $\det B\neq 0$ we have $\mathcal S_{Q,B}\neq\emptyset$ and $0\not\in\mathcal S_{Q,B}$.
Thus, minimizers $G\in\mathcal S_{Q,B}$ exist and correspond to non-flat, cylindrical surfaces.

\paragraph{Axial minimizers.} An important role for the upcoming discussion is played by $G\in\mathcal G$
with angle $\alpha(G)\in\{0,\frac\pi2\}$. Such $G$ correspond to a cylindrical surface with a principal direction 
that in local coordinates is parallel to one of the coordinate axes of $\R^2$ (Figure~\ref{fig:angle}). We call such matrices
\emph{axial}, and note that
\begin{equation*}
  G\in\mathcal G\setminus\{0\}\text{ is axial}
  \qquad \iff \qquad
  G\in\operatorname{span}\{G_1\}\cup\operatorname{span}\{G_2\}.
\end{equation*}
Since $Q$ is positive definite, the restrictions of $G\mapsto Q(G-B)$ to $\operatorname{span}\{G_1\}$
or $\operatorname{span}\{G_2\}$ are strictly convex. They thus admit unique minimizers, which can be
computed by elementary calculations. We obtain
\begin{equation}\label{eq:st:1004}
  \big\{G\in\mathcal S_{Q,B}\,:\,G\text{ is axial}\big\}
  \subseteq
  \Big\{\frac{2q_1\widehat B_1+q_{12}\widehat B_2}{2q_1} G_1,\,\frac{2q_{12}\widehat B_1+q_{2}\widehat B_2}{2q_2}G_2\Big\},
\end{equation}
where $\widehat B_1,\widehat B_2$ denote the coefficients of $B$ in the sense of \eqref{canonicalBasis}.

As we shall see, for most choices of $Q$ and $B$, every $G\in\mathcal S_{Q,B}$ is axial.
In this case, $\mathcal S_{Q,B}$ consists of at most two (global) minimizers, and to determine $\mathcal S_{Q,B}$
we only need to compare the energy values associated with the elements in the set on the
right-hand side in \eqref{eq:st:1004}.
However, we shall see that for certain choices of $Q$ and $B$, the set $\mathcal S_{Q,B}$ contains
two or even a one-parameter family of infinitely many non-axial minimizers.
In the following, we develop an algorithm to compute the set $\mathcal S_{Q,B}$ in this case.

\paragraph{Mirror symmetry of the solution set.}
Consider the bijective transformation
\begin{equation*}
  T:\mathcal G\to\mathcal G,\qquad
  \begin{pmatrix}
    a_1&a_3\\a_3&a_2
  \end{pmatrix}\mapsto   \begin{pmatrix}
    a_1&-a_3\\-a_3&a_2
  \end{pmatrix}.
\end{equation*}
Then by orthotropicity of $Q$ and diagonality of $B$ we have $Q(G-B)=Q(TG-B)$, and thus
\begin{equation}\label{eq:S:P1}
  G\in\mathcal S_{Q,B}\qquad \iff \qquad   TG\in\mathcal S_{Q,B}.
\end{equation}
Geometrically, the transformation $T$ is a reflection in the following sense: If $G$ describes a cylindrical surface with
curvature $\kappa$ and angle $\alpha\in(-\frac\pi2,\frac\pi2)$, then  $TG$ corresponds
to a cylindrical surface with the same curvature $\kappa$ but angle $-\alpha$.
From \eqref{eq:S:P1} we conclude that
\begin{align}
\label{eq:st:1000a}
  &\mathcal S_{Q,B}= \big\{G,TG\,:\,G\in\mathcal S^+_{Q,B} \big\},\\
\nonumber
  &\text{where }\mathcal S^+_{Q,B}\colonequals \argmin_{G\in\mathcal G^+} Q(G-B),
  \qquad
  \mathcal G^+\colonequals\big\{G\in\mathcal G\;\text{with}\;G:G_3\geq 0\big\}.
\end{align}
Moreover, we note that for all $G\in\mathcal G\setminus\{0\}$ we have $G=TG$, if and only if $G$ is axial.

\paragraph{Classification of minimizers.} The set $\mathcal G^+$ can be conveniently parametrized by the nonlinear, bijective transformation
\begin{align*}
    \Phi & :\mathcal{G}^+_{\R^2}\to\mathcal G^+,\quad \Phi(a_1,a_2)\colonequals a_1G_1+a_2G_2+\sqrt{2a_1a_2}G_3,\\
  \mathcal{G}^+_{\R^2} & \colonequals \big\{a=(a_1,a_2)\in\R^2\,:\,a_1a_2\geq 0\big\}.
\end{align*}
We remark that  the boundary $\partial \mathcal{G}^+_{\R^2}$ corresponds to axial $G\in\mathcal G^+$, i.e., $\Phi(\partial \mathcal{G}^+_{\R^2})=\{G\in\mathcal G^+\,:\,G\text{ is axial}\}$.
In order to express $\mathcal G^+\ni G\mapsto Q(G-B)$ in these coordinates, we introduce the quadratic function
\begin{equation*}
  \mathcal E_{Q,B}:\R^2\to\R,\qquad   \mathcal E_{Q,B}(a)\colonequals\frac12 a\cdot Ha-2a\cdot Ab,
\end{equation*}
where 
\begin{equation}\label{eq:reformulation}
  H\colonequals
  \begin{pmatrix}
    2q_1&q_{12}+2q_3\\
    q_{12}+2q_3&2q_2      
  \end{pmatrix},\qquad A\colonequals
  \begin{pmatrix}
    q_1&\frac{q_{12}}{2}\\
    \frac{q_{12}}{2}&q_2
  \end{pmatrix},\qquad b\colonequals(\widehat B_1,\widehat B_2)^\top.
\end{equation}
One can easily check that for all $a\in \mathcal{G}^+_{\R^2}$ we have $  \mathcal E_{Q,B}(a)=Q(\Phi(a)-B)+c$ for a constant $c$ that is independent of $a$. We thus conclude that
\begin{equation}\label{eq:st:1000b}
  \mathcal S_{Q,B}^+ = \Phi\bigg(\argmin_{a\in \mathcal{G}^+_{\R^2}}  \mathcal E_{Q,B}(a)\bigg).
\end{equation}
The problem on the right-hand side is a quadratic minimization problem subject to the nonlinear constraint $a_1a_2\geq 0$.
Since $q_1,q_2>0$, the quadratic part of $\mathcal E_{Q,B}$ is elliptic, parabolic, or hyperbolic, if and only if
$\det H>0$,  $\det H=0$, or $\det H<0$, respectively. In the case $\det H>0$ the minimizer of $\mathcal E_{Q,B}$ is unique and given by 
\begin{equation}\label{def:gstar}
  g_*\colonequals2H^{-1}Ab.
\end{equation}
We can always compute $g_*$ in closed form using Cramer's rule.

We obtain the following classification of the set of minimizers:
\begin{lemma}[Trichotomy of minimizers]\label{L:char}
  Let $Q:\R^{2\times 2}_{\sym}\to\R$ be a positive definite quadratic form that is
  orthotropic in the sense of Definition~\ref{def:orthotropicity}.  Let further $B \in \R^{2 \times 2}$
  be diagonal with $\det B \neq 0$. Then exactly one of the following three cases has to hold:
  \begin{enumerate}[(a)]
  \item\label{L:char:a} \emph{(Axial minimizers).} All minimizers $G\in\mathcal S_{Q,B}$ are axial. Furthermore, $\mathcal S_{Q,B}$ is characterized as follows:
    \begin{equation}
      \label{eq:form_min}
      \emptyset\neq \mathcal S_{Q,B}
      \subseteq
      \Big\{\frac{2q_1\widehat B_1+q_{12}\widehat B_2}{2q_1} G_1,\,\frac{2q_{12}\widehat B_1+q_{2}\widehat B_2}{2q_2}G_2\Big\},
    \end{equation}
    and equality holds if and only if $\frac{(\widehat B_2)^2}{q_1}=\frac{(\widehat B_1)^2}{q_2}$.
  \item\label{L:char:b} \emph{(Two non-axial minimizers).} We have $\det H>0$ and $g_*$ defined in \eqref{def:gstar} is an interior point of $\mathcal{G}^+_{\R^2}$. Furthermore, $\mathcal S_{Q,B}$ is characterized as follows:
    \begin{equation*}
      \mathcal S_{Q,B}=\{G,TG\},\qquad G\colonequals\Phi(g_*).
    \end{equation*}
  \item\label{L:char:c} (One-parameter family of minimizers). We have $\det H=0$ and $Ab\in\operatorname{range} H$. Furthermore, $\mathcal S_{Q,B}$ is characterized as follows:
    \begin{equation*}
      \mathcal S_{Q,B}
      =
      \big\{G,TG\,:\,G=\Phi(a)\text{ for all $a\in \mathcal{G}^+_{\R^2}$ s.t.\ $a\cdot q_*=s_*$}\big\},
    \end{equation*}
    where
    \begin{equation}
    \label{L:char:c:quantities}
      q_*\colonequals\frac{1}{\sqrt{(q_{12}+2q_3)^2+4q_2^2}}
      \begin{pmatrix}
        q_{12}+2q_3\\ 2q_2
      \end{pmatrix}
      \qquad\text{and}\qquad s_*\colonequals\frac{q_*\cdot Ab}{q_1+q_2}.
    \end{equation}
  \end{enumerate}
\end{lemma}
\reftoproof{SS:Formulas}

\begin{figure}
  \centering
  \begin{subfigure}[c]{\textwidth}
  \centering
 \includegraphics{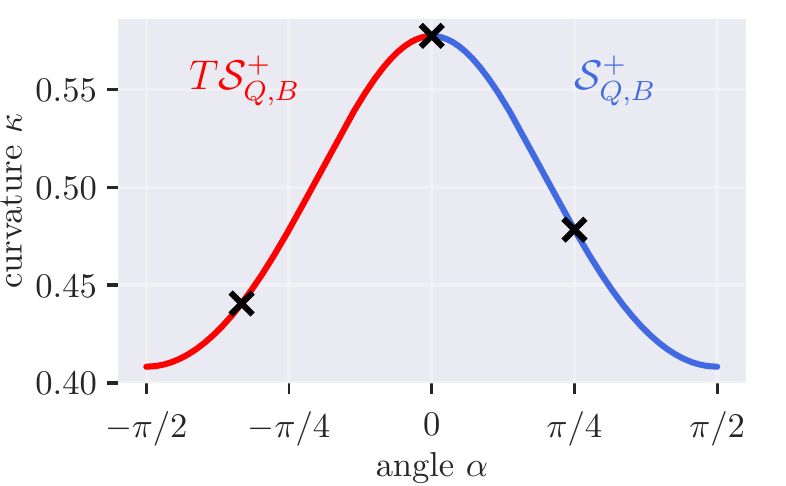}
 \subcaption{Plot of $\mathcal S_{Q,B}$ in $(\alpha,\kappa)$-coordinates.}
 \end{subfigure} \hfill%
\begin{subfigure}[c]{0.3\textwidth}
\includegraphics[width=\textwidth, trim={0 2cm 0 1.5cm},clip]{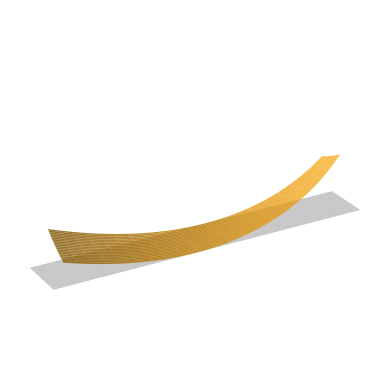}
\subcaption{$\alpha = -\frac{\pi}{3}, \kappa \approx 0.44$}
\end{subfigure} \hfill%
\begin{subfigure}[c]{0.3\textwidth}
\includegraphics[width=\textwidth, trim={0 2cm 0 1.5cm},clip]{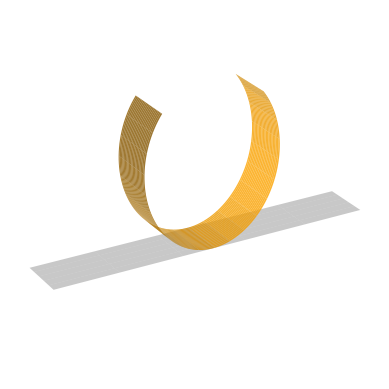}
\subcaption{$\alpha = 0, \kappa \approx 0.58$}
\end{subfigure} \quad
\begin{subfigure}[c]{0.3\textwidth}
\includegraphics[width=\textwidth, trim={0 2cm 0 1.5cm},clip]{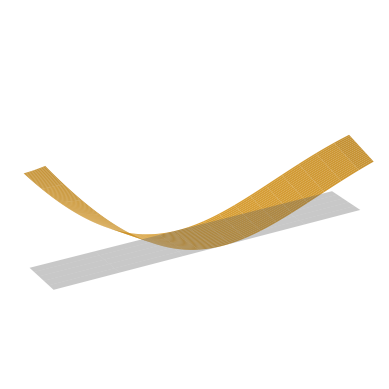}
\subcaption{$\alpha = \frac{\pi}{4}, \kappa \approx 0.48$}
\end{subfigure}
\caption{
  Visualization of the one-parameter family of minimizers $\mathcal S_{Q,B}$ in the case of Lemma~\ref{L:char}~\ref{L:char:c} with parameters $q_1=1,q_2=2,q_{12}=1/2$ and $q_3=\frac12(\sqrt{4q_1q_2}-q_{12})$. Furthermore, we chose $b= A^{-1}q_*$ in order to enforce the condition of Lemma~\ref{L:char}~\ref{L:char:c}. Subfigures (b), (c), (d) show the deformations corresponding to the points on $\mathcal S_{Q,B}$ marked in Subfigure (a).}
\label{fig:parabolic}
\end{figure}

\begin{figure}
    \centering
\begin{subfigure}[c]{0.49\textwidth}
\includegraphics[trim={0 0 0 1.5cm},clip]{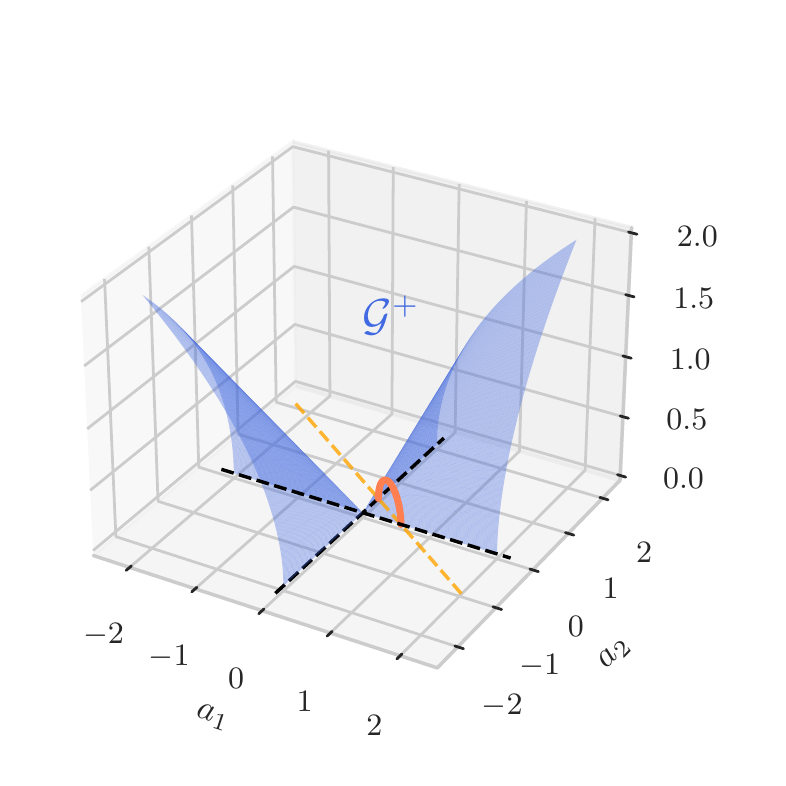} 
\subcaption*{$\mathcal S_{Q,B}\subset \mathcal{G}^+ \subset \R^3$}
\end{subfigure} \hfill \hfill%
\begin{subfigure}[c]{0.49\textwidth}
\includegraphics[trim={0 0 0 1.5cm},clip]{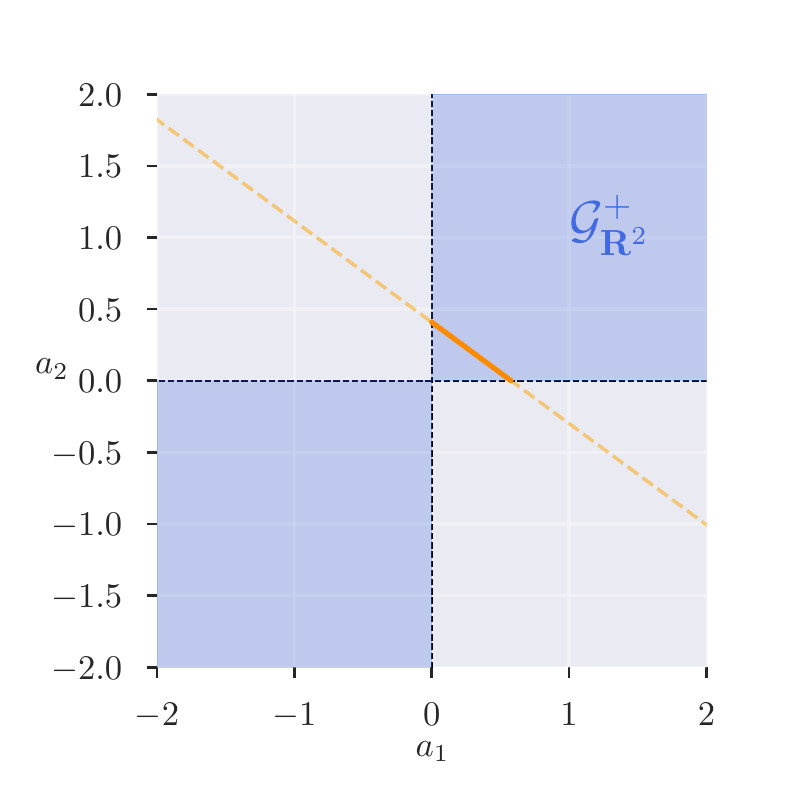}
\subcaption*{$\Phi^{-1}\left(\mathcal S_{Q,B}\right)\subset \mathcal{G}^+_{\R^2} \subset \R^2$}
\end{subfigure}
  \caption{The one-parameter family $\mathcal S_{Q,B}$ (orange) of Figure~\ref{fig:parabolic} as a subset of $\mathcal{G}^+$
    and of $\mathcal{G}^+_{\R^2}$}
  \label{fig:1-parFamG+}
\end{figure}

In Figure~\ref{Contour-Plot} we illustrate the three  Case of Lemma~\ref{L:char} by visualizing the level lines of $\mathcal E_{Q,B}$ and its minimizers in $\mathcal{G}^+_{\R^2}$ for different choices of the coefficients.
\begin{figure}[H]
    \centering
 \captionsetup{type=figure}\addtocounter{figure}{-1}
\begin{subfigure}[c]{0.49\textwidth}
 \ref{L:char:a}\\
\includegraphics[trim={0 0.5cm 0 0.5cm},clip]{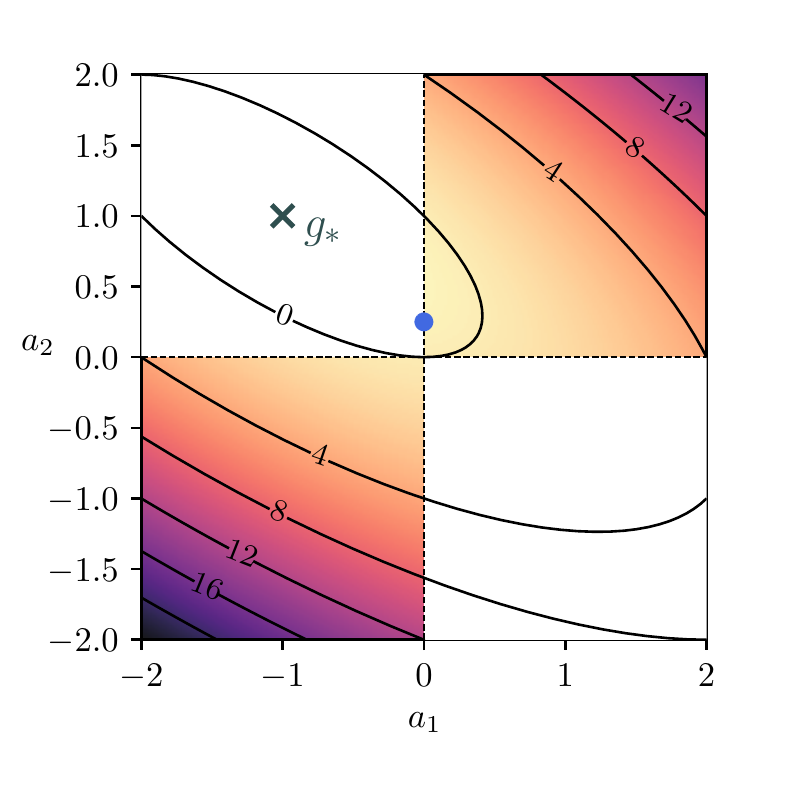}
\end{subfigure} \hfill%
\begin{subfigure}[c]{0.49\textwidth}
 \ref{L:char:b}\\
\includegraphics[trim={0 0.5cm 0 0.5cm},clip]{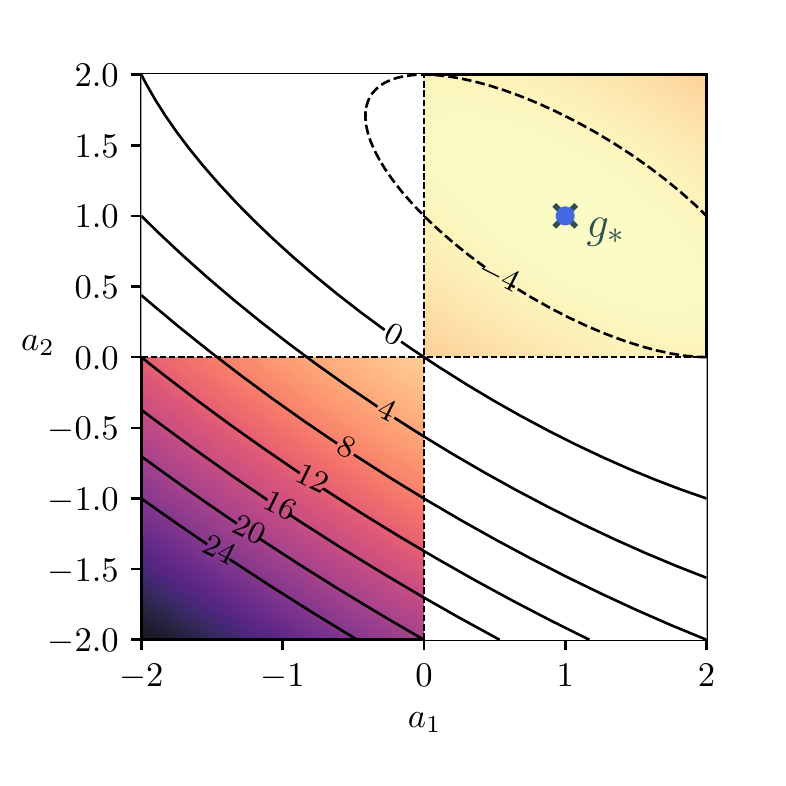}
\end{subfigure}\hfill%
\vspace{-\baselineskip}
\begin{subfigure}[c]{0.49\textwidth}
\ref{L:char:c}\\
\includegraphics[trim={0 0.5cm 0 0.5cm},clip]{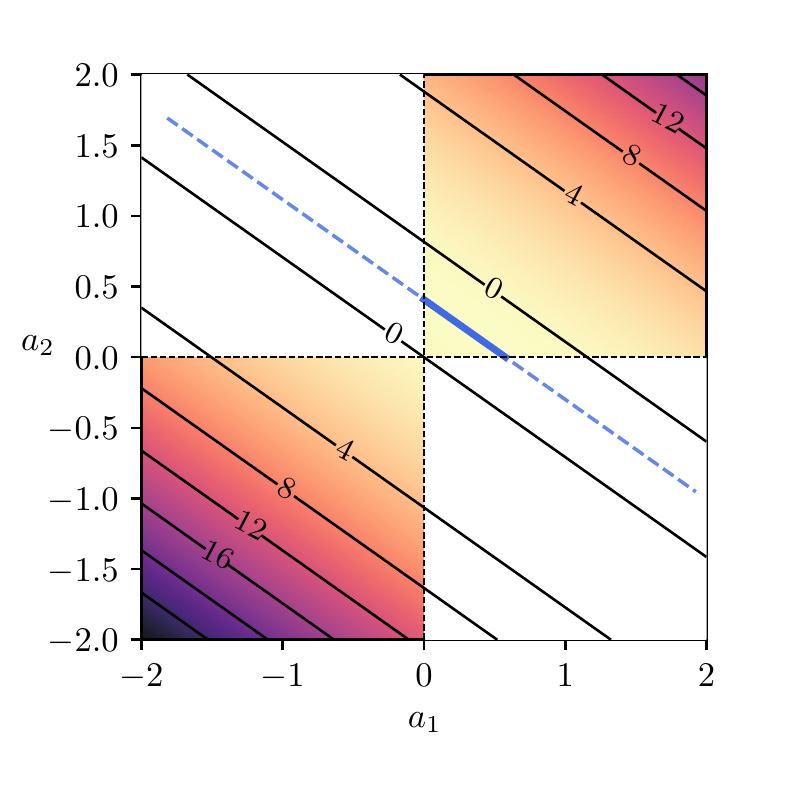}
\end{subfigure} \hfill%
\begin{minipage}[c]{0.48\textwidth}
\centering
    \begin{tabular}[c]{c|c|c|c|c|c|c}
	Case &$q_1$ & $q_2$ & $q_{12} $  & $q_3$  & $B^\gamma_{\eff,1}$ & $B^\gamma_{\eff,2}$ \\
	\hline 
	(a) & 1.0 	& 2.0 & 0 & 1.0	  &  0 & 0.5 \\
	(b) & 1.0 	& 2.0 & 0 & 1.0	  & 2.0 & 1.5\\
	(c) & 1.0   & 2.0 & 0.5 & 1.164 & 0.491 &  0.347    \\
    \end{tabular} \captionof{table}{Effective quantities rounded to 3 decimal places.} 
\captionof{figure}{Visualization of $\mathcal{G}^+_{\R^2}$ and the level lines of $\mathcal E_{Q,B}$ with each subfigure \ref{L:char:a}, \ref{L:char:b} and \ref{L:char:c} corresponding to one Case of Lemma~\ref{L:char}.
The set of minimizers of $\mathcal E_{Q,B}\vert_{\mathcal G^+_{\R^2}}$ is highlighted in blue.}
\label{Contour-Plot}    
 \end{minipage}
\end{figure}
\begin{remark}[Qualitative properties of $\mathcal S_{Q,B}$]\label{R:quali}\mbox{}
  \begin{enumerate}[(a)]
  \item \emph{(Sufficient and necessary conditions)}. A sufficient condition for the axiality
  of minimizers is $\det H = 4q_1q_2 - (q_{12}+2q_3)^2 < 0$, which is a condition only depending on $Q$.
  Likewise, a necessary condition for two non-axial minimizers is $\det H = 4q_1q_2 - (q_{12}+2q_3)^2 > 0$.

  \item \emph{(Stability with respect to perturbations of $Q$ and $B$).} In contrast to the case of a single axial minimizer and the case of Lemma~\ref{L:char}~\ref{L:char:b}, the condition in Lemma~\ref{L:char}~\ref{L:char:c}
    is very sensitive with respect to perturbations of $Q$ and $B$. In applications, we shall evaluate $Q$ and $B$ numerically and thus case~\ref{L:char:c} can typically be neglected.  The same holds for the case of two global axial minimizers.

  \item \emph{(One-parameter family).} Figure~\ref{fig:parabolic} visualizes the one-parameter family $\mathcal S_{Q,B}$ in the case of Lemma~\ref{L:char}~\ref{L:char:c}.   In Figure~\ref{fig:1-parFamG+} that very same one-parameter family is visualized as a subset of $\mathcal{G}^+$
    and of $\mathcal{G}^+_{\R^2}$.
  One can show that angles $\alpha(G)$ associated with $G\in\mathcal S_{Q,B}$ span the whole intervall $(-\pi/2,\pi/2]$. Furthermore, one can show that $\mathcal S_{Q,B}$ is a one-dimensional compact manifold.

  \item \emph{(Continuous dependence).} Lemma~\ref{L:char} reveals that minimizers may not depend
  continuously on the prestrain $B$. For instance, even in the stable case~\ref{L:char:a}, the global minimizer may jump from an axial minimizer with angle $\alpha=0$ to one with angle $\alpha=\frac\pi2$. On the other hand, in case~\ref{L:char:b}, the two global minimizers continuously depend on $b\in \frac12 A^{-1}H(\mathcal{G}^+_{\R^2}\setminus\partial \mathcal{G}^+_{\R^2})$.

  \item \emph{(Local minimizers).} The same techniques can be used to handle local minimizers.
In the case of Lemma~\ref{L:char} (a) both axial matrices on the right-hand side of \eqref{eq:form_min} turn out to be local minimizers in the case $\det H >0$. In the case of Lemma~\ref{L:char}~\ref{L:char:b} the two global minimizers are also the only local minimizers, and in the case of Lemma~\ref{L:char}~\ref{L:char:c}, the one-parameter families
  \end{enumerate}
\end{remark}

\begin{remark}[The case $q_{12}=0$ and the isotropic case]\label{R:simply}
  In the case $q_{12}=0$ we may simplify the statement of Lemma~\ref{L:char} further.
  In particular, $A$ is diagonal and the following holds:
  \begin{enumerate}[(a)]
  \item In the case of Lemma~\ref{L:char}~\ref{L:char:a} we have
    \begin{equation*}
      \mathcal S_{Q,B}\subseteq\{\widehat B_1G_1,\widehat B_2G_2\}
    \end{equation*}
  \item Formula \eqref{def:gstar} for $g_*$ simplifies to
    \begin{equation*}
      g_*=  \bigg(\frac{q_{1}q_{2}\widehat B_1 - q_{3}q_{2}\widehat B_2}{q_{1}q_{2} - q_{3}^{2}} ,
        \frac{q_{1}q_{2}\widehat B_2 - q_{3}q_{1}\widehat B_1}{q_{1}q_{2} - q_{3}^{2}} \bigg).
    \end{equation*}
  \end{enumerate}
  Furthermore, we note that if $Q$ is isotropic, i.e., $q_1=q_2=q_3$ and $q_{12}=0$,
  and $B$ is a multiple of the identity, we are always in the case of Lemma~\ref{L:char}~\ref{L:char:c} and $\mathcal S$ consists of all matrices of the form $G=\sqrt{\det B}(\minvec\otimes \minvec)$ with $\minvec\in \R^2$, $\abs{\minvec}=1$.
\end{remark}


\subsection{Microstructure--shape relation for the parametrized laminate}\label{S:param2}
We continue our study of the parametrized laminate considered in Section~\ref{S:ex:2} (shown in Figure~\ref{sketchofmicrostructure}), and now focus on the microstructure--shape relation. We  thus consider the algebraic minimization problem~\eqref{L:ex:1:eq1} with $Q^\gamma_{\hom}$ and $B^\gamma_{\eff}$ defined as in Lemma~\ref{S:ex:C3}, and study the dependence of the set of minimizers $\mathcal S$ on the parameters listed in Table~\ref{parList}. Throughout this section we set $\mu_1=\rho_1=1$; the set $\mathcal S$ for other values can then be easily obtained with the help of \eqref{eq:QBinvariance}.

\paragraph{Visualization of $\mathcal S$.} Before we start our exploration, we briefly comment on how we visualize the set $\mathcal S$. As we shall see, except for the special case of a homogeneous composite or for parameters $(\theta,\theta_\mu,\theta_\rho)$ in a small exceptional set, $\mathcal S$ consists of
\begin{equation}\label{eq:casegstar}
  \text{one unique axial minimizer \quad or \quad two non-axial minimizers}.
\end{equation}
Hence, in view of \eqref{eq:S:P1} there exists a unique $G_*\in\mathcal G^+$ such that $\mathcal S=\{G_*,TG_*\}$; note that $G_*=TG_*$ if and only if $G_*$ is axial.
This allows to visualize $\mathcal S$ by visualizing $G_*$. We use the angle--curvature parametrization introduced in Remark~\ref{R:kappaalpha}, i.e., we shall plot the angle $\alpha(G_*)$ and the curvature $\kappa(G_*)$ as functions of the parameters under consideration. 
In fact, the exceptional set of parameters $(\theta_\rho,\theta_\mu,\theta)$ that violate condition \eqref{eq:casegstar} are precisely those that lead to minimizers of Case~\ref{L:char:c} of Lemma~\ref{L:char}. In view of Remark~\ref{R:quali}, these exceptional parameters belong to a set of zero Lebesgue measure and thus can be neglected in the following presentation.
\begin{remark}[The homogeneous case: $\theta\in\{0,1\}$ or $\theta_\mu = 1$]\label{R:homcase}
  When $\theta\in\{0,1\}$ or $\theta_\mu=1$, Condition~\eqref{eq:casegstar} is not satisfied, and thus $\mathcal S$ cannot be visualized in terms of the angle--curvature parametrization. 
  %
   In the case $\theta\in\{0,1\}$, both the composite and the prestrain are homogeneous. In the case $\theta_\mu=1$ only the composite is  homogeneous. In both cases we conclude that
  \begin{equation*}
    q_1=q_2=q_3,\qquad\text{and}\qquad \widehat B_{\eff,1}^{\gamma}=\widehat B_{\eff,2}^{\gamma},
  \end{equation*}
  and we are thus in the isotropic case of Remark~\ref{R:simply}. In particular, we deduce that
  $\mathcal S = \{\kappa \minvec\otimes \minvec \,:\,\minvec\in \R^2,\;\abs{\minvec}=1\}$ with 
  \begin{equation*}
    \kappa \colonequals
    \begin{cases}
      \frac{3 }{2}\rho_{1}&\text{if }\theta=0,\\
      \frac{3}{2} \theta_{\rho} \rho_{1}&\text{if }\theta=1,\\
      \frac{3}{2} \rho_{1} \left( 1 - \theta (1+\theta_{\rho}\right)&\text{if }\theta_\mu=1.
    \end{cases}
  \end{equation*}
  \end{remark}

In the following discussion, we exclude the homogeneous case and focus on the dependence of $\mathcal S$ on parameters $(\theta_\rho,\theta_\mu,\theta)$ in the set
\begin{equation}\label{def:P}
\mathcal{P} \colonequals \left\{
\begin{pmatrix}
  \theta \\ \theta_\mu \\ \theta_\rho
\end{pmatrix} 
:  \theta_\mu > 1,\ \theta_\mu \neq 1 \text{ and } \theta \in (0,1) 
\right\} \subset \R^3.
\end{equation}
We note that parameters $(\theta,\theta_\mu,\theta_\rho)$ with $\theta_\mu\in(0,1)$ are also covered by the following discussion: In view of the symmetry of the compopsite, the associated set $\mathcal S$ is obtained by considering the parameters $(1-\theta,\tfrac1{\theta_\mu},\tfrac1{\theta_\rho})\in\mathcal P$.

%

\begin{figure}[H]
  \centering
  
 \tikzstyle{base} = [draw, thick, rounded corners, align=center]
 \tikzstyle{decision} = [base, diamond, fill=yellow!20, text width=7em, inner sep=0pt,aspect = 2]
\tikzstyle{startstop} = [rectangle, rounded corners, minimum width=3cm, minimum height=1cm ,text centered, draw=black, fill=orange!30,text width=5cm]
\tikzstyle{input} = [rectangle, minimum width=3cm, minimum height=1cm, text centered, text width=3cm, draw=black, fill=green!20!gray]
\tikzstyle{io} = [rectangle, minimum width=3cm, minimum height=1cm, text centered, text width=3cm, draw=black, fill=blue!20]
\tikzstyle{process} = [rectangle, minimum width=3cm, minimum height=1cm, text centered, text width=3cm, draw=black, fill=orange!30]
\tikzstyle{arrow} = [draw, -{Latex}]

\resizebox{\textwidth}{!}{%
\begin{tikzpicture}[%
    >=stealth,
    node distance=3.5cm,
    on grid,
    auto,
    scale=\textwidth,
    font=\normalsize,thick
  ] 

\node (in1) [input,text width=4.5cm, xshift=0cm, yshift=0cm]{\textbf{Input}: $\theta,\ \theta_\mu,\ \theta_\rho,\ \mu_1,\ \lambda_1,\ \rho_1,\ \gamma$};

\node (start) [startstop, right of=in1, text width=14cm,xshift=9cm,yshift=-2cm,align=left] {Compute via Lemma~\ref{S:ex:C3}: 
\begin{align*}
  &q_1 = \frac{1}{6} \langle \mu \rangle_h,\quad  q_2 =  \frac{1}{6} \Bar{\mu},\quad q_{12}=0,\\
 &q_3 =
   \begin{cases}
     q_2&\text{in the case }0<\gamma\ll 1,\\
     q_1&\text{in the case }\gamma\gg 1,\\
    \frac16\mu_\gamma&\text{else  (this involves solving \eqref{muGamma2}),}    
  \end{cases}\\
  &b = (\widehat B^\gamma_{\eff,1},\widehat B^\gamma_{\eff,2})^\top.
\end{align*}%
Assemble $H$ and $A$ according to \eqref{eq:reformulation}.
};

\node (dec1) [decision, left of=start,xshift=-7cm, yshift=-1cm] {$\operatorname{det} H > 0$?};
\node (pro1) [process, left of=dec1, xshift=-2cm]{Compute $g_* = 2H^{-1}Ab$};

\node (dec4) [decision, below of=pro1, yshift=-0.5cm] {$g_* \in \mathcal{G}^+_{\R^2}\setminus \partial \mathcal{G}^+_{\R^2}$?};
\node (pro4) [io, below of=dec4, yshift=-1.5cm, text width=5cm] {$ \mathcal S_{Q,B}=\{G,TG\},\  G=\Phi(g_*)$ };

\node (dec2) [decision, below of=dec1, yshift=-0.5cm] {$\operatorname{det} H = 0$?};
\node (pro2b) [process, right of=dec2, xshift=2.5cm] {Compute $q_*, s_*$ according to \eqref{L:char:c:quantities}};

\node (dec3) [decision, right of=pro2b, xshift=4cm] {$(Ab \cdot q_*)q_*= Ab$?};
\node (pro3) [io, below of=dec3,xshift=0.1cm, yshift=-1.5cm, text width=5cm ] {$\mathcal S_{Q,B}=\{G,TG\,:\,G=\Phi(a)$ \\ $\text{ with }a\in \mathcal{G}^+_{\R^2}\text{ s.t. }a\cdot q_*=s_*\}$};

\node (dec5) [decision, below of=dec2, yshift=-0.5cm, xshift=2cm, text width=3.5cm,  align=center] {$\frac{(\widehat B^\gamma_{\eff,2})^2}{q_1}=\frac{(\widehat B^\gamma_{\eff,1})^2}{q_2}$?};
\node (pro5c) [io, right of=dec5, yshift=-0.1cm, xshift=2.75cm, text width=3.8cm] {$\mathcal S_{Q,B} = \{\widehat B_{1}G_1, \widehat B_{2}G_2 \}$ };

\node (dec6) [decision, below of=dec5, yshift=-0.5cm, xshift=0cm, text width=4cm,  align=center] {$Q_{\hom}^\gamma( \widehat B_{\eff,1}^{\gamma}G_1{-}B_{\rm eff}^\gamma) $ \\ $>Q_{\hom}^\gamma(\widehat B_{\eff,2}^{\gamma}G_2{-}B_{\rm eff}^\gamma)?$};

\node (pro6a) [io, right of=dec6, yshift=-0.1cm, xshift=4cm, text width=3.5cm] {$\mathcal S_{Q,B} = \{\widehat B_{2}G_2 \}$ };
\node (pro6b) [io, left of=dec6, yshift=-0.1cm,xshift=-4cm, text width=3.5cm] {$\mathcal S_{Q,B} = \{\widehat B_{1}G_1 \}$ };
\draw [arrow] (in1) -- (0.0119,0); 
\draw [arrow] (0.0119,-0.0065) -- (dec1);
\draw [arrow] (dec1) -- node[pos=0.25,anchor=east,fill=white] {Yes} (pro1);
\draw [arrow] (dec1) -- node[pos=0.5,anchor=south,fill=white] {No}(dec2);
\draw [arrow] (pro1) -- (dec4);
\draw [arrow] (dec4) -- node[pos=0.5,anchor=south,fill=white] {Yes} (pro4);
\draw [arrow] (dec2.east) -- node[pos=0.25,anchor=west,fill=white] {Yes} (pro2b);
\draw [arrow] (pro2b) -- (dec3);
\draw [arrow] (dec4) -- node[pos=0.25,anchor=west,fill=white] {No} ++(0.007,0) |- (dec5.west);
\draw [arrow] (dec2) -- node[pos=0.75,anchor=south,fill=white] {No} ++(0.0,-0.0042) -| (dec5);
\draw [arrow] (dec5.east) -- node[pos=0.2,anchor=west,fill=white] {Yes} (pro5c);
\draw [arrow] (dec5) -- node[pos=0.75,anchor=south,fill=white] {No} (dec6);
\draw [arrow] (dec3.south)  --node[pos=0.85,anchor=south,fill=white] {Yes} (pro3.north);
\draw [arrow] (dec3.south) |- node[pos=0.75,anchor=east,fill=white] {No} (dec5.north);
\draw [arrow] (dec6.east)  -- node[pos=0.25,anchor=west,fill=white] {Yes}(pro6a.west);
\draw [arrow] (dec6.west) -- node[pos=0.25,anchor=east,fill=white] {No} (pro6b.east);
\end{tikzpicture} }  
  \caption{The algorithm for evaluating the set of minimizers $\mathcal S$ of the algebraic minimization problem \eqref{L:ex:1:eq1} in the case of the parametrized laminate of Lemma~\ref{S:ex:C3}. Note that in the cases $0<\gamma\ll1$ and $\gamma\gg 1$ we use the approximation $q_3=q_2$ and $q_3=q_1$, respectively, see~\eqref{S:ex:muGammaProp2} for the justification.}
  \label{fig:flowchart}
\end{figure}
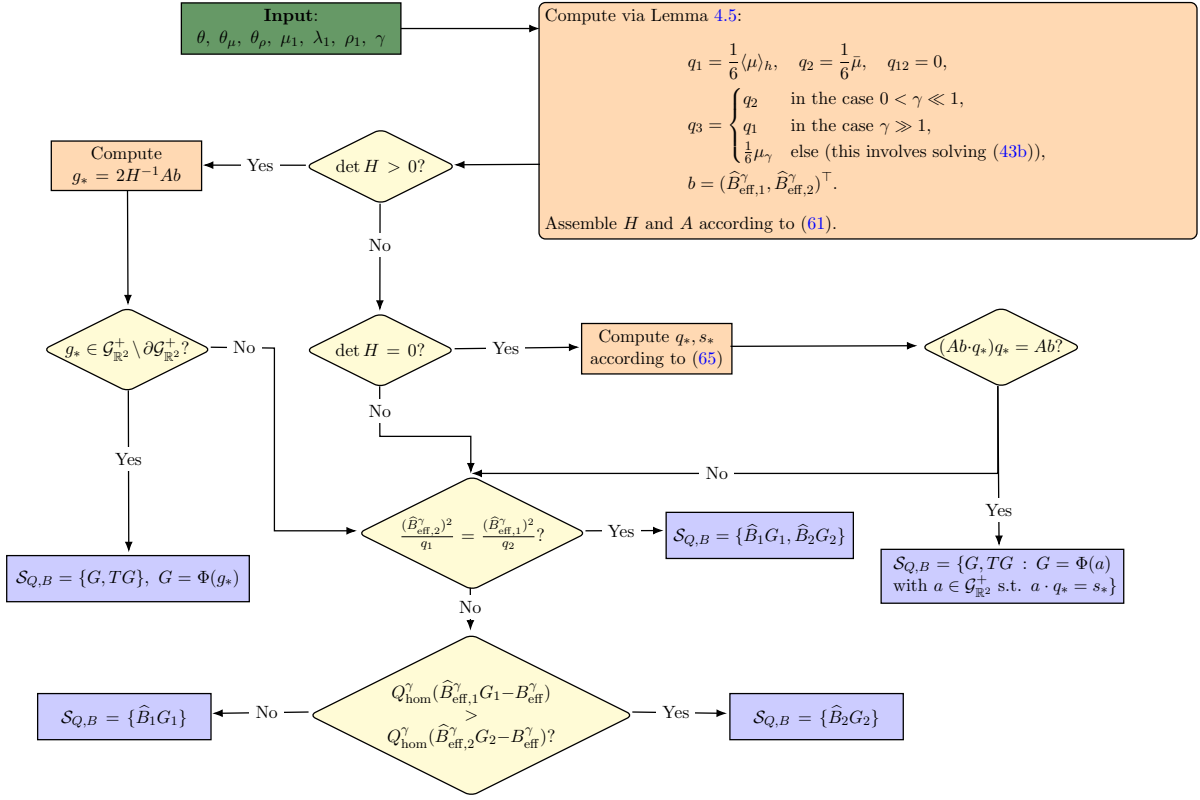

\paragraph{Computation of $\mathcal S$.}
The computation of 
the set of minimizers $\mathcal S$ associated with $Q^\gamma_{\hom}$ and $B^\gamma_{\eff}$ for prescribed parameters $(\theta, \theta_\mu,\theta_\rho)\in\mathcal P$  uses Lemma~\ref{L:char} and the simplifications in Remark~\ref{R:simply}.
The resulting algorithm is summarized as a flow-chart in Figure~\ref{fig:flowchart}. Recall that we denote by $q_1,q_2,q_{12}, q_3$ the coefficients of $Q^\gamma_{\hom}$ as in Definition~\ref{def:orthotropicity}, and by $\widehat B^\gamma_{\eff,1},\widehat B^\gamma_{\eff,2},\widehat B^\gamma_{\eff,3}$ the coefficients of $B^\gamma_{\eff}$ with respect to the basis in \eqref{canonicalBasis}.

\paragraph{Dependence of $\mathcal S$ on $\gamma$ and the critical value $\gamma_*$.}
According to Lemma~\ref{S:ex:ex1}, $B^\gamma_\text{eff}$ is independent of $\gamma$,
and the only coefficient of $Q^\gamma_{\hom}$ that depends on $\gamma$ is $q_3=q_3(\gamma)$.
By~\eqref{S:ex:muGammaProp1} we have $q_1< q_3(\gamma)< q_2$,
 and for $(\theta,\theta_\mu,\theta_\rho)\in\mathcal P$ the map,  $\gamma\mapsto q_3(\gamma)$
is continuous and strictly monotonically decreasing. Thus, in view of \eqref{S:ex:muGammaProp2}, we can find a unique value $\gamma_*\in(0,\infty)$ such that
\begin{equation*}
  q_3(\gamma_*)=q_{3}^{*} \colonequals \sqrt{q_{1} q_{2}}.
\end{equation*}
We call $\gamma_*$ the \emph{critical value} for the following reason: For $\gamma<\gamma_*$ we have $\det H<0$ and we are thus in Case~\ref{L:char:a} of Lemma~\ref{L:char}, which, in particular, means that all minimizers are axial. On the other hand, for $\gamma>\gamma_*$ we have $\det H>0$, which is a necessary condition for the case~\ref{L:char:b}. We conclude that non-axial minimizers may only  emerge for $\gamma>\gamma_*$. The critical condition $\gamma=\gamma_*$ is equivalent to $\det H=0$, and therefore necessary for the case Lemma~\ref{L:char}~\ref{L:char:c}, which is the only case where a one-parameter family  of minimizers may emerge.
\begin{figure}
  \centering
    \includegraphics[trim={0 0 0 1cm},clip]{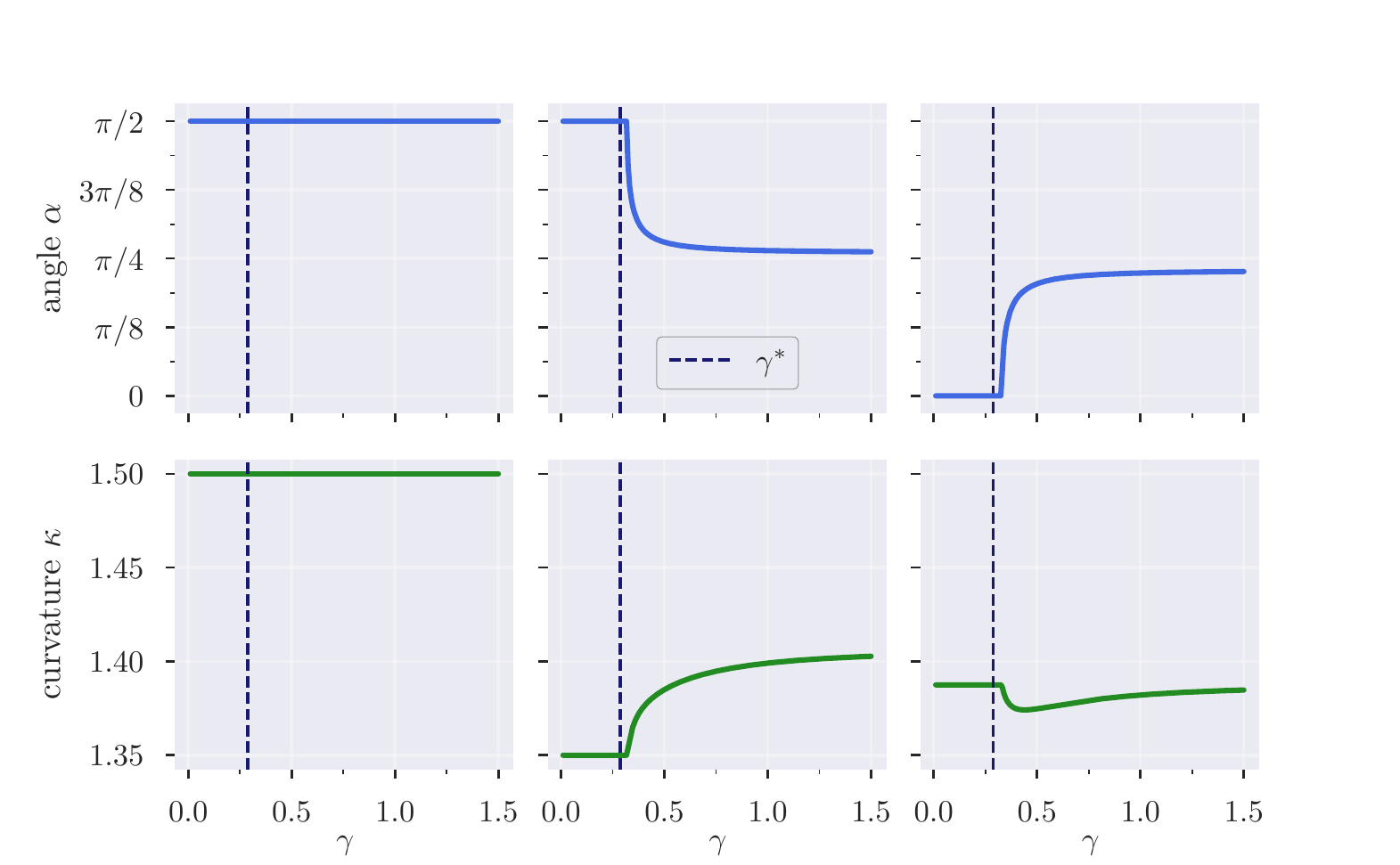}
     \caption{Angle $\alpha(G_*)$ (top) and curvature $\kappa(G_*)$ (bottom) for different values of $\gamma$ with $\theta_\mu=2$, $\theta=\frac14$. The three columns correspond to $\theta_\rho = -1.0$ (left), $\theta_\rho = -0.75$ (middle), and $\theta_\rho = -0.7$ (right). The dotted line shows the critical value $\gamma_*$}
    \label{anglevsq3}  
\end{figure} 

In order to demonstrate the dependence of $\mathcal S$ on the relative scaling parameter $\gamma$
we show in Figure~\ref{anglevsq3} the angle $\alpha(G_*)$ and the curvature $\kappa(G_*)$ as a function of $\gamma$.
The three columns of Figure~\ref{anglevsq3} feature different values of the prestrain ratio~$\theta_\rho$,
each resulting in qualitatively different behavior: The angle $\alpha(G_*)$ and the curvature $\kappa(G_*)$  may increase, decrease or remain constant in $\gamma$. Furthermore, Figure~\ref{anglevsq3} also includes a vertical dotted line to indicated the value of $\gamma^*$. Recall that $\gamma>\gamma^*$ is a necessary (but not sufficient) condition for the existence of non-axial minimizers, while $\gamma < \gamma^*$ is a sufficient condition for axiality of minimizers.
The curvature  $\kappa(G_*)$ remains constant if $\gamma< \gamma^*$.
Finally, note that both the angle $\alpha(G_*)$ and the curvature $\kappa(G_*)$ appear to
depend continuously on $\gamma$.%
\begin{figure}
  \begin{subfigure}[c]{0.49\textwidth}
    \includegraphics[width=\textwidth, trim={1cm 0cm 1cm 7cm },clip]{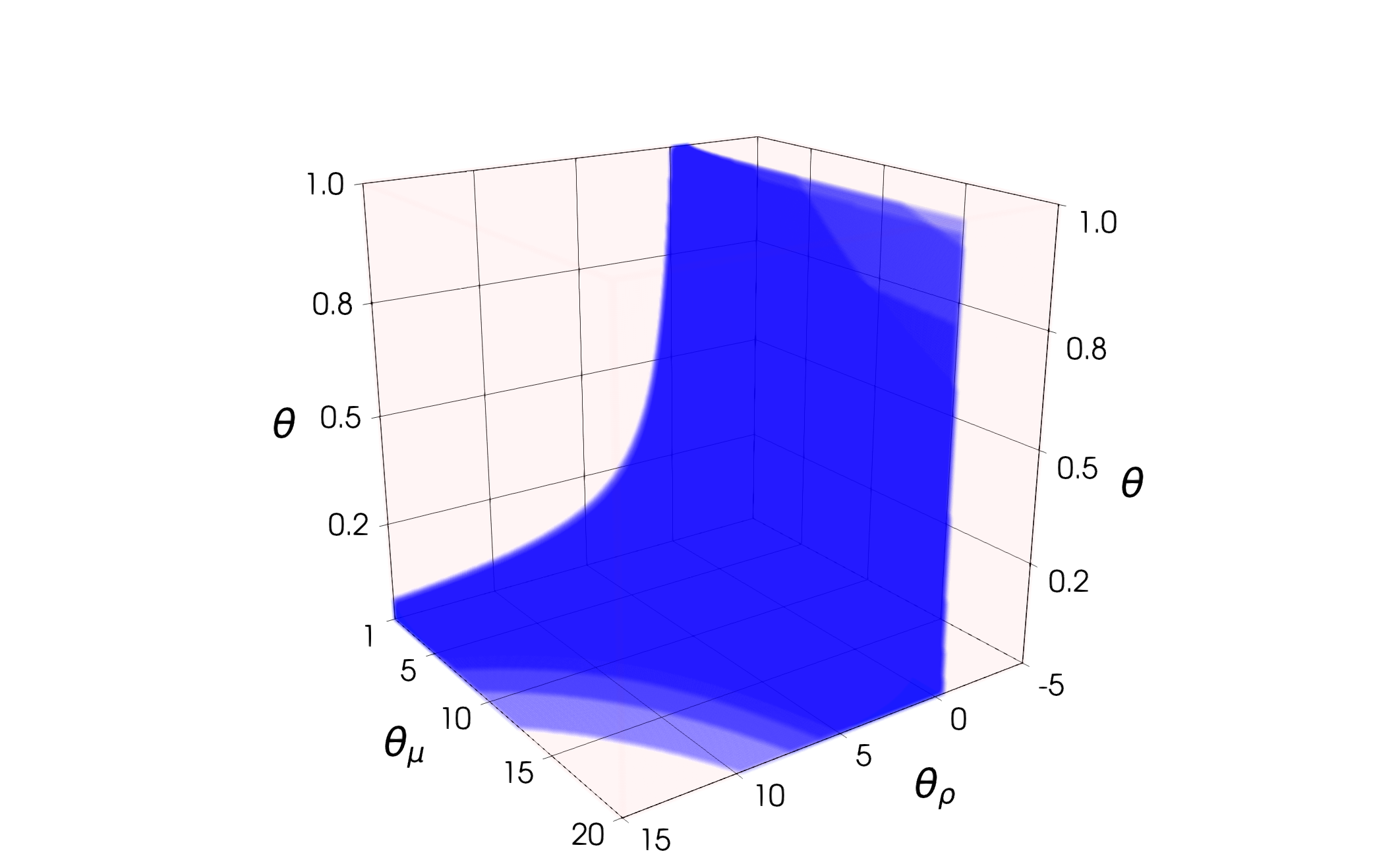}%
  \end{subfigure}%
  \begin{subfigure}[c]{0.49\textwidth}
    \includegraphics[width=\textwidth, trim={1cm 0cm 1cm 7cm },clip]{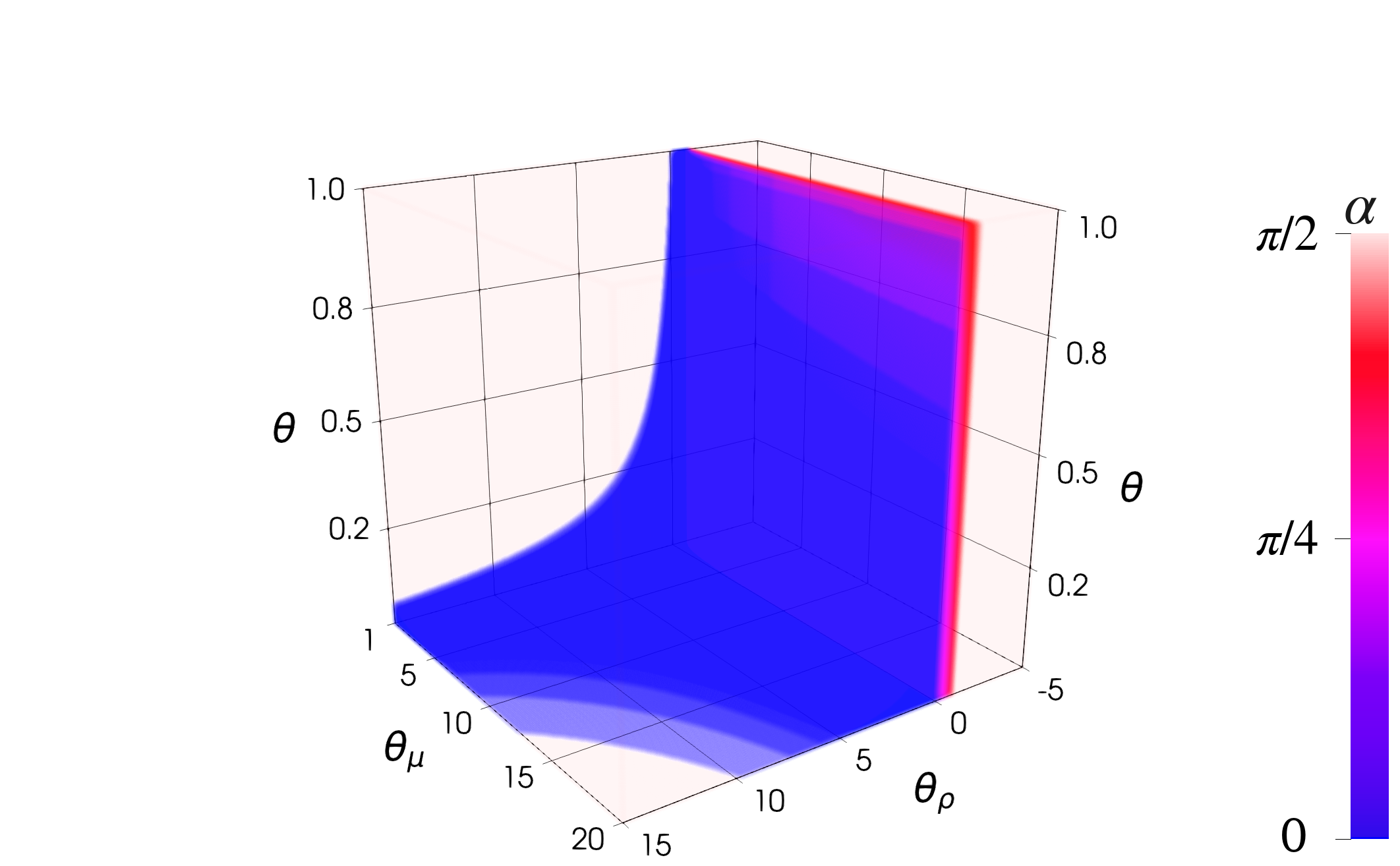}
  \end{subfigure}
  \begin{subfigure}[c]{0.49\textwidth}
    \includegraphics[width=\textwidth, trim={1cm 0cm 1cm 2cm },clip]{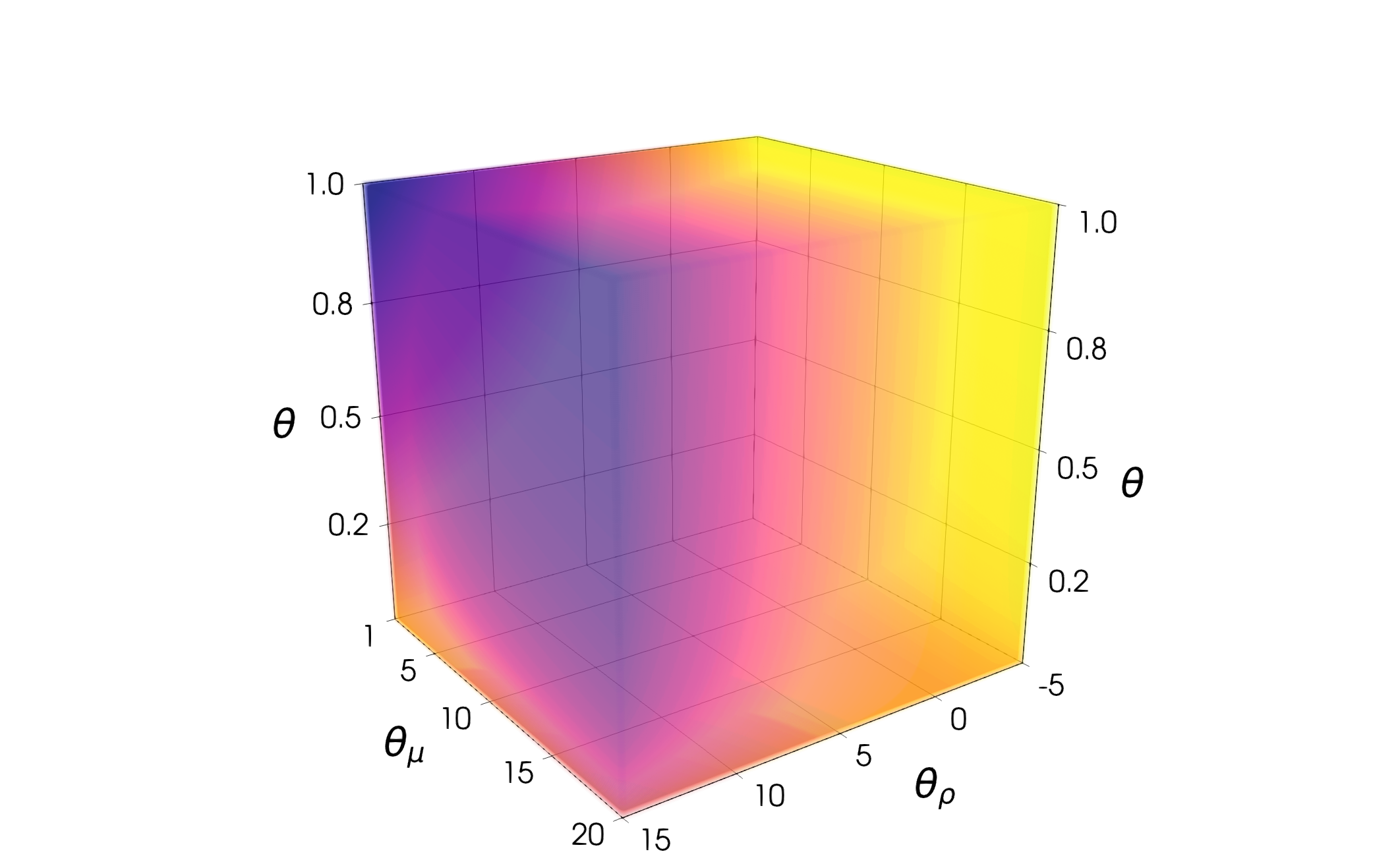}%
    \subcaption*{$0< \gamma \ll 1$ }
  \end{subfigure}%
  \begin{subfigure}[c]{0.49\textwidth}
    \includegraphics[width=\textwidth, trim={1cm 0cm 1cm 2cm },clip]{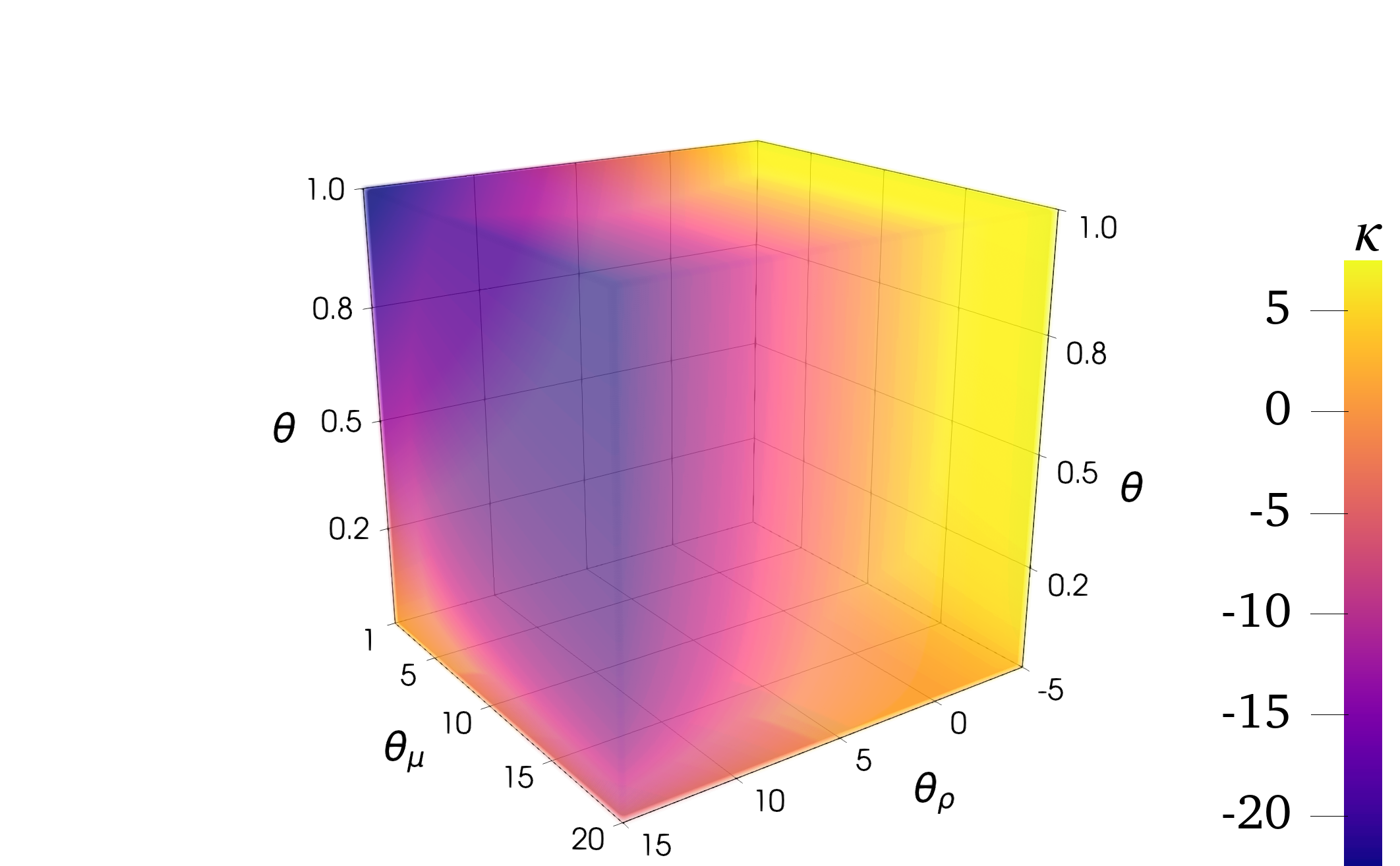}
    \subcaption*{$ \gamma \gg 1 $ }
  \end{subfigure}
    \caption{The angle $\alpha(G_*)$ (top) and curvature $\kappa(G_*)$ (bottom) as a functions of  $\theta_{\rho}$,\  $\theta_{\mu}$ and $\theta$ in the extreme regimes $0< \gamma \ll 1$ (left) and $ \gamma \gg 1 $ (right). In the plots for the angle, transparent regions correspond to the angle $\frac{\pi}{2}$.}
\label{3D-PhaseDiagram1}
\end{figure}
\paragraph{Dependence of $\mathcal S$ on $\theta,\theta_\mu$ and $\theta_\rho$.}
We investigate the dependence of $\mathcal S$ on the parameters $(\theta,\theta_\mu,\theta_\rho)\in\mathcal P$ (see \eqref{def:P}).
Since we shall consider a high number of sample points in the parameter set $\mathcal P$, we mostly focus
on the extreme regimes $0<\gamma \ll 1$ and $ \gamma \gg 1 $.
In these regimes, we can use the approximations $q_{3}=q_{2}$ for $0<\gamma\ll 1$ and $q_{3}=q_{1}$ for $\gamma\gg 1$ 
in order to avoid the high computational cost of solving \eqref{muGamma2} numerically, which would be required to evaluate $q_3(\gamma)$ for $\gamma\in(0,\infty)$.
Note that the regimes $0<\gamma \ll 1$ and $ 1 \ll \gamma$ illustrate the possible range of behaviors, since the behavior for $0 < \gamma < \infty$ is a middle ground between these extreme cases.
We illustrate this at the end of this section.
Figure~\ref{3D-PhaseDiagram1} visualizes $G_*$ in the extreme regimes $0<\gamma\ll 1$ and $\gamma\gg 1$
for $300^3$ uniformly spaced sample points in the parameter set $\mathcal{P}$. The transparent and blue regions in these plots corresponds to axial minimizers with $\alpha(G_*)=\tfrac{\pi}{2}$ and $\alpha(G_*)=0$, respectively.
The region with $\alpha(G_*)=\tfrac{\pi}{2}$ seems to consist of two connected components that are divided by a region where $\alpha(G_*)$ takes other values.
We refer to that later region as the \emph{transition region}.
\medskip

\begin{figure}[H]
  \centering
  \includegraphics[trim={0cm 0cm 0cm 0.25cm },clip]{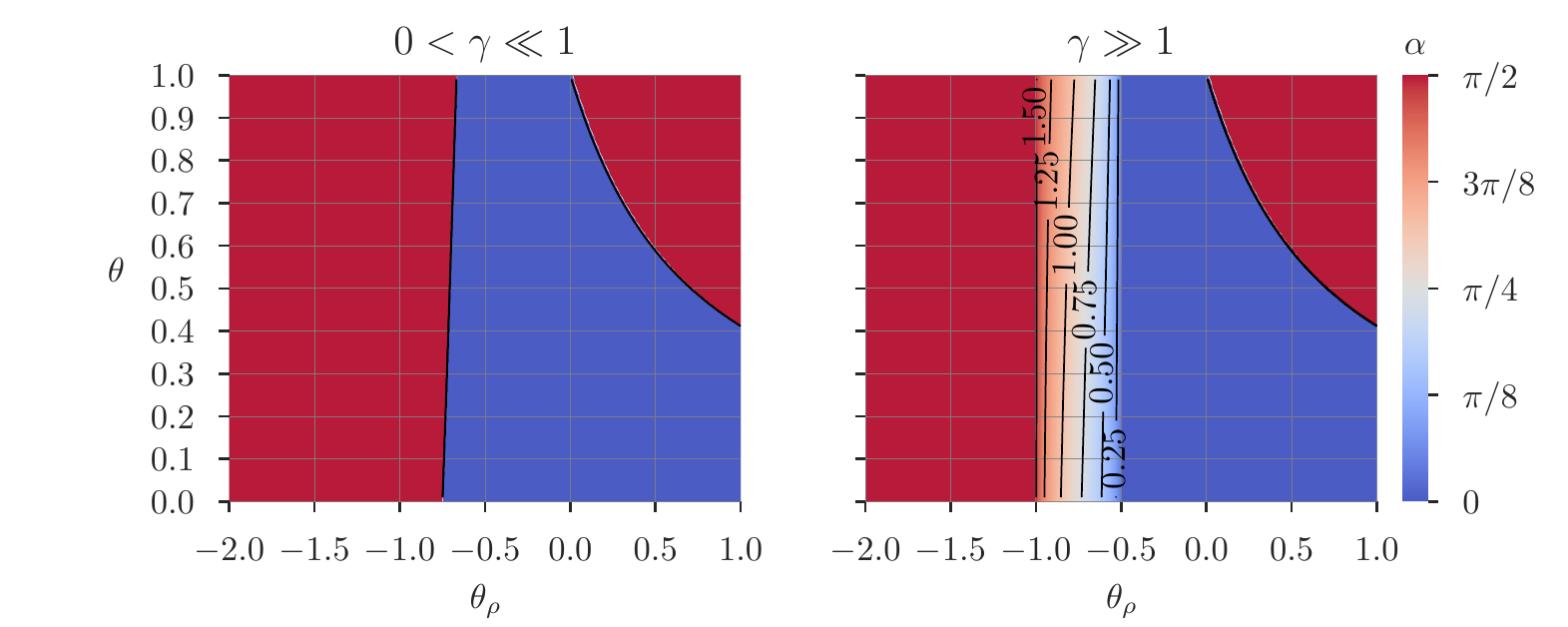}
  \includegraphics{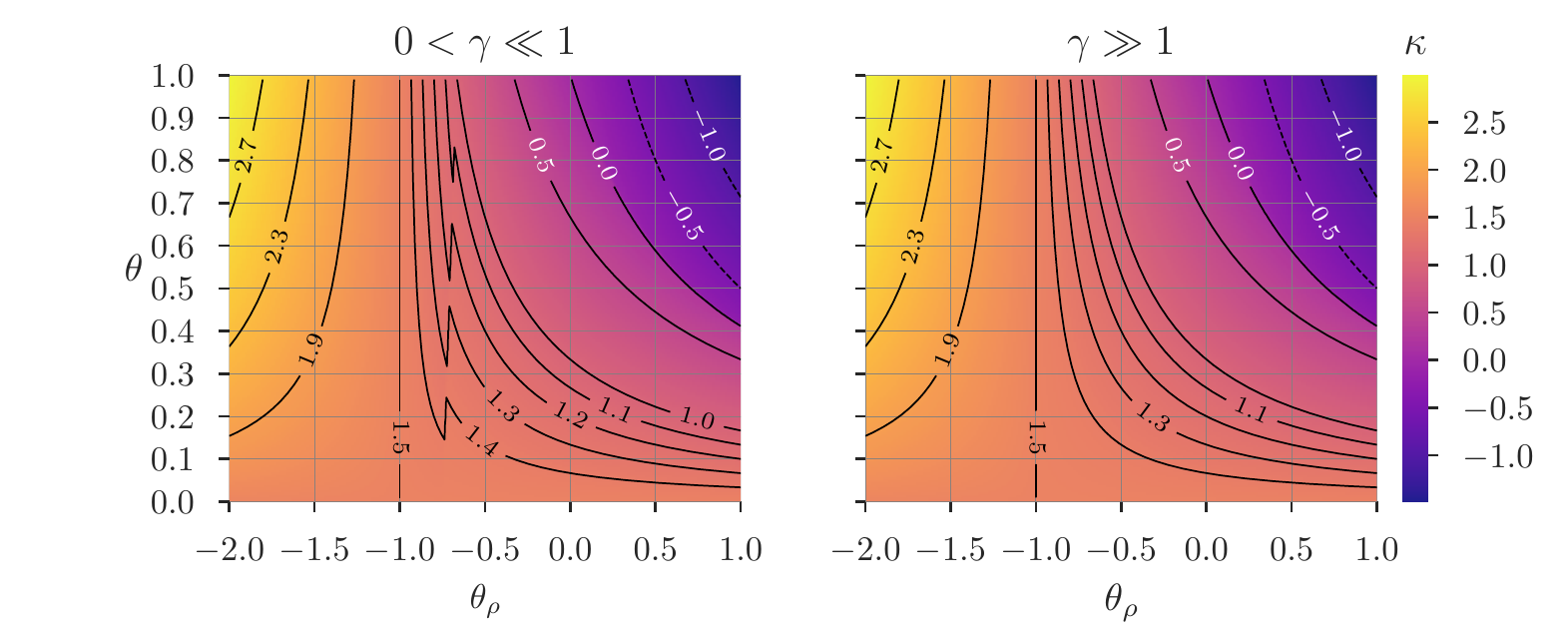}
  \caption{The angle $\alpha(G_*)$ (top) and curvature $\kappa(G_*)$ (bottom) in the case $0<\gamma\ll1$ (left) and $\gamma\gg1$ (right) for points in the parameter set $\mathcal P$, with $\theta_\mu = 2$}
  \label{2DPhaseDiagram}
\end{figure}

\begin{figure}[H]
    \centering
   \includegraphics{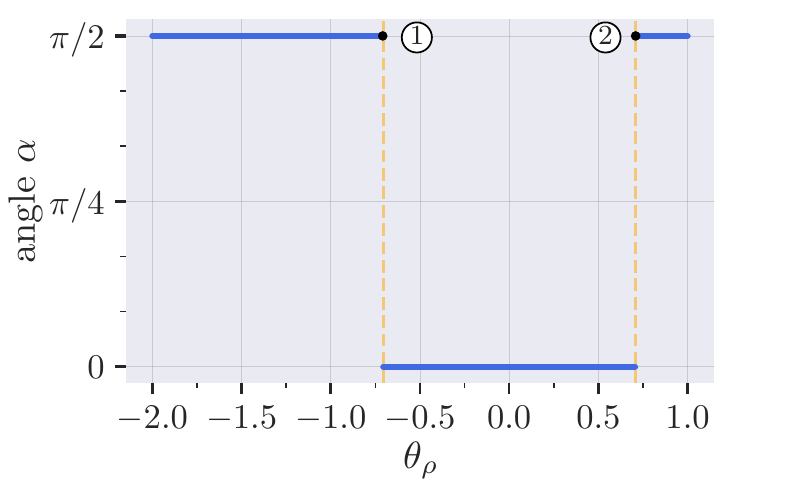}%
  \includegraphics{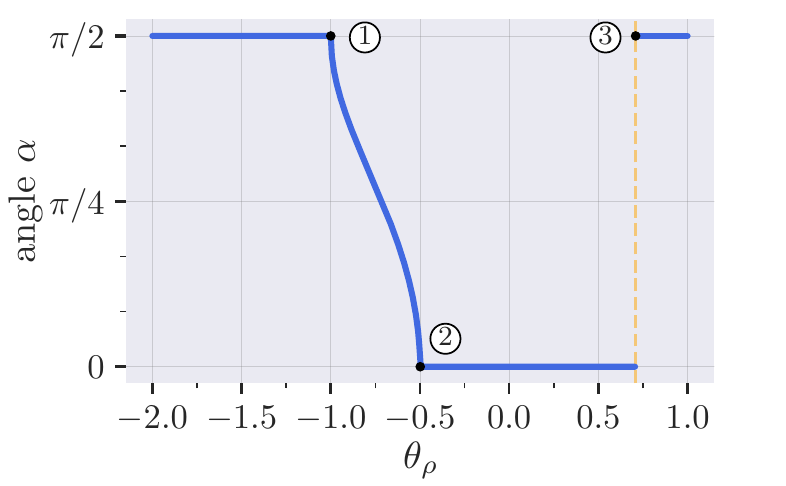}
  \begin{subfigure}[c]{0.49\textwidth}
    \includegraphics{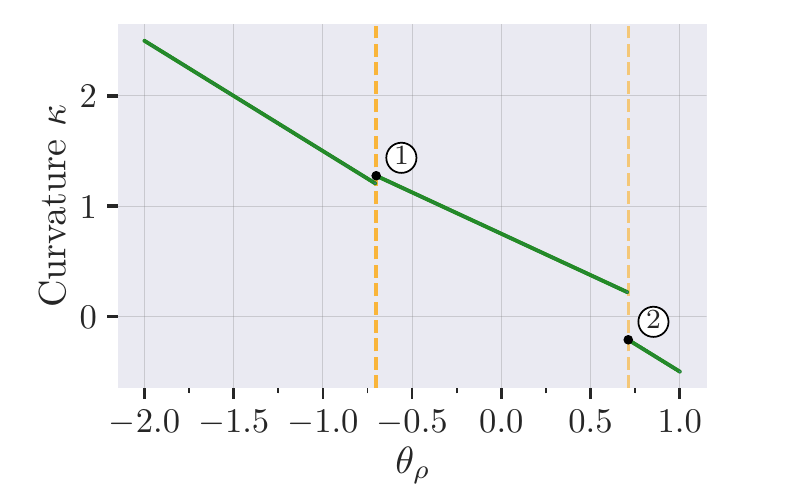}%
    \subcaption{$0< \gamma \ll 1$ }
  \end{subfigure}
  \begin{subfigure}[c]{0.49\textwidth}
    \includegraphics{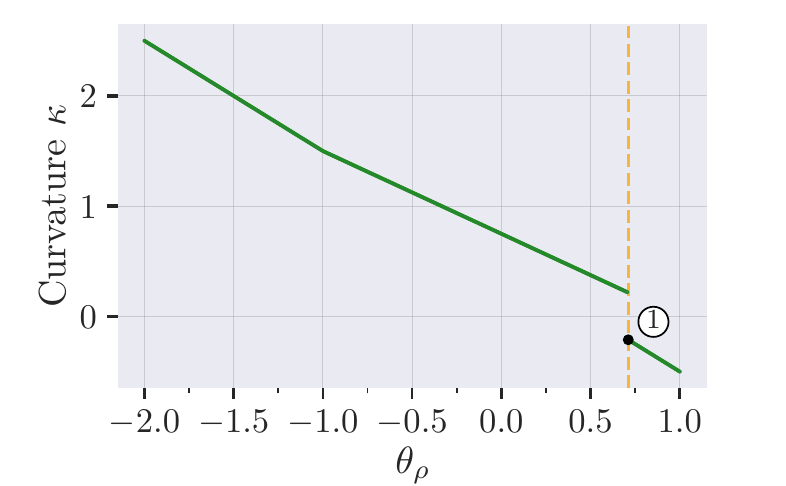}
    \subcaption{$\gamma \gg 1$ }
  \end{subfigure}
\centering
    \begin{tabular}[b]{c|c |c|c}
    Parameters &    $\theta_{\mu}$ & $\theta $     \\ \hline
  			   &   		2			&  $0.5$ \\
    \end{tabular}
    \caption{Angle $\alpha(G_*)$ (top) and curvature $\kappa(G_*)$ (bottom) as a function
     of the prestrain ratio~$\theta_\rho$ with $0< \gamma \ll 1$ (left) and $ \gamma \gg 1 $ (right). }
    \label{anglevstheta}
  \end{figure}

\begin{figure}[H]
\begin{subfigure}[c]{0.49\textwidth}
\includegraphics{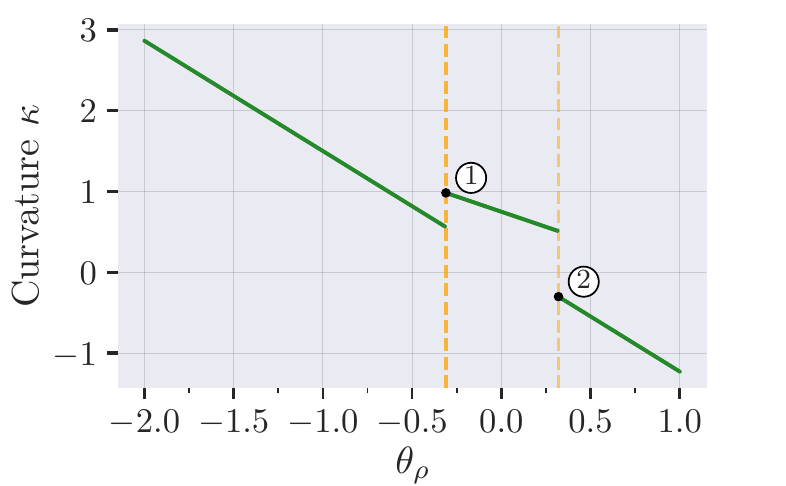}%
   	\subcaption{$0< \gamma \ll 1$ }
\end{subfigure}
\begin{subfigure}[c]{0.49\textwidth}
\includegraphics{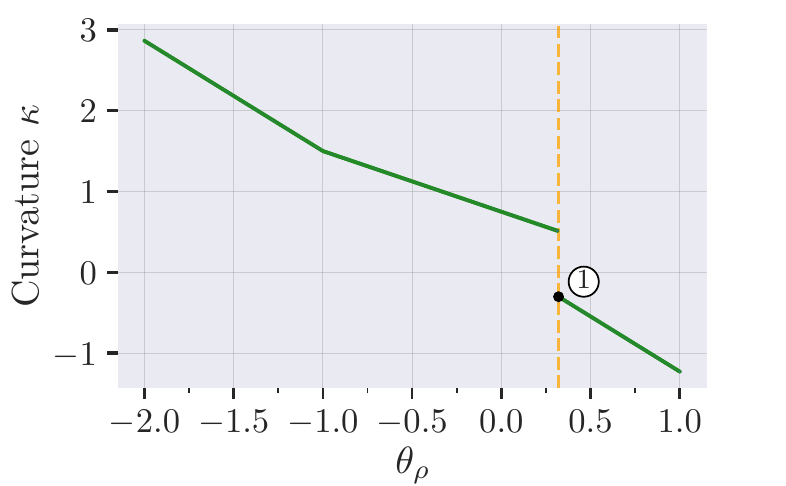}
    \subcaption{$ \gamma \gg 1 $ }
\end{subfigure}
\centering 
    \begin{tabular}[b]{c|c |c|c}
    Parameters &    $\theta_{\mu}$ & $\theta $  \\ \hline
  			   &   		10			&  $\frac{1}{2}$ \\
    \end{tabular}
    \caption{Curvature $\kappa(G_*)$ as a function of the prestrain ratio $\theta_\rho$.}
    \label{1dcurvatureplot(2)}
\end{figure}

\begin{figure}
    \centering
 	\includegraphics[trim={0 0 0 0.5cm},clip]{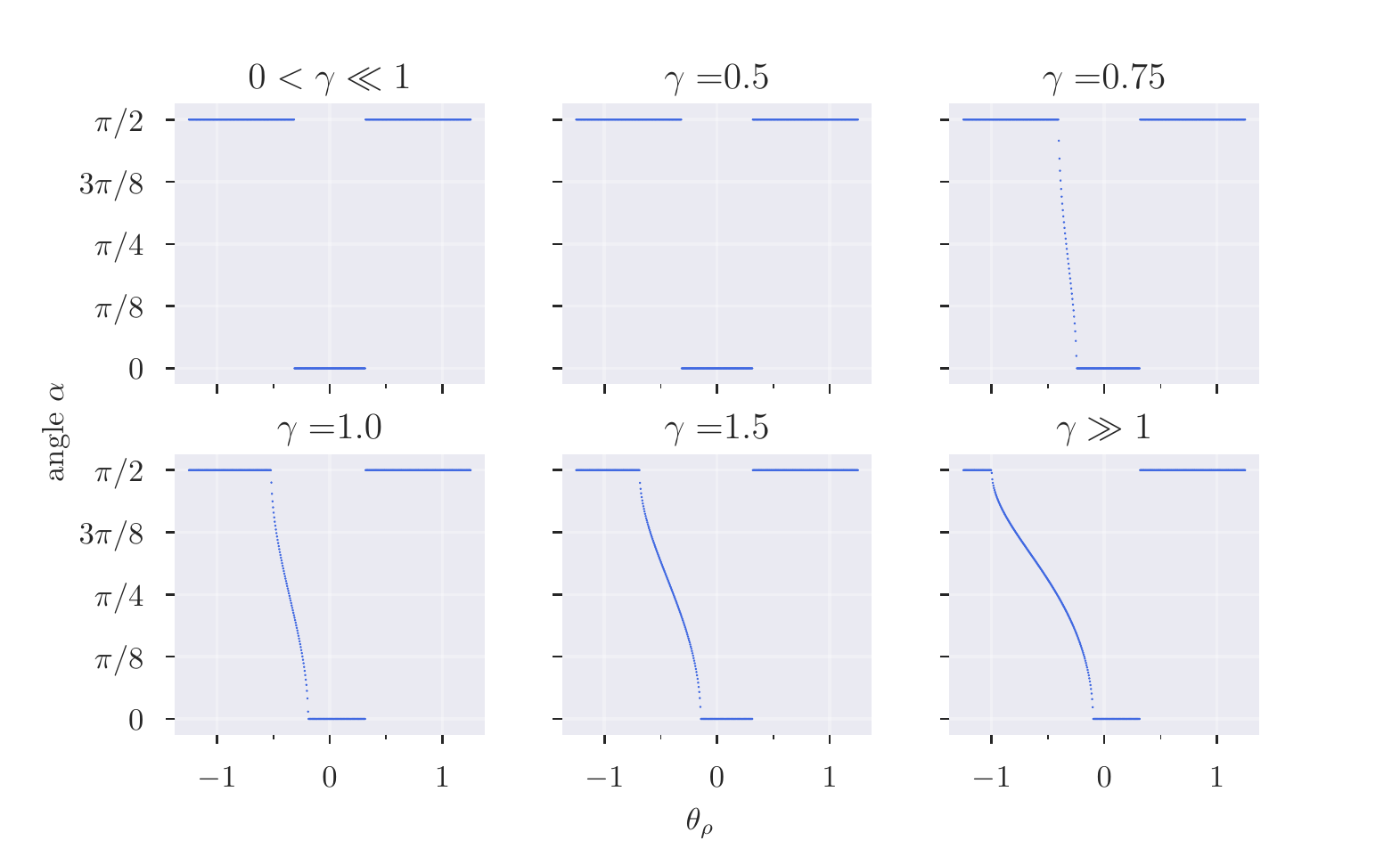}\hfill%
	\includegraphics[trim={0 0 0 0.5cm},clip]{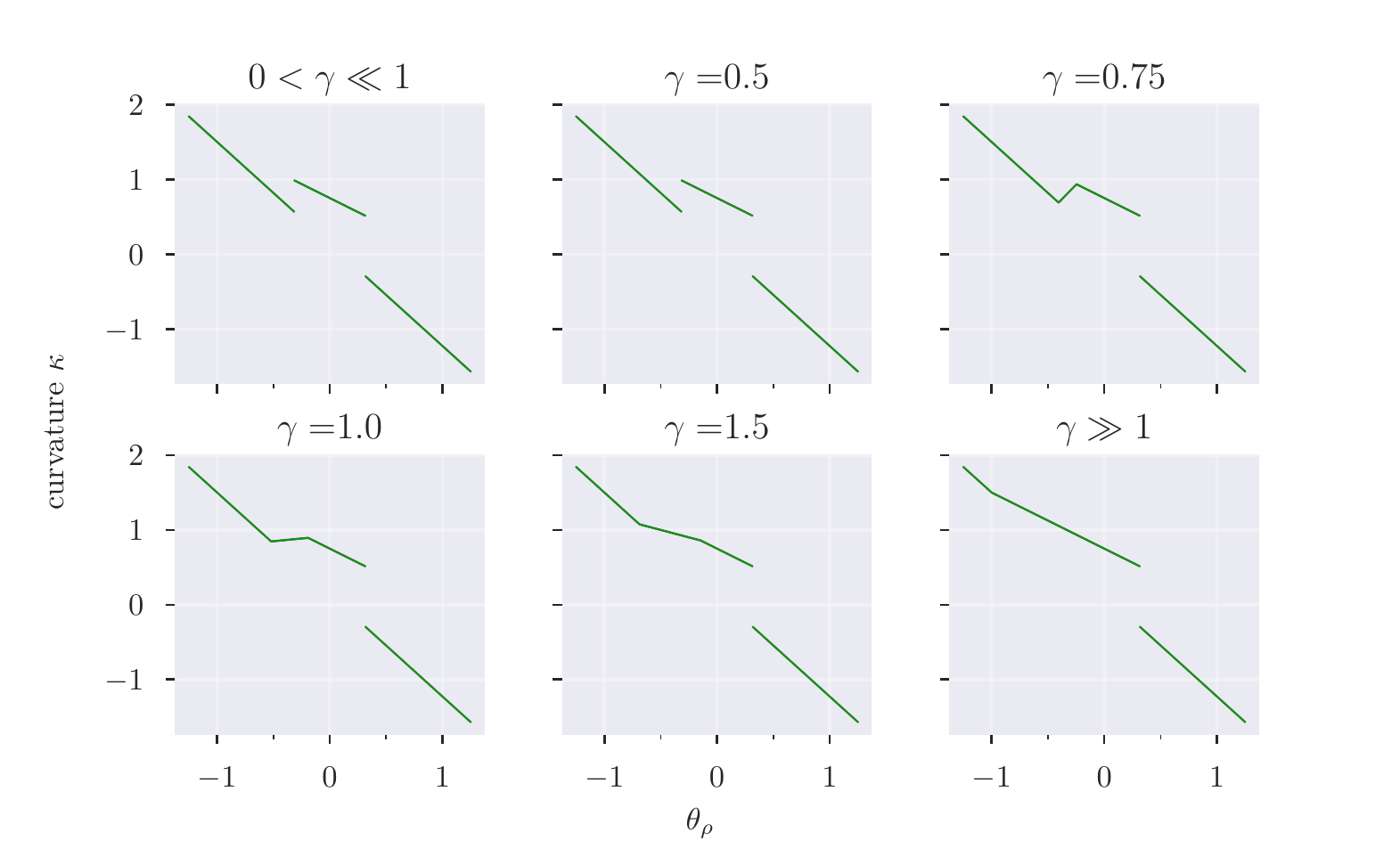}
	 \begin{tabular}[b]{c|c |c|c}
    Parameters &    $\theta_{\mu}$ & $\theta $  \\ \hline
  			   &   		10			&  $0.5$ \\
    \end{tabular}
    \caption{Angle $\alpha(G_*)$ (blue)  and curvature $\kappa(G_*)$ (green) as a function of the prestrain ratio $\theta_\rho$ for different values of $\gamma$.} 
\label{fig:midgamma} 
\end{figure}

We first discuss the regime $0<\gamma\ll1$.
In that case only angles $0$ and $\tfrac{\pi}{2}$ appear, which is in agreement with the trichotomy of minimizers, since $0<\gamma\ll 1$ is a sufficient condition for Lemma~\ref{L:char}~\ref{L:char:a}.
We see that the angle has a discontinuity at the boundary of the transition region.
This means that minimizers may flip from one  axial state to the other when parameters close to the boundary of the transition region are perturbed.
Figure~\ref{2DPhaseDiagram} (top, left) shows a cut through the diagram of Figure~\ref{3D-PhaseDiagram1} (top, left)  along the plane with $\theta_\mu = 2$.
We see that the boundary of the transition is tilted (left boundary) or curved (right boundary) with respect to the $\theta$-axis.
This means that for certain values of $\theta_\rho\in(-1,-0.5)\cup(0,\infty)$ the angle $\alpha(G_*)$ can be influenced by changing the volume fraction $\theta$.
Figure~\ref{anglevstheta} (top, left) visualizes the discontinuous dependence of the angle on the prestrain ratio $\theta_\rho$ for fixed parameters $\theta_\mu$ and $\theta$.
The marked points correspond to the boundary of the transition region.
The bottom left plots of Figures~\ref{3D-PhaseDiagram1}, \ref{2DPhaseDiagram}, and \ref{anglevstheta} visualize the curvature $\kappa(G_*)$.
Figure~\ref{2DPhaseDiagram} (bottom, left) suggests that the curvature is continuous as a function of $\theta$ in regions where the angle $\alpha(G_*)$ is constant, while we observe a jump in the curvature whenever the angle jumps.
It would further be natural to expect that the curvature is monotone in $\theta_\rho$, but Figure~\ref{anglevstheta} (bottom, left) shows that this is not the case: The curvature jumps upwards at the first point of discontinuity.
This observation can be seen more easily in Figure~\ref{1dcurvatureplot(2)} (left), where we consider  the larger value  $\theta_\mu=10$.
Next, we discuss the regime $\gamma\gg 1$.
The phase diagram in Figure~\ref{3D-PhaseDiagram1} (top, right) shows that in that regime the transition region features non-axial minimizers. Again, this is in agreement with the trichotomy of minimizers, since  $\gamma \gg 1$ is a necessary condition for Case~\ref{L:char:b} of Lemma~\ref{L:char}.
The cut shown in Figure~\ref{2DPhaseDiagram} (top, right) indicates that the transition on the left boundary of the transition region is continuous, while the transition on the right boundary is not.
In particular, the plot suggests that for a fixed volume fraction $\theta\in(0,1)$ we can obtain any angle in $[0,\tfrac{\pi}2]$ by choosing a suitable $\theta_\rho\in[-1,-0.5]$.
This is even more visible in Figure~\ref{anglevstheta} (top, right), which shows the $\theta_\rho$-dependence of the angle.
The marked points  \circled{1} and \circled{3} correspond to the boundary of the transition region. Furthermore, we observe a monotone behavior for $\theta_\rho\in[-1,-0.5]$.
Finally, a close look at the isolines of Figure~\ref{2DPhaseDiagram} (top, right) shows that for $\theta_{\rho}\in[-1,-0.5]$ the angle is monotone in the volume fraction $\theta$.
Figure~\ref{3D-PhaseDiagram1} (bottom, right) visualizes the curvature in the regime $\gamma \gg 1$. The phase diagram looks similar to the one in the regime $0<\gamma\ll1$.
However, the two are not identical, as becomes apparent by comparing the plots in Figure~\ref{anglevstheta} (bottom). 
The bottom right plots of Figures~\ref{3D-PhaseDiagram1}, \ref{2DPhaseDiagram}, and \ref{anglevstheta} suggest that also in the regime $\gamma\gg 1$, the curvature is continuous at points where the angle is continuous, and jumps if and only if the angle jumps.
In particular, Figure~\ref{anglevstheta} (bottom, right) shows that the curvature as a function of $\theta_\rho$ has only a single discontinuity at the marked point \circled{1}, which corresponds to the right boundary of the transition region.
In contrast to the regime $0<\gamma\ll1$, we observe in Figure~\ref{anglevstheta} (bottom, right) a monotone dependence on $\theta_\rho$.
\medskip
So far we have only analyzed the extreme regimes $0<\gamma \ll 1$ and $1 \ll \gamma $. 
We now briefly consider the intermediate regime $0<\gamma < \infty$.
We show these plots in Figure~\ref{fig:midgamma}. We note that the region of $\theta_{\rho}$ for which non-axial minimizers are observed is largest at $\gamma=\infty$. This region gets progressively smaller as $\gamma$ tends to $0$, and finally for values $ \gamma < \gamma^*$ it disappears and we observe a discontinuous transition between two axial minimizers. If $\theta_\mu = 10$, this region appears to lie always between $-1$ and $0$.

\paragraph{Discussion of the $\theta_\rho$-dependence.}
In the case $\theta_\rho=0$ the prestrain is only active in the top layer of the two-layer microstructure of Figure~\ref{sketchofmicrostructure}.
More precisely, the ``fibres'' in the top layer want to expand isotropically by a factor $\rho_1=1$ to gain equilibrium.
In view of this we expect that the plate bends downwards either in parallel or orthogonal to the fibres.
This corresponds to $\kappa(G_*)>0$ and $\alpha(G_*)\in\{0,\tfrac{\pi}{2}\}$.
Indeed, Figures~\ref{2DPhaseDiagram}, and \ref{fig:midgamma} show that this is indeed the case: Independently of $\gamma$, for $\theta_\mu\in\{2,10\}$ and all $\theta\in(0,1)$,  we have $\alpha(G_*)=0$ and $\kappa(G_*)>0$.
This means that the plate bends in a direction orthogonal to the fibres.
For $\theta_\rho\neq 0$ the prestrain is active in both layers. For $\theta_\rho<0$ the prestrain in the top and bottom layer have the opposite sign, and thus ``push'' into the same direction.
Hence, one could expect that the curvature monotonically increases when $\theta_\rho$ decreases.
But in the regime $0<\gamma\ll1$  we observe a downwards jump of the curvature, see Figure~\ref{anglevstheta} (bottom, left).
In the case $\theta_\rho>0$ the sign of the prestrain in the two layers is the same, and thus they are competing with regard to energy minimization.
We expect that for $\theta_\rho$ sufficiently large the bending induced by the fibres in the bottom layer dominate the behavior,
which would mean that the plate bends ``upwards''.
Indeed, as shown by Figure~\ref{2DPhaseDiagram}, for $\theta_\rho$ above the threshold (which decreases with the volume fraction $\theta$), the curvature changes its sign from negative to positive.
Surprisingly, at that point the angle jumps from $\alpha(G_*)=0$ to $\alpha(G_*)=\tfrac{\pi}{2}$.

\paragraph{Conclusion.}
Our analysis shows that the parametrized, prestrained laminate features a complex behavior.
In particular, we make the following observations:
\begin{enumerate}[(a)]
\item We observe a discontinuous dependence of the set of minimizers on the parameters.
\item We observe non-uniqueness of the global minimizers: In the regime $\gamma\gg 1$ we find parameters leading to non-axial minimizers. Those always come in pairs of the form $\{G_*,TG_*\}$. Likewise, for the special cases of Remark~\ref{R:homcase} we have a one-parameter family of minimizers. 
\item The mechanical system features a break of symmetry: For $\theta_\mu=1$, which corresponds to a homogeneous material, the set of minimizers is of the form $\{\kappa(\xi\otimes\xi)\,:\,\xi\in\mathbb R^2,\,|\xi|=1\}$ and \textit{rotationally symmetric}, while for $\theta_\mu\neq 1$ the set of minimizers degenerates to a set with one axial or two (possibly non-axial) minimizers.
\item The mechanical properties for $0<\gamma\ll1$ and $\gamma\gg 1$ are qualitatively and quantitatively different: In the former regime minimizers are always axial and $G_*$ as a function of $\theta_\rho$ is non-monotone and has two discontinuities, while the latter regime features minimizers of arbitrary angle and $G_*$ as a function of $\theta_\rho$ is monotone and only has one discontinuity.
\end{enumerate}


\section{Shape programming}\label{S:textures}
This section is devoted to the problem of shape programming: given \textit{a target shape} that can be parametrized by a bending deformation with second fundamental form $\II_*$, can we construct a \emph{simple}
prestrained composite such that minimizers of the corresponding three-dimensional elastic energy
approximate the target shape for $(\e,h)\to0$ ?  If the target shape is \textit{cylindrical} the problem can be solved relatively easily with the help of the parametrized laminate composite studied in Section~\ref{S:micro-shape}. In the general case, the problem is more complex and requires a composite whose microstructure changes with the in-plane position. We first discuss the simpler case of a cylindrical target shape in Section~\ref{S:textures1} and then turn to general target shapes in Section~\ref{S:textures2}. That section culminates in Theorem~\ref{T:shape} which provides our analytical answer to the problem of shape programming.
%

\subsection{The case of a cylindrical shape and the notion of a composite template}\label{S:textures1}
In the following we consider the problem of shape programming in the special case of a cylindrical target shape, that is, we are looking for a deformation with second fundamental form $\II_*=\kappa_*(\minvec_*\otimes \minvec_*)$ a.e.~in $S$ for some $\kappa_*\in\R$ and a unit vector $\minvec_*\in\R^2$. Our task is to identify a composite such that minimizers of the corresponding 3d-energy ``effectively'' describe a surface that coincides with the target shape. Here, ``effectively'' means that we consider accumulation points of sequences of minimizers of the 3d energy along the limit $(\e,h)\to 0$. As we explain next, for the cylindrical shapes the problem of shape programming can solved explicitly with the help of the parametrized laminate considered in Lemma~\ref{S:ex:C3}. Specifically, let $W_1,W_2$ be two isotropic stored energy functions satisfying \ref{item:nonlinear_material_w1} -- \ref{item:nonlinear_material_w4} with Lam\'e-parameters $\mu_1=1$, $\mu_2=2$, and $\lambda_1=\lambda_2=0$. We consider a $\Lambda \colonequals I_{2\times 2}$-periodic laminate that is constant in $e_2$-direction and composed of the materials described by $W_1$ and $W_2$. We denote by $\theta\in[0,1]$ the volume fraction of the material described by $W_2$, and we model the parametrized composite by a stored energy function $W(\theta;\cdot):\R\times\R^{3\times 3}\to[0,+\infty]$, $(y_1,F)\mapsto W(\theta;y_1,F)$. We assume that $W(\theta;y_1,F)$ is $1$-periodic in $y_1$, and
  \begin{equation*}
    W(\theta;y_1,F) \colonequals
    \begin{cases}
      W_1(F)&\text{if }|y_1|>\frac\theta 2,\\
      W_2(F)&\text{else,}
    \end{cases}\qquad \text{for }y_1\in(-\tfrac12,\tfrac12],\,F\in\R^{3\times 3}.
  \end{equation*}
  (By periodicity, $W(\theta;\cdot)$ is uniquely defined by this identity.) Let $Q(\theta;\cdot)$ be the quadratic form associated with $W(\theta;\cdot)$. Furthermore, for $\rho_1>0$ fixed, consider a prestrain~$B$ that is $\Lambda=I_{2\times 2}$-periodic and admissible in the sense of Definition~\ref{D:admissible}, and that satisfies
  \begin{equation*}
B(\theta;x_3,y) \colonequals
    \begin{cases}
      \rho_1I_{3\times 3}&\text{if }|y_1|>\frac\theta 2\text{ and }x_3>0,\\
      0&\text{else.}
    \end{cases}
  \end{equation*}
  Then $Q(\theta;\cdot)$ and $B(\theta;\cdot)$ coincide with the parametrized laminate considered in Lemma~\ref{S:ex:C3}  with parameters $\mu_1=1,\,\theta_\mu=2$, and $\theta_\rho=0$. By Lemma~\ref{S:ex:C3}, the coefficients of the associated homogenized quadratic form $Q^\gamma_{\hom}(\theta;\cdot)$ and effective prestrain $B^\gamma_{\eff}(\theta)$ are given by
  \begin{align*}
    &\widehat B_{\eff,1}^{\gamma}=\frac{3\rho_1}{2}(1-\theta),\qquad \widehat B_{\eff,2}^{\gamma}=\frac{3\rho_1}{2}\frac{1-\theta}{1+\theta},\qquad  \widehat B_{\eff,3}^{\gamma}=0,\\
    &q_1=\frac{1}{3(2-\theta)},\qquad q_2=\frac{1+\theta}{6},\qquad q_{12}=0.
  \end{align*}
  There is no explicit expression for $q_3(\gamma)$, but we know that $q_1\leq q_3(\gamma)\leq q_2$ by \eqref{S:ex:muGammaProp1}.
  In the following we fix an arbitrary value for $\gamma\in(0,\infty)$. For any $\theta\in(0,1]$ one can check (with help of a computer algebra system and only using $q_1\leq q_3(\gamma)\leq q_2$) that for this specific composite Case of \ref{L:char:a} of Lemma~\ref{L:char} always applies. In view of this, Lemma~\ref{L:char} implies that the algebraic minimization problem
  $$\min_{G\in\R^{2\times 2}\atop\det G=0} Q^{\gamma}_{\hom}\big(\theta;G-B^\gamma_{\hom}(\theta)\big)$$ admits a unique axial minimizer, which is given by
  \begin{equation*}
    \kappa(\theta)e_1\otimes e_1,\qquad\text{where }\kappa(\theta)=\frac{3\rho_1}{2}(1-\theta).
  \end{equation*}
  
  In the following let $\mathcal I^{\e,h}(\theta;\cdot)$ denote the 3d-energy defined by \eqref{def:ene} with $W_\e(x,F) \colonequals W(\theta;\frac{x_1}{\e},F)$ and $B_{\e,h}(x)=B(\theta;\frac{x_1}{\e})$. Also, let $\mathcal I^\gamma_{\hom}(\theta;\cdot)$ denote the  bending energy defined by \eqref{limit:eq1} with $Q^\gamma_{\hom}$ and $B^\gamma_{\eff}$ replaced by $Q^\gamma_{\hom}(\theta;\cdot)$ and $B^\gamma_{\eff}(\theta)$, respectively. Finally, let $(\deform_h)$ in $H^1(\Omega;\R^3)$ be an arbitrary sequence of almost minimizers for $\mathcal I^{\e(h),h}(\theta;\cdot)$. From Theorem~\ref{T1} we conclude that any accumulation point $\deform$ of $\big(\deform_h-\fint_\Omega \deform_h\big)$ for $h\to 0$ is a bending deformation that minimizes the energy $\mathcal I^\gamma_{\hom}(\theta;\cdot)$. Furthermore, with the help of Lemma~\ref{L:char:cylindrical} we conclude that $\II_\deform=\kappa(\theta)e_1\otimes e_1$. We may  summarize  this observation by saying that the composite described by $W(\theta;\cdot)$ and $B(\theta;\cdot)$ \emph{effectively programs} the target shape $\II_*=\kappa(\theta)(e_1\otimes e_1)$.

  Since the map $(0,1]\ni\theta\mapsto \kappa(\theta)\in[0,\frac{3\rho_1}{2})$ is a homeomorphism, we can \emph{program} any target shape $\II_*=\kappa_*(e_1\otimes e_1)$ with $\kappa_*$ in the range $[0,\frac{3\rho_1
  }{2})$---we just need to make an appropriate choice of $\theta_*\in(0,1]$. Morever, a non-positive target curvature $\kappa_*$ in the range $(-\frac{3\rho_1}{2},0]$ can be obtained by considering the reflected prestrain $(x_3,y)\mapsto B(\theta;-x_3,y)$. For a general unit vector $\minvec_*\in\R^2$ we can recover target shapes of the form $\II_*=\kappa_*(\minvec_*\otimes \minvec_*)$  by introducing a rotation of the laminate microstructure such that the direction of the rotated laminate is given by $\minvec_*$; this is made precise by the transformation rule of Corollary~\ref{C:transf} below.
  In conclusion, with help of the parametrized, two-phase laminate introduced above, we can program any cylindrical shape with curvature $\kappa_*\in \mathcal K \colonequals (-\frac{3\rho_1}{2},\frac{3\rho_1}{2})$ by an appropriate choice of the volume fraction $\theta$ and the laminate direction. Since the map $\theta\mapsto\kappa(\theta)$ is bijective, we may consider directly $\kappa\in\mathcal K$ (instead of $\theta$) as the design parameter of the composite.

\begin{figure}
    \begin{tabular}{ll}
      \begin{minipage}[t][4cm][t]{10.5cm}
        (a)%
        \begin{minipage}[t][1.2cm][b]{10cm}
          \includegraphics[width=10cm]{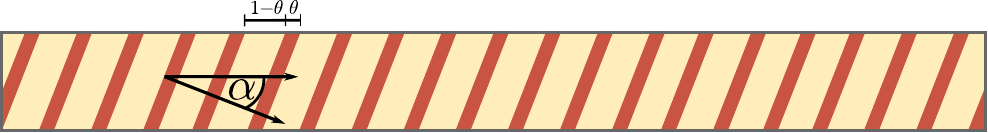}
        \end{minipage}\\[0.2cm]
        (b)%
        \begin{minipage}[t][2.8cm][b]{10cm}
          \includegraphics[width=10cm, trim={0 1cm 0 0cm},clip]{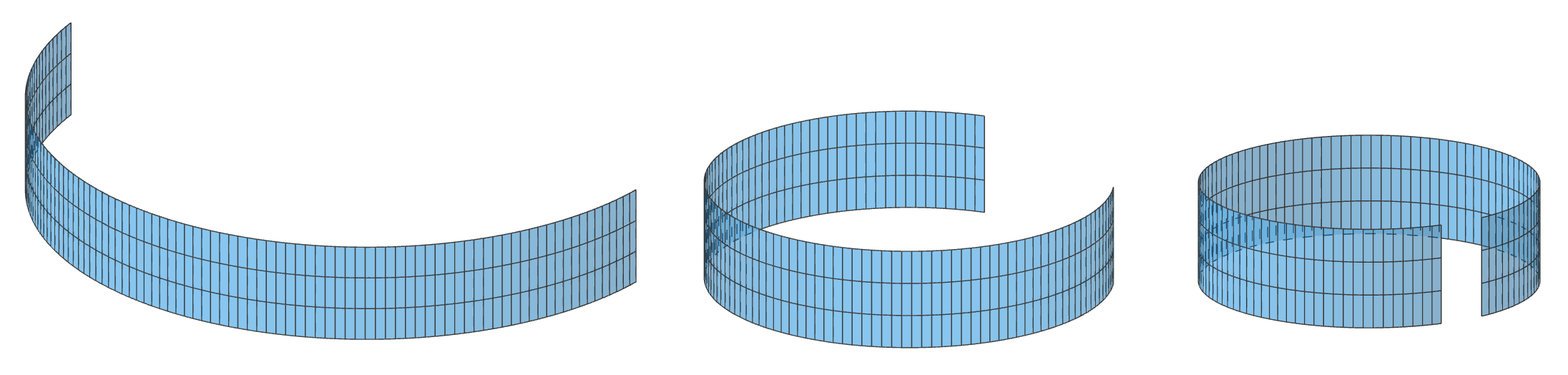}
          \vspace{-\baselineskip}
          \begin{center}
          \hspace{0.8cm} \footnotesize{$\kappa(\theta)=0.3$  \hspace{2cm} $\kappa(\theta)=0.5$ \hspace{1.4cm} $\kappa(\theta)=0.6$}
          \end{center}
        \end{minipage}%
      \end{minipage}
      &
        (c)%
        \begin{minipage}[t][4cm][b]{5cm}
          \includegraphics[width=4.5cm]{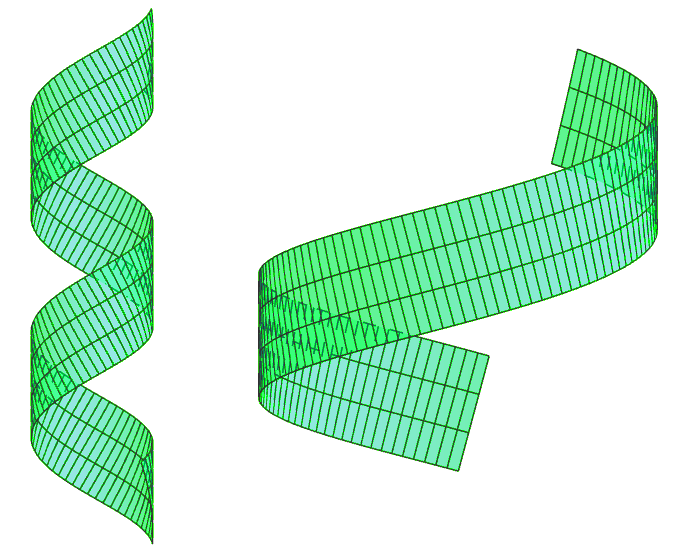}
        \end{minipage}
    \end{tabular}
     \vspace{0.7cm}
    \caption{Cylindrical shapes programmed via a two-phase laminaate composite. Figure (a) shows the reference configuration $S=(0,10)\times (0,1)$ of a composite plate. The plate is filled with the parametrized composite of Lemma~\ref{S:ex:C3} with fixed parameters $\mu_1=1$, $\theta_\mu=2$, $\lambda_1=\lambda_2=0$, $\theta_\rho=0$, and $\rho_1>0$. The volume fraction $\theta$ and the rotation angle $\alpha$ are considered as design parameters. We use the same color coding as in Figure~\ref{sketchofmicrostructure_a}:  yellow refers to the material with energy $W_1$ and no prestrain, while red refers to the prestrained material with energy $W_2$. Figure (b) shows global minimizers of the effective energy for the case when $\alpha=0$; from left to the right the value of $\theta$ decreases. Specifically, $\theta$ is chosen such that $\kappa(\theta)\in\{0.3,0.5,0.6\}$. Figure (c) shows (from left to right) the global minimizers of the effective energy for $\alpha\in\{15^\circ,30^\circ\}$ and a choice of $\theta$ with $\kappa(\theta)=1.5$.
}
    \label{fig:cylindrical_strip}
  \end{figure}

  Figure~\ref{fig:cylindrical_strip} illustrates some examples of cylindrical deformations of a strip obtained with a two-phase laminate. Similarly, Figure~\ref{fig:expstrip} shows cylindrical shapes which were experimentally observed for 3d--printed, two--layer composite plates with a laminar microstructure consisting of prestrained fibres. As in our analysis, the orientation and the volume fraction of the fibres are design parameters. Qualitatively, our analysis reproduces the relation between these parameters and the geometry of the strip.

  The two-phase laminate discussed above is an example of what we call a \textit{composite template}:
  \begin{definition}[Composite template]\label{D:bb}
  Let $\gamma\in(0,\infty)$, $\Lambda\in\R^{2\times 2}$ be invertible, and let $\mathcal K\subset\R$ be an open interval. The one-parameter family $(\bar W(\kappa,\cdot), \bar B(\kappa,\cdot))_{\kappa\in \mathcal K}$ is called a \emph{composite template}, if the following properties hold:
  \begin{enumerate}[(i)]
  \item \label{item:bb:i}
    For all $\kappa\in \mathcal K$ the energy density $\bar W(\kappa;\cdot):(-\frac12,\frac12)\times\R^2\times\R^{3\times 3}\to[0,+\infty]$ is a Borel function, and $\bar W(\kappa;x_3,y,\cdot)\in\mathcal W(\alpha,\beta,\rho,r)$ for a.e.~$(x_3,y)\in(-\frac12,\frac12)\times\R^2$. Furthermore, $\bar Q(\kappa;\cdot)$ and $\bar B(\kappa;\cdot)$ are $\Lambda$-periodic and admissible in the sense of Definition~\ref{D:admissible}, and
    \begin{equation*}
      \fint_{\Box_\Lambda}|\bar B(\kappa;x_3,y)|^2\dd(x_3,y)\leq C_B,
    \end{equation*}
    for a constant $C_B$ that is independent of $\kappa$.
    
  \item \label{item:bb:ii}
    For any sequence $\kappa_j\to\kappa$ in $\mathcal K$, we have
    \begin{alignat*}{2}
      \bar Q(\kappa_j;x_3,y,G) &\to \bar Q(\kappa;x_3,y,G) & \qquad & \text{for all $G\in\R^{3\times 3}$ and a.e.~}(x_3,y),\\
      \bar B(\kappa_j;\cdot)   &\to \bar B(\kappa;\cdot)   &        & \text{strongly in $L^2(\Box_\Lambda;\R^{3\times 3})$}.
    \end{alignat*}

  \item \label{item:bb:iii}For all $\kappa\in \mathcal K$ the matrix $\kappa(e_1\otimes e_1)$ is a unique global minimizer of the algebraic minimization problem
    \begin{equation*}
      \min_{G\in\R^{2\times 2}\atop\det G=0}\bar Q_{\hom}^{\gamma}\big(\kappa;G-\bar B^\gamma_{\eff}(\kappa)\big),
    \end{equation*}
    where $\bar Q^\gamma_{\hom}(\kappa;\cdot)$ and $\bar B^\gamma_{\eff}(\kappa)$ are the effective quantities associated with $\bar Q(\kappa;\cdot)$ and $\bar B(\kappa;\cdot)$
    by Definitions~\ref{def:Qgamma} and~\ref{D:eff}, respectively.
  \end{enumerate}
\end{definition}
With a composite template we can directly program shapes of the form $\II_*=\kappa(e_1\otimes e_1)$ for any $\kappa\in\mathcal K$. As already indicated, by combining the composite template with an in-plane rotation we can also program shapes with arbitrary bending direction. This can be seen from the following two transformation results:

\begin{lemma}[Transformation of microstructures]\label{L:transformation}
  Let $\Lambda,T\in\R^{2\times 2}$ be invertible, and let $(Q,B)$ be $\Lambda$-periodic and admissible in the sense of Definition~\ref{D:admissible}. 
  Consider the transformation
  \begin{equation*}
    \widetilde Q(x_3,y,\tilde G)\colonequals Q(x_3,T^{-1}y,\widehat T^{\top} \tilde G\widehat T),\qquad \tilde B(x_3,y)\colonequals\widehat T^{-\top}B(x_3,T^{-1}y)\widehat T^{-1},\qquad \widehat T\colonequals
    \begin{pmatrix}
      T&\\&1
    \end{pmatrix}.
  \end{equation*}
  Then $(\widetilde Q,\widetilde B)$ is $\widetilde\Lambda\colonequals T\Lambda$-periodic and admissible. Furthermore, for all $\tilde G\in\R^{2\times 2}_{\sym}$ we have
  \begin{equation*}
    \widetilde Q^\gamma_{\hom}(\tilde G)=Q^\gamma_{\hom}(T^{\top}\tilde GT)\text{ and }\tilde B^\gamma_{\rm eff}=T^{-\top}B^\gamma_{\rm eff}T^{-1}.
  \end{equation*}
  Above, $Q^\gamma_{\hom}$, $\widetilde Q^\gamma_{\hom}$. $B^\gamma_{\rm eff}$, and $\tilde B^\gamma_{\rm eff}$ are associated with $Q$,$\widetilde Q$, $B$, and $\tilde B$ via Definitions~\ref{def:Qgamma} and~\ref{D:eff}, respectively.
\end{lemma}
\reftoproof{SS:shapedesign}

As a corollary of Lemma~\ref{L:transformation} we obtain the following transformation rule for minimizers of the algebraic minimization problem \eqref{L:ex:1:eq1}. Roughly speaking, it states that the transformation by an in-plane rotation and the passage to the effective quantities commute.
\begin{corollary}[Rotation of microstructure of a composite template]\label{C:transf}
  Let $\big(\bar W(\kappa;\cdot)$, $\bar B(\kappa;\cdot)\big)_{\kappa\in \mathcal K}$ be a composite template
  as in Definition~\ref{D:bb}. For $R\in\SO 2$
  and $\widehat{R} \colonequals \operatorname{diag}(R, 1)$ consider the transformation
  \begin{align*}
    W(R,\kappa;x_3,y,F)\colonequals\bar W(\kappa;x_3,R^\top y,F\widehat R),
    \qquad
    B(R;\kappa;x_3,y)\colonequals\widehat R\bar B(\kappa;x_3,R^\top y)\widehat R^\top.
  \end{align*}
  Then:
  \begin{enumerate}[(a)]
  \item For all $R\in\SO 2$ the transformed composite $\big(W(R,\kappa;\cdot),B(R,\kappa;\cdot)\big)_{\kappa\in \mathcal K}$ satisfies Properties~\ref{item:bb:i} and~\ref{item:bb:ii}
  of Definition~\ref{D:bb}.

  \item For all $R\in\SO 2$ and $\kappa\in \mathcal K$ we have
    \begin{equation*}
      Q^\gamma_{\hom}(R,\kappa;G)=\bar Q^\gamma_{\hom}(\kappa; R^\top GR)\text{ for all }G\in\R^{2\times 2}_{\sym}\text{ and } B^\gamma_{\textnormal{eff}}(R,\kappa)=R\bar B^\gamma_{\rm eff}(\kappa)R^\top.
    \end{equation*}
    Here, $Q^\gamma_{\hom}(R,\kappa;\cdot)$ and $B^\gamma_{\rm eff}(R,\kappa)$ denote the effective quantities associated with\linebreak $\big(W(R,\kappa;\cdot),B(R,\kappa;\cdot)\big)$.
    Likewise, $\bar Q^\gamma_{\hom}(\kappa;\cdot)$ and $\bar B^\gamma_{\rm eff}(\kappa)$
    denote the effective quantities associated with $\big(\bar W(\kappa;\cdot),\bar B(\kappa;\cdot)\big)$.

  \item For all $G\in R^{2\times 2}_{\sym}$ the map
    \begin{equation*}
      \SO2\times \mathcal K\ni(R,\kappa)\mapsto \Big(Q^\gamma_{\hom}(R,\kappa;G),B^\gamma_{\rm eff}(R;\kappa)\Big)
    \end{equation*}
    is continuous.
  \item For all $R\in\SO 2$ and $\kappa\in \mathcal K$ the matrix $\kappa(Re_1\otimes Re_1)$ is the unique minimizer of the algebraic minimization problem
    \begin{equation*}
      \min_{G\in\R^{2\times 2}\atop\det G=0}Q^\gamma_{\hom}\big(R,\kappa;G-B^\gamma_{\rm eff}(R;\kappa)\big).
    \end{equation*}
  \end{enumerate}
\end{corollary}
\reftoproof{SS:shapedesign}

\subsection{The general case}\label{S:textures2}
We turn to the shape programming problem for surfaces parametrized by general bending deformations. More precisely, given a second fundamental form $\II_*$ associated with a bending deformation $\deform_*\in H^2_{\iso}(S;\R^3)$, the task is to identify a \textit{simple} composite such that any accumulation point $\deform\in H^2_{\iso}(S;\R^3)$ of a sequence of almost minimizers for the associated 3d-energy, parametrizes a surface that coincides with the target shape in the sense that $\II_\deform=\II_*$. The problem is trivial if we allow arbitrarily complex composites. Instead, we shall only consider a restricted class of composites that are \textit{simple} in the following sense:
\begin{itemize}
\item We assume that the 3d-plate is \emph{finitely structured} in the sense that
  $\Omega$ is partitioned into finitely many ``grains'' of the form $S_j\times(-\frac12,\frac12)$, $j\in J$ (cf.~Assumption~\ref{ass:pcperiodicity}).
\item A \textit{composite template} $\Big(W(\kappa;\cdot),B(\kappa;\cdot)\Big)_{\kappa\in\mathcal K}$ shall be given and each grain is filled with the composite template, possibly rotated in--plane.
\end{itemize}
We formalize this concept in the following definition:
\begin{definition}[Structured composite]\label{def:struct}
  Let $(\bar W(\kappa,\cdot), \bar B(\kappa,\cdot))_{\kappa\in \mathcal K}$ be a composite template. Let $\{S_j\}_{j\in J}$ denote a partition of $S$ as in Assumption~\ref{ass:pcperiodicity}. A \emph{structured composite} (based on the composite template and subordinate to the partition) is a pair $(W,B)$ consisting of Borel functions $W:\Omega\times\R^2\times\R^{3\times 3}\to[0,+\infty]$ and $B:\Omega\times\R^2\to\R^{3\times 3}$ such that for all $j\in J$ there exists a rotation $R_j\in \SO 2$ and a parameter $\kappa_j\in \mathcal K$ such that for all $x'\in S_j$,
  \begin{align*}
    W(x',x_3,y,F)
    & =
    \bar W(\kappa_j,x_3,R_j^\top y,F\widehat R_j), \\
    B(x',x_3,y)
    & =
    \widehat R_j\bar B(\kappa_j;x_3,R^\top_jy)\widehat R_j^\top,\qquad\text{where }\widehat R_j \colonequals
  \begin{pmatrix}
    R_j&\\&1
  \end{pmatrix}.
  \end{align*}
  If $J$ is finite, then $(W,B)$ is called a \emph{finitely structured composite}.
\end{definition}
\begin{remark}
Note that a structured composite is a special case of a prestrained composite satisfying Assumption~\ref{ass:W}. Specifically, if $(W,B)$ is a structured composite, then $W_\e(x,F) \colonequals W(x,\tfrac{x'}{\e},F)$ and $B_{\e(h),h}(x) \colonequals B(x,\tfrac{x'}{\e(h)})$ satisfy Assumption~\ref{ass:W}. Also note that once the composite template $(\bar W,\bar B)$  is fixed, the structured composite is fully characterized by the partition $\{S_j\}_{j\in J}$, the parameters $\{\kappa_j\}_{j\in J}\subseteq\mathcal K$, and the rotations $\{R_j\}_{j\in J}\subseteq \SO 2$.
\end{remark}
For an example of a composite template and a schematic picture of an associated finitely structured composite, see Figure~\ref{fig:cone_strip} below.
The following main result states that any target shape (satisfying a certain restriction on its curvature, see \eqref{eq:kappa} below) can effectively be approximated by almost--minimizers of a 3d-energy associated with a finitely structured composite.

\begin{theorem}[Shape programming]\label{T:shape}Let Assumption~\ref{A:gamma} be satisfied. Let $(\bar W(\kappa;\cdot), \bar B(\kappa;\cdot))_{\kappa\in \mathcal K}$ be a composite template,
  and let $\II_*$ be the second fundamental form of a bending deformation in $H^2_{\rm iso}(S;\R^3)$ such that
  \begin{equation}\label{eq:kappa}
    \II_*(x')\in\big\{\kappa (\minvec\otimes \minvec)\,:\,\kappa\in \mathcal K,\,\minvec\in\R^2,\,|\minvec|=1\big\}
    \qquad
    \text{for a.e.~$x'\in S$.}
  \end{equation}
  Then for all $\delta>0$, there exists a finitely structured composite $(W,B)$ (based on the composite template and subordinate to a finite partition of $S$ that depends on $\delta$)
  such that the following holds: Consider the functional
  \begin{equation*}
    \mathcal I^h(\deform)\colonequals
    \frac1{h^2}\int_{\Omega} W\Big(x,\tfrac{x'}{\e(h)},\nabla_h \deform(x)\big(I_{3\times 3}-hB\big(x,\tfrac{x'}{\e(h)}\big)\big)\Big)\,dx,
  \end{equation*}
  and let $(\deform_h)$ denote a sequence with $\fint_{\Omega}\deform_h\,dx=0$ that satisfies $\mathcal I^h(\deform_h)\leq \inf\mathcal I^h+h$. Then any accumulation point $\deform$ of $(\deform_h)$ (viewed as a sequence in $L^2(\Omega;\R^3)$) is a bending deformation, i.e., $\deform\in H^2_{\iso}(S;\R^3)$, and satisfies
  \begin{equation}\label{T:shape:1}
    \|\II_\deform-\II_*\|_{L^2(S)}\leq \delta.
  \end{equation}
\end{theorem}
\reftoproof{SS:shapedesign}

Condition~\eqref{eq:kappa} is  a restriction on the curvature of the the target shape. In view of Definition~\ref{D:bb} \ref{item:bb:iii}, the set   $\mathcal K$ is the curvature range that can be recovered by the composite template under consideration. For example, if we consider the model composite of Section~\ref{S:textures1} as a composite template, then \eqref{eq:kappa} takes the form of the restriction $|\II_*(x')|<\frac{3\rho_1}{2}$. We see that the range of admissible shapes increases for larger values of the parameter $\rho_1$, which in the model composite is the magnitude of the prestrain.

In a nutshell our result is:
For any composite template and any target shape satisfying the constraint on the curvature \eqref{eq:kappa}, we can find a partition and a subordinate, finitely structured composite that programs the target shape up to a tolerance that is quantified by the parameter $\delta$.
As shall become clear from the proof of Theorem~\ref{T:shape}, in order to decrease $\delta$ one needs to refine the partition.
We note that the proof of Theorem~\ref{T:shape} yields a partition of $S$ into dyadic squares and a boundary layer.
\medskip

An illustration of a finitely structured composite designed to approximate a conically deformed strip is shown in Figure~\ref{fig:cone_strip}. Figure~\ref{fig:cone_strip} (a) shows an isometrically deformed paper strip glued to the surface of a symmetric cone. The associated reference domain $S=(0,4)\times (0,1)$ is shown in Figure~\ref{fig:cone_strip} (a). The asymptotic lines of the deformed strip are indicated in blue. The deformation is affine along these lines.  Since the deformed strip lies on a cone, these lines intersect in a single point. The principal curvature curves of the strip are shown in red. In the conical setting, the curves lie on circles whose center is given by the intersection point of the blue lines. The mean curvature is constant along these red curves. Moreover, we note that the curvature increases when moving along the blue lines (from lower--left to upper--right).

As indicated in Figure~\ref{fig:cone_strip}~(b), the reference domain $S$ is equipartitioned into small squares $\{S_j\}_{j\in J}$. An example of a finitely structured composite can be defined by filling each 3d-cube $S_j\times(-\frac12,\frac12)$, $j\in J$, with a rotated instance of the composite template introduced in Section~\ref{S:textures1}. It is based on the parametrized laminate shown in Figure~\ref{sketchofmicrostructure} and features the volume fraction $\theta$ as a parameter. In order to obtain a finitely structured composite that leads to an approximation of the deformed strip, in each box $S_j\times(-\frac12,\frac12)$ the rotation of the composite template needs to be choosen according to the direction of the asymptotic (blue) line in $S_j$, and the parameter~$\theta$ needs to be chosen in dependence of the curvature in $S_j$. The two schematic drawings in Figure~\ref{fig:cone_strip}~(c) illustrate this: The left drawing is the parametrized composite in the lower--left square of $S$ and the right drawing is the parametrization in the upper--right square. In Figure~\ref{fig:cone_strip}~(c) the yellow color refers to the domain occupied by the material with energy $W_1$ (where no prestrain is present), while the red color refers to the prestrained material with energy $W_2$. We see that the laminate direction is aligned with the asymptotic line. Furthermore, we see that in the composite assigned to the lower-left square, the volume fraction of the prestrained material is smaller than in the composite assigned with the upper--right square. This corresponds to the fact that the curvature in the lower-left square is smaller than the curvature in the upper--right one.


\begin{figure}[h]
  (a)%
  \begin{minipage}[t][3.7cm][b]{5.5cm}
    \includegraphics[width=5cm]{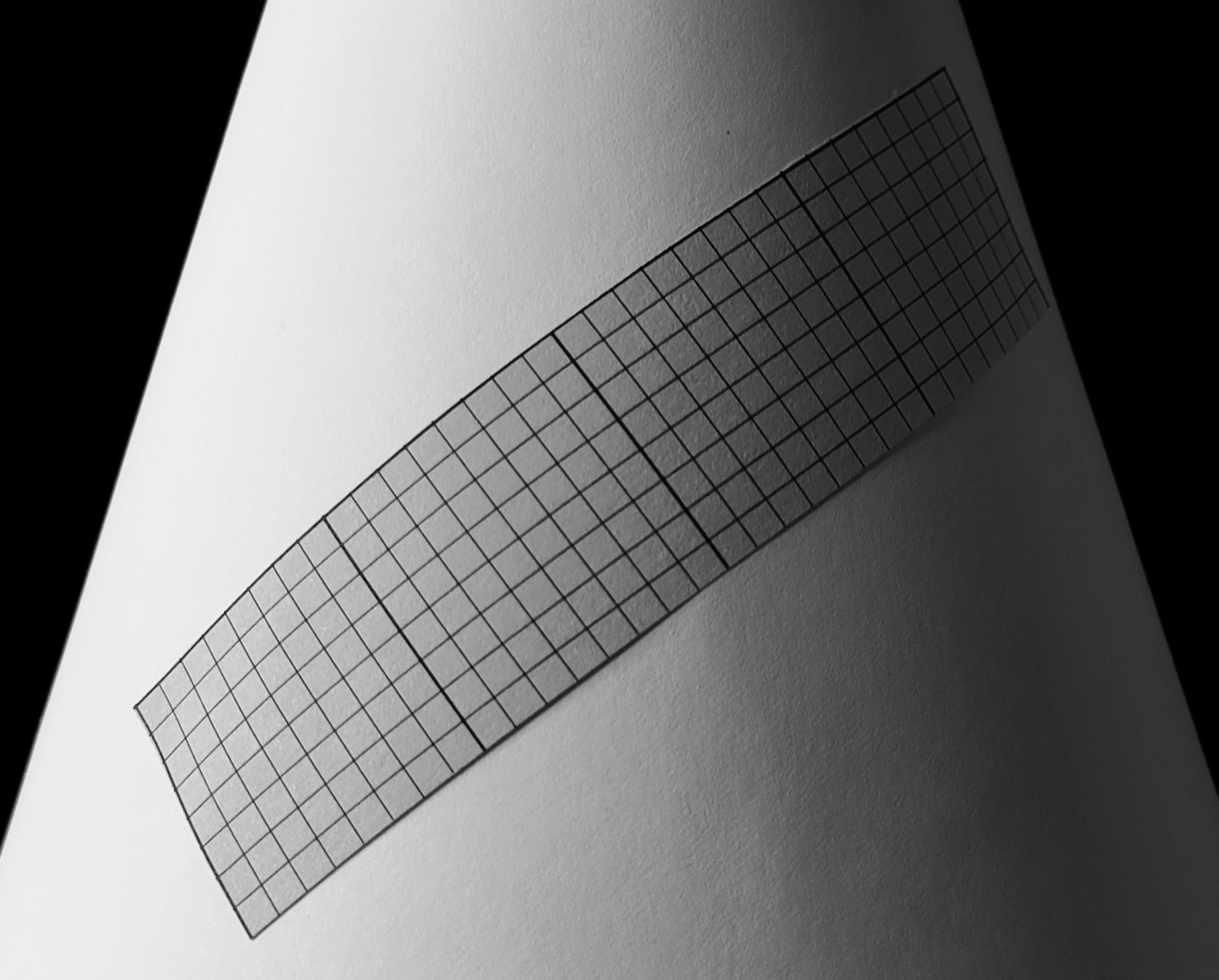}
  \end{minipage}%
  \begin{minipage}[t][3cm][b]{10cm}
    Image (a) shows a conically deformed strip. Figure (b) shows the reference domain $S=(0,4)\times (0,1)$ of the strip. The blue lines are asymptotic lines of the deformed strip and the red curves are integral curves of the principal direction along which the strip curved. Figure (c) illustrates the specific composite template that is assigned to the lower-left and upper-right subsquares of $S$.
\end{minipage}
\medskip

(b)%
  \begin{minipage}[t][1.8cm][b]{8.5cm}
    \includegraphics[width=8.5cm]{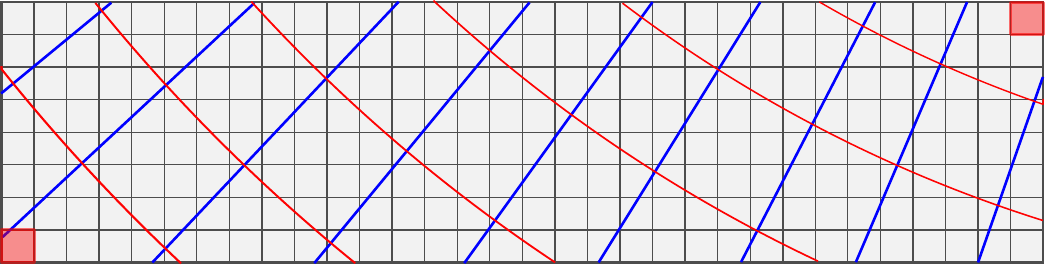}
  \end{minipage}%
  \ \ (c)%
  \begin{minipage}[t][1.8cm][b]{4.5cm}
    \includegraphics[width=4.5cm]{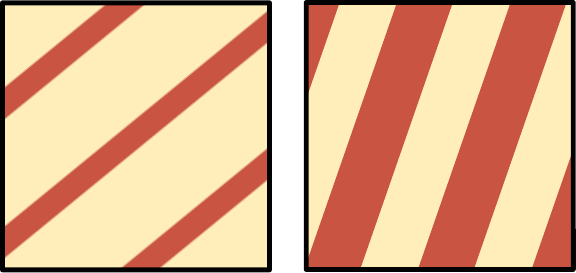}
  \end{minipage}%
  
    \caption{A finitely structured composite designed for a conically deformed strip}
    \label{fig:cone_strip}
\end{figure}

The idea of the proof of Theorem~\ref{T:shape} is as follows: Similarly to Section~\ref{S:textures1}, by the properties of the composite template and with the help of the transformation rule for in-plane rotations of Corollary~\ref{C:transf}, we find parameter fields $\kappa:S\to\mathcal K$ and $R:S\to\SO 2$ such that for each $x'\in S$ the composite template with parameter $\kappa(x')$ and in-plane rotation $R(x')$ leads to an algebraic minimization problem
\begin{equation*}
  \min_{G\in\R^{2\times 2}\atop\det G=0}Q^\gamma_{\hom}\Big(R(x'),\kappa(x');G-B^\gamma_{\eff}\big(R(x'),\kappa(x')\big)\Big)
\end{equation*}
which is uniquely minimized by $G=\II_*(x')$;
see Corollary~\ref{C:transf} for the definition of  $Q^\gamma_{\hom}(R,\kappa;\cdot)$ and $B^\gamma_{\eff}(R,\kappa)$.
In particular, we shall see that the associated bending energy is minimized by bending deformations $\deform\in H^2_{\iso}(S;\R^3)$ satisfying $\II_\deform=\II_*$.
In order to obtain a \textit{finitely} structured composite, we approximate the parameter fields $(R,\kappa)$ by fields $(\kappa_n,R_n)$ that are piecewise constant subordinate to a finite partition of $S$, and we shall prove that minimizers of the bending energy associated with $(\kappa_n,R_n)$ converge to a minimizer of the bending energy associated with $(\kappa,R)$. For this we use the following approximation result for bending energies. It shows that pointwise convergence of the quadratic form and of the prestrain implies $\Gamma$-convergence of the bending energy.
\begin{lemma}[$\Gamma$-convergence of parametrized bending energies]\label{L:gamma-spatial}
  For $n\in\N\cup\{\infty\}$ consider a Borel function $Q_n:S\times\R^{2\times 2}_{\sym}\to\R$ such that
  \begin{equation}\label{eq:Q2d}
    \alpha|G|^2\leq Q_n(x',G)\leq\beta|G|^2\qquad\text{for a.e.~}x'\in S\text{ and all }G\in \R^{2\times 2}_{\sym},
  \end{equation}
  where $0<\alpha\leq \beta$ are independent of $n$.
  Moreover, let $B_n\in L^2(S;\R^{2\times 2}_{\sym})$. Consider the functional
  \begin{equation*}
    \mathcal I_n:H^2_{\iso}(S;\R^3)\to\R,\qquad \mathcal I_n(\deform)\colonequals\int_SQ_n(x',\II(x')-B_n(x'))\,dx'.
  \end{equation*}
  Suppose that for $n\to \infty$, 
  \begin{alignat*}{2}
    Q_n(\cdot,G) & \to Q_\infty(\cdot,G) & \qquad & \text{a.e.~in $S$ and for all $G\in\R^{2\times 2}_{\sym}$},\\
    B_n & \to B_\infty & & \text{strongly in $L^2(S)$}.
  \end{alignat*}
  Then:
  \begin{enumerate}[(a)]
  \item \label{L:gamma-spatial:a} (Equicoercivity). Let $(\deform_n)$ have equibounded energy, i.e., $\limsup\limits_{n\to 0}\mathcal I_n(\deform_n)<\infty$. Then there exists $\deform\in H^2_{\iso}(S;\R^3)$ with $\fint_S \deform\,dx'=0$ such that
    \begin{equation*}
      \deform_n-\fint_S\deform_n\,dx'\wto \deform\qquad \text{weakly in $H^2(S;\R^3)$},
    \end{equation*}
    for a subsequence.
  \item \label{L:gamma-spatial:b} (\,$\Gamma$-convergence). For $n\to\infty$, the functional  $\mathcal I_n$ $\Gamma$-converges to $\mathcal I_\infty$ with respect to weak convergence in $H^2(S;\R^3)$.
  \item \label{L:gamma-spatial:c} (Strong convergence of minimizers). Let $(\deform_n)\subset H^2_{\iso}(S;\R^3)$ be a sequence such that for all $n\in\N$, $\deform_n$ is a minimizer of $\mathcal I_n$ subject to $\fint_S \deform_n=0$. Then, up to extraction of a subsequence, we have $\deform_n\to \deform_\infty$ strongly in $H^2(S;\R^3)$ where $\deform_\infty$ is a minimizer of $\mathcal I_\infty$.
  \end{enumerate}
\end{lemma}
\reftoproof{SS:shapedesign}

\begin{remark}
  Convergence of minimizers in the weak topology of $H^2(S;\R^3)$ follows directly from the general theory of $\Gamma$-convergence. The upgrade to strong convergence requires an additional argument that exploits the fact that $\mathcal I_n(\deform)$ is quadratic when seen as a function $\II_\deform$, see Step~3 in the proof of Lemma~\ref{L:gamma-spatial}. 
\end{remark}


\section{Proofs}
\label{S:proofs}
In the following we present the proofs of our results. We start in Section~\ref{SS:homogenization} to prove results concerning the definition of the effective quantities $(Q^{\gamma}_{\hom},B^\gamma_{\eff})$ and their representation via correctors. In Sections~\ref{SS:compactness} -- \ref{SS:upperbound} are devoted to the proofs of Theorems~\ref{T1} and \ref{T2}. Our results on the two-scale structure of the nonlinear strain are proven in Section~\ref{SS:chracterization}. In Section~\ref{SS:strongconvergence} we present the argument for the  results on the microstructure--properties--shape relation and in Section~\ref{SS:shapedesign} we prove the results on shape programming.

\subsection{Homogenization formula: Proofs of Lemmas~\ref{L:E}, \ref{L:3.3}, \ref{L:corrector},  Proposition~\ref{P:1}, and Lemma~\ref{L:cont} }\label{SS:homogenization}

\begin{proof}[Proof of Lemma~\ref{L:E}]
  We first note that since $Q\in\mathcal Q(\alpha,\beta)$ (see Definition~\ref{D1}), it suffices to prove
  \begin{equation*}
    \alpha\fint_{\Box_\Lambda}|\iota(x_3 G)|^2\leq\|P^{\gamma,\perp}_{{\rm rel},\Lambda}(\iota(x_3G))\|_{\Lambda}^2\leq\beta\fint_{\Box_\Lambda}|\iota(x_3 G)|^2.
  \end{equation*}
  Since $P^{\gamma,\perp}_{{\rm rel},\Lambda}$ is a projection, we have
  \begin{equation*}
    \|P^{\gamma,\perp}_{{\rm rel},\Lambda}(\iota(x_3G))\|_{\Lambda}^2\leq\|\iota(x_3 G)\|^2_{\Lambda},
  \end{equation*}
  and thus the upper bound follows. For the lower bound note that by definition of $\HH^\gamma_{{\rm rel},\Lambda}$ we have
  \begin{equation}\label{L:E:eq1}
    \|P^{\gamma,\perp}_{{\rm rel},\Lambda}(\iota(x_3G))\|_{\Lambda}^2=\inf_{M,\varphi}\|\iota(x_3G)+\iota(M)+\sym\nabla_\gamma\varphi\|_{\Lambda}^2,
  \end{equation}
  where the infimum is taken over all $M\in\R^{2\times 2}_{\sym}$ and $\varphi\in H^1_\gamma(\Box_\Lambda;\R^3)$. Since $Q\in\mathcal Q(\alpha,\beta)$, we have for all $M\in\R^{2\times 2}_{\sym}$ and $\varphi\in H^1_\gamma(\Box_\Lambda;\R^3)$,
  \begin{equation}\label{L:E:eq2}
    \|\iota(x_3G)+\iota(M)+\sym\nabla_\gamma\varphi\|_{\Lambda}^2\geq \alpha \fint_{\Box_\Lambda}|\iota(x_3G)+\iota(M)+\sym\nabla_\gamma\varphi|^2\geq\alpha \fint_{\Box_\Lambda}|\iota(x_3 G)|^2,
  \end{equation}
  where the last inequality holds thanks to the orthogonality
  \begin{equation*}
    \fint_{\Box_\Lambda}\iota(x_3G):(\iota(M)+\sym\nabla_\gamma\varphi)=0,
  \end{equation*}
  which can be checked by a direct computation. In combination with \eqref{L:E:eq1}, the lower bound follows.
\end{proof}

\begin{proof}[Proof of Lemma~\ref{L:3.3}]
  For the argument, fix $j\in J$, $x'\in S_j$, and note that for all $\chi\in\HH^\gamma_{\rm rel,\Lambda_j}$ we have
  \begin{eqnarray*}
    &&\fint_{\Box_{\Lambda_j}}Q\Big(x',x_3,y,\iota(x_3\II_\deform(x')))+\chi-B(x',x_3,y)\Big)\dd (x_3,y)\\
    &=&\|\iota(x_3\II_\deform(x'))+\chi-P^\gamma_{\Lambda_j}(\sym B(x',\cdot))\|^2_{\Lambda_j}+\|(I-P^\gamma_{\Lambda_j})(\sym B((x',\cdot))\|^2_{\Lambda_j},
  \end{eqnarray*}
  where $I$ denotes the identity operator on $\HH_{\Lambda_j}$.
  Moreover, by Lemma~\ref{L:E}, Definition~\ref{D:eff}, and Definition~\ref{def:Qgamma},
  \begin{eqnarray*}
    \inf_{\chi\in\HH^\gamma_{\rm rel,\Lambda_j}}\|\iota(x_3\II_\deform(x'))+\chi-P^\gamma_{\Lambda_j}(\sym B(x',\cdot)))\|^2_{\Lambda_j}
    &=&\|P^{\gamma,\perp}_{{\rm rel},\Lambda_j}\Big(\iota(x_3(\II_\deform(x')-B_{\eff}^\gamma(x')))\Big)\|^2_{\Lambda_j}\\
    &=&Q_{\hom}^\gamma(x',\II_\deform(x')-B_{\eff}^\gamma(x')).
  \end{eqnarray*}
  Now, the claim follows by combining the previous two identities, integration over $S_j$ and summation in $j\in J$.
\end{proof}

\begin{proof}[Proof of Lemma~\ref{L:corrector} and Proposition~\ref{P:1}]
  \step 1 {Representation of $\HH^{\gamma}_{{\rm rel},\Lambda}$}
  We claim that for all $\chi\in \HH^{\gamma}_{{\rm rel},\Lambda}$ there exists a unique pair $(M,\varphi)$ with $M\in \R^{2\times 2}_{\sym}$ and $\varphi\in H^1_\gamma(\Box_\Lambda;\R^3)$ such that
  \begin{equation}\label{pf:P:1:1}
    \chi=\iota(M)+\sym\nabla_\gamma\varphi\qquad\text{and}\qquad\fint_{\Box_\Lambda}\varphi=0,
  \end{equation}
  and
  \begin{equation}\label{pf:P:1:1b}
    |M|^2+\fint_{\Box_\Lambda}|\varphi|^2+|\nabla_\gamma\varphi|^2\leq C\|\chi\|_{\Lambda}^2,
  \end{equation}
  for a constant $C$ only depending on $\alpha, \gamma$ and  $C_{\Lambda}$.

  Here comes the argument: The existence of $(M,\varphi)$ is clear by definition.  Moreover we have  $\fint_{\Box_\Lambda}\iota(M):\sym\nabla_\gamma\varphi=0$, since $\varphi$ is in-plane periodic. Combined with the Korn inequality in Lemma~\ref{L:korn-lambda} and the Poincaré--Wirtinger inequality, the bound \eqref{pf:P:1:1b} follows. Uniqueness of the representation is a consequence of \eqref{pf:P:1:1b}.

  \step 2 {Proof of Lemma~\ref{L:corrector}}
  Let $G\in \R^{2\times 2}_{\sym}$. Let $\chi_G\in \HH^{\gamma}_{{\rm rel},\Lambda}$ denote the unique minimizer to
  \begin{equation*}
    \|\iota(x_3 G)+\chi_G\|_{\Lambda}^2=\min_{\chi\in \HH^{\gamma}_{{\rm rel},\Lambda}}\|\iota(x_3 G)+\chi\|^2_\Lambda,
  \end{equation*}
  and note that  $\chi_G$ is characterized by the variational problem
  \begin{equation}\label{pf:P:1:2}
    \big(\iota(x_3 G)+\chi_G,\chi\big)_\Lambda=0\qquad\text{for all }\chi\in \HH^{\gamma}_{{\rm rel},\Lambda},
  \end{equation}
  with $\big(\cdot,\cdot\big)_{\Lambda}$ denoting the inner product of $\HH_{\Lambda}$.
  Since projections are contractions and $Q\in\mathcal Q(\alpha,\beta)$, we have
  \begin{equation}\label{pf:P:1:3}
    \|\chi_G\|_{\Lambda}\leq \|\iota(x_3 G)\|_{\Lambda}\leq \sqrt{\tfrac{\beta}{12}}|G|.
  \end{equation}
  For future reference, we also note that we have the identity
  \begin{equation}\label{pf:P:1:4}
    \iota(x_3 G)+\chi_G=P^{\gamma,\perp}_{\rm rel,\Lambda}(\iota(x_3 G))=\mathbf E^\gamma_{\Lambda}(G).
  \end{equation}
  Now, let $(M_G,\varphi_G)$ denote the unique pair associated with $\chi_G$ via \eqref{pf:P:1:1}. Then, we deduce from \eqref{pf:P:1:2} that $(M_G,\varphi_G)$ is characterized by \eqref{eq:corrector_equation}, and we obtain \eqref{eq:corrector_apriori} from \eqref{pf:P:1:3} and \eqref{pf:P:1:1b}.

  \step 3 {Proof of Proposition~\ref{P:1}~\ref{P:1:b}}
  First, we recall that the definition of $Q_{\hom}^\gamma$ can be rephrased as $Q_{\hom}^\gamma(G)=\|P^{\gamma,\perp}_{\rm rel,\Lambda}(\iota(x_3 G))\|_{\Lambda}$ for all $G\in\R^{2\times 2}_{\sym}$.
  Using the orthonormal basis $G_1,G_2,G_3$, we write $G=\widehat G_iG_i$, where here and below we appeal to Einstein's summation convention. By \eqref{pf:P:1:4} and linearity we have $P^{\gamma,\perp}_{\rm rel,\Lambda}(\iota(x_3G))=\widehat G_i(\iota(x_3G_i)+\chi_{G_i})$, and thus
  \begin{equation*}
        Q_{\hom}^\gamma(G)=\widehat G_i\widehat G_k\Big(\iota(x_3 G_i)+\chi_{G_i},\,\iota(x_3G_k)+\chi_{G_k}\Big)_{\Lambda}=\widehat G_i\widehat G_k\Big(\iota(x_3 G_i)+\chi_{G_i},\,\iota(x_3G_k)\Big)_{\Lambda}=\widehat G_i\widehat G_k\widehat Q_{ik},
  \end{equation*}
  where the last two identities hold thanks to \eqref{pf:P:1:2} and the definition of $\widehat Q$, respectively. The identity shows that $\widehat Q$ is symmetric.
  The bound \eqref{P:1:coercivityQhat} is a direct consequence of \eqref{eq:bounds}.

 \step 4 {Proof of Proposition~\ref{P:1}~\ref{P:1:c}}
 Set $\tilde B\colonequals\sum_{i=1}^3\big(\widehat Q^{-1}\widehat B\big)_iG_i$. In view of the definition of $B^\gamma_{\eff}$, cf.~\eqref{eq:def:beff}, and since $\mathbf E^\gamma_{\Lambda}$ is injective, it suffices to show that
 \begin{equation}\label{pf:P1:4}
   \mathbf E^\gamma_{\Lambda}(\tilde B)=P^{\gamma,\perp}_{\rm rel,\Lambda}(\sym B).
 \end{equation}
 Since $P^{\gamma,\perp}_{\rm rel,\Lambda}$ is a projection onto $\HH^{\gamma,\perp}_{\rm rel,\Lambda}$, and since the latter is spanned by the fields $\iota(x_3G_i)+\chi_{G_i}$, the above identity is equivalent to
 \begin{equation}\label{pf:P1:5}
   \big(\mathbf E^\gamma_{\Lambda}(\tilde B),\iota(x_3G_i)+\chi_{G_i}\big)_{\Lambda}=\big(\sym B,\iota(x_3G_{G_i})+\chi_{G_i}\big)_{\Lambda}
 \end{equation}
 for all $i=1,2,3$. By definition, the right-hand side is equal to $\widehat B_i$. On the other hand, by \eqref{pf:P:1:4} and linearity, we have $\mathbf E^\gamma_{\Lambda}(\tilde B)=\sum_{k=1}^3\big(\widehat Q^{-1}\widehat B\big)_k(\iota(x_3 G_k)+\chi_{G_k})$. Thus, \eqref{pf:P1:5} is equivalent to $\big(\widehat Q^{-1}\widehat B\big)_i=\widehat B_i$, and we conclude \eqref{pf:P1:4}.
\end{proof}

\begin{proof}[Proof of Lemma~\ref{L:cont}]
  It suffices to prove \eqref{pf:L:cont:eq1} and \eqref{pf:L:cont:eq2}. The other statements then follow with help of Proposition~\ref{P:1}.

  We first note that there exists $\tilde M$ and $\tilde\varphi$ such that for a subsequence (not relabeled) we have
  \begin{equation}\label{pf:L:cont:eq3}
    M_{n,i}\to \tilde M\qquad\text{and}\qquad \varphi_{n,i}\wto\tilde\varphi\text{ weakly in }H^1_{\gamma}(\Box_\Lambda;\R^3).
  \end{equation}
  Indeed, this follow from \eqref{eq:corrector_apriori}, the variant of Korn's inequality of Lemma~\ref{L:korn-lambda} and weak compactness of bounded sequences in Hilbert spaces. Thanks to the assumed convergence of $(Q_n)$, we may pass to the limit in the corrector equation \eqref{eq:corrector_equation} for $(M_{n,i},\varphi_{n,i})$ and deduce that
    \begin{equation*}
      \int_{\Box_\Lambda}\mathbb L_\infty\big(\iota(x_3 G_i+\tilde M)+\sym(\nabla_\gamma\tilde \varphi)\big):\big(\iota(M')+\sym(\nabla_\gamma\varphi')\big)=0
    \end{equation*}
    for all test functions $M'$ and $\varphi'$.
    Above, $\mathbb L_\infty$ denotes the fourth-order tensor associated with $Q_\infty$ via the polarization identity \eqref{eq:polarization}.
    By uniqueness of the corrector (see Lemma~\ref{L:corrector}), we conclude that $\tilde M=M_{\infty,i}$ and $\tilde\varphi=\varphi_{\infty,i}$, and thus the convergence \eqref{pf:L:cont:eq3} holds for the entire sequence. It remains to show that $\varphi_{n,i}\to\varphi_{\infty,i}$ strongly in $H^1$. In view of Korn's inequaltiy, cf.~Lemma~\ref{L:korn-lambda}, it suffices to show that $\|\sym\nabla_\gamma\varphi_{n,i}-\sym\nabla_\gamma\varphi_{\infty,i}\|_{L^2(\Box_\Lambda)}\to 0$.
    Since $Q_n\in\mathcal Q(\alpha,\beta)$ and by appealing to the corrector equations \eqref{eq:corrector_equation} for $\varphi_{n,i}$ and $\varphi_{\infty,i}$, we have
    \begin{eqnarray*}
      &&\alpha\fint_{\Box_\Lambda}|\sym\nabla_\gamma\varphi_{n,i}-\sym\nabla_\gamma\varphi_{\infty,i}|^2\\
      &\leq& \fint_{\Box_\Lambda}\mathbb L_n\Big(\iota(M_n-M_\infty)+\sym\nabla_\gamma(\varphi_{n,i}-\varphi_{\infty,i})\Big):\Big(\iota(M_n-M_\infty)+\sym\nabla_\gamma(\varphi_{n,i}-\varphi_{\infty,i})\Big)\\
      &=& \fint_{\Box_\Lambda}(\mathbb L_\infty-\mathbb L_n)\Big(\iota(x_3G_i+M_\infty)+\sym\nabla_\gamma\varphi_{\infty,i})\Big):\Big(\iota(M_\infty-M_n)+\sym\nabla_\gamma(\varphi_{n,i}-\varphi_{\infty,i})\Big).
    \end{eqnarray*}
    Note that the last integral converges to $0$, since we have $(\mathbb L_\infty-\mathbb L_n)\Big(\iota(x_3G_i+M_\infty)+\sym\nabla_\gamma\varphi_{\infty,i})\Big)\to 0$ strongly in $L^2$, and  $\iota(M_\infty-M_n)+\sym\nabla_\gamma(\varphi_{n,i}-\varphi_{\infty,i})\wto 0$ weakly in $L^2$.
\end{proof}

\subsection{Compactness:  Proof of Theorem~\ref{T1} (a)}\label{SS:compactness}

\begin{proof}[Proof of Theorem~\ref{T1}~\ref{item:T1:compactness}]
By the triangle inequality and Young's inequality, for all $F\in\R^{3\times3}$ and $h>0$, we have
\begin{equation}\label{eq:comp}
  \frac12\dist(F,\SO3)\leq \dist\Big(F(I_{3\times 3}-hB_{\e(h),h}),\SO 3\Big)+\frac12h^2|B_{\e(h),h}|^2+h|B_{\e(h),h}|
\end{equation}
almost everywhere in $\Omega$. Together with Assumption~\ref{ass:W}~\ref{ass:W:4} this yields the implication
\begin{equation*}
  \limsup\limits_{h\to 0}\mathcal I^{\e(h),h}(\deform_h)<\infty\quad\Rightarrow\quad  \limsup\limits_{h\to 0}\frac1{h^2}\int_\Omega \dist^2(\nabla_h \deform_h(x),\SO 3)\dd x<\infty.
\end{equation*}
Hence, the statement of Theorem~\ref{T1}~\ref{item:T1:compactness} follows from  \cite[Theorem~4.1]{friesecke2002theorem}.
\end{proof}

\subsection{Lower bound: Proof of Theorem~\ref{T1} (b)}\label{SS:lowerbound}
We only need to consider the case when \eqref{eq:equibounded} is satisfied since otherwise, the claim is trivial. We note that by Theorem~\ref{T1}~\ref{item:T1:compactness}  we have $\deform\in H^2_{\iso}(S;\R^3)$. In view of Lemma~\ref{L:3.3} it suffices to prove that
\begin{equation*}
  \liminf\limits_{h\to 0}\mathcal I^{\e(h),h}(\deform_h)\geq   \widetilde{\mathcal I}^\gamma(\deform),
\end{equation*}
where $\widetilde{\mathcal I}^\gamma$ is defined in \eqref{def:Igamma}.
For the proof we modify the argument in \cite[Proof of Theorem~3.3 (Lower bound)]{neukamm2012rigorous} as follows:  Without loss of generality, we may assume that
\begin{equation*}
  \liminf\limits_{h\to 0}\mathcal I^{\e(h),h}(\deform_h)=\lim\limits_{h\to 0}\mathcal I^{\e(h),h}(\deform_h)<\infty.
\end{equation*}
Hence, by Theorem~\ref{T1}~\ref{item:T1:compactness}, $(\deform_h)$ is a sequence of finite bending energy in the sense of \eqref{eq:FBE}. By appealing to Proposition~\ref{P:chartwoscale}~\ref{item:char_twoscale_a} we can pass to a subsequence  (that we do not relabel) such that
\begin{equation}
  \label{pf:eq1}
  E_h(\deform_h)\wtto E\colonequals\iota(x_3\II_\deform+M)+\sym\nabla_\gamma\varphi,
\end{equation}
for some $M\in L^2(S;\R^{2\times 2}_{\sym})$ and $\varphi\in L^2(S;H^1_{\gamma,{\rm uloc}})$ that is locally periodic in the sense of \eqref{eq:locper}. Set
\begin{equation}\label{eq:thetah}
  \theta_h(x)\colonequals
  \begin{cases}
    1&\text{if }\dist(\nabla_h\deform_h(x),\SO 3)\leq h^{1/2},\\
    0&\text{else,}
  \end{cases}
\end{equation}
and note that \eqref{eq:FBE} yields
\begin{equation}\label{eq:thetabound}
  \limsup\limits_{h\to 0}\frac{1}{h}\int_\Omega |1-\theta_h|\leq \limsup\limits_{h\to \infty}\frac{1}{h^2}\int_\Omega\dist^2(\nabla_h\deform_h(x),\SO 3)\dd x<\infty.
\end{equation}
By polar factorization (of $\nabla_h\deform_h(x)$ for $x\in\{\theta_h=1\}$) and by Assumption~\ref{ass:W}~\ref{ass:W:4}, there exists a rotation field $R_h:\Omega\to\SO 3$ such that $\theta_h  \nabla_h\deform_h=\theta_hR_h\big(I_{3\times 3}+hE_h(\deform_h))$. Thus, $\theta_h\nabla_h\deform_h(I_{3\times 3}-hB_{\e(h),h})=\theta_h R_h(I_{3\times 3}+hG_h)$ where
\begin{equation}\label{pf:lb:unifgh}
  G_h\colonequals\frac{1}{h}\Big((I_{3\times 3}+h \theta_hE_h(\deform_h))(I_{3\times 3}-h  \theta_hB_{\e(h),h})\,-I_{3\times 3}\Big)\to 0\text{ uniformly in }\Omega.
\end{equation}
Frame-indifference~\ref{item:nonlinear_material_w1} and the natural state condition~\ref{item:nonlinear_material_w3} thus yield
\begin{equation*}
  \frac{\theta_h}{h^2}W_{\e(h)}(x,\nabla_h\deform_h(I_{3\times 3}-h B_{\e(h),h}))=\frac{1}{h^2}W_{\e(h)}(x,I_{3\times 3}+h G_h(x)).
\end{equation*}
Combined with~\ref{item:nonlinear_material_w4} and \eqref{pf:lb:unifgh},  we get
\begin{equation}\label{pf:lb:taylor}
  \begin{aligned}
  \limsup\limits_{h\to 0}\Big|\frac{1}{h^2}\int_\Omega \theta_hW_{\e(h)}(x,\nabla_h\deform_h(I_{3\times 3}-h B_{\e(h),h}))\dd x\\
  -\int_\Omega Q_{\e(h)}\Big(x,\theta_h\big(E_h(\deform_h)-B_{\e(h),h}\big)\Big)\dd x \Big|=0,
\end{aligned}
\end{equation}
and thus (by non-negativity of $W_{\e(h),h}$),
\begin{equation}\label{pf:lb:2}
  \liminf\limits_{h\to 0}\mathcal I^{\e(h),h}(\deform_h)\geq
  \liminf\limits_{h\to 0}\int_\Omega Q_{\e(h)}\Big(x,\theta_h\big(E_h(\deform_h)-B_{\e(h),h}\big)\Big)\dd x.
\end{equation}
By appealing to Assumption~\ref{ass:W}~\ref{ass:W:1}  and Assumption~\ref{ass:W}~\ref{ass:W:4} , the right-hand side is equal to

\begin{equation}\label{pf:eq2}
  \liminf\limits_{h\to 0}\int_\Omega Q\Big(x,\tfrac{x'}{\e(h)},\theta_h\big(E_h(\deform_h)-B(x,\tfrac{x'}{\e(h)})\big)\Big)\dd x.
\end{equation}
We claim that
\begin{equation}\label{pf:eq3}
  \theta_h\big(E_h(\deform_h)-B(x,\tfrac{x'}{\e(h)})\big)\wtto E-B\qquad\text{weakly two-scale in }L^2,
\end{equation}
where $E$ is defined in \eqref{pf:eq1}.
For the argument, we first note that $E_h(\deform_h)-B(x,\tfrac{x'}{\e(h)})$ weakly two-scale converges to the right-hand side thanks to \eqref{pf:eq1} and Assumption~\ref{ass:W}~\ref{ass:W:4}. The left-hand side of \eqref{pf:eq3} converges to the same limit, since $(\theta_h)$ is bounded in $L^\infty$ and $\theta_h\to 1$ in $L^1$ by \eqref{eq:thetabound}.
\smallskip

In view of \eqref{pf:lb:2}, \eqref{pf:eq2}, \eqref{pf:eq3}, and \eqref{pf:eq1}, the two-scale lower semicontinuity result Lemma~\ref{L:twoscale:lsc} yields the lower bound
\begin{equation*}
  \liminf\limits_{h\to 0}\mathcal I^{\e(h),h}(\deform_h)\geq \sum_{j\in J}\int_{S_j}\fint_{\Box_{\Lambda_j}}Q(x,y,\iota(x_3\II_\deform+M)+\sym\nabla_\gamma\varphi-B)\dd (x_3,y)\dd x'.
\end{equation*}
In view of the definition of $\widetilde I^\gamma$ (see \eqref{L:3.3}) the claim follows by taking the infimum on the right-hand side over $M$ and $\varphi$.
\qed

\subsection{Recovery sequence: Proofs of Theorem~\ref{T1} (c) and Theorem~\ref{T2}}
\label{SS:upperbound}
We only discuss the proof of Theorem~\ref{T2}, since, in comparison to Theorem~\ref{T1}~\ref{item:T1:recovery_sequence}, it features the additional difficulty to take care of the clamped boundary conditions \eqref{eq:BC3d}.

\begin{proof}[Proof of Theorem~\ref{T2}]

  \step 1 {A priori estimate on $M$ and choice of $\delta$}
Let $\tilde \deform\in H^2_{\iso}(S;\R^3)$. Set $\tilde\II\colonequals\II_{\tilde \deform}$, and denote by $(\tilde M,\tilde \varphi)$ the associated pair of correctors satisfying
  \begin{equation*}
    \mathcal I^\gamma_{\hom}(\tilde \deform)=\sum_{j=1}^\infty\int_{S_j}\fint_{\Box_{\Lambda_j}}Q(x',x_3,y,\iota(x_3\tilde\II+\tilde M)+\sym\nabla_\gamma\tilde \varphi-\widetilde B)\dd (x_3,y)\dd x',
  \end{equation*}
  where
  \begin{equation}
    \label{eq:12314}
    \widetilde B(x',\cdot)\colonequals P^\gamma_{\Lambda_j}(\sym B(x',\cdot))\qquad\text{for $x'\in S_j$ and $j\in J$.}
  \end{equation}
  We claim that there exists a constant $C=C(\alpha,\beta,\gamma,C_{J},C_B)$ such that
  \begin{equation}\label{pf:recov:eq:0}
    |\tilde M(x')|^2\leq C(|\tilde\II(x')|^2+1).
  \end{equation}
  Before we prove this estimate, we note that this implies that
  \begin{equation}\label{pf:recov:eq:step1}
    \tilde M(x')+\delta I_{2\times 2}\geq 0\text{ a.e. in }\{\tilde\II=0\},
  \end{equation}
  for a constant $\delta=\delta(\alpha,\beta,\gamma,C_{J},C_B)\geq 0$ that is idependent of the specific isometry $\tilde \deform$.
  We now prove \eqref{pf:recov:eq:0} and first argue that for all $j\in J$,
  \begin{equation}\label{pf:recov:eq:1}
    \fint_{\Box_{\Lambda_j}}|\widetilde B(x',x_3,y)|^2\dd (x_3,y)\lesssim 1\qquad\text{for a.e. }x'\in S_j;
  \end{equation}
  where here and below, $\lesssim$ means $\leq$ up to a multiplicative constant that only depends  $\alpha,\beta,\gamma,C_{J},C_B$.
  Indeed, since $Q(x,y,\cdot)\in\mathcal Q(\alpha,\beta)$, and since $P^\gamma_{\Lambda_j}$ is a projection (adapted to $Q(x',\cdot)$), we have
  \begin{eqnarray*}
    \fint_{\Box_{\Lambda_j}}|\widetilde B(x',x_3,y)|^2\dd (x_3,y)
    &\lesssim& \fint_{\Box_{\Lambda_j}}Q(x',\widetilde B(x',x_3,y))\dd (x_3,y)\\
    &\lesssim& \fint_{\Box_{\Lambda_j}}Q(x',\sym B(x',x_3,y))\dd (x_3,y)\\
    &\lesssim&\fint_{\Box_{\Lambda_j}}|\sym B(x',x_3,y)|^2\dd (x_3,y),
  \end{eqnarray*}
  and thus the claim follows from \eqref{ass:unifB}.
  Next, we note that by orthogonality,
  \begin{eqnarray*}
    |\tilde M(x')|^2&\leq&     \fint_{\Box_{\Lambda_j}}|\iota(x_3\tilde \II+\tilde M)+\sym\nabla_\gamma\tilde\varphi|^2\\
                                                  &\leq&     2\fint_{\Box_{\Lambda_j}}|\iota(x_3\tilde\II+\tilde M)+\sym\nabla_\gamma\tilde\varphi-\widetilde B|^2+2\fint_{\Box_{\Lambda_j}}|\widetilde B|^2.
  \end{eqnarray*}
  With \eqref{pf:recov:eq:1} and since $Q\in\mathcal Q(\alpha,\beta)$, we conclude
  \begin{eqnarray*}
    |\tilde M(x')|^2&\lesssim& \fint_{\Box_{\Lambda_j}}Q(x',x_3,y,\iota(x_3\tilde \II+\tilde M)+\sym\nabla_\gamma\tilde \varphi-\widetilde B)\dd (x_3,y)+1.
  \end{eqnarray*}
  Since $(\tilde M(x'),\tilde \varphi(x',\cdot))$ is the corrector pair, we conclude that
  \begin{eqnarray*}
    |\tilde M(x')|^2&\lesssim& Q^\gamma_{\rm \hom}(x',\tilde \II(x'))+1\lesssim |\tilde\II(x')|^2+1,
  \end{eqnarray*}
  as claimed.

  \step 2 {Construction in the smooth case}
  Let $\delta\geq 0$ be as in \eqref{pf:recov:eq:step1} and denote by
  \begin{equation*}
    \mathcal A_{BC}^h\colonequals\{\deform_h\in H^1(\Omega;\R^3)\,:\,\deform_h\text{ satisfies \eqref{eq:BC3d}}\},
  \end{equation*}
  the space of deformations that satisfy the 3d-boundary conditions. We claim that for all $\rho>0$ there exists $\deform_\rho\in \mathcal A_{BC}\cap C^\infty(\bar S;\R^3)$ with
  \begin{equation}\label{pf:recov:step2:1}
    \|\deform_{\rho}-\deform\|_{L^2(S)}+|\mathcal I^\gamma_{\hom}(\deform_\rho)-\mathcal I^\gamma_{\hom}(\deform)|\leq\rho,
  \end{equation}
  and a recovery sequence $(\deform_{h,\rho})$ with $\deform_{h,\rho}\in\mathcal A_{BC}^h$ satisfying the claim of Theorem~\ref{T2} with $\deform$ replaced by $\deform_{\rho}$. Here comes the argument: First, since $\mathcal I^\gamma_{\hom}$ is continuous w.r.t.~strong convergence in $H^2(S;\R^3)$ and in view of the approximation result \eqref{ass:BC:dense}, we can find $\deform_{\rho}$ such that \eqref{pf:recov:step2:1} holds. Let $(M_\rho,\varphi_\rho)$ denote the pair of correctors associated with $\deform_\rho$, that is,
  \begin{equation*}
    \mathcal I^\gamma_{\hom}(\deform_\rho)=\sum_{j=1}^\infty\int_{S_j}\fint_{\Box_{\Lambda_j}}Q(x',x_3,y,\iota(x_3\II_{\deform_\rho}+M_\rho)+\sym\nabla_\gamma\varphi_{\rho}-\widetilde B)\dd (x_3,y)\dd x',
  \end{equation*}
  where $\widetilde B$ is defined in \eqref{eq:12314}.
  By \eqref{pf:recov:eq:step1} we have $M_\rho+\delta I_{2\times 2}\geq 0$ on $\{\II_{\deform_\rho}=0\}$. Therefore, we can apply Proposition~\ref{P:chartwoscale} to obtain a recovery sequence $\deform_{h,\rho}\in \mathcal A_{BC}^h$ such that $\deform_{h,\rho}\to \deform_\rho$ in $L^2(\Omega)$,
  \begin{equation*}
    E_h(\deform_{h,\rho})\stto \iota(x_3\II_{\deform_\rho}+M_\rho)+\sym\nabla_\gamma\varphi_\rho\text{ strongly two-scale,}
  \end{equation*}
  and in addition
  \begin{equation*}
    \limsup\limits_{h\to 0}h\|E_h(\deform_{h,\rho})\|_{L^{\infty}}=\,0\quad\text{ and }\quad\lim_{h\to 0}\|\det(\nabla_h\deform_{h,\rho})-1\|_{L^{\infty}}=0.
  \end{equation*}
  Note that the latter implies that (for $h$ sufficiently small) there exists a rotation field $R_{h,\rho}\in L^\infty(\Omega;\SO 3)$ such that
  \begin{equation*}
    \nabla_h \deform_{h,\rho}(I_{3\times 3}-h B_{\e(h),h})=R_{h,\rho}(I_{3\times 3}+h E_h(\deform_{h,\rho}))(I_{3\times 3}-h B_{\e(h),h}).
  \end{equation*}
  Moreover, with Assumption~\ref{ass:W}~\ref{ass:W:4},
  \begin{equation*}
    \nabla_h \deform_{h,\rho}(I_{3\times 3}-h B_{\e(h),h})=R_{h,\rho}(I_{3\times 3}+h G_h+ho_h),\qquad G_h=E_h(\deform_{\rho,h})-B_{\e(h),h},
  \end{equation*}
  for some remainder $o_h:\Omega\to\R^{3\times 3}$ and $h(\|G_h\|_{L^\infty(\Omega)}+\|o_h\|_{L^\infty(\Omega)})\to 0$. We thus conclude by frame-indifference and a Taylor expansion at $I_{3\times 3}$ of $W_{\e(h)}$ that
  \begin{equation*}
    \limsup\limits_{h\to\infty}\mathcal I^{\e(h),h}(\deform_{\rho,h})=\limsup\limits_{h\to\infty}\int_\Omega Q_{\e(h)}(x,G_h(x))\dd x.
  \end{equation*}
  In view of Assumption~\ref{ass:W}~\ref{ass:W:2}, we have
  \begin{equation}\label{pf:recov:step2:eq3}
    \limsup\limits_{h\to\infty}\int_\Omega Q_{\e(h)}(x,G_h(x))\dd x= \limsup\limits_{h\to\infty}\int_\Omega Q(x',x_3,\tfrac{x'}{\e(h)},G_h(x))\dd x,
  \end{equation}
  and in view of Assumption~\ref{ass:W}~\ref{ass:W:4}, we have
  \begin{equation*}
    G_h\stto \iota(x_3\II_{\deform_\rho}+M_\rho)+\sym\nabla_\gamma\varphi_\rho-B\qquad\text{strongly two-scale in }L^2.
  \end{equation*}
  By the continuity of the convex functional on the right-hand side of \eqref{pf:recov:step2:eq3} w.r.t.~strong two-scale convergence, cf.~Lemma~\ref{L:twoscale:lsc}, we conclude that
  \begin{align*}
    &\limsup\limits_{h\to\infty}\mathcal I^{\e(h),h}(\deform_{\rho,h})\\
    =\,&\,\sum_{j=1}^\infty\int_{S_j}\fint_{\Box_{\Lambda_j}}Q(x',x_3,y,\iota(x_3\II_{\deform_\rho}+M_\rho)+\sym\nabla_\gamma\varphi_{\rho}-\sym B)\dd (x_3,y)\dd x'\\
    =\,&\,\mathcal I^\gamma_{\hom}(\deform_\rho)+\mathcal I^\gamma_{\rm res}(B),
  \end{align*}
  see Lemma~\ref{L:3.3}.
  This concludes the argument.

  \step 3 {Extraction of a diagonal sequence}
  Consider
  \begin{equation*}
    c_{h,\rho}\colonequals|\mathcal I^{\e(h),h}(\deform_{h,\rho})-(\mathcal I^\gamma_{\hom}(\deform)+\mathcal I^\gamma_{\rm res}(B))|+\int_\Omega|\deform_{h,\rho}-\deform|^2.
  \end{equation*}
  Step~2 yields $\limsup\limits_{h\to 0}\limsup\limits_{\rho\to 0}c_{h,\rho}=0$, and thus there exists a diagonal sequence $h\mapsto \rho(h)$ such that $c_{h,\rho(h)}\to 0$. Hence, $\deform_h\colonequals \deform_{h,\rho(h)}$ defines the recovery sequence $(\deform_h)$ we were looking for.
\end{proof}

\subsection{Characterization of the limiting strain: Proof of Proposition~\ref{P:chartwoscale} and Lemma~\ref{L:3.3}}
\label{SS:chracterization}

We start with the proof of Proposition~\ref{P:chartwoscale}~\ref{item:char_twoscale_a}, which is almost a direct consequence of \cite[Proposition 3.2]{hornung2014derivation}, where the ``single-grain'' case with $\Lambda_j=I_{2\times 2}$ is considered.

\begin{proof}[Proof of Proposition~\ref{P:chartwoscale}~\ref{item:char_twoscale_a}]
  By Lemma~\ref{L:twoscale:prop}~\ref{L:twoscale:prop:a} we may pass to a subsequence (not relabled) such that
  \begin{equation*}
    E_h(\deform_h)\wtto E\qquad\text{weakly two-scale in }L^2,
  \end{equation*}
  where $E\in L^2(\Omega;L^2_{\uloc}(\R^2;\R^{3\times 3}_{\sym}))$. In view of Lemma~\ref{L:twoscale:parad}, it suffices to show that for all $j\in J$, we have (after possibly passing to a further subsequence)
  \begin{equation*}
    E_h(\deform_h)\wtto \iota(x_3\II+M_j)+\sym\nabla_\gamma\varphi_j,
  \end{equation*}
  weakly two-scale in the sense of \eqref{eq:L:twoscale:parad} for some $M_j\in L^2(S_j;\R^{2\times 2}_{\sym})$ and $\varphi_j\in L^2(S_j;H^1_\gamma(\Box_{\Lambda_j};\R^3))$. In the case $\Lambda_j=I_{2\times 2}$, \eqref{eq:L:twoscale:parad} directly follows from \cite[Proposition 3.2]{hornung2014derivation}. With help of Proposition~\ref{P:twoscale:gradrecov} and the variant of Korn's inequality Lemma~\ref{L:korn-lambda}, the argument of \cite[Proposition 3.2]{hornung2014derivation} extends verbatim to the case of general invertible $\Lambda_j\in\R^{2\times 2}$.
\end{proof}

We turn to the proof of Proposition~\ref{P:chartwoscale}~\ref{item:char_twoscale_b}. We only discuss the construction in the case with prescribed boundary data, since this adds an additional layer of difficulties. The construction is based on the density of $\mathcal A_{BC}\cap C^\infty(\bar S;\R^3)$ in $\mathcal A_{BC}$, see~\eqref{ass:BC:dense}, and a ``single-scale'' approximation result that we recently obtained together with M.~Griehl et al. in \cite{GNPG1}. We recall it in the following (slightly weaker) form:
\begin{proposition}[{cf.~\cite[Proposition~4.4]{GNPG1}}]\label{P:general_construction}
  Let $\deform\in H^2_\mathrm{ iso}(S;\R^3)\cap C^\infty(\overline S;\R^3)$, $M\in L^2(S;\R^{2\times 2}_{\sym})$, $d\in L^2(\Omega;\R^3)$ and assume that there exists $\delta\geq 0$ such that
  \begin{equation*}
    M+\delta I_{2\times 2}\geq 0\text{ a.e.~in }\{\II_\deform=0\}.
  \end{equation*}
  Then there exists a sequence $(\deform_h)\subseteq C^{\infty}(\bar\Omega;\R^3)$ such that
  \begin{equation}\label{eq:recov:onescale}
    \begin{aligned}
      \deform_h\to\, &\deform\qquad\text{uniformly in }\Omega,\\
      E_h(\deform_h)\,\to\,&\,\iota(x_3 \II+M)+\sym(d\otimes e_3)\quad\text{strongly in }L^2(\Omega),\\
      \limsup\limits_{h\to 0}h\|E_h(\deform_h)\|_{L^\infty(\Omega)}=\,&0\quad\text{ and }\quad\lim_{h\to 0}\|\det(\nabla_h\deform_h)-1\|_{L^\infty(\Omega)}=0.
    \end{aligned}
  \end{equation}
  Furthermore, if Assumption~\ref{ass:BC} is satisfied and $\deform\in\mathcal A_{BC}\cap C^\infty(\overline S;\R^3)$, then we may additionally enforce the clamped, affine boundary condition \eqref{eq:BC3d}.
\end{proposition}
\begin{proof}[Proof of Proposition~\ref{P:chartwoscale}~\ref{item:char_twoscale_b}]
  By \eqref{ass:BC:dense}, for any $\rho$ there exists $\deform_\rho\in\mathcal A_{BC}\cap C^\infty(\bar S;\R^3)$ such that
  \begin{equation*}
    \int_S|\deform-\deform_\rho|^2+|\II_\deform-\II_{\deform_\rho}|^2+|(\nabla'\deform,b_\deform)-R_\rho|^2\leq \rho,\qquad\text{where }R_\rho\colonequals(\nabla'\deform_\rho,b_{\deform_\rho}).
  \end{equation*}
  We claim that there exists a sequence $(\tilde \deform_{h,\rho})\subset C^1(\bar\Omega;\R^3)$ that satisfies the clamped boundary conditions \eqref{eq:BC3d} and
  \begin{subequations}
    \begin{align}    \label{eq:Pid:2a}
      \tilde \deform_{h,\rho}\to\, &\deform_\rho\qquad\text{uniformly in }\Omega,\\
      \label{eq:Pid:2b}
      E_h(\tilde \deform_{h,\rho})\,\stto\,&\,\iota(x_3 \II_{\deform_\rho}+M)+\sym\nabla_\gamma\varphi\quad\text{strongly two-scale in }L^2(\Omega),\\
      \label{eq:Pid:2c}
      \limsup\limits_{h\to 0}h\|E_h(\tilde \deform_{h,\rho})\|_{L^\infty(\Omega)}=\,&0\quad\text{ and }\quad\lim_{h\to 0}\|\det(\nabla_h\tilde \deform_{h,\rho})-1\|_{L^\infty(\Omega)}=0.
    \end{align}
  \end{subequations}
  Before we prove this claim we note that with the sequence $(\tilde \deform_{h,\rho})$ at hand, the searched for sequence is obtained (similar to Step~3 in the proof of Theorem~\ref{T2}) as a diagonal sequence , i.e., $\deform_h=\tilde \deform_{h,\rho(h)}$ for a suitable function $h\mapsto\rho(h)$ (note that the notion of strong two-scale convergence is metrizable). We leave the details to the reader and turn to the construction of $(\tilde \deform_{h,\rho})$:
  By Proposition~\ref{P:general_construction} we can associate with $\deform_\rho$ a sequence $(\deform_{h,\rho})\subseteq C^{1}(\bar\Omega;\R^3)$ that satisfies \eqref{eq:recov:onescale} with $d=0$ and the clamped coundary conditions \eqref{eq:BC3d}. Furthermore, by Proposition~\ref{P:twoscale:gradrecov} there exists a sequence $(\varphi_h)\subseteq C^\infty_c(\Omega)$ such that $\nabla_h\varphi_h\stto \nabla_\gamma(R_\rho\varphi)$, and
  \begin{equation*}
    h\|\varphi_h\|_{L^\infty(\Omega)}+h\|\nabla_h\varphi_h\|_{L^\infty(\Omega)}\to 0.
  \end{equation*}
  We now define
  \begin{equation*}
    \tilde \deform_{h,\rho}\colonequals \deform_{h,\rho}+h\varphi_h.
  \end{equation*}
  Obviously, $\tilde \deform_{h,\rho}$ satisfies  the clamped coundary conditions \eqref{eq:BC3d} (since $\varphi_h$ is compactly supported in $\Omega$). Statements \eqref{eq:Pid:2a} and \eqref{eq:Pid:2c} are also obvious by construction. It remains to prove  \eqref{eq:Pid:2b}. To that end, we first note that by \eqref{eq:recov:onescale}, for $h>0$ sufficiently small, we have $\det(\nabla_h \deform_{h,\rho}(x))>0$ for all $x\in\Omega$. Thus, by the polar factorization, there exists a unique rotation field $R_{h,\rho}:\Omega\to\SO 3$ such that
  \begin{equation}\label{eq:Pid:1}
    \nabla_h\deform_{h,\rho}=R_{h,\rho}(I_{3\times 3}+h E_h(\deform_{h,\rho})),
  \end{equation}
  and thus
  \begin{eqnarray*}
    \nabla_h\tilde \deform_{h,\rho}=R_{h,\rho}\Big(I_{3\times 3}+hG_h\Big),\qquad G_h\colonequals E_h(\deform_{h,\rho})+R_{h,\rho}^\top\nabla_h\varphi_h.
  \end{eqnarray*}
  By construction, we have $h\|G_h\|_{L^\infty(\Omega)}\to 0$, and thus a Taylor expansion yields
  \begin{equation*}
    \sqrt{(\nabla_h\tilde \deform_{h,\rho})^\top \nabla_h\tilde \deform_{h,\rho}}=\sqrt{I_{3\times 3}+2h\sym G_h+h^2G_h^\top G_h}=I_{3\times 3}+h\sym G_h+ho_h
  \end{equation*}
  with a remainder $o_h:\Omega\to\R^{3\times 3}$ satisfying $\|o_h\|_{L^\infty(\Omega)}\to 0$. Thus, in order to prove \eqref{eq:Pid:2b}, we only need to show that
  \begin{equation*}
    \sym G_h=E_h(\deform_{h,\rho})+\sym\big(R_{h,\rho}^\top\nabla_h\varphi_h\big)\stto \iota(x_3 \II_{\deform_\rho}+M)+\sym\nabla_\gamma\varphi\text{ strongly two-scale in }L^2.
  \end{equation*}
  In view of \eqref{eq:recov:onescale}, it suffices to show that $R_{h,\rho}^T\nabla_h\varphi_h\stto \nabla_\gamma\varphi$. We first note that $R_{h,\rho}\to R_{\rho}$ in $L^2(\Omega)$, since the left-hand side of \eqref{eq:Pid:1} converges in $L^2(\Omega)$  to $R_\rho$. We therefore conclude with help of Lemma~\ref{L:twoscale:prop}~\ref{L:twoscale:prop:c} that $R_{h,\rho}^\top\nabla_h\varphi_h\wtto R^\top_\rho(R_\rho\nabla_\gamma\varphi)=\nabla_\gamma\varphi$ weakly two-scale. On the other hand, the norm converges:
  \begin{equation*}
    \|R^\top_{h,\rho}\nabla_h\varphi_h\|_{L^2(\Omega)}^2=    \|\nabla_h\varphi_h\|_{L^2(\Omega)}^2\to \sum_{j\in J}\int_{S_j}\fint_{\Box_{\Lambda_j}}|\nabla_\gamma\varphi|^2\dd (x_3,y)\dd x',
  \end{equation*}
  since  $|R^\top_{h,\rho}(x')\nabla_h\varphi_h(x)|=|\nabla_h\varphi_h(x)|$ (recall $R^\top_{h,\rho}(x)\in\SO 3$), and because $\nabla_h\varphi_{h,\rho}\stto R_\rho\nabla_\gamma\varphi$. This completes the proof of \eqref{eq:Pid:2b}.
\end{proof}

\subsection{Strong two-scale convergence of the nonlinear strain: Proof of Proposition \ref{strongconvergenceforalmostminimizers}}\label{SS:strongconvergence}

For brevity we only present the proof in the single-grain case, i.e. we assume that $Y_j=Y_1$ for all $j\in J$. The argument extends without larger modifications to the general case.

\step 1 {Reduction of the problem}
We first note that
\begin{eqnarray*}
  \mathcal I^\gamma_{\hom}(\deform_*)+\mathcal I^\gamma_{\rm res}(B)&=&\int_S\fint_{\Box_{\Lambda_1}}Q(x,y,E_*(x,y))\dd (x_3,y)\dd x'\\
                                                              &=&
                                                                     \min_{(M,\varphi)}
                                                                       \int_{S}\fint_{\Box_{\Lambda_1}}Q(x,y,\iota(x_3\II_{\deform_*}+M)+\sym\nabla_\gamma\varphi-B)\dd (x_3,y)\dd x',
\end{eqnarray*}
where the minimum runs over all $M\in L^2(S;\R^{2\times 2}_{\sym})$ and $\varphi\in L^2(S;H^1_{\gamma}(\Box_{\Lambda_1};\R^3))$, see Lemma~\ref{Igamma}.
With help of  Theorem~\ref{T1}~\ref{item:T1:recovery_sequence}  we choose a recovery sequence $(\tilde \deform_h)$ for $\deform_*$ satisfying
\begin{equation*}
\lim\limits_{h\to\infty}\mathcal I^{\e(h),h}(\tilde \deform_h)= \mathcal I^\gamma_{\hom}(\deform_*)+\mathcal I^\gamma_{\rm res}(B).
\end{equation*}
Note that the recovery sequence constructed in  Theorem~\ref{T1}~\ref{item:T1:recovery_sequence} additionally satisfies
\begin{equation}\label{pf:strongconv:1a}
  E_h(\tilde \deform_h)\stto E_*\qquad\text{strongly two-scale in }L^2.
\end{equation}
We claim that the statement of the Proposition~\ref{strongconvergenceforalmostminimizers} follows from the following two claims:
\begin{align}\label{pf:strongconv:1}
  &\limsup\limits_{h\to 0}  \int_\Omega Q\Big(x,\frac{x'}{\e(h)},\theta_hE_h(\deform_h)-E_h(\tilde \deform_h)\Big)\dd x=0,\\
    \label{pf:strongconv:2}
  &\limsup\limits_{h\to 0}  \int_\Omega (1-\theta_h)|E_h(\deform_h)|^2\dd x=0,
\end{align}
where $\theta_h$ is the indicator function defined in \eqref{eq:thetah}. Indeed,  since $E_h(\deform_h)$ and $E_h(\tilde \deform_h)$ are symmetric  and since $Q\in\mathcal Q(\alpha,\beta)$, the combination of both statements, implies that $F_h\colonequals E_h(\deform_h)-E_h(\tilde \deform_h)\to 0$ strongly in $L^2(\Omega)$, and thus also strongly two-scale in $L^2$. In view of \eqref{pf:strongconv:1a}, this implies that $E_h(\deform_h)=E_h(\tilde \deform_h)+F_h\stto E_*$ strongly two-scale in $L^2$, which is the claim of the Proposition~\ref{strongconvergenceforalmostminimizers}.

We shall prove \eqref{pf:strongconv:1} and \eqref{pf:strongconv:2} in Steps~3 and 4, respectively.

\step 2 {Weak two-scale  convergence towards $E_*$}
For brevity set $B_h(x)\colonequals B(x,\tfrac{x'}{\e(h)})$. We claim that
\begin{align}\label{pf:strongconv5c}
  &E_h(\deform_h)\wtto E_*\qquad\text{weakly two-scale in }L^2,\\\label{pf:strongconv5a}
  &\theta_h(E_h(\deform_h)-B_h)\wtto E_*-\sym B\qquad\text{weakly two-scale in }L^2,\\\label{pf:strongconv5}
  &\lim\limits_{h\to 0}\frac1{h^2}\int_{\Omega}\theta_h W_{\e(h)}(x,\nabla_h\deform_h(I_{3\times 3}-h B_{\e(h),h}))\dd x\\\notag
  &\qquad =\,\lim\limits_{h\to 0}\int_{\Omega}Q\Big(x,\tfrac{x'}{\e(h)},\theta_h(E_h(\deform_h)-B_h)\Big)\dd x\\\notag
  &\qquad =\, \mathcal I^\gamma_{\hom}(\deform_*)+\mathcal I^\gamma_{\rm res}(B).
\end{align}
Here comes the argument: By compactness, see Proposition~\ref{P:chartwoscale}~\ref{item:char_twoscale_a}, we have (for a subsequence)
\begin{equation*}
  E_h(\deform_h)\wtto \iota(x_3\II_{\deform_*})+M+\sym\nabla_\gamma\varphi\qquad\text{weakly two-scale in }L^2,
\end{equation*}
for some corrector pair $(M,\varphi)$. Moreover, we may additionally assume that $\int_{S}|\int_{\Box_{\Lambda_1}}\varphi|=0$. Next, we recall from Step~1 in the proof of Theorem~\ref{T1} part \ref{item:T1:lower_bound} (cf.~\eqref{pf:lb:taylor}, \eqref{pf:lb:2} and \eqref{pf:eq2}) that
\begin{align*}
  \liminf\limits_{h\to 0}\mathcal I^{\e(h),h}(\deform_h)\,\geq\,&
  \liminf\limits_{h\to 0}\frac1{h^2}\int_{\Omega}\theta_h W_{\e(h)}(x,\nabla_h\deform_h(I_{3\times 3}-h B_{\e(h),h}))\dd x\\
  \geq\,&\liminf\limits_{h\to 0}\int_{\Omega}Q\Big(x,\tfrac{x'}{\e(h)},\theta_h(E_h(\deform_h)-B_h)\Big)\dd x\\
  \geq\,&
          \int_{S}\fint_{\Box_{\Lambda_1}}Q(x',x_3,y,\iota(x_3\II_{\deform_*}+M)+\sym\nabla_\gamma\varphi-B)\dd (x_3,y)\dd x'\\
  \geq\,&\mathcal I^\gamma_{\hom}(\deform_*)+\mathcal I^\gamma_{\rm res}(B).
\end{align*}
In view of the convexity of $Q$ and assumption~\eqref{eq:strongconvergenceforalmostminimizers:1}, we can upgrade this lower bound to identity \eqref{pf:strongconv5} as claimed.
Furthermore, with Lemma~\ref{L:3.3} we conclude that $(M,\varphi)$ are minimizers of the functional in the second last line, and we conclude that $M=M_*$ and $\varphi=\varphi_*$ (that is, $E=E_*$). Thus, \eqref{pf:strongconv5c} holds not only for the subsequence, but the entire sequence. Moreover, \eqref{pf:strongconv5a} then follows from \eqref{pf:eq3}.

\step 3 {Proof of \eqref{pf:strongconv:2}}
For brevity set $w_h(x)\colonequals\frac1{h^2}W_{\e(h)}(x,\nabla_h\deform_h(x)(I_{3\times 3}-h B_{\e(h),h}(x)))$. By assumption we have $\int_\Omega w_h(x)\dd x=\mathcal I^{\e(h),h}(\deform_h)\to \mathcal I^\gamma_{\hom}(\deform_*)+\mathcal I^{\gamma}_{\rm res}(B)$, and by \eqref{pf:strongconv5} we have $\int_\Omega\theta_hw_h\dd x\to \mathcal I^\gamma_{\hom}(\deform_*)+\mathcal I^{\gamma}_{\rm res}(B)$. Hence, since $w_h$ is non-negative, we conclude that
\begin{equation*}
  \int_{\Omega}(1-\theta_h)|w_h|=  \int_{\Omega}w_h-  \int_{\Omega}\theta_hw_h\to 0.
\end{equation*}
On the other hand, from the elementary, pointwise bound $|E_h(\deform_h)|\leq\frac{1}{h^2}\dist^2(\nabla_h\deform_h,\SO 3)$, \eqref{eq:comp}, and \eqref{item:nonlinear_material_w1}, we conclude that
\begin{equation*}
  \alpha\int_\Omega(1-\theta_h)|E_h(\deform_h)|^2\leq   \int_{\Omega}(1-\theta_h)\big(|w_h|+\frac12h^2|B_{\e(h),h}(x)|^4+|B_{\e(h),h}|^2\big).
\end{equation*}
Thus, it remains to show that $\int_\Omega(1-\theta_h)|B_{\e(h),h}|^2\to 0$. The latter is a consequence of $  \int_\Omega(1-\theta_h)|B_{\e(h),h}|^2\leq \|1-\theta_h\|_{L^1(\Omega)}\|B_{\e(h),h}\|_{L^\infty(\Omega)}^2$, \eqref{eq:thetabound}, and Assumption~\ref{ass:W}~\ref{ass:W:4}.

\step 4 {Proof of \eqref{pf:strongconv:1}}
By adding and subtracting $B_h\colonequals B(x,\tfrac{x'}{\e(h)})$, and by  expanding the square, we have
\begin{eqnarray*}
  &&\int_\Omega Q\Big(x,\tfrac{x'}{\e(h)},\theta_hE_h(\deform_h)-E_h(\tilde \deform_h)\Big)\dd x\\
  &=&
      \int_\Omega Q\Big(x,\tfrac{x'}{\e(h)},\theta_hE_h(\deform_h)-B_h\Big)\dd x-\int_\Omega Q\Big(x,\tfrac{x'}{\e(h)},E_h(\tilde \deform_h)-B_h\Big)\dd x\\
  &&+2 \int_\Omega \mathbb L(x,\tfrac{x'}{\e(h)})(E_h(\tilde \deform_h)-B_h):(E_h(\tilde \deform_h)-\theta_h E_h(\deform_h))\dd x.
\end{eqnarray*}
It suffices to show that the right-hand side converges to $0$ for $h\to 0$. The first integral converges by \eqref{pf:strongconv5}, while the second integral converges to the same limit, since $\tilde \deform_h$ is a recovery sequence. Thus, it remains to show that
the third integral vanishes. Note that $E_h(\tilde \deform_h)-\theta_h E_h(\deform_h)\wtto 0$ weaky in $L^2$, since \eqref{pf:strongconv5a} and since $\tilde \deform_h$ is a recovery sequence. Since $\mathbb L(x,\tfrac{x'}{\e(h)})(E_h(\tilde \deform_h)-B_h)$ is strongly two-scale conergent in $L^2$, we may pass to the limit by appealing to Lemma~\ref{L:twoscale:prop}~\ref{L:twoscale:prop:c}. We thus conclude that
\begin{equation*}
  \lim\limits_{h\to 0}\int_\Omega \mathbb L(x,\tfrac{x'}{\e(h)})(E_h(\tilde \deform_h)-B_h):(E_h(\tilde \deform_h)-\theta_h E_h(\deform_h))\dd x=0.
\end{equation*}
\qed

\subsection{Formulas for the orthotopic case: Proofs of Lemmas~\ref{L:orth1}, \ref{S:ex:ex1},  \ref{S:ex:C3}}\label{SS:Formulas}
Throughout this section we use the shorthand notation $\Box\colonequals(-\frac12,\frac12)\times (-\frac12,\frac12)^2$.

\begin{proof}[Proof of Lemma~\ref{L:orth1}]

  \step 1 {Symmetry properties}
  For $i=1,2,3$ define $P^{(i)}\in\R^{3\times 3}$ and $\pi^{(i)}:\Box\to\Box$ by
  \begin{equation*}
    P^{(i)}\colonequals
    \begin{cases}
      \operatorname{\rm diag}(-1,1,1)&i=1,\\
      \operatorname{\rm diag}(1,-1,1)&i=2,\\
      \operatorname{\rm diag}(1,1,-1)&i=3,
    \end{cases}\qquad
    \pi^{(i)}(x_3,y_1,y_2)\colonequals\begin{cases}
      (x_3,-y_1,y_2)&i=1,\\
      (x_3,y_1,-y_2)&i=2,\\
      (-x_3,y_1,y_2)&i=3,
    \end{cases}
  \end{equation*}
  and remark that thanks to the symmetries and isotropy of $Q$ we have
  \begin{align}\label{Q:symmetry}
    Q(x_3,y,F)=&Q\big(\pi^{(i)}(x_3,y),F\big)=Q\big(x_3,y,P^{(i)}FP^{(i)}\big),
  \end{align}
  for all $F\in\R^{3\times 3}_{\sym}$, a.e.~$(x_3,y)\in\Box$, and $i=1,2,3$.
  As a consequence of these symmetries, we obtain for all $G\in\R^{2\times 2}_{\sym}$, $M\in\R^{2\times 2}_{\sym}$, and $\varphi\in H^1_{\gamma}(\Box;\R^3)$ the identity
  \begin{equation}\label{fundamentalsymmetry}
    \begin{aligned}
    &\int_{\Box} Q\big(x_3,y,\iota(x_3 G)+\iota(M)+\sym\nabla_\gamma\varphi\big)\\
    &\qquad=
    \left\{\begin{aligned}
        &\int_{\Box} Q\big(x_3,y,P^{(i)}\iota(x_3G)P^{(i)}+\iota(\widetilde M^{(i)})+\sym\nabla_\gamma\widetilde \varphi^{(i)}\big)&&\text{for }i=1,2,\\
      &\int_{\Box} Q\big(x_3,y,\iota(-x_3G)+\iota(M)+\sym\nabla_\gamma\widetilde \varphi^{(3)}\big)&&\text{for }i=3,
    \end{aligned}\right.
  \end{aligned}
  \end{equation}
  where
  \begin{equation}\label{trafoMvarphi}
    \widetilde\varphi^{(i)}\colonequals P^{(i)}\varphi\circ\pi^{(i)},\qquad \widetilde M^{(i)}
    \colonequals
    \begin{cases}
      \operatorname{\rm diag}(-1,1)M\operatorname{\rm diag}(-1,1)&\text{for }i=1,\\
      \operatorname{\rm diag}(1,-1)M\operatorname{\rm diag}(1,-1)&\text{for }i=2.
    \end{cases}
  \end{equation}
  Indeed, this follows by a direct calculation using \eqref{Q:symmetry} and  the identities
  \begin{equation*}
    \nabla_\gamma\tilde\varphi^{(i)}=P^{(i)}\big(\nabla_\gamma\varphi\circ\pi^{(i)}\big)P^{(i)},\qquad P^{(3)}\iota(x_3G+M)P^{(3)}=\iota(x_3G+M).
  \end{equation*}

\step 2 {Symmetries of the corrector}
Let $G\in\R^{2\times 2}_{\sym}$ and let $(M_G,\varphi_G)$ denote the associated corrector in the sense of Lemma~\ref{L:corrector}.
We claim that
\begin{align}\label{symetryforgeneralG}
  M_G=0\quad\text{and}\quad\varphi_G=
  \begin{cases}
    P^{(i)}\varphi_G\circ\pi^{(i)}&\text{if }i=1,2\text{ and $G$ is diagonal},\\
    -P^{(i)}\varphi_G\circ\pi^{(i)}&\text{if }i=1,2\text{ and $G$ vanishes on the diagonal},\\
    -P^{(3)}\varphi_G\circ\pi^{(3)}&\text{if }i=3.
  \end{cases}
\end{align}
The proof is as follows: Identity \eqref{fundamentalsymmetry} with $i=3$ in combination with the characterization of the corrector as a minimizer (see Remark~\ref{R:corrector}) implies that $(-M_G, -P^{(3)}\varphi_G\circ\pi^{(3)})$ is also a corrector associated with $G$. By uniqueness of the corrector, we especially conclude that $M_G=-M_G$. Hence, $M_G=0$ and \eqref{L:orth1:a} follow.
Likewise, we conclude that $\varphi_G=-P^{(3)}\varphi_G\circ\pi^{(3)}$. By the same argument and by using that for $i=1,2$ and all $G\in\R^{2\times 2}_{\sym}$ we have
\begin{equation*}
  P^{(i)}\iota(x_3G)P^{(i)}=
  \begin{cases}
    \iota(x_3G)&\text{if $G$ is diagonal},\\
    \iota(-x_3G)&\text{if $G$ vanishes on the diagonal},
  \end{cases}
\end{equation*}
the remaining identities follow. Since $Q$ is independent of $y_2$ we conclude that $\partial_{y_2}\varphi_G=0$  and $\varphi_G \circ \pi_{2} = \varphi_G$. Combined with \eqref{symetryforgeneralG}, we deduce that
  \begin{equation}\label{prop:1}
    \begin{aligned}
      &\varphi_G=(\varphi_{G,1},0,\varphi_{G,3})&&\text{if $G$ is diagonal},\\
      &\varphi_G = (0,\varphi_{G,2},0)&&\text{if $G$ vanishes on the diagonal}.
    \end{aligned}
  \end{equation}
  In particular, \eqref{L:orth1:c} follows.

  \step 3 {Conclusion}
  To prove that $Q_{\hom}^\gamma$ is orthotopic, by Proposition~\ref{P:1} and since $M_{G_3}=0$,  it suffices to show that for $i=1,2$ we have
  \begin{equation}\label{eq:euler22}
    \int_\Box\mathbb L(x_3,y_1)(\iota(x_3G_3)+\sym\nabla_\gamma\varphi_{G_3}):\iota(x_3G_i)=0.
  \end{equation}
  For the argument, note that in view of \eqref{formulaforQ} we have
  \begin{equation*}
    \mathbb L(x_3,y_1)G=\lambda(x_3,y_1)({\rm tr}G)I_{3\times 3}+2\mu(x_3,y_1)\sym G.
  \end{equation*}
  Furthermore, from \eqref{prop:1} and the fact that $\varphi_{G_3}$ is independent of $y_2$, we conclude that all diagonal entries of $\iota(x_3G_3)+\sym\nabla_\gamma\varphi_{G_3}$ vanish a.e.~in $\Box$. We conclude \eqref{eq:euler22} (and thus orthotropicity) and the identity \eqref{eq:trace}.
\end{proof}

\begin{proof}[Proof of Lemma~\ref{S:ex:ex1}]

Throughout the proof, we will use the shorthand notation $\varphi_{i,3} = \partial_{3} \varphi\cdot e_i$ and $\varphi_{i,j} = \partial_{y_{j}} \varphi\cdot e_i$ for $i=1,2,3$ and $j =1,2$. Since the setting under consideration is a special case of Lemma~\ref{L:orth1}, we conclude that $Q^\gamma_{\hom}$ is orthotropic. Hence, it suffices to prove the formulas for the coefficients $q_1,q_2,q_3,q_{12}$ and the properties of the map $\gamma\mapsto\mu_\gamma$.

\step 1 {Correctors for $i=1,2$ and coefficients $q_1,q_2,q_{12}$}
Consider the following subspace of $H^1_{\gamma}(\Box;\R^3)$,
\begin{align*}
  \mathcal H_{\rm diag} \colonequals &\,\Big\{\varphi\in H^1_\gamma(\Box;\R^3)\,:\,\varphi_{1,2}=\varphi_{3,2}=\varphi_2=0\text{ a.e.~in~$\Box$, and }\int_{\Box}\varphi_1=\int_{\Box}\varphi_3=0\,\Big\},
\end{align*}
and note that by \eqref{L:orth1:c} we have $\varphi_{G_1},\varphi_{G_2}\in\mathcal H_{\rm diag}$. Since $M_{G_i}=0$ by Lemma~\ref{L:orth1}, and in view of the characterization of the corrector as a minimizer, we obtain for $i=1,2$ the identity
\begin{align*}
  q_1=Q^\gamma_{\hom}(G_i)=\,&\fint_\Box Q\big(x_3,y_1,\iota(x_3 G_i)+\sym\nabla_\gamma\varphi_{G_i}\big)\\
  =\,&\min_{\varphi\in \mathcal H_{\rm diag}}\fint_{\Box}\mu(\tfrac{1}{\gamma}\varphi_{1,3}+\varphi_{3,1})^2+2\mu(\tfrac{1}{\gamma}\varphi_{3,3})^2+2\mu(\delta_{i1}x_3+\varphi_{1,1})^2+2\mu\delta_{i2}x_3^2,
\end{align*}
where $\delta_{ij}=1$ for $i=j$ and $\delta_{ij}=0$ for $i\neq j$. Based on this identity, one can check that for $i=1$ we have $\varphi_{G_1}=(x_3w,0,W)$ with
\begin{align*}
  w(y_1)\colonequals\int_{-\frac12}^{y_1}\big(\frac{\meanh\mu}{\mu}-1\big)\dd s,\qquad W(y_1)\colonequals-\frac1\gamma\big(\int_{-\frac12}^{y_1}w(s)\dd s-\int_{-\frac12}^{\frac12}\int_{-\frac12}^{t}w(s)\dd s\dd t\big).
\end{align*}
We compute
\begin{equation}\label{SN:corr1}
  \nabla_\gamma\varphi_{G_1}=\big(\frac{\meanh\mu}{\mu}-1\big)\iota(x_3 G_1)
\end{equation}
and conclude that  $q_1=Q^\gamma_{\hom}(G_1)= \frac{1}{6} \meanh{\mu}$ as claimed. Similarly, for $i=2,$ we get  $\varphi_{G_2}=0$ and thus $q_2=Q^\gamma_{\hom}(G_2)= \frac{1}{6} \overline{\mu}$. Combined with \eqref{SN:corr1} and Proposition~\ref{P:1}, we further conclude that
\begin{align*}
  q_{12}=\,&\fint_\Box\mathbb L(\iota(x_3G_1+M_{G_1})+\sym\nabla_\gamma\varphi_{G_1}):\iota(x_3 G_2)\\
  =\,&
  \fint_{\Box}2\mu\Big(\iota(x_3G_1)+\big(\frac{\meanh\mu}{\mu}-1\big)\iota(x_3 G_1)\Big):\iota(x_3G_2)=0,
\end{align*}
since $\iota(x_3G_1):\iota(x_3G_2)=0$.

\step 2 {Corrector for $i=3$, coefficient $q_3$, and properties of $\mu_\gamma$}
For the argument it is convenient to introduce for $\gamma\in(0,\infty)$ the functional
\begin{equation*}
  \mathcal E_\gamma:\mathcal H\to\R,\qquad \mathcal E_\gamma(w)\colonequals\fint_{\Box}\mu\Big((\sqrt{12} x_3+\partial_{y_{1}}w)^2+(\frac{1}{\gamma}\partial_3w)^2\Big).
\end{equation*}
We note that  by construction, we have $\mu_\gamma=\min\mathcal E_\gamma$ where $\mu_\gamma$ is defined by \eqref{muGamma2}.
Moreover, by \eqref{L:orth1:c} in Lemma~\ref{L:orth1} we have $\varphi_{G_3}=(0,w_*,0)$ for some $w_*\in\mathcal H$. As in the previous step we thus conclude that
\begin{equation}\label{problemforQ33}
  q_3=Q^\gamma_{\hom}(G_3)=\min_{w\in \mathcal H} \fint_{\Box}\mu\Big((\sqrt 2 x_3+\partial_{y_{1}}w)^2+(\frac{1}{\gamma}\partial_3w)^2\Big)=\frac16\min_{w\in\mathcal H}\mathcal E_\gamma(w),
\end{equation}
and thus the identity for $q_3$ in \eqref{eq:q_hom_laminate_coefficients} follows.

\substep{2.1 (Proof of  \eqref{S:ex:muGammaProp1})} Taking $w=0$ as a test function in \eqref{problemforQ33}, we obtain $\frac16\mu_\gamma\leq  \fint_{\Box}2\mu x_3^2= \frac{1}{6} \overline{\mu}$. On the other hand,
\begin{equation*}
  \begin{aligned}
    \meanh\mu &= \min_{w\in\mathcal H}\fint_{\Box}\mu\big(\sqrt{12} x_3+\partial_{y_{1}}w\big)^2 \leq  \min_{w \in\mathcal H}\fint_{\Box}\mu\Big((\sqrt{12} x_3+\partial_{y_{1}}w)^2+(\frac{1}{\gamma}\partial_3w)^2\Big)=\mu_\gamma,
\end{aligned}
\end{equation*}
and thus \eqref{S:ex:muGammaProp1} follows.

\substep{2.2 (Continuity of $\mu_\gamma$)}
Let $(\gamma_n)$ be a sequence that converges to $\gamma\in(0,\infty)$.
It is straightforward to check that $\mathcal E_{\gamma_n}$ $\Gamma$-converges to $\mathcal E_{\gamma}$ w.r.t.~weak convergence in $\mathcal H$, where we consider $\mathcal H$ as a closed subspace of $H^1_\gamma(\Box;\R)$.
Hence, standard arguments from the theory of $\Gamma$-convergence show that $\min\mathcal E_{\gamma_n}\to\min\mathcal E_{\gamma}$ and we conclude that $\gamma\mapsto\mu_\gamma$ is continuous.

\substep{2.3 (Asymptotic behavior of $\mu_\gamma$)}
We show that $\lim_{\gamma\to \infty}\mu_\gamma=\meanh\mu$.
Since the map $\gamma\mapsto\mu_\gamma$ is monotone and bounded by \eqref{S:ex:muGammaProp1}, the limit $\lim_{\gamma\to \infty}\mu_\gamma$ exists.
To identify the limit, consider $w_\infty(x_3,y)\colonequals\sqrt{12}x_3\int_{-\frac12}^{y_1}(\frac{\meanh\mu}{\mu(s)}-1)\dd s$ and note that we have $w_\infty\in\mathcal H$ and $\meanh\mu=\lim_{\gamma\to\infty}\mathcal E_\gamma(w_\infty)$.
Combined with the lower bound $\mathcal E_\gamma(w_\infty)\geq\mu_\gamma\geq\meanh\mu$, we get $\meanh\mu\geq\lim\limits_{\gamma\to\infty}\mu_\gamma\geq \meanh\mu$,
and thus the claim follows.

Next, we prove $\lim\limits_{\gamma\to 0}\mu_\gamma=\overline\mu$.
As above, by monotonicity, the limit exists. Let $w_\gamma\in\mathcal H$ be a minimizer of $\mathcal E_\gamma$.
By \eqref{S:ex:muGammaProp1} we have $\meanh\mu\leq \mathcal E_\gamma(w_\gamma)\leq\overline\mu$, and thus $\limsup_{\gamma\to 0}\fint_\Box|\partial_{y_1}w_\gamma|^2+|\frac1\gamma\partial_3 w_\gamma|^2<\infty$.
Hence, $(w_\gamma)$ is bounded in $\mathcal H$ and we can pass to a subsequence (not relabeled) such that $w_\gamma\wto w_0$ weakly in $\mathcal H$, where $w_0$ satisfies $\partial_3w_0=0$.
By weak lower semicontinuity of convex functionals and since $\fint_\Box\mu\big(\sqrt{12}x_3+\partial_{y_1}w_\gamma\big)^2\leq\mathcal E_\gamma(w_\gamma)=\mu_\gamma\leq\overline\mu$, we conclude
\begin{equation*}
  \fint_\Box\mu\big(\sqrt{12}x_3+\partial_{y_1}w_0\big)^2\leq \liminf\limits_{\gamma\to 0}  \fint_\Box\mu\big(\sqrt{12}x_3+\partial_{y_1}w_\gamma\big)^2\leq \lim\limits_{\gamma\to 0}\mathcal E_\gamma(w_\gamma)\leq\overline\mu.
\end{equation*}
On the other hand, since $\fint_\Box\mu\sqrt{12}x_3\partial_{y_1}w_0=0$, we have
\begin{equation}\label{eq:SN:orthogonal}
  \fint_\Box\mu\big(\sqrt{12}x_3+\partial_{y_1}w_0\big)^2=  \fint_\Box\mu(\sqrt{12}x_3)^2+\fint_\Box\mu(\partial_{y_1}w_0)^2\geq \overline\mu.
\end{equation}
We conclude that $\lim_{\gamma\to 0}\mu_\gamma=\overline\mu$.

\substep{2.4 (Strict monotonicity of $\mu_\gamma$)} It suffices to prove that if $\mu_\gamma$ is not strictly monotone, then $\mu$ is constant. For the argument, suppose that $\mu_\gamma$ is not strictly monotone.
Since $\mu_\gamma$ is monotone, there exist $\gamma_1<\gamma_2$ such that $(\gamma_1,\gamma_2)\ni\gamma\mapsto\mu_\gamma$ is constant.
Hence $\frac{d}{d\gamma}\mathcal E_\gamma(w_\gamma)=0$ for all $\gamma\in(\gamma_1,\gamma_2)$, where $w_\gamma$ denotes the minimizer of $\mathcal E_\gamma$.
By combining the latter with the Euler-Lagrange equation for $w_\gamma$, i.e.,
\begin{equation}\label{eq:EL}
  \fint_\Box\mu\Big((\sqrt12 x_3+\partial_{y_1}w_\gamma)\partial_{y_1}\eta+\frac1{\gamma^2}\partial_3w_\gamma\partial_3\eta\Big)=0\qquad\text{for all }\eta\in\mathcal H,
\end{equation}
we find that $-2\gamma^{-3}\fint_\Box\mu(\partial_3w_\gamma)^2=0$. Hence, $w_\gamma$ is independent of $x_3$, and thus $\mathcal E_\gamma(w_\gamma)=\overline\mu+\fint_\Box\mu(\partial_{y_1}w_\gamma)^2$. Since $w_\gamma$ minimizes $\mathcal E_\gamma$, we conclude that $w_\gamma=0$ and thus \eqref{eq:EL} reduces to $\fint_\Box\mu(\sqrt12 x_3)\partial_{y_1}\eta=0$, which especially holds for test-functions of the form $\eta(x_3,y_1)=x_3\tilde\eta(y_1)$. We thus conclude that $\mu=\fint_\omega\mu$ a.e.~in $\Box$ and the claim follows.

\step 3 {Evaluation of effective prestrain}
From Proposition~\ref{P:1} \ref{P:1:c} we recall the definition of $\widehat B\in\R^3$, which we combine with \eqref{SN:corr1}, the identity $\varphi_{G_2}=0$, and \eqref{eq:trace}. We get
\begin{equation*}
  \widehat B_{1}=2\meanh\mu\fint_{\Box}\rho x_3\dd (x_3,y),\qquad
  \widehat B_{2}=2\fint_{\Box}\mu\rho x_3\dd (x_3,y),\qquad \widehat B_{3}^\gamma=0.
\end{equation*}
Furthermore, we recall from Proposition~\ref{P:1} that $\widehat B^\gamma_{\eff}=\widehat Q^{-1}\widehat B$ and note that we have $\widehat Q=\frac16\operatorname{\rm diag}(\meanh\mu,\overline\mu,\mu_\gamma)$ by Step~1 and Step~2. Now, \eqref{eq:Beff} follows from  a short calculation.
\end{proof}

\begin{proof}[Proof of Lemma~\ref{S:ex:C3}]
This follows from Lemma~\ref{S:ex:ex1} and direct computations.
\end{proof}

\begin{proof}[Proof of Lemma~\ref{L:char}]
  The proof is based on case discrimination depending on the sign of $\det H$.

  \step 1 {Case~\ref{L:char:a} and its logical complement}
  By the invariance property \eqref{eq:S:P1} and identity \eqref{eq:st:1000b} we have
  \begin{equation}\label{L:char:st:P1}
    \text{Case~\ref{L:char:a}}\qquad\Leftrightarrow\qquad \mathop{\argmin}_{a\in \mathcal{G}^+_{\R^2}}\mathcal E_{Q,B}(a)\subset\partial \mathcal{G}^+_{\R^2}.
  \end{equation}
  Suppose \ref{L:char:a}. Then the inclusion of $\mathcal S_{Q,B}$ claimed in \ref{L:char:a} follows from \eqref{eq:st:1004}. On the other hand, if \ref{L:char:a} does not hold, then by \eqref{L:char:st:P1}, $\mathcal E_{Q,B}$ admits a local minimizer in the interior of $\mathcal{G}^+_{\R^2}$, and thus there exists $a\in \mathcal{G}^+_{\R^2}\setminus\partial \mathcal{G}^+_{\R^2}$ such that $\nabla \mathcal E_{Q,B}(a)=0$ and $\nabla^2 \mathcal E_{Q,B}(a)$ is positive semi-definite. Since $\nabla^2 \mathcal E_{Q,B}=H$ and  $\mathop{\operatorname{trace}}(H)=2(q_1+q_2)>0$ (by positive definiteness of $Q$), we have $\det H\geq 0$. Thus, to conclude the trichotomy, it suffices to show
  \begin{align}\label{L:char:Pb}
    \big(\,\neg\ref{L:char:a}\text{ and }\det H>0\,\big)\qquad\Leftrightarrow\qquad \ref{L:char:b},\\
    \label{L:char:Pc}
    \big(\,\neg\ref{L:char:a}\text{ and }\det H=0\,\big)\qquad\Leftrightarrow\qquad \ref{L:char:c}.
  \end{align}

  \step 2 {Case~\ref{L:char:b} and proof of \eqref{L:char:Pb}}
  Suppose $\det H>0$. Then $\mathcal E_{Q,B}$ is strictly convex and $g_*$ is the unique minimizer of $\mathcal E_{Q,B}$ on $\R^2$ as can be seen by checking the identity $\nabla \mathcal E_{Q,B}(g_*)=0$. In view of \eqref{L:char:st:P1}, the equivalence \ref{L:char:b} follows. Furthermore, if \ref{L:char:b} holds, then $\mathop{\argmin}_{\mathcal{G}^+_{\R^2}}\mathcal E_{Q,B}=\{g_*\}$ and the characterization of $\mathcal S_{Q,B}$ claimed in \ref{L:char:b} follows from \eqref{eq:st:1000a}.

  \step 3 {Case~\ref{L:char:c} and proof of \eqref{L:char:Pc}}
  Let $\det H=0$. Since $\mathop{\operatorname{trace}}(H)=2(q_1+q_2)>0$,  $H$ is positive semi-definite and thus $\mathcal E_{Q,B}$ is convex.
  Therefore, any critical point of $\mathcal E_{Q,B}$ is a minimizer of $\mathcal E_{Q,B}$ (in $\R^2$). Furthermore, we have $H=2(q_1+q_2)q_*\otimes q_*$.

  Now assume that \ref{L:char:a} is not valid. Then by the argument of Step~1, $\mathcal E_{Q,B}$ admits a critical point $a\in \mathcal{G}^+_{\R^2}\setminus\partial \mathcal{G}^+_{\R^2}$, and thus  $0=\nabla \mathcal E_{Q,B}(a)=Ha-2Ab$.
    We conclude that $2Ab\in\mathop{\operatorname{range}} H$, and thus \ref{L:char:c} holds.

    On the other hand, if \ref{L:char:c} holds, then $2Ab=(2Ab\cdot q_*)q_*$, and thus $H(s_*q_*)=2Ab$. Hence, $g_*=s_*q_*$ is a critical point. In fact, every point on the line $\mathcal L\colonequals\{a\in\R^2\,:\,a\cdot q_*=s_*\}$ is a critical point and thus a minimizer of $\mathcal E_{Q,B}$. In order to conclude $\neg$\ref{L:char:a} (and thus the validity of \eqref{L:char:Pc}), we only need to show that $\mathcal L\cap (\mathcal{G}^+_{\R^2}\setminus\partial \mathcal{G}^+_{\R^2})\neq\emptyset$. The latter can be seen as follows: We parametrize $\mathcal L$ by  $\ell:\R\to\R^2$, $\ell(t)\colonequals s_* q_*+t(-q_{*,2},q_{*,1})$ and consider the function $\varphi(t) \colonequals (s_*q_{*,1}-tq_{*,2})(s_*q_{*,2}+tq_{*,1})$. By construction we have $\ell(t)\in \mathcal{G}^+_{\R^2}\setminus\partial \mathcal{G}^+_{\R^2}$ if and only if $\varphi(t)>0$. It is easy to see that $\varphi$ is concave. We claim that
    \begin{equation}\label{eq:st:1000g}
      \varphi(0)=s_*^2q_{*,1}q_{*,2}>0.
    \end{equation}
    Indeed, combining \eqref{eq:st:1001} with the identity $0=\det H=4q_1q_2-(q_{12}+2q_3)^2$ yields $q_{12}+q_3>0$ and thus $2q_2(q_{12}+2q_3)>0$.
    In view of the definition of $q_*$ and the property $s_*^2>0$ (which follows from $b\neq 0$), we obtain \eqref{eq:st:1000g}.
    Since $\varphi$ is concave and continuous, we conclude from \eqref{eq:st:1000g} that $\{t\in\R\,:\,\varphi(t)\geq 0\}$ is an interval with non-empty interior,
    and thus $\mathcal L\cap \mathcal{G}^+_{\R^2}\setminus\partial \mathcal{G}^+_{\R^2}\neq \emptyset$ follows.
    This completes the argument for \eqref{L:char:Pb}.
    Furthermore, the above argument proves the characterization of $\mathcal S_{Q,B}$ stated in \ref{L:char:c}.
\end{proof}

  \subsection{Shape programming: Proofs of Lemma~\ref{L:transformation}, Corollary~\ref{C:transf}, Lemma~\ref{L:gamma-spatial}, and  Theorem~\ref{T:shape}}\label{SS:shapedesign}

  \begin{proof}[Proof of Lemma~\ref{L:transformation}]
    For convenience set $\Box\colonequals\Box_\Lambda$ and $\widetilde\Box\colonequals\Box_{\tilde\Lambda}$.

    \step 1 {Reduction using linearity}
     We first note that the map
    \begin{equation*}
      \Phi\,:\,\mathbf H_{\widetilde\Lambda}\to\mathbf H_{\Lambda},\qquad \Phi(\tilde F)\colonequals F,\qquad\text{where }F(x_3,y)\colonequals\widehat T^{\top}\tilde F(x_3,Ty)\widehat T,
    \end{equation*}
    is an isometric isomorphism. In particular, for all $\tilde F,\tilde G\in\mathbf H_{\tilde\Lambda}$ we have
    \begin{align*}
      \fint_{\tilde\Box}\widetilde Q(x_3,y,\tilde F)=      \fint_{\Box}Q(x_3,y,\Phi(\tilde F)),\qquad      \fint_{\tilde\Box}\widetilde{\mathbb L}\tilde F:\tilde G=      \fint_{\Box}\mathbb L\Phi(\tilde F):\Phi(\tilde G).
    \end{align*}
    Furthermore, the map
    \begin{equation*}
      H^1_\gamma(\Box_\Lambda;\R^3)\to H^1_\gamma(\Box_{\tilde\Lambda};\R^3),\qquad \varphi\mapsto\tilde\varphi,\qquad \tilde\varphi(x_3,y)=\widehat T^{-\top}\varphi(x_3,T^{-1} y),
    \end{equation*}
    defines an isomorphism satisfying $\Phi\big(\sym\nabla_\gamma\tilde\varphi\big)=\sym\nabla_\gamma\varphi$. With help of this identity it is easy to check that $\Phi(\mathbf H_{\rm rel, \widetilde\Lambda}^\gamma)=\mathbf H_{\rm rel,\Lambda}^\gamma$ and $\Phi(\mathbf H_{\widetilde\Lambda}^\gamma)=\mathbf H_{\Lambda}^\gamma$.

  \step 2 {Reduction using projections}
  We claim that
    \begin{equation}
      \label{eq:pf:transf1}
      \Phi(P^\gamma_{\widetilde\Lambda}\tilde B)=P^\gamma_{\Lambda}B.
    \end{equation}
    Indeed,  by the definition of $\tilde B$ and the properties of $\Phi$, we have for all $F\in\mathbf H^\gamma_{\Lambda}$,
    \begin{eqnarray*}
      \fint_{\Box}\mathbb L\big(\Phi(P^\gamma_{\tilde\Lambda}(\sym\tilde B)-\sym B\big):F=\fint_{\tilde\Box}\widetilde{\mathbb L}\big(P^\gamma_{\tilde\Lambda}(\sym\tilde B)-\sym \tilde B\big):\Phi^{-1}(F)=0,
    \end{eqnarray*}
    where the last identity holds, since $\Phi^{-1}(F)\in\mathbf H^\gamma_{\tilde\Lambda}$ and because $P^\gamma_{\tilde\Lambda}$ is the orthogonal projection onto $\mathbf H^\gamma_{\tilde\Lambda}$.

    \step 3 {Conclusion}
     Let $\tilde G\in \R^{2\times 2}_{\sym}$. Then by Step~1,
    \begin{eqnarray*}
      \widetilde Q^\gamma_{\hom}(\tilde G)&=&\inf_{\tilde\chi\in\mathbf H^\gamma_{\rm rel,\tilde\Lambda}}\fint_{\tilde\Box}\widetilde Q(x_3,y,\iota(x_3 \tilde G)+\tilde\chi)\,=\,\inf_{\tilde\chi\in\mathbf H^\gamma_{\rm rel,\tilde\Lambda}}\fint_{\Box}Q\Big(x_3,y,\Phi\big(\iota(x_3 \tilde G)+\tilde\chi\big)\Big)\\
                                          &=&\inf_{\chi\in\mathbf H^\gamma_{\rm rel,\Lambda}}\fint_{\Box}Q\Big(x_3,y,\iota(x_3 T^\top \tilde G T)+\chi\Big)\,=\,Q^\gamma_{\hom}(T^\top\tilde GT),
    \end{eqnarray*}
    as claimed. Similarly,
    \begin{eqnarray*}
      &&\widetilde Q^\gamma_{\hom}(\tilde G-\tilde B^\gamma_{\rm eff})+\fint_{\tilde\Box}\widetilde Q(x_3,y,(I-P^\gamma_{\tilde\Lambda})(\sym\tilde B))\\
      &=&\inf_{\tilde\chi\in\mathbf H^\gamma_{\rm rel,\tilde\Lambda}}\fint_{\tilde\Box}\widetilde Q(x_3,y,\iota(x_3 \tilde G)+\tilde\chi-\tilde B)
      =\inf_{\chi\in\mathbf H^\gamma_{\rm rel,\Lambda}}\fint_{\Box}Q\Big(x_3,y,\iota(x_3 T^\top \tilde G T)+\chi-B\Big)\\&=&Q^\gamma_{\hom}(T^\top\tilde GT-B^\gamma_{\rm eff})+\fint_{\Box}Q(x_3,y,(I-P^\gamma_{\Lambda})(\sym B)).
    \end{eqnarray*}
    From \eqref{eq:pf:transf1} and the properties of $\Phi$, we deduce that on both sides the second integrals are equal, and thus
    \begin{equation*}
      \widetilde Q^\gamma_{\hom}(\tilde G-\tilde B^\gamma_{\rm eff})=Q^\gamma_{\hom}(T^\top\tilde GT-B^\gamma_{\rm eff}).
    \end{equation*}
    Since this is true for arbitrary $\tilde G\in\R^{2\times 2}_{\sym}$, we conclude that $\tilde B^\gamma_{\rm eff}=T^{-\top}B^\gamma_{\rm eff}T^{-1}$ as claimed.
\end{proof}

\begin{proof}[Proof of Corollary~\ref{C:transf}]
  Note that for all $G\in\R^{3\times 3}$ we have by frame indifference
  \begin{equation*}
    W(R,\kappa;x_3,y,I_{3\times 3}+G)=\bar W(\kappa;x_3,R^\top y,I_{3\times 3}+R^\top GR).
  \end{equation*}
  Thus, the quadratic forms associated with $W(R,\kappa;\cdot)$ and $\bar W(\kappa;\cdot)$ transform as in Lemma~\ref{L:transformation}, that is,
  \begin{equation*}
    Q(R,\kappa;x_3,y,G)=\bar Q(\kappa;x_3,R^\top y,R^\top GR).
  \end{equation*}
  With help of Lemma~\ref{L:transformation} it is straightforward to check the claims of the corollary.
\end{proof}

  \begin{proof}[Proof of Lemma~\ref{L:gamma-spatial}]
  \step 1 {Proof of ~\ref{L:gamma-spatial:a}}
   W.l.o.g.~we may assume $\int_S\deform_n=0$. From \eqref{eq:Q2d} we get
  \begin{equation*}
    \frac{\alpha}{2}\int_S|\II_{\deform_n}|^2\leq\mathcal I_n(\deform_n)+\alpha\|B_n\|^2_{L^2(S)},
  \end{equation*}
  and thus we have $\|\II_{\deform_n}\|^2_{L^2(S)}<\infty$. Since $|\II_{\deform_n}|^2=|\nabla'\nabla'\deform_n|^2$, the Poincar\'e--Wirtinger inequality implies that $\limsup_{n\to\infty}\|\deform_n\|_{H^2(S;\R^3)}<\infty$. Therefore, we may pass to a subsequence that weakly converges in $H^2(S;\R^3)$ to some limit $\deform$. Since $H^2_{\iso}(S;\R^3)$ is weakly closed, we conclude that $\deform\in H^2_{\iso}(S;\R^3)$.

  \step 2 {Proof of ~\ref{L:gamma-spatial:b}}
   We first note that the construction of a recovery sequence $(\deform_n)$ for $\deform\in H^2_{\iso}(S;\R^3)$ is trivial, since we just may consider the constant sequence $\deform_n=\deform$. To prove the lower bound part of $\Gamma$-convergence, consider a sequence $(\deform_n)\subset H^2_{\iso}(S;\R^3)$ that weakly converges in $H^2$ to some $\deform_\infty\in H^2(S;\R^3)$. We need to show that
  \begin{equation*}
    \liminf\limits_{n\to\infty}\mathcal I_n(\deform_n)\geq \mathcal I_\infty(\deform_\infty).
  \end{equation*}
  For the argument set $F_n\colonequals\II_{\deform_n}-B_n$ (for $n\in\N\cup\{\infty\}$) and note that we have
  \begin{equation}\label{eq:texture:1}
    F_n\wto F_\infty\qquad\text{weakly in }L^2.
  \end{equation}
  By expanding the square and by positivity of $Q_n$, we have
  \begin{eqnarray*}
    \int_S  Q_n(x',F_n)
    &=&
        \int_S Q_n(x',F_\infty)+Q_n(x',F_n-F_\infty)+\mathbb L_n(x') F_\infty:(F_n-F_\infty)\\
    &\geq&
           \int_S    Q_n(x',F_\infty)+\mathbb L_n(x') F_\infty :(F_n-F_\infty).
  \end{eqnarray*}
  Thanks to the convergence of $Q_n$ we have $\int_S    Q_n(x',F_\infty)\to\int_SQ_\infty(x',F_\infty)=\mathcal I_\infty(\deform_\infty)$. Moreover, the assumption on $Q_{n}$ implies that $\mathbb L_n(x') F_\infty\to \mathbb L_\infty(x') F_\infty$ strongly in $L^2(S)$. Together with \eqref{eq:texture:1} this yields $\int_S\mathbb L_n(x') F_\infty:(F_n-F_\infty)\to 0$.

  \step 3 {Proof of ~\ref{L:gamma-spatial:c}}
  Standard arguments from the theory of $\Gamma$-convergence imply that we have  $\deform_n\wto \deform_\infty$ weakly in $H^2$ (up to extraction of a subsequence) where $\deform_\infty\in H^2_{\iso}(S;\R^3)$ is a minimizer of $\mathcal I^\gamma_{\hom,\infty}$. In the following we explain how to upgrade this to $\deform_n\to \deform_\infty$ strongly in $H^2$. Note that it suffices to show
  \begin{equation}\label{L:gamma-spatial:1}
    \|\II_{\deform_n}-\II_{\deform_\infty}\|_{L^2(S)}\to 0.
  \end{equation}
  Indeed, thanks to the identity $|\II_{\deform_n}|^2=|\nabla'\nabla'\deform_n|^2$, \eqref{L:gamma-spatial:1} implies $\|\nabla'\nabla'\deform_n\|_{L^2(S)}\to \|\nabla'\nabla' \deform\|_{L^2(S)}$, and thus $\nabla'\nabla'\deform_n\to \nabla'\nabla'\deform$ strongly in $L^2$. To see \eqref{L:gamma-spatial:1} we argue that
  \begin{equation*}
    \int_SQ_n(x',\II_{\deform_n}-\II_{\deform_\infty})\to 0.
  \end{equation*}
  The latter can be seen as follows: By adding and substracting $B_n$, and by expanding the square, we get
  \begin{eqnarray*}
    &&\int_SQ_n(x',\II_{\deform_n}-\II_{\deform_\infty})\\
    &=&\int_SQ_n(x',\II_{\deform_n}-B_n)-\int_SQ_n(x',\II_{\deform_\infty}-B_n)+2\int_S\mathbb L_n(x')(\II_{\deform_\infty}-B_n):(\II_{\deform_\infty}-\II_{\deform_n}).
  \end{eqnarray*}
  Since $(\deform_n)$ is a sequence of minimizers, and since $\deform_\infty$ is a minimizer of the $\Gamma$-limit, we have $\int_SQ_n(x',\II_{\deform_n}-B_n)\to \int_SQ_\infty(x',\II_{\deform_\infty}-B_\infty)$. Similarly, the second term $\int_SQ_n(\II_{\deform_\infty}-B_n)$ converges to the same limit, since $Q_n$ and $B_n$ converge by assumption. It remains to show that the last term converges to $0$. This is the case, since $\mathbb L_n(\II_{\deform_\infty}-B_n)$ strongly converges in $L^2$, and $\II_{\deform_n}\wto \II_{\deform_\infty}$ weakly in $L^2$.
\end{proof}

  \begin{proof}[Proof of Theorem~\ref{T:shape}]
  \step 1 {Approximation of $\II_*$}
  We claim that there exist sequences $(\kappa_n)\subset L^2(S)$ and $(R_n)\subset L^2(S;\SO 2)$ such that
  \begin{align*}
    &\kappa_n(R_ne_1\otimes R_ne_1)\to\II_*\text{ strongly in }L^2,\\
    &(\kappa_n), (R_n)\text{ strongly converge in $L^2$},
  \end{align*}
  and for each $n\in\N$, the functions $\kappa_n,R_n$ are piecewise constant functions subject to a partition of $S$ (up to a null-set) into finitely many open, mutually disjoint sets; (in fact, the argument below yields a specific partition of $S$ into equally sized lattice-squares contained in $S$ and a remaining set close to the boundary of $S$). We prove this claim in two steps. In the first step, we establish the representation $\II_*=\kappa(Re_1\otimes Re_1)$ with $\kappa\in L^2(S)$ and $R:S\to\SO 2$ measurable. In the second step, we approximate $\kappa$ and $R$ by functions $\kappa_n:S\to\R$ and $R_n:S\to\SO 2$ that are piecewise constant subordinate to a finite partition of $S$. The approximation $\kappa_n$ can be easily obtained, e.g., by an $L^2$-projection of $\kappa$ onto the piecewise constant functions. On the other hand, the construction of $R_n$ is a bit more delicate, since the target manifold $\SO 2$ is a non-convex subset of $\R^{2\times 2}$. Therefore, in the first step we additionally show that the $R$ in the representation of $\II_*$ can be choosen to be locally Lipschitz. For the construction of $(\kappa_n,R_n)$ it is convenient to consider a dyadic decomposition of $S$: Below, we denote by $\mathcal Q_{\rm dyadic}\colonequals\{2^{-n}(\Z^2+[-\frac12,\frac12)^2)\,:\,n\in\N\}$ the set of dyadic squares in $\R^2$, and by $x_\Box\in 2^{-n}\Z^2,n\in\N$, the center point of the cube $\Box\colonequals x_\Box+2^{-n}[-\frac12,\frac12)^2\in\mathcal Q_{\rm dyadic}$.

  \smallskip

  {\it Substep 1.1. Representation of $\II_*$}

  We claim that there exist $\kappa\in L^2(S)$, $R:S\to\SO 2$ measurable, and a partition $\{\Box\}_{\Box\in\mathcal Q}$ of $S$ into (at most countably many) dyadic, mutually disjoint squares $\Box\in\mathcal Q\subset\mathcal Q_{\rm dyadic}$ such that
  \begin{equation}\label{eq:repr}
    \II_*=\kappa(R e_1\otimes Re_1)\text{ almost everywhere in }S,
  \end{equation}
  and for all $\Box\in\mathcal Q$ the restriction $R\vert_\Box$ is Lipschitz with
  \begin{equation}\label{eq:Lipsch}
    \operatorname{Lip}(R\vert_\Box)\leq C\operatorname{dist}(\Box,\R^2\setminus S)^{-1},
  \end{equation}
  for a universal constant $C$. For the argument let $\mathcal Q$ denote a dyadic Whitney covering, that is, $\mathcal Q$ is a countable family of mutually disjoint, dyadic squares of the form $\Box\in 2^{-n}(\Z^2+[-\frac12,\frac12)^2)$, $n\in\N$, such that
  \begin{equation*}
    S=\bigcup_{\Box\in\mathcal Q}\Box,\qquad \forall \Box\in\mathcal Q\,:\,\frac1{30}\dist(\Box,\R^2\setminus S)\leq\operatorname{diam}(\Box)\leq \frac1{10}\dist(\Box,\R^2\setminus S).
  \end{equation*}
  Let $\Box\in\mathcal Q$. By construction there exists a ball $B$ such that $\Box\subset B\subset 2B\subset S$. Thus, we may apply \cite[Lemma 8]{neukamm2015homogenization} to obtain a Lipschitz field $\xi:\Box\to S^1$ such that $\II_*(x')\in\operatorname{span}(\xi(x')\otimes \xi(x'))$ for a.e.~$x'\in \Box$, and $\operatorname{Lip}(\xi)\leq\frac{2}{\operatorname{diam}(\Box)}$. Since $\SO 2\ni R\mapsto Re_1\in S^1$ is a diffeomorphism, we find a Lipschitz field $R:\Box\to\SO 2$ such that $\xi=Re_1$ and $\operatorname{Lip}(R)\leq C\dist(\Box,\R^2\setminus S)^{-1}$. The representation \eqref{eq:repr} (on $\Box$) then follows by setting $\kappa\colonequals\II_*:(Re_1\otimes Re_1)$. Since this can be done for all $\Box\in\mathcal Q$, the claim follows.

  {\it Substep 1.2. Conclusion}

  We only explain the construction of $(R_n)$, since the ``linear'' construction of $\kappa_n$ is easier. Here comes the argument:  For $n\in\N$ set $S_n\colonequals\bigcup\{\Box\,:\,\Box\in\mathcal Q\text{ with }\operatorname{diam}(\Box)\geq 2^{-n}\,\}$. By construction, this is a finite union of squares in $\mathcal Q$ and we have $S_n\uparrow S$ for $n\to\infty$. We partition $S_n$ into (finer) dyadic squares with diameter $2^{-2n}$. To that end set $\mathcal Q_n\colonequals\{\Box\in\mathcal Q_{\rm dyadic}\,:\,\operatorname{diam}(\Box)=2^{-2n}\text{ and }\Box\subseteq S_n\,\}$.
  We now define the piecewise constant approximation $R_n:S\to\SO 2$ as
  \begin{equation*}
    R_n(x')\colonequals
    \begin{cases}
      R(x_\Box)&\text{if }x'\in \Box\text{ for some }\Box\in\mathcal Q_n,\\
      I_{2\times 2}&\text{else.}
    \end{cases}
  \end{equation*}
  (Recall that $x_\Box$ denotes the center point of $\Box$). We claim that $R_n(x')\to R(x)$ for all $x'\in S$ (and thus strongly in $L^2$ by the dominated convergence theorem). For the argument let $x'\in S_n$ and choose $n\in\N$ large enough such that $x'\in S_n$. Then there exists $\Box\in\mathcal Q_n$ with $x'\in \Box$ and $\Box'\in\mathcal Q$ with $\Box\subset \Box'$ and $\operatorname{diam}(\Box')\geq 2^{-n}$. Hence, by \eqref{eq:Lipsch}
  \begin{equation*}
    |R_n(x')-R(x')|=|R(x_\Box)-R(x')|\leq \operatorname{diam}(\Box)\operatorname{Lip}(\kappa\vert_{\Box'})\leq C2^{-2n}2^n,
  \end{equation*}
  and thus $R_n(x')\to R(x')$ for $n\to\infty$.

  \step 2 {Definition of a sequence of structured composites and its limit}
  For $n\in\N$ we define, as in Corollary~\ref{C:transf},
  \begin{eqnarray*}
    W_n(x,y,F)&\colonequals&\bar W(\kappa_n(x');x_3,R_n(x')^\top y,F\widehat R_n(x')),\\
    B_n(x,y)&\colonequals&\widehat R_n(x')\bar B(\kappa_n(x');x_3,R_n(x')^\top y)\widehat R_n^\top(x'),
  \end{eqnarray*}
  and note that $(W_n,B_n)$ is a structured composite. (We remark that $(W_n,B_n)$ is Borel measurable, since $(\bar W,\bar B)$ is measurable and $(\kappa_n,R_n)$ is piecewise constant subordinate to a finite partition of $S$ into measurable sets). By Theorem~\ref{T1}, the functionals $\mathcal I_n^h:H^1(\Omega;\R^3)\to[0,\infty]$,
  \begin{equation*}
    \mathcal I_n^h(\deform)\colonequals\frac1{h^2}\int_{\Omega} W_n(x,\tfrac{x'}{\e(h)},\nabla_h \deform(x)(I_{3\times 3}-hB_n(x,\tfrac{x'}{\e(h)})))\dd x.
  \end{equation*}
  $\Gamma$-converge for $h\to 0$ to a functional of the form $\mathcal I^\gamma_{\hom,n}(\cdot)+\mathcal I^\gamma_{\rm res}(B_n)$ with
  \begin{equation*}
    \mathcal I^\gamma_{\hom,n}(\deform)=
    \begin{cases}
      \int_S Q^\gamma_{\hom,n}(x',\II_{\deform}-B^\gamma_{\eff,n}(x'))\dd x'&\text{for }\deform\in H^2_{\iso}(S;\R^3),\\
      \infty&\text{else.}
    \end{cases}
  \end{equation*}
  Since $(W_n,B_n)$ is defined by locally rotating the building block $(\bar W,\bar B)$, we deduce with help of Corollary~\ref{C:transf} that
  \begin{equation*}
    Q^\gamma_{\hom,n}(x',G)=Q^\gamma_{\hom}(\kappa_n(x');R_n(x')^\top GR_n(x')),\qquad B^\gamma_{\hom,n}(x')=R_n(x')\bar B^\gamma_{\eff}(\kappa_n(x'))R_n(x')^\top.
  \end{equation*}
  Furthermore, from the pointwise convergence of $(\kappa_n,R_n)$ we infer that for a.e.~$x'\in S$ and $G\in\R^{2\times 2}_{\sym}$ we have
  \begin{eqnarray*}
    Q^\gamma_{\hom,n}(x',G)&\to& Q^\gamma_{\hom}(x',G)\colonequals Q^\gamma_{\hom}(\kappa(x');R(x')^\top GR(x')),\\
    B^\gamma_{\eff,n}(x')&\to& B^\gamma_{\eff}(x')\colonequals R(x') \bar B^\gamma_{\eff}(\kappa(x'))R(x')^\top,
  \end{eqnarray*}
  as $n\to \infty$. Hence, Lemma~\ref{L:gamma-spatial} implies that $\mathcal I^\gamma_{\hom,n}$ $\Gamma$-converges for $n\to\infty$ to the functional
  \begin{equation*}
    \mathcal I^\gamma_{\hom}(\deform)=\int_S Q^\gamma_{\hom}(x',\II_{\deform}(x')-B^\gamma_{\eff}(x'))\dd x'.
  \end{equation*}
  In view of the definition of $(Q^\gamma_{\hom},B^\gamma_{\eff})$ and the identity $\II_*=\kappa(Re_1\otimes Re_1)$, the algebraic minimization problem $\R^{2\times 2}_{\sym}\ni G\mapsto Q^\gamma_{\hom}(x',G-B^\gamma_{\eff}(x'))$ subject to $\det G=0$ is uniquely minimized by $G=\II_*(x')$ for a.e. $x'\in S$. We thus conclude that the minimizer of $\mathcal I^\gamma_{\hom}$ is unique in the sense that any minimizer $\deform_*\in H^2_{\iso}(S;\R^3)$ of $\mathcal I^\gamma_{\hom}$ satisfies
  \begin{equation}\label{T:shape:eq2}
    \II_{\deform_*}=\II_*\text{ a.e. in }S.
  \end{equation}

  \step 3 {Convergence of minimizers as $n\to\infty$}
  We consider
  \begin{equation*}
    \Delta_n\colonequals\sup\{\|\II_{\deform_*}-\II_*\|_{L^2(S)}\,:\,\deform_*\in\argmin\mathcal I^\gamma_{\hom,n}\},
  \end{equation*}
  and claim that $\lim\limits_{n\to\infty}\Delta_n=0$.
  For the argument, we choose $\deform_n\in \argmin\mathcal I^\gamma_{\hom,n}$ with $\fint_S \deform_n=0$ such that
  \begin{equation}\label{T:shape:eq1}
    \Delta_n\leq \|\II_{\deform_n}-\II_*\|_{L^2(S)}+\frac1n.
  \end{equation}
  By passing to a subsequence (not relabeled) we may assume that $\limsup\limits_{n\to\infty}\Delta_n=\lim\limits_{n\to\infty}\Delta_n$ and that $\deform_n\wto \deform_\infty$ weakly in $H^2(S;\R^3)$ for some $\deform_\infty\in H^2_{\iso}(S;\R^3)$. An application of Lemma~\ref{L:gamma-spatial} shows that $\deform_\infty\in\argmin\mathcal I^\gamma_{\hom}$, and thus $\II_{\deform_\infty}=\II_*$ in view of \eqref{T:shape:eq2}. Moreover, the lemma implies that $\deform_n\to \deform_\infty$ strongly in $H^2(S;\R^3)$, and thus we conclude that $\II_{\deform_n})\to \II_{\deform_\infty}=\II_*$ strongly in $L^2(S)$, which in combination with \eqref{T:shape:eq1} implies $\Delta_n\to 0$.

  \step 4 {Conclusion}
  By Step~3 we may choose $n$ large enough such that
  \begin{equation}\label{T:shape:eq4}
    \sup\{\|\II_{\deform_*}-\II_*\|_{L^2(S)}\,:\,\deform_*\in\argmin\mathcal I^\gamma_{\hom,n}\}\leq \delta.
  \end{equation}
  We claim that $(W_n,B_n)$ is the sought after finitely structured composite. Indeed, if $(\deform_h)\subset H^1(\Omega;\R^3)$ is an almost minimizing sequence with $\fint_\Omega \deform_h=0$, then Theorem~\ref{T1}~\ref{item:T1:compactness} implies that from any subsequence we can extract a subsequence that converges in the sense of \eqref{T1:conv} to a limit $\deform_*\in H^2_{\iso}(S;\R^3)$. Since $\mathcal I^h_n\to \mathcal I_{\hom,n}^\gamma+\mathcal I^\gamma_{\rm res}(B_n)$ in the sense of $\Gamma$-convergence, we conclude that $\deform_*$ is a minimizer of the limiting energy and thus of $\mathcal I_{\hom,n}^\gamma$. Finally, $\|\II_{\deform}-\II_*\|_{L^2(S)}\leq \delta$ follows from \eqref{T:shape:eq4}.

\end{proof}

\subsection*{Acknowledgments}
The authors acknowledge support by the German Research Foundation (DFG) via the
research unit FOR 3013 ``Vector- and tensor-valued surface PDEs'' (grant no.~NE2138/3-1).

\section{Appendix}

\subsection{Spaces of \texorpdfstring{$\Lambda$}{Lambda}-periodic functions}\label{A:lambda}
In this section we introduce various function spaces for periodic functions.
Unless stated otherwise, in this section, $\Lambda$ denotes a matrix in $\R^{2\times 2}$ such that $\Lambda^\top\Lambda>0$.
The column--vectors of $\Lambda$ generate a Bravais lattice and we consider functions that are periodic w.r.t.~this lattice.
From Definition~\ref{ass:pcperiodicity} \ref{ass:pcperiodicity:item1} we recall that $u:\R^2\to\R$ is called $\Lambda$-periodic, if $u(y+\tau)=u(y)$ for all $\tau\in\Lambda\Z^2$ and almost every $y\in\R^2$.
We associate with $\Lambda$ the parallelogram $Y_\Lambda\colonequals\Lambda[-\frac12,\frac12)^2$. It represents the \emph{reference cell of $\Lambda$-periodicity}. We denote by
    \begin{equation*}
      C(\Lambda)\colonequals\{u\in C(\R^2)\,:\,u\text{ is }\Lambda\text{-periodic.}\,\},
    \end{equation*}
    the space of continuous, $\Lambda$-periodic functions; it is a Banach space when endowed with the norm $\|u\|_{L^\infty(Y_\Lambda)}$.
    We denote by
    \begin{equation*}
      L^p(\Lambda)\colonequals\{u\in L^p_{\rm loc}(\R^2)\,:\,u\text{ is }\Lambda\text{-periodic.}\,\},
    \end{equation*}
    the space of $\Lambda$-periodic $L^p$-functions; it is a Banach space when endowed with the norm of $L^p(Y_\Lambda)$.
    We denote by $$H^1(\Lambda)\colonequals\{u\in H^1_{\rm loc}(\R^2)\,:\,u\text{ is }\Lambda\text{-periodic.},\}$$ the space of $\Lambda$-periodic $H^1$-functions; it is a Hilbert space when endowed with the norm of $H^1(Y_\Lambda)$.

    In this paper, we simultaneously have to deal with $L^2$-functions that are periodic w.r.t.~different lattices. Therefore, for the analysis, it is convenient to introduce the common superspace
    \begin{equation}\label{D:uloc}
      L^2_{\uloc}(\R^2)\colonequals\{u\in L^{2}_{\rm loc}(\R^2)\,:\,\|u\|_{L^2_{\uloc}}<\infty\},\qquad \|u\|_{L^2_{\uloc}}\colonequals\sup_{z\in\R^2}\|u(\cdot+z)\|_{L^2((-\frac12,\frac12)^2)},
    \end{equation}
    which is a Banach space.
    We  note that $L^2(\Lambda)$ is a closed subspace of $L^2_{\uloc}(\R^2)$. In particular, given a family  $\{\Lambda_j\}_{j\in J}\subset\R^{2\times 2}$ satisfying \eqref{ass:pcperidicity:eq}, then there exists a constant $C$ only depending on $C_\Lambda$ of \eqref{ass:pcperidicity:eq}, such that for all $j\in J$ and $u\in L^2(\Lambda_j)$ we have
    \begin{equation*}
      \frac{1}{C}\|u\|_{L^2(\Lambda_j)}\leq  \|u\|_{L^2_\uloc}\leq C\|u\|_{L^2(\Lambda_j)}.
    \end{equation*}

\subsection{Two-scale convergence for grained structures}\label{A:twoscale}

  In this section discuss the notion of two-scale convergence \emph{for grained microstructures} introduced in Definition~\ref{D:twoscale}. It is a slight variant of the notion introduced in \cite{neukamm2010homogenization}. Throughout the section we suppose that $\{\Lambda_j\}_{j\in J}$ and $\{S_j\}_{j\in J}$ are as in Assumption~\ref{ass:pcperiodicity}, and  $h\mapsto\e(h)$ satisfies Assumption~\ref{A:gamma}. The notion of two-scale convergence of Definition~\ref{D:twoscale} is tailored made to keep track of locally periodic structures that frequently appear in this paper. In particular, within this notion, sequences of the form
  \begin{equation*}
    \varphi_h(x)\colonequals\varphi(x,\tfrac{x'}{\e(h)}),
  \end{equation*}
  strongly two-scale converge to $\varphi$, if the function $\varphi(x,y)$ (next to integrability properties) is \emph{locally periodic} in the sense of \eqref{eq:locper}.
  Note that the elastic tensor $\mathbb L$ associated with $Q$ in Assumption~\ref{ass:W} and the prestrain tensor $B$ of Assumption~\ref{ass:W} (iii) satisfy this periodicity assumption.
  In the single grain case (when the lattice of periodicity is everywhere the same in $S$, i.e.~$\Lambda_j=\Lambda$), the two-scale limit of a sequence in $L^2$ is given by a function in $L^2(\Omega;L^2(\Lambda))$. In the multiple grain case, the lattice of periodicity changes from grain to grain, and thus the family of spaces $L^2(\Lambda_j)\subset L^2_{\uloc}(\R^2)$, $j\in J$, needs to be considered.
  A simple, yet useful, observation is that a sequence $(u_h)\subset L^2(\Omega)$ two-scale converges in the sense of Definition~\ref{D:twoscale}, if and only if for all $j\in J$ the restrictions $u_h\vert_{S_j}$ two-scale converge in the \emph{usual} sense:
\begin{lemma}\label{L:twoscale:parad}
  Let $(u_h)$ be a sequence in $L^2(\Omega)$ and $u\in L^2(\Omega;L^2_{\uloc}(\R^2))$. The following are equivalent
  \begin{enumerate}[(a)]
  \item $u_h\wtto u$ weakly two-scale in $L^2$ (in the sense of Definition~\ref{D:twoscale}).
  \item For all $j\in J$ we have $u\vert_{S_j\times(-\tfrac12,\tfrac12)\times\R^2}=u_j$ where $u_j\in L^2(S_j\times(-\tfrac12,\tfrac12);L^2(\Lambda_j))$ is the weak two-scale limit of $u_{h,j}\colonequals u_h\vert_{S_j\times(-\tfrac12,\tfrac12)}$ in the usual sense, that is, $(u_j)$ is bounded in $L^2(S_j\times(-\tfrac12,\tfrac12))$ and
    \begin{equation}\label{eq:L:twoscale:parad}
      \int_{S_j\times(-\tfrac12,\tfrac12)}u_{h,j}(x)\varphi(x,\tfrac{x'}{\e(h)})\dd x\to \int_{S_j}\fint_{\Box_{\Lambda_j}}u_j(x,y)\varphi(x,y)\dd(x_3,y)\dd x',
    \end{equation}
    for all $\varphi\in C^\infty_c(S_j\times(-\tfrac12,\tfrac12);C(\Lambda_j))$.
  \end{enumerate}
\end{lemma}
With help of this lemma we can lift various properties of two-scale convergence (in the usual sense) to the variant introduced in Definition~\ref{D:twoscale}:
\begin{lemma}\label{L:twoscale:prop}
  \begin{enumerate}[(a)]
  \item \label{L:twoscale:prop:a} (Compactness). Let $(u_h)\subset L^2(\Omega)$ be a bounded sequence. Then there exists a subsequence and $u\in L^2(\Omega;L^2_{\uloc}(\R^2))$ satisfying \eqref{eq:locper} such that $u_h\wtto u$ weakly two-scale in $L^2$.
  \item \label{L:twoscale:prop:b}  (Lower semicontinuity of the norm). Let $(u_h)\subseteq L^2(\Omega)$ weakly two-scale converge to $u$. Then
    \begin{equation*}
      \liminf\limits_{h\to 0}\int_{\Omega}|u_h|^2\dd x\geq\sum_{j\in J}\int_{S_j}\fint_{\Box_{\Lambda_j}}|u|^2\dd(x_3,y)\dd x_3.
    \end{equation*}
  \item \label{L:twoscale:prop:c}  (Strong times weak). Let $(u_h)\subseteq L^2(\Omega)$ and $u\in L^2(\Omega;L^2_{\uloc}(\R^2))$. Then  $u_h\stto u$ strongly two-scale, if and only if
        \begin{equation*}
          \int_\Omega u_h(x)\varphi_h(x)\dd x\to \sum_{j\in J}\int_{S_j}\fint_{\Box_{\Lambda_j}}u(x,y)\varphi(x,y)\dd(x_3,y)\dd x_3
        \end{equation*}
        for all weakly two-scale converging sequences $(\varphi_h)\subset L^2(\Omega)$ with limit $\varphi$.
  \item \label{L:twoscale:prop:d} (Approximation of two-scale limits). Let $u\in L^2(\Omega;L^2_{\uloc}(\R^2))$ satisfy \eqref{eq:locper}. Then there exists a sequence $(u_h)\subset C^\infty_c(\Omega)$ such that $u_h\stto u$ strongly two-scale in $L^2$.
  \end{enumerate}
\end{lemma}
For the proof we refer to \cite{neukamm2010homogenization, neukamm2012rigorous}, where the single-grain case (i.e., $\Lambda_j=I_{2\times 2}$) is discussed. The extension to the above setting is obvious and left to the reader.
\medskip

In the context of dimension reduction we are especially interested in two-scale limits of sequences of scaled gradients $\nabla_hu_h$ of displacements $(u_h)\subset H^1(\Omega;\R^3)$. As we shall see, a two-scale limit in $L^2$ of such a sequence is a vector field $F\in L^2(\Omega;L^2_{\uloc}(\R^2);\R^{3\times 3})$ that satisfies \eqref{eq:locper} and that can be represented with a help of potential $\varphi$ in the sense that
\begin{equation*}
  F=\nabla_\gamma\varphi,\qquad \nabla_\gamma=(\nabla_y,\frac{1}{\gamma}\partial_{3}).
\end{equation*}
For the precise statement we need to introduce the Banach space
\begin{equation}\label{def:Hgamma}
  \begin{aligned}
  &H^1_{\gamma,\uloc}\colonequals\big\{u\in H^{1}_{\loc}((-\tfrac12,\tfrac12)\times\R^2;\R^3)\,:\,\|u\|_{H^1_{\gamma,\uloc}}<\infty\big\},\\
  & \|u\|_{H^1_{\gamma,\uloc}}^2\colonequals \sup_{z\in\R^2}\int_{(-\frac12,\frac12)^3}|u(x_3,y+z)|^2+|\nabla_\gamma u(x_3,y+z)|^2\dd(x_3,y),
\end{aligned}
\end{equation}
and the subspace of $\Lambda_j$-periodic functions,
\begin{equation}\label{def:Hgammaper}
  \begin{aligned}
  &H^1_{\gamma}(\Box_{\Lambda_j};\R^3)\colonequals H^{1}_{\gamma,\uloc}\cap L^2((-\tfrac12,\tfrac12);L^2(\Lambda_j;\R^3)),\\
  & \|u\|_{H^1_{\gamma}(\Box_{\Lambda_j})}^2\colonequals\fint_{\Box_{\Lambda_j}}|u|^2+|\nabla_\gamma u|^2\dd(x_3,y).
\end{aligned}
\end{equation}
Note that $H^1_{\gamma}(\Box_{\Lambda_j})$ is a closed subspace of $H^1_{\gamma,\uloc}((-\tfrac12,\tfrac12)\times\R^2)$, and we have
\begin{equation}\label{app:embed:locper}
  \frac{1}{C}\|\cdot\|_{H^1_{\gamma}(\Box_{\Lambda_j})}\leq \|\cdot\|_{H^1_{\gamma,\uloc}}\leq C  \|\cdot\|_{H^1_{\gamma}(\Box_{\Lambda_j})},
\end{equation}
for some $C=C(\gamma,C_\Lambda$).
\begin{proposition}[Two-scale limit of scaled gradients]\label{P:twoscale:gradrecov}
  \begin{enumerate}[(a)]
  \item Let $(u_h)\subset H^1(\Omega;\R^3)$ be a sequence and suppose that $(\nabla_hu_h)\subseteq L^2(\Omega;\R^{3\times 3})$ is bounded. Then there exists $u\in H^1(S;\R^3)$ and $\varphi\in L^2(S;H^1_{\gamma,\uloc})$ that satisfies \eqref{eq:locper} such that
    \begin{equation*}
      \nabla_hu_h\wtto (\nabla'u,0)+\nabla_\gamma\varphi\qquad\text{weakly two-scale in $L^2$}.
    \end{equation*}
  \item Let $\varphi\in L^2(S;H^1_{\gamma,\uloc})$ satisfy \eqref{eq:locper}. Then there exists a sequence $(\varphi_h)\subset C^\infty_c(\Omega;\R^3)$  such that
    \begin{equation*}
      \nabla_h\varphi_h\stto \nabla_\gamma\varphi\qquad\text{strongly two-scale in $L^2$},
    \end{equation*}
    and
    \begin{equation*}
      \lim\limits_{h\to 0}\sqrt{h}(\|\varphi_h\|_{L^\infty(\Omega)}+\|\nabla_h\varphi_h\|_{L^\infty(\Omega)})=0.
    \end{equation*}
  \end{enumerate}
\end{proposition}
\begin{proof}With help of Lemma~\ref{L:twoscale:parad} the argument of \cite[Section 6.3]{neukamm2010homogenization} easily extends to the case of grained two-scale convergence, see also \cite{neukamm2012rigorous}.
\end{proof}

\begin{lemma}[(Lower semi-)continuity of convex functionals]\label{L:twoscale:lsc}
  Let $Q$ be as in Assumption~\ref{ass:W} and let $G\in L^2(\Omega; L^2_{\uloc}(\R^2;\R^{3\times 3}))$. Then:
  \begin{enumerate}[(a)]
  \item Let $(G_h)\subset L^2(\Omega;\R^{3\times 3})$ be a sequence that weakly two-scale converges to $G$. Then
    \begin{equation*}
      \liminf\limits_{h\to 0}\int_{\Omega} Q(x,\tfrac{x'}{\e(h)},G_h(x))\dd x\geq \sum_{j\in J}\int_{S_j}\fint_{\Box_{\Lambda_j}}Q(x,y,G(x,y))\dd(x_3,y)\dd x'.
    \end{equation*}
  \item Let $(G_h)\subset L^2(\Omega;\R^{3\times 3})$ be a sequence that strongly two-scale converges to $G$. Then
    \begin{equation*}
      \lim\limits_{h\to 0}\int_{\Omega} Q(x,\tfrac{x'}{\e(h)},G_h(x))\dd x= \sum_{j\in J}\int_{S_j}\fint_{\Box_{\Lambda_j}}Q(x,y,G(x,y))\dd(x_3,y)\dd x'.
    \end{equation*}
  \end{enumerate}
\end{lemma}
\begin{proof}With help of Lemma~\ref{L:twoscale:parad} the argument of \cite[Section 3.2]{neukamm2010homogenization} easily extends to the case of grained two-scale convergence.
\end{proof}

\begin{lemma}[Korn's inequality for scaled gradient]\label{L:korn-lambda}
  Let $\gamma\in(0,\infty)$ and $\bar C>0$. Then there exists a constant $C=C(\gamma,\bar C)$ such that for all $\Lambda\in\R^{2\times 2}$ with $\frac{1}{\bar C}\leq\Lambda^\top\Lambda\leq\bar C$ the following Korn's inequality holds:
  For all $\varphi\in H^1_{\gamma}(\Box_\Lambda;\R^3)$, we have
  \begin{equation}\label{eq:korn:1}
    \fint_{\Box_\Lambda}|\nabla_\gamma\varphi|^2\leq C    \fint_{\Box_\Lambda}|\sym \nabla_\gamma\varphi|^2.
  \end{equation}
\end{lemma}
\begin{proof}
  In the following we write $\lesssim$ if $\leq$ holds up to a multiplicative constant only depending on $\gamma$ and $\bar C$. We consider the scaled function
  \begin{equation*}
    \tilde\varphi(x_3,y)\colonequals\frac{1}{\gamma}\varphi(x_3,\gamma y),\qquad \tilde\Lambda\colonequals\frac{1}{\gamma}\Lambda,\qquad \tilde Y\colonequals Y_{\tilde\Lambda}.
  \end{equation*}
  Then $\tilde\varphi\in H^1_\gamma((-\frac12,\frac12)\times\tilde\Lambda;\R^3)$ and $\nabla\tilde\varphi(x_3,y)=\nabla_\gamma\varphi(x_3,\gamma y)$, and thus \eqref{eq:korn:1} is equivalent to
  \begin{equation}\label{eq:korn:2}
    \fint_{(-\tfrac12,\tfrac12)\times Y_{\tilde\Lambda}}|\nabla\tilde\varphi|^2\leq C    \fint_{(-\tfrac12,\tfrac12)\times Y_{\tilde\Lambda}}|\sym\nabla\tilde\varphi|^2.
  \end{equation}
  Thus, it suffices to prove \eqref{eq:korn:2} for all $\tilde\varphi\in H^1_\gamma((-\frac12,\frac12)\times\tilde\Lambda;\R^3)$. In order to see that the constant can be choosen only depending on $\gamma$ and $\bar C$, we first note that there exists two concentric cubes $Q_i\colonequals\ell_i[-\frac12,\frac12)^2$ such that
  \begin{equation*}
    Q_1\subset \tilde Y\subset Q_2\qquad\text{and}\qquad\tilde C\colonequals\frac{|Q_2|}{|Q_1|}\lesssim 1.
  \end{equation*}
  Thus, we have
  \begin{equation}\label{pf:korn:3}
    \fint_{(-\frac12,\frac12)\times \tilde Y}|F|^2\lesssim    \fint_{(-\frac12,\frac12)\times Q_2}|F|^2,
  \end{equation}
  for all $F\in L^2((-\frac12,\frac12)\times Q_2;\R^{3\times 3})$. Furthermore, by exploiting $\tilde\Lambda$-periodicity, we see that
  \begin{equation}\label{pf:korn:4}
    \fint_{(-\frac12,\frac12)\times 2Q_2}|F|^2\lesssim    \fint_{(-\frac12,\frac12)\times \tilde Y}|F|^2
  \end{equation}
  for all $\tilde\Lambda$-periodic functions $F\in L^2((-\frac12,\frac12)\times \tilde\Lambda;\R^{3\times 3})$.

  By the standard Korn's inequality there exists $K\in\R^{3\times 3}_{\rm skw}$ such that
  \begin{equation}\label{eq:korn:5}
    \fint_{(-\frac1 2,\frac1 2)\times 2Q_2}|\nabla\psi-K|^2\lesssim\fint_{(-\tfrac 1 2,\tfrac1 2)\times 2Q_2}|\sym\nabla\psi|^2,
  \end{equation}
  for all $\psi\in H^1((-\tfrac1,\tfrac12)\times 2Q_2;\R^3)$. In particular, we may apply this estimate to the scaled function $\tilde\varphi$. Combined with \eqref{pf:korn:3} and \eqref{pf:korn:4} we get
  \begin{eqnarray*}
    \fint_{(-\tfrac12,\tfrac12)\times \tilde Y}|\nabla\tilde\varphi-K|^2
    &\lesssim & \fint_{(-\frac12,\frac12)\times 2Q_2}|\nabla\tilde\varphi-K|^2\,\lesssim\,\fint_{(-\frac12,\frac12)\times 2Q_2}|\sym\nabla\tilde\varphi|^2\\
    &\lesssim & \fint_{(-\frac12,\frac12)\times \tilde Y}|\sym\nabla\tilde\varphi|^2.
  \end{eqnarray*}
  In order to conclude \eqref{eq:korn:2}, we need to estimate $|K|$. To that end let $\tilde\Lambda_1$ and $\tilde\Lambda_2$ denote the first and second column vector of $\tilde\Lambda$, respectively. Consider
  \begin{equation*}
    u(x_3,y)\colonequals\tilde\varphi(x_3,y)-K
    \begin{pmatrix}
      y\\x_3
    \end{pmatrix}.
  \end{equation*}
  Then, by $\tilde\Lambda$-periodicity of $\tilde\varphi$, we have for a.e. $(x_3,y)\in (-\frac12,\frac12)\times\R^2$,
  \begin{equation*}
    K
    \begin{pmatrix}
      \tilde\Lambda_i\\0
    \end{pmatrix}=u(x_3,y+\tilde\Lambda_i)-u(x_3,y)=\int_0^1\nabla' u(x_3,y+s\tilde\Lambda_i)\cdot\tilde\Lambda_i\dd s.
  \end{equation*}
  Taking the square, integration w.r.t.~$y$, and Jensen's inequality yield
  \begin{eqnarray*}
    |K
    \begin{pmatrix}
      \tilde\Lambda_i\\0
    \end{pmatrix}|^2&\leq& |\tilde\Lambda_i|^2\int_0^1\fint_{(-\frac12,\frac12)\times Q_2}|\nabla' u(x_3,y+s\Lambda_i)|^2\dd y\dd x_3\dd s\\
    &\lesssim&\fint_{(-\frac12,\frac12)\times 2Q_2}|\nabla u(x_3,y)|^2\dd y\dd x_3,
  \end{eqnarray*}
  where for the last estimate, we used that $y+s\Lambda_i\in 2Q_2$ for all $y\in Q_2$ and $s\in(0,1)$.
  Hence, combined with \eqref{eq:korn:5} and the definition of $u$, we get
  \begin{eqnarray*}
    \sum_{i=1,2}|K
    \begin{pmatrix}
      \tilde\Lambda_i\\0
    \end{pmatrix}|^2&\lesssim&\fint_{(-\frac12,\frac12)\times 2Q_2}|\nabla\tilde\varphi-K|^2\lesssim \fint_{(-\frac12,\frac12)\times \tilde Y}|\sym\nabla\tilde\varphi|^2.
  \end{eqnarray*}
  Moreover, by skew symmetry, we have $|K|^2\lesssim \sum_{i=1,2}|K
  \begin{pmatrix}
    \tilde\Lambda_i\\0
  \end{pmatrix}|^2$, and thus,
  \begin{equation*}
    \fint_{(-\tfrac12,\tfrac12)\times \tilde Y}|\nabla\tilde\varphi|^2\lesssim |K|^2+\fint_{(-\tfrac12,\tfrac12)\times \tilde Y}|\nabla\tilde\varphi-K|^2\lesssim \fint_{(-\frac12,\frac12)\times \tilde Y}|\sym\nabla\tilde\varphi|^2.
  \end{equation*}
\end{proof}

\bibliographystyle{plainnat}
  \bibliography{bibliography.bib}


\end{document}